\documentclass{article}
\newtheorem{Theorem}{Theorem}
\newtheorem{Lemma}[Theorem]{Lemma}
\newtheorem{Corollary}[Theorem]{Corollary}
\newtheorem{Proposition}[Theorem]{Proposition}
\newtheorem{Proof of Theorem}{Proof of Theorem}
\newtheorem{Definition}[Theorem]{Definition}
\newtheorem{Example}[Theorem]{Example}

\newtheorem{Claim}[Theorem]{Claim}
\begin{document}
\title{Uniform Random Planar Graphs with~Degree~Constraints}
\author{Chris Dowden}
\date{August 11, 2008}
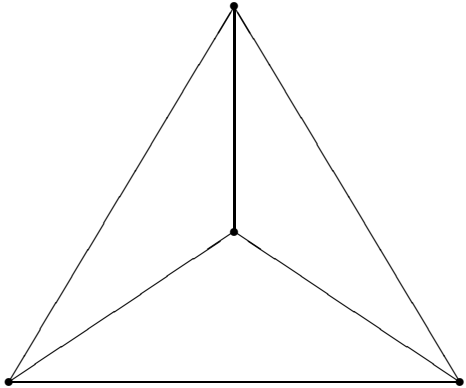
\maketitle
\setlength{\unitlength}{1cm}
\begin{picture}(20,7)
\put(2.5,0){\line(1,0){6}}
\put(2.5,0){\line(3,2){3}}
\put(2.5,0){\line(3,5){3}}
\put(5.5,2){\line(0,1){3}}
\put(8.5,0){\line(-3,2){3}}
\put(8.5,0){\line(-3,5){3}}
\put(2.5,0){\circle*{0.1}}
\put(8.5,0){\circle*{0.1}}
\put(5.5,2){\circle*{0.1}}
\put(5.5,5){\circle*{0.1}}

\end{picture}
\thispagestyle{empty}

\newpage

Uniform Random Planar Graphs with Degree Constraints 

Chris Dowden (Merton College) 

D.Phil.~thesis, Trinity 2008

\begin{abstract}
Random planar graphs have been the subject of much recent work.
Many basic properties of the standard uniform random planar graph $P_{n}$,
by which we mean a graph chosen uniformly at random 
from the set of all planar graphs with vertex set $ \{ 1,2, \ldots, n \}$,
are now known,
and variations on this standard random graph are also attracting interest.

Prominent among the work on $P_{n}$ have been asymptotic results
for the probability that $P_{n}$ will be connected or contain given components/ subgraphs.
Such progress has been achieved through a combination of counting arguments \cite{mcd}
and a generating function approach \cite{gim}.

More recently,
attention has turned to $P_{n,m}$,
the graph taken uniformly at random from the set of all planar graphs on $\{ 1,2, \ldots, n \}$ with exactly $m(n)$ edges
(this can be thought of as a uniform random planar graph with a constraint on the average degree).
In \cite{ger} and \cite{gim},
the case when $m(n) =~\!\lfloor qn \rfloor$ for fixed $q \in (1,3)$ has been investigated,
and results obtained for the events that 
$P_{n, \lfloor qn \rfloor}$ will be connected
and that $P_{n, \lfloor qn \rfloor}$ will contain given subgraphs.

In Part I of this thesis,
we use elementary counting arguments to extend the current knowledge of $P_{n,m}$.
We investigate the probability that $P_{n,m}$ will contain given components,
the probability that $P_{n,m}$ will contain given subgraphs,
and the probability that $P_{n,m}$ will be connected,
all for general $m(n)$,
and show that there is different behaviour
depending on which `region' the ratio $\frac{m(n)}{n}$ falls into.
In Part II,
we investigate the same three topics for a uniform random planar graph
with constraints on the maximum and minimum degrees. 
\end{abstract}

\newpage
\section*{Acknowledgements}

I would like to thank the University of Oxford, Merton College, the Statistics Department and the Mathematical Institute
for the use of various resources,
and the Engineering and Physical Sciences Research Council for providing me with funding.
I would also like to thank Louigi Addario-Berry for his work on Section~\ref{4reg},
and Stefanie Gerke and Alex Scott for refereeing the thesis.

On a personal note,
I am grateful to Dominic Welsh and Colin McDiarmid for taking the time to reply 
to my enquiries about coming to Oxford,
and also to the latter for supervising me.
Most importantly,
I would like to thank my parents for their love and support.

\newpage

\begin{displaymath}
\end{displaymath} \\

\phantom{p} 

\tableofcontents

\newpage
\section{Introduction} \label{aintro}

\subsection*{Random Graphs}

Over the last 50 years, 
random graphs 
have been the subject of much activity.
Two main types of random graph have been studied in particular --- 
the random graph with edge probability $p$ and the uniform random graph with $m$ edges
(where we take `graph' to mean `simple labelled graph',
as throughout this thesis).

The \textit{random graph with edge probability $p$} is the random graph on the vertex set $\{ 1,2, \ldots, n \}$
where each of the $\left( ^{n}_{2} \right)$ possible edges occur independently with probability $p=p(n)$.
Thus, we would expect the number of edges to be around~$p \left( ^{n}_{2} \right)$,
but any number is possible.
By contrast, the \textit{uniform random graph with $m$ edges} is the graph taken uniformly at random
from among all the graphs on $\{ 1,2, \ldots, n \}$ that have \textit{exactly} $m=m(n)$ edges.
Thus, each edge occurs with probability $\frac{m}{ \left( ^{n}_{2} \right) }$,
but not independently of the other edges.

An alternative way to define the uniform random graph with $m$ edges is via the \textit{random graph process},
where one starts at stage zero with an empty graph and inserts all $\left( ^{n}_{2} \right)$ edges
one by one in a (uniformly) random order.
The $m$th stage of this random graph process is a uniform random graph with $m$ edges.

Much is known about how the properties of our two types of random graph depend on the functions $p(n)$ and $m(n)$.
For example, it is known that $n^{-1} \log n$ is a `threshold' for the property that
a random graph with edge probabilty $\!p$ will be connected,
meaning that
\begin{displaymath}
\mathbf{P} [\textrm{connected}] \to
\left\{ \begin{array}{ll}
0 & \textrm{if $\frac{np(n)}{\log n} \to 0 \textrm{ as } n \to \infty$} \\
1 & \textrm{if $\frac{np(n)}{\log n} \to \infty \textrm{ as } n \to \infty$.} 
\end{array} \right. 
\end{displaymath} 
Thresholds (for both types of random graphs) are also known for properties such as 
`containing a subgraph isomorphic to $H$' or `containing a component isomorphic to $H$'
and it has become customary to refer to this development of random graphs as $p$ or $m$ grows as 
the `evolution' of the random graph.

An important tool in the investigation of the random graph with edge probabilty $p$ is the use of probabilistic methods.
For example, if $X$ is a function of a graph 
(such as `the number of subgraphs isomorphic to $H$ in a graph')
and we wish to know bounds for $\mathbf{P}[X=0]$,
then it is often helpful to work out $\mathbf{E}(X)$ and $\textrm{var}(X)$,
the calculation of which is facilitated by the fact that each edge occurs independently of all the others.
Certain results can then be transferred over to uniform random graphs,
since it can be shown that a close relationship exists between the random graph with edge probabilty $p$
and the uniform random graph with $m=p \left( ^{n}_{2} \right)$ edges. \\
\\

\subsection*{Planar Graphs} 

A graph is said to be \textit{planar} if it is possible to draw it in the plane 
(or, equivalently, on the sphere) in such a way that the edges do not cross, meeting only at vertices.
For example, $K_{4}$ is planar, since it can be drawn as on the title page.
In the rest of this subsection,
we shall collect together (without proof) various basic properties of planar graphs that will be useful later.

Throughout this thesis,
we shall use $\mathcal{P}(n)$ to denote the class of all planar graphs on the vertex set $\{ 1,2, \ldots, n \}$.
It can be shown that, for $n \geq 3$, 
the maximum number of edges of a graph in $\mathcal{P}(n)$
is $3n-6$, so $K_{5}$, for example, cannot be planar, since it has $10$ edges.
The maximum size 
(where we use `size' to mean the number of edges)
is achieved if and only if a graph is a \textit{triangulation},
i.e.~if and only if it is possible to draw the graph in the plane
in such a way that each face, including the outside face, is a triangle.
It is always possible to extend a planar graph to a triangulation by inserting extra edges.

Another useful property of $\mathcal{P}(n)$ is \textit{edge-addability}: \label{planar}
for each graph $G$ in $\mathcal{P}(n)$, 
any graph that is obtained from $G$ by adding an edge between two vertices in different components 
is also in $\mathcal{P}(n)$.
This follows from the important fact that we may draw $G$ in such a way that any given face is on the outside. \\ 
\\

\subsection*{Random Planar Graphs}

Since they were first investigated in \cite{den},
random planar graphs have generated much interest.
There are three main models --- 
the random planar graph with edge probability $p$,
the random planar graph process,
and the uniform random planar graph with $m$ edges. \\

The \textit{random planar graph with edge probability $p$} is defined to be the graph obtained by 
repeatedly sampling a (general) random graph with edge probability~$p$ until we find one that is planar.
Note that this planarity condition distorts the randomness in such a way that the probabilistic methods 
used for the original case are no longer helpful,
making this model difficult to study.

However, if we let $p=\frac{1}{2}$ then we obtain $P_{n}$,
the graph taken uniformly at random from the set of all graphs in $\mathcal{P}(n)$.
In \cite{mcd}, the limiting probability (as $n \to \infty$) was found for the event that
$P_{n}$ will contain a component isomorphic to an arbitrary fixed connected planar graph.
This result was proven using `counting' methods:
to calculate the proportion of graphs with a particular property,
we think of a way to construct graphs with this property,
we count how many graphs we can create this way,
and we count how many times each graph will be constructed (i.e.~the amount of double-counting).
An inherent difficulty with this approach is that we need to be able to count with accuracy.

A vital tool in the aforementioned result of \cite{mcd} was a rather precise asymptotic estimate of $|\mathcal{P}(n)|$,
which was proven in \cite{gim} using the concept of generating functions.
The \textit{exponential generating function} for a class $\mathcal{A}$
is defined to be $A(x) = \sum_{n \geq 0} \frac{a_{n}}{n!} x^{n}$,
where $a_{n}$ denotes the number of elements of $\mathcal{A}$ that have parameter $n$
(e.g.~$\mathcal{A}$ could be the class of planar graphs and $a_{n}$ could be $|\mathcal{P}(n)|$).
The structure of $\mathcal{A}$ can sometimes be used to produce algebraic equations involving $A(x)$
and its derivatives,
which can then be solved to find $A(x)$ explicitly.
By analysing the singularities of $A(x)$,
one may then be able to derive the asymptotic behaviour of $a_{n}$.
Unfortunately, the equations involved may be extremely difficult to solve. \\

The \textit{random planar graph process}
is a random graph process equipped with an additional acceptance test:
before we insert an edge,
we check whether the resulting graph would be planar and, if not, 
we reject the edge (and never look at it again).
Properties such as connectivity and containing given subgraphs were studied for this model in \cite{ger4},
using counting methods. \\

The \textit{uniform random planar graph with $m$ edges}, $P_{n,m}$, is
the graph taken uniformly at random from the set of all graphs in $\mathcal{P}(n)$
with exactly $m$ edges.
We shall use $\mathcal{P}(n,m)$ to denote this set.
Thus, the probability that $P_{n,m}$ will have a particular property is simply equal to 
the proportion of graphs in $\mathcal{P}(n,m)$ that have that property.
Unlike with general random graphs, it turns out that $P_{n,m}$ is not equivalent to the graph
obtained by the random planar graph process after $m$ edges have been accepted.

It is known from general uniform random graph theory 
(see, for example, Theorem 5.5 of \cite{jan})
that, \textit{asymptotically almost surely}
(a.a.s., that is, with probability tending to $1$ as $n$ tends to infinity),
a graph taken uniformly at random from the class of all size $m$ graphs on $\{1,2, \ldots, n \}$
will be planar if $\frac{m(n)}{n}$ is bounded above by $C < \frac{1}{2}$.
Thus, for this region of $m$,
uniform random planar graphs behave in the same way as general uniform random graphs.
Since the latter have already been extensively investigated,
the interest for planar graphs lies with the case when $m \geq Cn$.
Also, recall that (for $n \geq 3$) we must have $m \leq 3n-6$ for planarity to be possible.

Generating functions were used in \cite{gim} to give rather precise asymptotic expressions for both
$|\mathcal{P}(n, \lfloor qn \rfloor)|$ and $|\mathcal{P}_{c}(n, \lfloor qn \rfloor)|$ for fixed $q \in (1,3)$,
where $\mathcal{P}_{c}(n,m)$ denotes the class of connected graphs in $\mathcal{P}(n,m)$,
and these results can be combined to give an expression for 
$\mathbf{P}\left[P_{n,\lfloor qn \rfloor} \textrm{ will be connected}\right]$
for~$q \in~\!(1,3)$.
Also,  
$\mathbf{P}\left[P_{n,\lfloor qn \rfloor} \textrm{ will contain a given subgraph}\right]$
has been investigated in \cite{ger} using counting methods.
A more detailed summary of all such results shall be given in Section~\ref{previous}. \\ 
\\

\subsection*{Overview of Thesis}

This thesis shall be split into two distinct parts.
These shall deal, respectively,
with uniform random planar graphs with $m$ edges,
adding to the existing literature,
and then with the unexplored topic of random planar graphs with bounds on the minimum and maximum degrees.

In Part~\ref{I}, we shall use counting methods to investigate the behaviour of~$P_{n,m}$,
the uniform random planar graph with $n$ vertices and $m$ edges, as $m(n)$ varies.
We will extend the connectivity and subgraph results of \cite{ger} and \cite{gim} to cover~the case when 
$m$ is not of the form $\lfloor qn \rfloor$, 
as well as investigating the largely uncharted region of $m \leq (1+o(1))n$.
We shall also examine the thresholds~for~$P_{n,m}$ containing given components.
A summary of our main~results~is~given~on~page~\pageref{main cpt}.

Clearly, $P_{n,m}$ can be thought of as a random planar graph with constraints~on the average degree.
In Part~\ref{II},
we shall instead look at random planar graphs with constraints on the minimum and maximum degrees.
Again, we shall use counting methods to investigate the typical properties of such graphs 
and~see~how these vary with our constraints.
A summary of these results is given on page~\pageref{sum}.

\newpage
\part{The Evolution of Uniform Random Planar Graphs} \label{I}
\section{Outline of Part~\ref{I}} \label{intro}

As already mentioned in Section~\ref{aintro},
we shall now investigate how the properties of planar graphs change depending on the number of edges.
In particular, we will look at a graph $P_{n,m(n)}$ 
taken uniformly at random from the set of all labelled planar graphs on $\{ 1,2, \ldots, n \}$
with $m(n)$ edges,
and see how the behaviour of~$P_{n,m}$
varies with the ratio $\frac{m(n)}{n}$
(as mentioned in Section~\ref{aintro},
we may assume that $\liminf_{n \to \infty} \frac{m}{n} > 0$,
since otherwise our planarity condition has no impact).
We shall focus on three topics:
the probability that $P_{n,m}$ will contain given components,
the probability that $P_{n,m}$ will be connected,
and the probability that $P_{n,m}$ will contain given subgraphs.
Our objective will be to show exactly when,
in terms of $\frac{m(n)}{n}$,
these probabilities converge to $0$ ,
converge to $1$,
or are bounded away from both $0$ and $1$.

We shall start in Section~\ref{previous}
by giving a detailed summary of the state of knowledge of $P_{n,m}$ prior to this thesis.
We will also include results for $P_{n}$,
the graph taken uniformly at random from the set of all labelled planar graphs on~$\{ 1,2, \ldots, n \}$,
since this model is closely related to $P_{n,m}$.

In Sections~\ref{pen} and \ref{add} of this thesis, 
we shall then do the groundwork for our later theorems by investigating the number of pendant edges
(i.e.~the number of edges incident to a vertex of degree $1$)
and the number of addable edges
(i.e.~the number of edges that can be added individually to a planar graph without destroying planarity).
These results will be important ingredients for our later counting arguments.

We shall then start to use our ingredients to examine
the probability that $P_{n,m}$ will contain a component isomorphic to a given graph $H$,
the results for which are summarised on page~\pageref{main cpt}.
First, in Section~\ref{cptlow},
we will produce various lower bounds for this probability
(splitting into different cases depending on both $e(H)-|H|$ and $m(n)$),
where by `lower bound' we mean results such as 
$\mathbf{P}[P_{n,m} \textrm{ will contain a component isomorphic to }H] \to 1$
as $n \to \infty$ 
or 
$\liminf_{n \to \infty} 
\mathbf{P}[P_{n,m} \textrm{ will contain a component isomorphic to }H] > 0$,
rather than precise figures.
In Sections~\ref{conkappa} and~\ref{cyccpt},
we shall then obtain exactly complementary upper bounds.

The upper bounds of Section~\ref{conkappa} will, in fact,
be obtained through achieving another of our objectives
by producing an account of the probability that $P_{n,m}$ will be connected
(in which case it clearly won't contain any component of order~$<n$),
and a summary of the results for this topic
is also given on page~\pageref{main cpt}.
As a spin-off,
we shall obtain (in the second half of Section~\ref{conkappa}) 
some results on the total number of components in $P_{n,m}$,
which (although not one of our primary themes)
is quite an interesting subject in its own right.

In Sections~\ref{sub}--\ref{multsub},
we will turn our attention to 
$\mathbf{P}[P_{n,m} \textrm{ will contain a given}$ \textit{subgraph}],
again dealing separately with different cases
depending on the number of edges of the subgraph.
These results are again summarised on page~\pageref{main cpt},
for the simplified case when the subgraph is connected. 

Often we may prove slightly stronger results than those stated on page~\pageref{main cpt}.
For example, we might show that the probability that
$P_{n,m}$ has at least $t$ \textit{vertex-disjoint induced order-preserving} copies of $H$
converges to $1$,
rather than just
$\mathbf{P}\left[P_{n,m} \textrm{ has a copy of } H \right] \to 1$,
or
$\mathbf{P}\left[P_{n,m} \textrm{ has a component isomorphic to } H \right] > 1-e^{- \Omega (n)}$
rather than just
$\mathbf{P}\left[P_{n,m} \textrm{ has a component isomorphic to } H \right] \to 1$.
Typically, these extensions will not alter the idea of a proof,
but may make some of the details slightly more complicated. 

Throughout,
we shall use `Lemma', `Theorem' and `Corollary' in the usual way,
with `Proposition' reserved for those results that are at a tangent to
our three main objectives
(such as when we look at the total number of components in $P_{n,m}$)
and also for results that were given in other papers.

To aid the reader,
a diagram will be given at the start of each section to illustrate how it is structured.
Arrows are used to show the relationship between the results,
with those of that section highlighted in bold,
and the main theorems circled.

\newpage
\label{main cpt}
\subsection*{Summary of Component Results}

For a connected planar graph $H$,
let $\mathbf{P} := \mathbf{P}[P_{n,m} \textrm{ will have a component } \cong H]$. \\
\\ 
\begin{tabular}{c|c|c|c}
&
$e(H)<|H|$ &
$e(H)=|H|$ &
$e(H) > |H|$ \\
\hline
$0 < \textrm{\small{$\underline{\lim}$} } \frac{m}{n}$ &
$\mathbf{P} \to 1$ (Cor.~\ref{tree6}) &
\small{$\underline{\lim}$} $\mathbf{P} > 0$ (T\ref{gen4}) &
$\mathbf{P} \to 0$ (Thm.~\ref{cyc41}) \\
\& $\frac{m}{n} \leq 1+o(1)$ &
&
\small{$\overline{\lim}$} $\mathbf{P} < 1$ (L\ref{cyc52}) &
\\
\hline
$1 < \textrm{\small{$\underline{\lim}$} } \frac{m}{n}$ &
\small{$\underline{\lim}$} $\mathbf{P} > 0$ (T\ref{gen3}) &
\small{$\underline{\lim}$} $\mathbf{P} > 0$ (T\ref{gen3}) &
\small{$\underline{\lim}$} $\mathbf{P} > 0$ (T\ref{gen3}) \\
\& $\textrm{\small{$\overline{\lim}$} } \frac{m}{n} < 3$ &
\small{$\overline{\lim}$} $\mathbf{P} < 1$ (L\ref{conn1}) &
\small{$\overline{\lim}$} $\mathbf{P} < 1$ (L\ref{conn1}) &
\small{$\overline{\lim}$} $\mathbf{P} < 1$ (L\ref{conn1}) \\
\hline
$\frac{m}{n} \to 3$ &
$\mathbf{P} \to 0$ (Cor.~\ref{conn5}) &
$\mathbf{P} \to 0$ (Cor.~\ref{conn5}) &
$\mathbf{P} \to 0$ (Cor.~\ref{conn5})
\end{tabular}
\begin{displaymath}
\end{displaymath}

\subsection*{Summary of Connectivity Results}

Let $\mathbf{P}_{c} := \mathbf{P}[P_{n,m} \textrm{ will be connected}]$. \\
\\
\begin{tabular} {c|c}
$\frac{m}{n} \leq 1+o(1)$ &
$\mathbf{P}_{c} \to 0$ (Corollary~\ref{tree6}) \\
\hline
$1 < \underline{\lim}~\frac{m}{n}$ \& $\overline{\lim}~\frac{m}{n} < 3$ &
$\underline{\lim}~\mathbf{P}_{c} > 0$ (Lemma~\ref{conn1}) \\
&
$\overline{\lim}~\mathbf{P}_{c} < 1$ (Theorem~\ref{gen4}) \\
\hline
$\frac{m}{n} \to 3$ &
$\mathbf{P}_{c} \to 1$ (Corollary~\ref{conn5})
\end{tabular}
\begin{displaymath}
\end{displaymath}

\subsection*{Summary of Subgraph Results}

For a connected planar graph $H$,
let $\mathbf{P}_{s} := \mathbf{P}[P_{n,m} \textrm{ will have a copy of } H]$. \\
\\
\begin{tabular}{c|c|c|c}
&
$e(H)<|H|$ &
$e(H)=|H|$ &
$e(H) > |H|$ \\
\hline
$0 < \underline{\lim}~\frac{m}{n}$ &
$\mathbf{P}_{s} \to 1$ (Cor.~\ref{tree6}) &
\small{$\underline{\lim}$} $\mathbf{P}_{s} > 0$ (T\ref{gen4}) &
$\mathbf{P}_{s} \to 0$ (C\ref{msub3}) \\
\& $\overline{\lim}~\frac{m}{n} < 1$ &
&
\small{$\overline{\lim}$} $\mathbf{P}_{s} < 1$ (T\ref{cyc54}) &
\\
\hline
$\frac{m}{n} \to 1$ &
$\mathbf{P}_{s} \to 1$ (Cor.~\ref{tree6}) &
$\mathbf{P}_{s} \to 1$ (Lem.~\ref{unisub1}) &
Unknown \\
&
&
&
\small{(see Section~\ref{multsub})} \\
\hline
$\underline{\lim}~\frac{m}{n} > 1$ &
$\mathbf{P}_{s} \to 1$ (Thm.~\ref{sub4}) &
$\mathbf{P}_{s} \to 1$ (Thm.~\ref{sub4}) &
$\mathbf{P}_{s} \to 1$ (T\ref{sub4})
\end{tabular}

\newpage
\section{Previous Results} \label{previous}

Recall that $P_{n}$ is the graph taken uniformly at random from the set of all graphs in $\mathcal{P}(n)$
and that $P_{n,m}$ is the graph taken uniformly at random from the set of all graphs in $\mathcal{P}(n,m)$.
In this section, 
a detailed account of the existing results on $P_{n}$ and $P_{n,m}$ will be given,
with sketch-proofs of some of the major theorems. \\

We shall begin by looking at $P_{n}$.
We will see results on the number of components in $P_{n}$ that are isomorphic to a given $H$ (Proposition~\ref{mcd 5.6}),
the number of special copies of $H$ (called `appearances') in $P_{n}$ (Proposition~\ref{gim T4}),
the probability that $P_{n}$ is connected (Propositions~\ref{mcd 2.1} \&~\ref{gim T6}),
and the total number of components in $P_{n}$ (also Propositions~\ref{mcd 2.1} \&~\ref{gim T6}).
In addition, 
we shall see (in Proposition~\ref{gim T1}) a precise asymptotic estimate for $|\mathcal{P}(n)|$,
the number of planar graphs of order $n$.

We will then summarise results about $P_{n,m}$.
Again, we shall see results concerning appearances (Proposition~\ref{ger T3.1}),
connectivity (Proposition~\ref{gim T3}),
and the total number of components (Propositions~\ref{ger L2.6} \&~\ref{ger3 6.6}),
as well as estimates for $|\mathcal{P}(n,m)|$ (Propositions~\ref{gim T3} \&~\ref{ger T2.1}). 
\begin{displaymath}
\end{displaymath} \\

We start with an estimate of how many planar graphs there are:

\begin{Proposition}[\cite{gim}, Theorem 1] \label{gim T1}
\begin{eqnarray*}
|\mathcal{P}(n)| \sim g \cdot n^{-7/2} \gamma_{l}^{n}n! ,
\end{eqnarray*}
where $g \approx 0.4260938569 \cdot 10^{-5}$
and $\gamma_{l} \approx 27.2268777685$ are constants given by explicit analytic expressions. 
\end{Proposition}
\textbf{Sketch of Proof}
Since a planar graph is a set of connected planar graphs,
the~e.g.f.~(exponential generating function)
for the class of planar graphs can be expressed in terms of the e.g.f.~for the class of connected planar graphs.
Similarly, a connected planar graph may be decomposed into $2$-connected components,
and so this enables us to relate the e.g.f.~for the class of connected planar graphs 
to the e.g.f.~for the class of $2$-connected planar graphs.
This latter class has already been analysed in \cite{ben},
via a further decomposition into $3$-connected components
(which have a unique embedding in the sphere~\cite{whi})
and the use of known results on planar map enumeration.
Hence, we are then able to use the e.g.f.~equations 
(together with a large amount of algebraic manipulation)
to obtain our asymptotic estimate for $|\mathcal{P}(n)|$.
\phantom{qwerty}
\setlength{\unitlength}{.25cm}
\begin{picture}(1,1)
\put(0,0){\line(1,0){1}}
\put(0,0){\line(0,1){1}}
\put(1,1){\line(-1,0){1}}
\put(1,1){\line(0,-1){1}}
\end{picture} \\
\\

Note that Proposition~\ref{gim T1} implies that 
$\left(\frac{|\mathcal{P}(n)|}{n!}\right)^{1/n} \to \gamma_{l} \textrm{ as } n \to \infty$.
Thus, we call $\gamma_{l}$ the \textit{labelled planar graph growth constant}. \\
\begin{displaymath}
\end{displaymath}

The precise nature of Proposition~\ref{gim T1} enables the structure of $P_{n}$ to be investigated in detail.
The main theorem is:

\begin{Proposition}[implicit in \cite{mcd}] \label{mcd 5.6}
Let $H_{1},\ldots ,H_{r}$ denote a fixed collection of pairwise non-isomorphic connected planar graphs,
and let $X_{n}^{(i)}$ denote a random variable which counts 
the number of components isomorphic to $H_{i}$ in $P_{n}$.
Then 
$\left( X_{n}^{(1)},\ldots,X_{n}^{(r)} \right) \stackrel{d}{\to} \left( Z_{1},\ldots,Z_{r} \right)$, 
where $Z_{i} \in \emph{Poi} \left( \left(|\emph{Aut}(H_{i})| \gamma _{l}^{|H_{i}|}\right)^{-1} \right) $ are independent. 
\end{Proposition}
\textbf{Sketch of Proof} 
We construct graphs of order $n$ with at least $k_{i}$ components isomorphic to $H_{i}$, for all $i$, 
and find that we have built 
$|\mathcal{P}(n)| \mathbf{E} 
\left[ \prod_{i \leq r} \left( (X_{n}^{(i)})_{k_{i}} \right) \right]$
graphs in total,
where $(X)_{k} = X(X-1)\cdots(X-k+1)$ denotes the $k$th factorial moment.
Thus, we obtain a formula for 
$\mathbf{E} \left[ \prod_{i \leq r} \left( (X_{n}^{(i)})_{k_{i}} \right) \right]$, 
which turns out to simplify in terms of $I_{n}$, the expected number of isolated vertices in $P_{n}$.
But 
$I_{n} = n\mathbf{P}[v_{n} \textrm{ is isolated in }P_{n}] = 
\frac{n|\mathcal{P}(n-1)|}{|\mathcal{P}(n)|} \to \gamma_{l}^{-1}$, 
by Theorem~\ref{gim T1}.
A standard result on the factorial moments of the Poisson distribution then completes the proof. 
$\phantom{qwerty}$ 
\begin{picture}(1,1)
\put(0,0){\line(1,0){1}}
\put(0,0){\line(0,1){1}}
\put(1,1){\line(-1,0){1}}
\put(1,1){\line(0,-1){1}}
\end{picture} \\
\\

As a corollary to Theorem~\ref{mcd 5.6}, 
we can see that the limiting probability for $P_{n}$ having a component isomorphic to $H$ is 
$1-e^{-\frac{1}{\left( \left| \scriptsize{\textrm{Aut}}(H) \right| \gamma_{l} ^{|H|}\right)}}$. \\
\begin{displaymath}
\end{displaymath}

It is easier to investigate the number of components isomorphic to $H$ than the number of subgraphs,
since components do not interfere with one another.
However, subgraphs may be approached via the concept of `appearances':

\begin{Definition} \label{defapps}
Let $H$ be a graph on the vertex set $\{1,2,\ldots,h\},$ and let $G$ be a graph on the vertex set $\{1,2,\ldots,n\}$, 
where $n>h$.
Let $W \subset V(G)$ with $|W|=h$, and let the `root' $r_{W}$ denote the least element in $W$.
We say that $H$ \textbf{\emph{appears}} at~$W$ in $G$ if 
(a) the increasing bijection from ${1,\ldots,h}$ to $W$ gives 
an isomorphism between $H$ and the induced subgraph $G[W]$ of $G$;
and (b) there is exactly one edge in $G$ between $W$ and the rest of $G$, 
and this edge is incident with the root $r_{W}$. 
We let $\mathbf{f_{H}(G)}$ denote the number of appearances of $H$ in $G$, 
that is the number of sets $W \subset V(G)$ such that $H$ appears at $W$ in $G$.
\end{Definition}

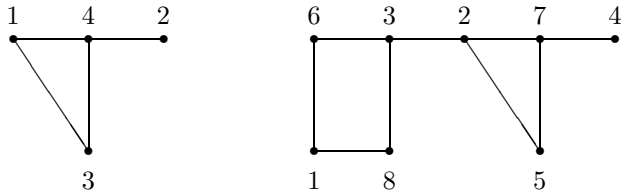
\begin{figure} [ht]
\setlength{\unitlength}{1cm}
\begin{picture}(10,2.2)(0,0)
\put(2,1.5){\line(1,0){2}}
\put(1.91,1.7){$1$}
\put(2,1.5){\circle*{0.1}}
\put(2.91,1.7){$4$}
\put(3,1.5){\circle*{0.1}}
\put(3.91,1.7){$2$}
\put(4,1.5){\circle*{0.1}}
\put(3,1.5){\line(0,-1){1.5}}
\put(2,1.5){\line(2,-3){1}}
\put(2.91,-0.5){$3$}
\put(3,0){\circle*{0.1}}

\put(6,1.5){\line(1,0){4}}
\put(6,0){\line(1,0){1}}
\put(6,1.5){\line(0,-1){1.5}}
\put(7,1.5){\line(0,-1){1.5}}
\put(5.91,1.7){$6$}
\put(6,1.5){\circle*{0.1}}
\put(6.91,1.7){$3$}
\put(7,1.5){\circle*{0.1}}
\put(5.91,-0.5){$1$}
\put(6,0){\circle*{0.1}}
\put(6.91,-0.5){$8$}
\put(7,0){\circle*{0.1}}

\put(7.91,1.7){$2$}
\put(8,1.5){\circle*{0.1}}
\put(8.91,1.7){$7$}
\put(9,1.5){\circle*{0.1}}
\put(9.91,1.7){$4$}
\put(10,1.5){\circle*{0.1}}
\put(9,1.5){\line(0,-1){1.5}}
\put(8,1.5){\line(2,-3){1}}
\put(8.91,-0.5){$5$}
\put(9,0){\circle*{0.1}}

\end{picture}
\caption{A graph $H$ and an appearance of $H$.} 
\end{figure} 

\phantom{p}

The clean structure of an appearance makes it a suitable candidate for the generating function approach.
Hence, we may obtain:

\begin{Proposition}[\cite{gim}, Theorem 4] \label{gim T4}
Let $H$ be a fixed connected planar graph on the vertex set $\{ 1,\ldots,h \}$.
Then $f_{H}(P_{n})$ is asymptotically Normal, and the mean~$\mu_{n}$ and variance $\sigma_{n}^{2}$ satisfy
$\mu_{n} \sim \frac{n}{\gamma_{l}^{h}h!} \textrm{ and } \sigma_{n}^{2} \sim \frac{n}{\gamma_{l}}$.  
\end{Proposition} 
\phantom{p}

Thus, a.a.s.~$P_{n}$ will contain at least linearly many copies of any given connected planar graph. \\
\begin{displaymath}
\end{displaymath}

One further topic that has been looked at is that of connectivity 
and the number of components of $P_{n}$.
Recall (from page~\pageref{planar}) that $\mathcal{P}(n)$ is edge-addable.
By a counting argument based on this,
it is shown in \cite{mcd} that we have:

\begin{Proposition}[\cite{mcd}, 2.1] \label{mcd 2.1}
The random number $\kappa (P_{n})$ of components of $P_{n}$ 
is stochastically dominated 
\footnote{We say the distribution $(p_{j})$ is \textit{stochastically dominated} by the distribution $(q_{j})$ 
to mean that we have 
$\Sigma_{k \leq l} \phantom{p} p_{k} \phantom{p} \geq 
\phantom{p} \Sigma_{k \leq l} \phantom{p} q_{k} \phantom{p} \forall l$.}
by $1+X$, 
where $X$ has the Poisson distribution with mean $1$.
In particular, 
\begin{displaymath}
\mathbf{P} [P_{n} \textrm{ connected} ] \geq 1/e \textrm{ and } \mathbf{E}[\kappa (P_{n})] \leq 2. 
\end{displaymath}
\end{Proposition} 
\phantom{p}

These bounds hold for all $n$,
but they may be significantly improved in the asymptotic case by using generating functions:

\begin{Proposition}[\cite{gim}, Theorem 6 \& Corollary 1] \label{gim T6}
Asymptotically, $\kappa (P_{n}) -~\!1$ is distributed like a Poisson law of parameter $\upsilon$,
where $\upsilon \approx 0.0374393660$ is a constant given by an explicit analytic expression.
In particular, 
\begin{eqnarray*}
\mathbf{P}[P_{n} \textrm{ connected}] & \to & e^{-\upsilon} \approx 0.9632528217 \\
\textrm{ and } \mathbf{E}[\kappa (P_{n})] & \to & 1+ \upsilon \approx 1.0374393660. 
\end{eqnarray*} 
\end{Proposition}

We now turn our attention to planar graphs with $n$ vertices and $m$ edges.
Again, precise estimates for $|\mathcal{P}(n,m)|$ and $|\mathcal{P}_{c}(n,m)|$,
the number of connected graphs in $\mathcal{P}(n,m)$,
are given in \cite{gim}:

\begin{Proposition}[\cite{gim}, implicit in Theorem 3] \label{gim T3}
Let $q \in (1,3)$.
Then
\begin{eqnarray*}
|\mathcal{P}(n,\lfloor qn \rfloor)| & \sim & 
g(q) \cdot \left( u(q) \right) ^{qn - \lfloor qn \rfloor} \cdot n^{-4} \gamma (q)^{n}n! \\ 
\textrm{and } \phantom{w} |\mathcal{P}_{c}(n,\lfloor qn \rfloor)| & \sim & 
g_{c}(q) \cdot \left( u(q) \right) ^{qn - \lfloor qn \rfloor} \cdot n^{-4} \gamma (q)^{n}n!,
\end{eqnarray*}
where $g(q)$, $g_{c}(q)$, $u(q)$ and $\gamma (q) >0$ are computable analytic functions.
\end{Proposition}
\textbf{Sketch of Proof}
We define the \textit{bivariate exponential generating function} for a class of graphs, $\mathcal{A}$,
to be 
$A(x,y) = \sum_{n,m \geq 0} \frac{a_{n,m}}{n!} x^{n}y^{m}$,
where $a_{n,m}$ denotes the number of graphs in $\mathcal{A}$ that have order $n$ and size $m$.
Hence, we may again use the decomposition of the proof of Proposition~\ref{gim T1} 
to obtain the given asymptotic estimates.
$\phantom{qwerty}$ 
\begin{picture}(1,1)
\put(0,0){\line(1,0){1}}
\put(0,0){\line(0,1){1}}
\put(1,1){\line(-1,0){1}}
\put(1,1){\line(0,-1){1}}
\end{picture} \\

Note this provides an expression for 
$\lim_{n \to \infty} \mathbf{P}\left[P_{n,\lfloor qn \rfloor} \textrm{ will be connected}\right]$.
As it turns out,
this limit is strictly between $0$ and $1$ for all $q \in (1,3)$. \\

Analogously to with $\gamma_{l}$,
we call $\gamma (q)$ the \textit{growth function for $q$}.
A fairly detailed picture of this function is given in \cite{ger} and \cite{gim}:

\begin{Proposition}[\cite{ger} and \cite{gim}] \label{ger T2.1}
For each $q \in [0,3]$,
there is a finite constant $\gamma(q) \geq 0$ such that, as $n \to \infty$,
if $0 \leq q < 3$ then both 
$(|\mathcal{P}_{c}(n, \lfloor qn \rfloor)|/n!)^{1/n}$ and $(|\mathcal{P} (n, \lfloor qn \rfloor)|/n!)^{1/n}$
tend to $\gamma(q)$,
and $(|\mathcal{P}(n,3n-6)|/n!)^{1/n}$ tends to $\gamma(3)$.

The function $\gamma(q)$ is equal to $0$ for $q \in [0,1)$
and is then unimodal and uniformly continuous on $[1,3]$,
satisfying
$\gamma(1) = e$, $\gamma(3) = 256/27$
and achieving a maximum value of $\gamma_{l}$
at $q = \mathbf{E} \left( e(P_{n}) \right) \approx 2.2132652385$. \\
\end{Proposition}

It is also shown in \cite{ger} that there is \textit{uniform} convergence to $\gamma(q)$,
in the following sense:

\begin{Proposition} [\cite{ger}, Lemma 2.9] \label{ger L2.9} 
Let $a \in (1,3)$ and let $\eta>0$.
Then there exists $n_{0}$ such that for all $n \geq n_{0}$ and all $s \in [an,3n-6]$ we have
\begin{displaymath}
\left| \left( \frac{|\mathcal{P}(n,s)|}{n!} \right)^{1/n} - \gamma \left( \frac{s}{n} \right) \right| < \eta.
\end{displaymath} \\
\end{Proposition}
\phantom{p}

No analogue to Proposition~\ref{mcd 5.6} has been given for $P_{n, \lfloor qn \rfloor}$.
However, via counting arguments we do have a (less precise) version of Proposition~\ref{gim T4}:

\begin{Proposition}[\cite{ger}, Theorem 3.1] \label{ger T3.1}
Let $q \in [1,3)$ and let $H$ be a fixed connected planar graph on the vertices $\{1, \ldots ,h\}$, 
with the extra condition that $H$ is a tree if $ q=1$.
Then there exists $\alpha= \alpha(H,q)>0$ such that 
\begin{displaymath}
\mathbf{P} \left[f_{H} \left( P_{n, \lfloor qn \rfloor} \right) \leq \alpha n \right] = e^{- \Omega (n)}. 
\end{displaymath} 
\end{Proposition} 
\phantom{p}

Thus, similarly to with $P_{n}$, we know that (a.a.s.)~$P_{n, \lfloor qn \rfloor}$ 
will contain at least linearly many copies of any given connected planar graph, if $q>1$.

\begin{displaymath}
\end{displaymath} \\

We finish with two upper bounds for the number of components in $P_{n,m}$.
The first deals with when $m = \lfloor qn \rfloor$ for $q \geq 1$,
and the second with when $m$ is slightly below $n$: \\

\begin{Proposition}[\cite{ger}, Lemma 2.6] \label{ger L2.6}
Let $q \in [1,3)$
and let $c > \ln \frac{\gamma_{l}}{\gamma(q)}$.
Then
\begin{displaymath}
\mathbf{P} \left[ \kappa \left( P_{n, \lfloor qn \rfloor} \right) > \lceil cn/ \ln n \rceil \right] = e^{-\Omega(n)}.
\end{displaymath} 
\end{Proposition}
\phantom{p}

\begin{Proposition}[\cite{ger3}, Lemma 6.6] \label{ger3 6.6}
Let $\beta >0$ be fixed,
and let $m=m(n)= n - (\beta +o(1))(n/ \ln n)$.
Let the constant $c>0$ satisfy $c > \beta + \ln \gamma_{l} - 1$.
Then
\begin{displaymath}
\mathbf{P} \left[ \kappa \left( P_{n,m} \right) > cn/ \ln n \right] = e^{-\Omega(n)}.
\end{displaymath} \\
\end{Proposition}

It is the aim of this project to expand on the current state of knowledge of~$P_{n,m}$.

\newpage
\section{Pendant Edges} \label{pen}

In this section, we will do some groundwork by investigating the number of pendant edges in random planar graphs,
which will be an important ingredient in later sections. 
The key result here will be Theorem~\ref{pen5}, where we shall see that (a.a.s.)
$P_{n,m}$ will have linearly many pendant edges, provided that $\frac{m}{n}$ is bounded away from both $0$ and $3$. 

Clearly, it suffices to show that (a.a.s.) $P_{n,m}$ will have linearly many vertices of degree $1$.
We will first see (in Corollary~\ref{pen3}) 
that the part of this result dealing with the case when $\frac{m}{n} \in [b,B]$,
where $b>1$ and $B<3$,
may be deduced fairly easily by analysing an important result from Section~\ref{previous}
concerning the concept of appearances.
We shall then (in Lemma~\ref{pen2}) also prove our result for smaller values of $\frac{m}{n}$, using counting arguments.

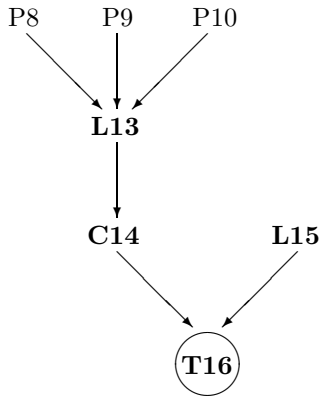
\begin{figure} [ht]
\setlength{\unitlength}{1cm}
\begin{picture}(10,4.5)(-0.6,0.6)
\put(3.35,4.9){P\ref{ger T2.1}}
\put(4.6,4.9){P\ref{ger L2.9}}
\put(5.8,4.9){P\ref{ger T3.1}}

\put(4.8,4.8){\vector(0,-1){1}}
\put(3.6,4.8){\vector(1,-1){1}}
\put(6,4.8){\vector(-1,-1){1}}

\put(4.45,3.45){\textbf{L\ref{gen2}}}
\put(4.8,3.35){\vector(0,-1){1}}
\put(4.4,2){\textbf{C\ref{pen3}}}
\put(4.8,1.9){\vector(1,-1){1}}
\put(6.85,2){\textbf{L\ref{pen2}}}
\put(7.2,1.9){\vector(-1,-1){1}}
\put(5.65,0.275){\textbf{T\ref{pen5}}}
\put(6,0.4){\circle{0.8}}
\end{picture}
\caption{The structure of Section~\ref{pen}.} 
\end{figure} 

\phantom{p}

We start by recalling our result on appearances in $P_{n, \lfloor qn \rfloor}$: \\
\\
\textbf{Proposition~\ref{ger T3.1} (\cite{ger}, Theorem 3.1)} 
\textit{Let $q \in [1,3)$ and let $H$ be a fixed connected planar graph on the vertices $\{1, \ldots ,h\}$, 
with the extra condition that $H$ is a tree if $ q=1$.
Then there exists $\alpha= \alpha(H,q)>0$ such that} 
\begin{displaymath}
\mathbf{P} \left[f_{H} \left( P_{n, \lfloor qn \rfloor} \right) \leq \alpha n \right] = e^{- \Omega (n)}. 
\end{displaymath}

By analysing the proof of Proposition~\ref{ger T3.1}, we may in fact obtain the following useful improved version:

\begin{Lemma} \label{gen2}
Let $H$ be a fixed connected planar graph on the vertices $\{1,2, \ldots ,h\}$,
let $b>1$ and $B<3$ be given constants,
and let $m=m(n) \in [bn,Bn]$ $\forall n$.
Then there exist constants $N(H,b,B)$ and $\alpha(H,b,B) >0$ such that
\begin{displaymath}
\mathbf{P}[f_{H}(P_{n,m}) \leq \alpha n] < e^{- \alpha n} \textrm{ } \forall n \geq N.
\end{displaymath}
\end{Lemma}
\textbf{Proof}
In the proof of Proposition~\ref{ger T3.1} in \cite{ger},
it is actually shown that $\exists N(H,q)$ and $\alpha (H,q) >0$ such that 
$\mathbf{P}[f_{H}(P_{n, \lfloor qn \rfloor}) \leq \alpha n] < e^{- \alpha n} \textrm{ } \forall n \geq N.$
It will suffice to show that this holds uniformly,
in the sense that $\exists N(H,b,B)$ and $\alpha(H,b,B)>0$
such that, \textit{for all $q \in [b,B]$},
we have
$\mathbf{P}[f_{H}(P_{n, \lfloor qn \rfloor}) \leq \alpha n] < e^{- \alpha n} \textrm{ } \forall n \geq N$.

The proof of Proposition~\ref{ger T3.1} given in \cite{ger} implicitly provides us with the value 
$\alpha (H,q) = \frac{1}{9e^{2}(\gamma (q))^{h+x}(h+x+2)(h+x)!}$,
where $x$ is the least integer such that 
$\frac{x+2}{x+1} \leq~q$ and $\frac{3x-5}{x+1} > q$.
Recall from Proposition~\ref{ger T2.1} 
that $\gamma (q) \in [e, \gamma_{l}]$~$\forall q \in [1,3]$.
Thus, it follows that we may take $\alpha$, \textit{independently of $q$}, 
to be $\frac{1}{9e^{2}\gamma_{l}^{h+z}(h+z+2)(h+z)!}$, 
where $z$ is the least integer such that $\frac{z+2}{z+1} \leq b$ and $\frac{3z-5}{z+1} > B$.

The value of $N(H,q)$ provided by the proof of Proposition~\ref{ger T3.1} given in \cite{ger} depends on~$q$ 
only in that we must have 
\begin{displaymath}
(1-\epsilon)^{n}n!(\gamma (q))^{n} \leq |\mathcal{P}(n, \lfloor qn \rfloor)| \leq 
(1+\epsilon)^{n}n!(\gamma (q))^{n} \textrm{ }\forall n \geq N,
\end{displaymath}
where $\epsilon$ is defined by the equation $\left( \frac{1}{9} \right)^{\alpha} = 1-3\epsilon$.
Thus, since our $\epsilon$ will only be a function of $H,b$ and $B$,
it suffices for us to prove $\exists N(b,B,\epsilon)$ such that, \emph{for all} $q \in [b,B]$, we have 
\begin{displaymath}
(1-\epsilon)^{n}n!(\gamma (q))^{n} \leq |\mathcal{P}(n, \lfloor qn \rfloor)| \leq 
(1+ \epsilon)^{n}n!(\gamma (q))^{n} \textrm{ }\forall n \geq N.
\end{displaymath}

By uniform convergence (Proposition~\ref{ger L2.9}),
we know that $\exists N_{1} (b, \epsilon)$ such that, for all $q \in [b,3)$, 
we have
$\left|\left(\frac{|\mathcal{P}(n, \lfloor qn \rfloor)|}{n!} \right)^{1/n} - 
\gamma \left(\frac{\lfloor qn \rfloor}{n} \right)\right| < \frac{\epsilon e}{2} $ 
$\forall n \geq N_{1}$.
But since $\gamma (q)$ is uniformly continuous on $[1,3]$
(Proposition~\ref{ger T2.1}), 
we also know 
$\exists N_{2} (\epsilon)$ such that, for all $q \in [1,3]$, we have
$\left|\gamma \left(\frac{\lfloor qn \rfloor}{n} \right)- \gamma (q)\right| < 
\frac{\epsilon e}{2} $ $\forall n \geq N_{2}$.
Thus, for all $q \in [b,3)$, we have 
\begin{displaymath}
\left|\left(\frac{|\mathcal{P}(n, \lfloor qn \rfloor)|}{n!} \right)^{1/n} - \gamma (q)\right| < 
\epsilon e \leq \epsilon \gamma (q) \textrm{ } \forall n \geq N(b,B, \epsilon) = \max \{ N_{1}, N_{2} \}. 
\end{displaymath}
Hence, we are done.
$\phantom{qwerty}$ 
\setlength{\unitlength}{.25cm}
\begin{picture}(1,1)
\put(0,0){\line(1,0){1}}
\put(0,0){\line(0,1){1}}
\put(1,1){\line(-1,0){1}}
\put(1,1){\line(0,-1){1}}
\end{picture} \\
\\

The first part of our pendant edges/vertices of degree $1$ result now follows:

\begin{Corollary} \label{pen3}
Let $b > 1$ and $B < 3$ be given constants
and let $m = m(n) \in [bn,Bn]$~$\forall n$.
Then there exist constants $N(b,B)$ and $\beta(b,B) >0$ such that
\begin{displaymath}
\mathbf{P}[P_{n,m} \textrm{ will have $< \beta n$ vertices of degree }1] < 
e^{- \beta n} \textrm{ } \forall n \geq N.
\end{displaymath}
\end{Corollary}
\textbf{Proof}
This follows from Lemma~\ref{gen2}, with $H$ as an isolated vertex. 
$\phantom{qwerty}$ 
\setlength{\unitlength}{.25cm}
\begin{picture}(1,1)
\put(0,0){\line(1,0){1}}
\put(0,0){\line(0,1){1}}
\put(1,1){\line(-1,0){1}}
\put(1,1){\line(0,-1){1}}
\end{picture} \\
\\

We are left with proving the result for the case when $\frac{m}{n}$ is small:

\begin{Lemma} \label{pen2}
Let $c>0$ and $\delta<1/8$ be given constants and let $m=m(n) \in [cn, (1+\delta)n]$ for all large $n$. 
Then there exists a constant $\beta(c,\delta) >0$ such that 
\begin{displaymath}
\mathbf{P}[P_{n,m} \textrm{ will have } < \beta n \textrm{ vertices of degree } 1] < e^{-\beta n} 
\textrm{ for all large } n.
\end{displaymath} 
\end{Lemma} 
\textbf{Sketch of Proof}
We suppose, hoping for a contradiction, that there will 
be a decent proportion of graphs in $\mathcal{P}(n,m)$ with
only `a few' vertices of degree~$1$ and
we consider separately the cases when 
(a)~there are also `many' isolated vertices and
(b) there are not `many' isolated vertices.

For case (a), we construct new graphs in $\mathcal{P}(n,m)$ by turning
some of the isolated vertices into vertices of degree~$1$ 
(and deleting some edges elsewhere)
and find that we can construct so many graphs that we obtain our desired contradiction.

For case (b), we note that we must have `lots' of vertices of degree $2$.
In fact, we must have `lots' of vertices of degree $2$ that are adjacent only to other vertices of degree $2$.
Such a vertex must belong either to a component that is a triangle or to a larger component.
In both cases, we construct new graphs by turning the chosen vertex into a vertex of degree $1$
(and inserting an edge elsewhere).
Again, we find that we can construct so many graphs that we obtain our desired contradiction.

In all cases of the proof, 
the key idea is to construct our graphs in such a way that there is not much double-counting.
This is helped by our assumption that we started with only `a few' vertices of degree $1$. \\
\\
\textbf{Full Proof}
Choose $\beta >0$ and suppose $\exists$ arbitrarily large $n$ such that 
\begin{displaymath}
\mathbf{P}[P_{n,m} \textrm{ will have $< \beta n$ vertices of degree $1$}] \geq e^{-\beta n}.
\end{displaymath}
Consider one of these $n$ and
let $\mathcal{G}_{n}$ denote the set of graphs in $\mathcal{P}(n,m)$ with 
$< \beta n$ vertices of degree $1$. 
Thus, $|\mathcal{G}_{n}| \geq e^{-\beta n} |\mathcal{P}(n,m)|$.

Choose $\epsilon >0$, let $\mathcal{H}_{n}$ denote the set of graphs in $\mathcal{P}(n,m)$ with 
$< \beta n$ vertices of degree~$1$ 
and with $> \epsilon n$ vertices of degree $0$,
and let
$\mathcal{J}_{n}$ denote the set of graphs in $\mathcal{P}(n,m)$ with 
$< \beta n$ vertices of degree $1$ and with $\leq \epsilon n$ vertices of degree $0$.
Then $|\mathcal{H}_{n}| + |\mathcal{J}_{n}| = |\mathcal{G}_{n}| \geq e^{-\beta n} |\mathcal{P}(n,m)|$.
Thus, either $|\mathcal{H}_{n}| \geq \frac{e^{-\beta n}}{2} |\mathcal{P}(n,m)|$ or 
$|\mathcal{J}_{n}| \geq \frac{e^{-\beta n}}{2} |\mathcal{P}(n,m)|$. 

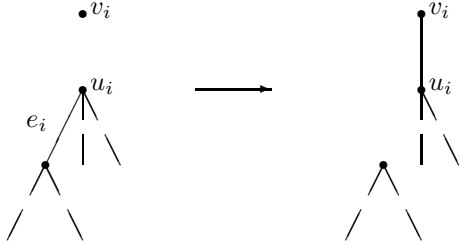
\begin{figure} [ht]
\setlength{\unitlength}{1cm}
\begin{picture}(10,2.5)(0,0.5)
\put(4,3){\circle*{0.1}}
\put(4.1,3){$v_{i}$}
\put(4,2){\circle*{0.1}}
\put(3.5,1){\circle*{0.1}}
\put(3.25,1.5){$e_{i}$}
\put(4,2){\line(-1,-2){0.5}}
\put(4,2){\line(0,-1){0.4}}
\put(4,1.4){\line(0,-1){0.4}}
\put(4,2){\line(1,-2){0.2}}
\put(4.3,1.4){\line(1,-2){0.2}}
\put(3.5,1){\line(1,-2){0.2}}
\put(3.8,0.4){\line(1,-2){0.2}}
\put(3.5,1){\line(-1,-2){0.2}}
\put(3.2,0.4){\line(-1,-2){0.2}}
\put(4.1,2){$u_{i}$}

\put(5.5,2){\vector(1,0){1}}

\put(8.5,3){\circle*{0.1}}
\put(8.6,3){$v_{i}$}
\put(8.5,3){\line(0,-1){1}}
\put(8.5,2){\circle*{0.1}}
\put(8,1){\circle*{0.1}}
\put(8.5,2){\line(0,-1){0.4}}
\put(8.5,1.4){\line(0,-1){0.4}}
\put(8.5,2){\line(1,-2){0.2}}
\put(8.8,1.4){\line(1,-2){0.2}}
\put(8,1){\line(1,-2){0.2}}
\put(8.3,0.4){\line(1,-2){0.2}}
\put(8,1){\line(-1,-2){0.2}}
\put(7.7,0.4){\line(-1,-2){0.2}}
\put(8.6,2){$u_{i}$}
\end{picture}

\caption{Constructing a vertex of degree $1$ in Stage $2$.} \label{penfig}
\end{figure}

In a moment,
we shall split our proof into two cases based on the observation of the previous sentence.
First, though,
we should just note that there will be several places later where we shall require
that our choices of $\beta$, $\epsilon$ and~$n$ satisfy various inequalities.
Formally, this can be done by assuming that we first chose~$\epsilon$ to be `sufficiently' small
(depending on $c$ and $\delta$),
that we then chose $\beta$ to be `sufficiently' small
(depending on $c$, $\delta$ and $\epsilon$),
and that we finally chose $n$ to be `sufficiently' large
(depending on $c$, $\delta$, $\epsilon$ and $\beta$). \\
\\
\\
\underline{Case(a)} \\
Suppose $|\mathcal{H}_{n}| \geq \frac{e^{-\beta n}}{2} |\mathcal{P}(n,m)|$ and consider a graph $G \in \mathcal{H}_{n}$.
Using $G$, we shall construct a graph (in $\mathcal{P}(n,m)$)
with $\geq \beta n$ vertices of degree $1$ as follows: \\
\\
Stage $1$: \\
Choose $\lceil \beta n \rceil$ isolated vertices 
(we have $\left(^{\phantom{p}d_{0}}_{\lceil \beta n \rceil}\right) > 
\left(^{\phantom{p}\epsilon n}_{\lceil \beta n \rceil}\right)$ 
choices 
\footnote{We use the definition 
$\left(^{\phantom{p}\epsilon n}_{\lceil \beta n \rceil}\right) :=
\frac{\epsilon n \cdot (\epsilon n -1) \cdots (\epsilon n - \lceil \beta n \rceil + 1)}
{\lceil \beta n \rceil !}$
if $\epsilon n$ is non-integral,
as is standard.}
for these,
where $d_{0}$ denotes the number of vertices of degree $0$).
Let us denote the chosen vertices, in order of their labels, as $v_{1}, v_{2}, \ldots, v_{\lceil \beta n \rceil}$,
and let $i=1$. \\
\\
Stage $2$: \\
Choose a vertex $u_{i}$ that was non-isolated in $G$ and that was also not 
incident to any of the edges $e_{1},e_{2}, \ldots, e_{i-1}$ defined in previous iterations
(we have at least $n-d_{0}-2(i-1)$ choices,
and we may assume that $\beta>0$ is sufficiently small that this is positive $\forall i \leq \lceil \beta n \rceil$,
since by planarity $n-d_{0} > \frac{m}{3} \geq \frac{cn}{3}$);
delete an edge $e_{i}$ incident to $u_{i}$; and join $u_{i}$ to $v_{i}$ (see Figure~\ref{penfig}). 
\\
Stage $3$: \\
If $i = \lceil \beta n \rceil$, then we terminate the algorithm.
Otherwise, we increase $i$ by $1$ and return to Stage $2$. \\
\\

By considering all possible initial graphs $G \in \mathcal{H}_{n}$,
it is clear to see that 
the number of ways to
build a graph with $\geq \beta n$ vertices of degree $1$
is at least
$\textrm{\small{$
\left(^{\phantom{p}\epsilon n}_{\lceil \beta n \rceil}\right) 
\left( \prod_{i=1}^{\lceil \beta n \rceil -1} (n-d_{0}-2(i-1)) \right) |\mathcal{H}_{n}|
\geq \left(^{\phantom{p}\epsilon n}_{\lceil \beta n \rceil}\right) 
\left( \frac{m}{3}-2\lceil \beta n \rceil \right)^{\lceil \beta n \rceil}  \frac{e^{-\beta n}}{2}|\mathcal{P}(n,m)|
.$}}$

We will now consider the amount of double-counting,
i.e.~how many times each new graph will be constructed.
We know that each new graph will contain at most $3 \lceil \beta n \rceil$ vertices of degree $1$
(since there were $< \beta n$ to begin with;
we have deliberately added $\lceil \beta n \rceil$;
and we may have created at most one extra one each time we deleted an edge),
so we will have at most $\left(^{3 \lceil \beta n \rceil} _{\phantom{l} \lceil \beta n \rceil} \right)$
possibilities for which were our chosen (originally isolated) vertices
$v_{1},v_{2}, \ldots, v_{\lceil \beta n \rceil}$.
We will then have $<(n- \lceil \beta n \rceil)^{\lceil \beta n \rceil}$ possibilities for where the edges 
attached to these vertices were originally.
Thus, we will build each graph at most $\left(^{3 \lceil \beta n \rceil} _{\phantom{l} \lceil \beta n \rceil} \right) 
(n-\lceil \beta n \rceil)^{\lceil \beta n \rceil}$ times.

Therefore, the number of distinct graphs in $\mathcal{P}(n,m)$ with $\geq \beta n$ vertices of degree~$1$ is at least 
\begin{eqnarray*}
& & \frac{
\left(^{\phantom{p} \epsilon n}_{\lceil \beta n \rceil}\right) 
\left( \frac{m}{3}-2\lceil \beta n \rceil \right)^{\lceil \beta n \rceil} 
\frac{e^{-\beta n}}{2} |\mathcal{P}(n,m)|}
{\left(^{3 \lceil \beta n \rceil} _{\phantom{l} \lceil \beta n \rceil} \right)  
(n- \lceil \beta n \rceil)^{\lceil \beta n \rceil}
} \\
& \geq & \left( \frac{\epsilon n}{3 \lceil \beta n \rceil} \right) ^{\lceil \beta n \rceil}
\left( \frac{
\frac{m}{3}-2\lceil \beta n \rceil}
{n- \lceil \beta n \rceil} \right)^{\lceil \beta n \rceil} 
\frac{e^{-\beta n}}{2} |\mathcal{P}(n,m)| \\
& &
\textrm{ since } \frac{ \left(^{x}_{y}\right) }{ \left(^{z}_{y}\right) } 
= \frac{x(x-1) \cdots (x-y+1)}{z(z-1) \cdots (z-y+1)}
\geq \left( \frac{x}{z} \right)^{y}
\textrm{ if } z \leq x \\
& & \textrm{ and here we may assume }
3 \lceil \beta n \rceil \leq \epsilon n \\
& > & \left( \left( \frac{\epsilon}{4 \beta} \right) \left( \frac{\frac{c}{3}-3 \beta}{1-\beta} \right) \right)
^{\lceil \beta n \rceil} \frac{e^{- \beta n}}{2} |\mathcal{P}(n,m)| 
\phantom{www} \textrm{ for sufficiently large $n$} \\
& > & 3^{\lceil \beta n \rceil} \frac{e^{- \beta n}}{2} |\mathcal{P}(n,m)|
\phantom{www} \textrm{ for sufficiently small $\beta$}, \\
& & \textrm{ since } 
\left( \frac{\epsilon}{4\beta} \right) \left( \frac{\frac{c}{3}-3\beta}{1-\beta} \right) \to \infty
\textrm{ as } \beta \to 0 \\
& > & |\mathcal{P}(n,m)|
\phantom{www} \textrm{ for sufficiently large $n$ (since $3>e$),} 
\end{eqnarray*}
which is a contradiction. \\
\\
\\
\underline{Case(b)} \\
Suppose instead that $|\mathcal{J}_{n}| \geq \frac{e^{-\beta n}}{2} |\mathcal{P}(n,m)|$ 
and consider a graph $G \in \mathcal{J}_{n}$.
We shall start by showing the intuitive fact that $G$ must contain many vertices of degree $2$. 

Let $d_{i}$ denote the number of vertices of degree $i$ in $G$.
Then 
\begin{eqnarray}
d_{1} + 2d_{2} + 3 \sum_{i \geq 3}d_{i} 
& \leq & \sum_{i \geq 1}id_{i} \nonumber \\
& = & 2m \nonumber \\
& \leq & 2n +2\delta n \nonumber \\
& = & 2d_{0} + 2d_{1} + 2d_{2} + 2 \sum_{i \geq 3}d_{i} + 2\delta n. \label{eq:label1}
\end{eqnarray}
Thus, 
\begin{eqnarray}
\sum_{i \geq 3}d_{i} 
& \leq & 2d_{0} + d_{1} + 2\delta n \nonumber \\
& < & 2 \epsilon n + \beta n + 2\delta n \label{eq:label2}
\end{eqnarray}
and so
\begin{eqnarray*}
d_{2} 
& = & n - \sum_{i \geq 3}d_{i} - d_{1} - d_{0} \\
& > & n - 2 \epsilon n - \beta n - 2 \delta n - \beta n - \epsilon n \\
& = & \left( 1-3\epsilon - 2\beta -2\delta \right)n. 
\end{eqnarray*}

We shall now see that, in fact,
$G$ must contain many vertices of degree~$2$ that are adjacent only to other vertices of degree $2$. 
Recall from (\ref{eq:label1}) that we have
$\sum_{i \geq 1}id_{i} \leq 2d_{0} + 2d_{1} + 2d_{2} + 2 \sum_{i \geq 3}d_{i} +~2\delta n$.
Thus, 
\begin{eqnarray*}
\sum_{i \geq 3}id_{i} 
& \leq & 2d_{0} + d_{1} + 2\sum_{i \geq 3}d_{i} + 2\delta n \\
& < & 2 \epsilon n + \beta n + 2(2 \epsilon n + \beta n + 2\delta n) + 2 \delta n \phantom{www}
\textrm{ by (\ref{eq:label2})} \\
& = & (6 \epsilon + 3 \beta + 6 \delta)n.
\end{eqnarray*}
Therefore, at most $(6 \epsilon + 3 \beta + 6 \delta)n$ of the degree $2$ vertices 
will be adjacent to a vertex of degree $\geq 3$.
Similarly, at most $d_{1} < \beta n$ of the degree $2$ vertices will be adjacent to a vertex of degree $1$.
Hence, at least 
$d_{2} - \beta n - (6 \epsilon +3 \beta + 6 \delta)n
> (1- 3 \epsilon - 2 \beta - 2 \delta)n - \beta n - (6 \epsilon + 3 \beta + 6 \delta)n
= (1-9 \epsilon - 6 \beta - 8\delta)n$
of the degree~$2$ vertices will be adjacent only to other degree $2$ vertices. 

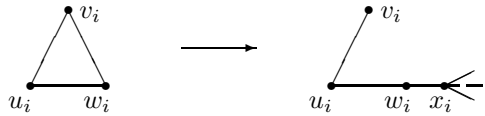
\begin{figure} [ht]
\setlength{\unitlength}{1cm}
\begin{picture}(20,0.5)(0,0.25)
\put(3.5,1){\line(-1,-2){0.5}}
\put(3.5,1){\line(1,-2){0.5}}
\put(3.5,1){\circle*{0.1}}
\put(3.65,0.9){$v_{i}$}
\put(3,0){\line(1,0){1}}
\put(3,0){\circle*{0.1}}
\put(4,0){\circle*{0.1}}
\put(2.7,-0.3){$u_{i}$}
\put(3.7,-0.3){$w_{i}$}

\put(5,0.5){\vector(1,0){1}}

\put(7.5,1){\line(-1,-2){0.5}}
\put(7.5,1){\circle*{0.1}}
\put(7.65,0.9){$v_{i}$}
\put(7,0){\line(1,0){1.5}}
\put(7,0){\circle*{0.1}}
\put(8,0){\circle*{0.1}}
\put(8.5,0){\circle*{0.1}}
\put(6.7,-0.3){$u_{i}$}
\put(7.7,-0.3){$w_{i}$}
\put(8.3,-0.3){$x_{i}$}

\put(8.5,0){\line(1,0){0.2}}
\put(8.8,0){\line(1,0){0.2}}
\put(8.5,0){\line(2,1){0.4}}
\put(8.5,0){\line(2,-1){0.4}}

\end{picture}
\caption{Constructing a vertex of degree $1$ in Stage $3a$.} \label{fig3a}
\end{figure}

Let $A$ denote the set of vertices of degree $2$ that are adjacent only to other degree $2$ vertices.
Using $G$, we shall construct a graph (in $\mathcal{P}(n,m)$)
with $ \geq \beta n$ vertices of degree $1$ by the following algorithm: \\
\\
Stage $1$: \\
Let $B_{0} = \emptyset$ and let $i=1$. \\
\\
Stage $2$: \\
Choose a vertex, $v_{i}$, in $A-B_{i-1}$ 
(it will become clear that $v_{i}$ will still be a degree $2$ vertex that is adjacent only to other degree $2$ vertices).
Let the vertices adjacent to $v_{i}$ be denoted by $u_{i}$ and $w_{i}$. \\
\\
Stage $3a$ (If $u_{i}w_{i}$ $\in E(G)$): \\
If $u_{i}w_{i}$ is an edge in $G$ (in which case $\{ u_{i}, v_{i}, w_{i} \}$ form a component that is a triangle, 
since $d(u_{i}) = d(v_{i}) = d(w_{i}) =2$), 
then delete the edge $v_{i}w_{i}$ and join $w_{i}$ to a vertex $x_{i} \notin (B_{i-1} \cup \{ u_{i}, v_{i} \} )$
(see Figure~\ref{fig3a}).
This is possible if $i \leq \lceil \beta n \rceil$,
since it will become clear that $|B_{i-1}| + 2 \leq 6(i-1) +2 <6( \lceil \beta n \rceil -1) +2 < n$
if $\beta < \frac{1}{6}$.
Planarity will be maintained, since $w_{i}$ and $x_{i}$ were in separate components. 

If $d(x_{i}) \neq 2$, let $B_{i} = B_{i-1} \cup \{ u_{i}, v_{i}, w_{i} \}$.
If $d(x_{i}) = 2$, let the vertices adjacent to $x_{i}$ be denoted by $y_{i}$ and $z_{i}$ and 
let $B_{i} = B_{i-1} \cup \{ u_{i}, v_{i}, w_{i}, x_{i}, y_{i}, z_{i} \}$. 
Thus, as already mentioned,
it is clear that $|B_{i}|$ increases by at most $6$ in each iteration
and that $A-B_{i}$ still just contains degree $2$ vertices that are adjacent only to other degree $2$ vertices.
\\ Stage $3b$ (If $u_{i}w_{i}$ $\notin E(G)$): \\
If $u_{i}w_{i}$ is not an edge in $G$, then insert an edge between $u_{i}$ and $w_{i}$
(this can be done arbitrarily close to the edges $u_{i}v_{i}$ and $v_{i}w_{i}$, so planarity is maintained).
Delete the edge $v_{i}w_{i}$.
Let $x_{i}$ denote the neighbour of $u_{i}$ that is not $v_{i}$ or $w_{i}$ and 
let $B_{i} = B_{i-1} \cup \{ u_{i}, v_{i}, w_{i}, x_{i} \}$. 
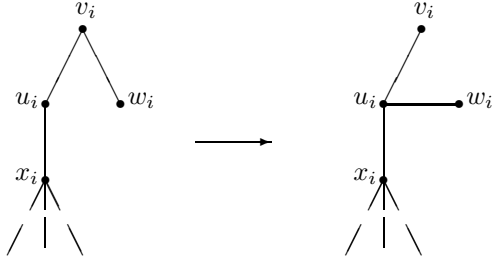
\begin{figure} [ht]
\setlength{\unitlength}{1cm}
\begin{picture}(20,2.75)(0,0.5)

\put(4,3){\line(-1,-2){0.5}}
\put(4,3){\line(1,-2){0.5}}
\put(3.5,2){\line(0,-1){1}}
\put(4,3){\circle*{0.1}}
\put(3.9,3.2){$v_{i}$}
\put(3.5,2){\circle*{0.1}}
\put(3.1,2){$u_{i}$}
\put(4.5,2){\circle*{0.1}}
\put(4.6,2){$w_{i}$}
\put(3.5,1){\circle*{0.1}}
\put(3.1,1){$x_{i}$}
\put(3.5,1){\line(0,-1){0.4}}
\put(3.5,0.5){\line(0,-1){0.4}}
\put(3.5,1){\line(1,-2){0.2}}
\put(3.8,0.4){\line(1,-2){0.2}}
\put(3.5,1){\line(-1,-2){0.2}}
\put(3.2,0.4){\line(-1,-2){0.2}}

\put(5.5,1.5){\vector(1,0){1}}

\put(8.5,3){\line(-1,-2){0.5}}
\put(8,2){\line(0,-1){1}}
\put(8,2){\line(1,0){1}}
\put(8.5,3){\circle*{0.1}}
\put(8.4,3.2){$v_{i}$}
\put(8,2){\circle*{0.1}}
\put(7.6,2){$u_{i}$}
\put(9,2){\circle*{0.1}}
\put(9.1,2){$w_{i}$}
\put(8,1){\circle*{0.1}}
\put(7.6,1){$x_{i}$}
\put(8,1){\line(0,-1){0.4}}
\put(8,0.5){\line(0,-1){0.4}}
\put(8,1){\line(1,-2){0.2}}
\put(8.3,0.4){\line(1,-2){0.2}}
\put(8,1){\line(-1,-2){0.2}}
\put(7.7,0.4){\line(-1,-2){0.2}} 

\end{picture}
\caption{Constructing a vertex of degree $1$ in Stage $3b$.} 
\end{figure} 
\\ Stage $4$: \\
If $i = \lceil \beta n \rceil$, then we terminate the algorithm.
Otherwise, we increase $i$ by $1$ and return to Stage $2$. \\

Note that we have $d(v_{i}) = 1$ $\forall i$ in our new graph.
At each iteration,
we had $|A-B_{i-1}| \geq (1 - 9 \epsilon - 6 \beta - 8 \delta)n - 6 \beta n
=(1- 9 \epsilon - 12 \beta - 8 \delta)n$ 
choices for $v_{i}$.
Thus, by considering all possible initial graphs $G \in \mathcal{J}_{n}$,
we have at least 
$(1- 9 \epsilon - 12 \beta - 8 \delta)^{\lceil \beta n \rceil} n^{\lceil \beta n \rceil} 
\frac{e^{-\beta n}}{2} |\mathcal{P}(n,m)|$ ways to construct a graph with $\geq \beta n$ vertices of degree $1$,
and so it remains only to consider the amount of double-counting.

Each new graph will contain at most $3 \lceil \beta n \rceil$ vertices of degree $1$
(since there were $< \lceil \beta n \rceil$ to begin with; we have deliberately created $\lceil \beta n \rceil$; 
and we may have created at most one extra one each time we used Stage $3a$, if $x_{i}$ was an isolated vertex).
Thus, we have $\leq 3 \lceil \beta n \rceil$ possibilities for which vertex was $v_{\lceil \beta n \rceil}$, 
which will now be adjacent only to $u_{\lceil \beta n \rceil}$.

If $u_{\lceil \beta n \rceil}$ now has degree $2$, then we must have used Stage $3a$ in the final iteration and 
hence $w_{\lceil \beta n \rceil}$ is the other neighbour of $u_{\lceil \beta n \rceil}$ 
and $x_{\lceil \beta n \rceil}$ is the remaining neighbour of $w_{\lceil \beta n \rceil}$.
Thus, we know how the graph changed during the final iteration.

If $u_{\lceil \beta n \rceil}$ now has degree $3$, 
then we must have used Stage $3b$ in the final iteration and 
hence we have two possibilities for which vertex was $w_{\lceil \beta n \rceil}$ 
(since it must be one of the other neighbours of $u_{\lceil \beta n \rceil}$ in the new graph).
Thus, we have two possibilities for how the graph changed during the final iteration.

Therefore, 
we have at most $6 \lceil \beta n \rceil$ possibilities in total for which vertex was $v_{\lceil \beta n \rceil}$
and how the graph looked before the final iteration.
Repeating this argument,
we find that we have at most $\left( 6 \lceil \beta n \rceil \right)^{\lceil \beta n \rceil}$ possibilities
for what the original graph was and which vertices were $v_{1},v_{2}, \ldots, v_{\lceil \beta n \rceil}$
(in order).
Hence, we have built each of our new graphs at most 
$(6 \lceil \beta n \rceil)^{\lceil \beta n \rceil}$ times and,
therefore, the number of graphs in $\mathcal{P}(n,m)$ 
with $\geq \beta n$ vertices of degree $1$ is at least 
$ \left ( \frac{(1-9 \epsilon - 12 \beta - 8 \delta)n}{6 \lceil \beta n \rceil} \right) ^{\lceil \beta n \rceil} 
\frac{e^{-\beta n}}{2} |\mathcal{P}(n,m)|$.

Recall that $\delta < 1/8$.
Thus, since we were free to choose $\epsilon$ and $\beta$ arbitrarily small,
we may assume that
$ \left ( \frac{1-9 \epsilon - 12 \beta - 8 \delta}{6 \beta} \right) >3$.
Therefore, for sufficiently large~$n$, we have
$ \left ( \frac{(1-9 \epsilon - 12 \beta - 8 \delta)n}{6 \lceil \beta n \rceil} \right) ^{\lceil \beta n \rceil} 
\frac{e^{-\beta n}}{2} |\mathcal{P}(n,m)| 
> 3^{\lceil \beta n \rceil} \frac{e^{-\beta n}}{2} |\mathcal{P}(n,m)| > |\mathcal{P}(n,m)|$,
which provides us with a contradiction. \\

Thus, we get a contradiction whether $\frac{|\mathcal{H}_{n}|}{|\mathcal{P}(n,m)|} \geq \frac{e^{-\beta n}}{2}$ or 
$\frac{|\mathcal{J}_{n}|}{|\mathcal{P}(n,m)|} \geq \frac{e^{-\beta n}}{2}$.~
\setlength{\unitlength}{.25cm}
\begin{picture}(1,1)
\put(0,0){\line(1,0){1}}
\put(0,0){\line(0,1){1}}
\put(1,1){\line(-1,0){1}}
\put(1,1){\line(0,-1){1}}
\end{picture} \\
\\

By combining Corollary~\ref{pen3} and Lemma~\ref{pen2}, we obtain our main result of this section:

\begin{Theorem} \label{pen5}
Let $b>0$ and $B<3$ be given constants and 
let $m(n) \in [bn,Bn]$ for all large $n$.
Then there exists a constant $\alpha(b,B) > 0$ such that
\begin{displaymath}
\mathbf{P}[P_{n,m} \textrm{ will have } < \alpha n \textrm{ pendant edges }] < e^{- \alpha n} 
\textrm{ for all large } n.
\end{displaymath} 
\end{Theorem}
\textbf{Proof}
The number of vertices of degree $1$ is at most twice the number of edges incident to a vertex of degree $1$.
$\phantom{qwerty}$ 
\setlength{\unitlength}{.25cm}
\begin{picture}(1,1)
\put(0,0){\line(1,0){1}}
\put(0,0){\line(0,1){1}}
\put(1,1){\line(-1,0){1}}
\put(1,1){\line(0,-1){1}}
\end{picture}

\newpage
\section{Addable Edges} \label{add}

In this section, 
we will continue to lay the groundwork for later counting arguments.
In future sections,
we shall often wish to choose an edge to insert into a graph without violating planarity,
and so our focus here will be to examine how many choices we have. 

\begin{Definition}
Given a planar graph $G$, we call a non-edge $e$ \emph{\textbf{addable}} in $G$ 
if the graph $G+e$ obtained by adding $e$ as an edge is still planar.
We let \emph{\textbf{add(}$\mathbf{G}$\textbf{)}} denote the set of addable non-edges of $G$ 
(note that the graph obtained by adding all the edges in \emph{add}$(G)$ may well not be planar)
and we let \emph{\textbf{add(}$\mathbf{n,m}$\textbf{)}} 
denote the minimum value of $|\emph{add}(G)|$ over all graphs $G \in \mathcal{P}(n,m)$.
\end{Definition}

In later sections,
we will require bounds (both lower and upper) on add$(n,m)$ for the case when $m \leq (1+o(1))n$,
and so that is our main purpose here.~
However,
as an interesting aside,
we shall also provide results for larger values of $m$
(recall that we use `Proposition' for such asides).
Hence,
we will in fact build up a fairly complete description,
showing that there are four main results: 
\begin{table} [ht] 
\begin{displaymath}
\textrm{add}(n,m(n)) =
\left\{ \begin{array}{ll}
\Theta(dn) & \textrm{if $d=n-m>0$ is such that $d=\Omega(n^{1/2})$} \\
& \textrm{(Theorems~\ref{add53} \&~\ref{add55})}\\
\Theta(n^{3/2}) & \textrm{if $|m-n|=O(n^{1/2})$ (Theorems~\ref{add3} \&~\ref{add21})}\\
\Theta \left( \frac{n^{2}}{d} \right) 
& \textrm{if $d=m-n>0$ is such that} \\
& \textrm{$d=\Omega(n^{1/2})$ and 
\small{$\limsup$} } \frac{d}{n} < 2
\textrm{ (Thms.~\ref{add54} \&~\ref{add56})} \\
\Theta(3n\!-\!m) & \textrm{if $m = \Omega(n)$ (Props.~\ref{add712} \&~\ref{add713} and Thm.~\ref{add56}})
\end{array} \right. 
\end{displaymath} 
\caption{The four main addability results.} \label{addsum1}
\end{table}

The bounds of Table~\ref{addsum1} will be sufficient for the rest of this thesis.
However,
as the topic of addable edges is quite interesting,
we shall flesh out our account by also giving more detailed results for seven special cases: 
\begin{table} [ht] 
\begin{displaymath}
\textrm{add}(n,m(n)) = 
\left\{ \begin{array}{ll}
(1+o(1)) \frac{(1-A)(1+A)}{2}n^{2} & \textrm{if $m = An + o(n)$, for $A<1$} \\
& \textrm{(Propositions~\ref{add2} \&~\ref{add61})} \\
(1+o(1)) dn & \textrm{if $d=n-m>0$ is such that} \\
& d=\omega(n^{1/2}) \textrm{ and } o(n) \\
& \textrm{(Propositions~\ref{add51} \&~\ref{add52})} \\
(1+o(1)) (2 - \lambda) n^{3/2} 
& \textrm{if $m=n+ \lambda n^{1/2} +o(n^{1/2})$,} \\
& \textrm{for $\lambda \leq 1$ \phantom{qq} (Thms.~\ref{add3} \&~\ref{add21})} \\
(1+o(1)) \frac{1}{\lambda} n^{3/2} & 
\textrm{if $m=n+ \lambda n^{1/2} +o(n^{1/2})$,} \\
& \textrm{for $\lambda \geq 1$ \phantom{qq} (Thms.~\ref{add3} \&~\ref{add21})} \\
\frac{n^{2}}{d} + O(n) & \textrm{if $d=m-n>0$ is such that} \\
& d=\omega(n^{1/2}) \textrm{ and } o(n) \\
& \textrm{(Propositions~\ref{add11} \&~\ref{add22})} \\
\mu(c)n+O(1) \textrm{ (for known $\mu$)} 
& \textrm{if $m=cn+O(1)$, for $c \in (1,3]$} \\
& \textrm{(Proposition~\ref{ger2 1.2}(\cite{ger2}, 1.2))} \\
\left \lceil \frac{3}{2} (3n-6-m) \right \rceil 
& \textrm{if $n \geq 6$ and $m \geq 2n-3$} \\
& \textrm{(Propositions~\ref{add712} \&~\ref{add713})}
\end{array} \right. 
\end{displaymath} 
\caption{Seven secondary addability results.} \label{addsum2}
\end{table}

The results of Table~\ref{addsum1},
together with the regions detailed in Table~\ref{addsum2},
are illustrated pictorially in Figure~\ref{addpic}:

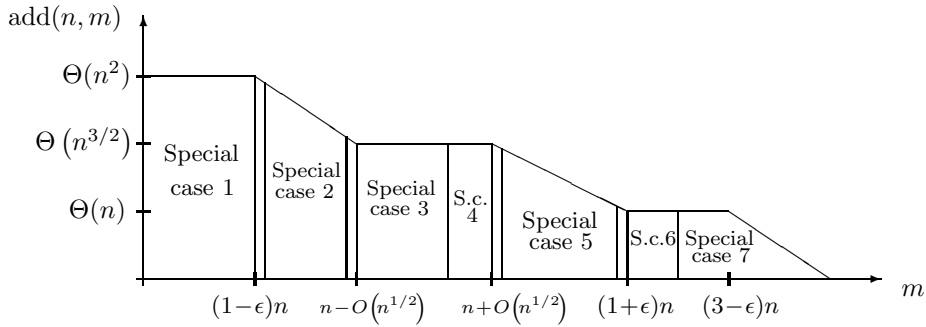
\begin{figure} [h!] 
\setlength{\unitlength}{0.9cm}
\begin{picture}(2,4.15)(0,2)
\put(0,6){\textrm{add$(n,m)$}}
\put(0.9,3.15){$\Theta(n)$}
\put(0.4,4.15){$\Theta\left(n^{3/2}\right)$}
\put(0.8,5.15){$\Theta(n^{2})$}
\put(3,1.75){\small{$(1\!-\!\epsilon)n$}}
\put(4.6,1.75){\scriptsize{$n\!-\!O\!\left(\!n^{1/2}\!\right)$}}
\put(6.7,1.75){\scriptsize{$n\!+\!O\!\left(\!n^{1/2}\!\right)$}}
\put(8.7,1.75){\small{$(1\!+\!\epsilon)n$}}
\put(10.25,1.75){\small{$(3\!-\!\epsilon)n$}}
\put(2,2.15){\vector(0,1){4}}
\put(1.9,2.25){\vector(1,0){11}}
\put(3.65,2.15){\line(0,1){3.1}}
\put(5.15,2.15){\line(0,1){2.1}}
\put(7.15,2.15){\line(0,1){2.1}}
\put(9.15,2.15){\line(0,1){1.1}}
\put(10.65,2.15){\line(0,1){0.2}}
\put(1.9,5.25){\line(1,0){1.75}}

\put(1.9,5.25){\line(1,0){1.75}}
\put(3.65,5.25){\line(3,-2){1.5}}
\put(5.15,4.25){\line(1,0){2}}
\put(7.15,4.25){\line(2,-1){2}}
\put(9.15,3.25){\line(1,0){1.5}}
\put(10.65,3.25){\line(3,-2){1.5}}

\put(3.8,2.25){\line(0,1){2.9}}
\put(5,2.25){\line(0,1){2.1}}
\put(6.5,2.25){\line(0,1){2}}
\put(7.3,2.25){\line(0,1){1.925}}
\put(9,2.25){\line(0,1){1.075}}
\put(9.9,2.25){\line(0,1){1}}

\thicklines
\put(3.65,2.15){\line(0,1){0.2}}
\put(5.15,2.15){\line(0,1){0.2}}
\put(7.15,2.15){\line(0,1){0.2}}
\put(9.15,2.15){\line(0,1){0.2}}
\put(10.65,2.15){\line(0,1){0.2}}
\put(1.9,3.25){\line(1,0){0.2}}
\put(1.9,4.25){\line(1,0){0.2}}
\put(1.9,5.25){\line(1,0){0.2}}

\put(2.3,4){\textrm{\small{Special}}}
\put(2.4,3.5){\textrm{\small{case $1$}}}
\put(3.9,3.75){\textrm{\footnotesize{Special}}}
\put(3.95,3.45){\textrm{\footnotesize{case $2$}}}
\put(5.3,3.5){\textrm{\footnotesize{Special}}}
\put(5.35,3.2){\textrm{\footnotesize{case $3$}}}
\put(6.58,3.35){\textrm{\footnotesize{S.c.}}}
\put(6.78,3.1){\textrm{\footnotesize{$4$}}}
\put(7.6,3){\textrm{\small{Special}}}
\put(7.7,2.7){\textrm{\small{case $5$}}}
\put(9.19,2.75){\textrm{\footnotesize{S.c.$6$}}}
\put(10,2.75){\textrm{\footnotesize{Special}}}
\put(10.1,2.5){\textrm{\footnotesize{case $7$}}}

\put(13.2,2){$m$}

\end{picture}

\caption{Summary of results on add$(n,m)$.} \label{addpic}
\end{figure}

We shall prove the lower bounds first.
We will start (in Lemma~\ref{add33}) with a simple argument that gives us a useful result for when $m<n$,
and from this we shall derive Theorem~\ref{add53} and Propositions~\ref{add51}--\ref{add2}.
We will then prove the remaining lower bounds (in Theorem~\ref{add3} to Proposition~\ref{add712})
by using some very helpful detailed results from \cite{ger2}.

In the second half of this section (Theorem~\ref{add55} to Proposition~\ref{add713})
we will prove the upper bounds,
by copying a construction used in \cite{ger2}.

\begin{figure} [ht]
\setlength{\unitlength}{1cm}
\begin{picture}(20,3.6)(0.25,0)

\put(4.875,0){\textbf{Tables 1 \& 2}}
\put(6,0.125){\oval(2.5,0.5)}
\put(5.325,1.575){\textbf{T\ref{add54}}}
\put(6,2.925){\vector(1,-3){0.333}}
\put(6,2.925){\vector(-1,-3){0.333}}
\put(4.9,3.025){\cite{ger2}(3.5 \& 3.10)}
\put(4.575,1.575){\textbf{T\ref{add3}}}
\put(6.825,1.575){\textbf{P\ref{add712}}}
\put(6.075,1.575){\textbf{P\ref{add11}}}
\put(6,2.925){\vector(1,-1){1}}
\put(6,2.925){\vector(-1,-1){1}}
\put(8.525,1.7){\vector(1,0){0.5}}
\put(7.725,1.575){\textbf{P\ref{ger2 1.2}}}
\put(9.125,1.575){\textbf{C\ref{add41}}}
\put(10.8,1.575){\textbf{T\ref{add55}--P\ref{add713}}}
\put(11.5,1.5){\vector(-4,-1){4.2}}
\put(5,1.475){\vector(1,-2){0.5}}
\put(7,1.475){\vector(-1,-2){0.5}}
\put(6,2.925){\vector(2,-1){2}}
\put(8,1.5){\vector(-1,-1){1}}
\put(5.6,1.475){\vector(1,-4){0.25}}
\put(6.4,1.475){\vector(-1,-4){0.25}}

\put(2.65,2.925){\vector(0,-1){1}}
\put(2.65,2.925){\vector(1,-1){1}}
\put(2.65,2.925){\vector(-1,-1){1}}
\put(2.325,3.025){\textbf{L\ref{add33}}}
\put(1.375,1.575){\textbf{T\ref{add53}}}
\put(2.375,1.575){\textbf{P\ref{add51}}}
\put(3.375,1.575){\textbf{P\ref{add2}}}
\put(2,1.475){\vector(2,-1){2.6}}
\put(2.65,1.475){\vector(2,-1){2.1}}
\put(3.65,1.475){\vector(3,-2){1.5}}

\end{picture}

\caption{The structure of Section~\ref{add}.}
\end{figure}

\begin{displaymath}
\end{displaymath} \\

As mentioned, we shall now start with lower bounds and,
in particular, with a result which will be useful for when $m<n$:

\begin{Lemma} \label{add33}
Let $m=m(n)$. 
Then 
\begin{displaymath}
\emph{add}(n,m) \geq \frac{(n+m)(n-m-1)}{2}.
\end{displaymath}
\end{Lemma}
\textbf{Proof} 
Let $G \in \mathcal{P}(n,m)$.
Clearly, $G$ must have at least $n-m$ components,
and we know that any non-edge between two vertices in different components is addable.

Note that 
the number of possible edges between disjoint sets $X$ and $Y$ is $|X||Y|$ and
that if $|X| \leq |Y|$ then $|X||Y| > (|X|-1)(|Y|+1)$.
Hence, it follows that the number of addable edges between vertices in different components is minimized when
we have $\kappa(G)-1$ isolated vertices and one component of $n- \kappa(G)+1$ vertices.
Thus, 
\begin{eqnarray*}
|\textrm{add}(G)| & \geq & \left(^{\kappa(G)-1}_{\phantom{ww}2}\right) + (\kappa(G)-1)(n-\kappa(G)+1) \\
& = & \frac{(\kappa(G)-1)(\kappa(G)-2)}{2} + (\kappa(G)-1)(n- \kappa(G)+1) \\
& = & (\kappa(G)-1) \left( n - \frac{\kappa(G)}{2} \right).
\end{eqnarray*}

By differentiation, it can be seen that taking $\kappa(G)$ to be $n-m$ 
minimizes $(\kappa(G)-1) \left( n - \frac{\kappa(G)}{2} \right)$
in the region $\kappa(G) \in [n-m,n]$.
Thus, $|\textrm{add}(G)| \geq (n-m-1) \left( n - \frac{n-m}{2} \right) = \frac{(n+m)(n-m-1)}{2}$.
\phantom{qwerty}
\setlength{\unitlength}{.25cm}
\begin{picture}(1,1)
\put(0,0){\line(1,0){1}}
\put(0,0){\line(0,1){1}}
\put(1,1){\line(-1,0){1}}
\put(1,1){\line(0,-1){1}}
\end{picture} \\
\\

Lemma~\ref{add33} provides us with the lower bound for the first of our four main results from Table~\ref{addsum1}:

\begin{Theorem} \label{add53}
Let $d=d(n)>0$ be such that $d= \Omega(n^{1/2})$.
Then 
\begin{displaymath}
\emph{add}(n,n-d) =\Omega(dn).
\end{displaymath} 
\end{Theorem}
\phantom{p}

(Note that we could improve the $d=\Omega(n^{1/2})$ condition to just $d>1$,
directly from Lemma~\ref{add33},
or even beyond if we altered the proof of the lemma slightly,
but we won't bother with this here 
as we shall see a stronger result for when $|m-n|=O(n^{1/2})$ in Theorem~\ref{add3}). \\
\\

Similarly, 
Lemma~\ref{add33} also provides us with the lower bounds for the first two of our six secondary results
from Table~\ref{addsum2}:

\begin{Proposition} \label{add51}
Let $d=d(n)>0$ be such that $d=\Omega(n^{1/2})$ and $d=o(n)$.
Then 
\begin{displaymath}
\emph{add}(n,n-d) \geq (1+o(1)) dn. 
\end{displaymath} 
\end{Proposition} 
\phantom{p}

\begin{Proposition} \label{add2}
Let $A < 1$ be a fixed constant and let $m(n)$ be such that~$m \leq~\!An$ for all large $n$.
Then 
\begin{displaymath}
\emph{add}(n,m) \geq (1+o(1)) \frac{(1-A)(1+A)}{2} n^{2}. 
\end{displaymath}
\end{Proposition} 
\textbf{Proof}
Note that add$(n,m)$ is clearly monotonic in $m$,
so it suffices to consider the case when $m=An$.
$\phantom{qwerty}$ 
\setlength{\unitlength}{.25cm}
\begin{picture}(1,1)
\put(0,0){\line(1,0){1}}
\put(0,0){\line(0,1){1}}
\put(1,1){\line(-1,0){1}}
\put(1,1){\line(0,-1){1}}
\end{picture} \\

\begin{displaymath}
\end{displaymath} 

The remainder of this section relies heavily on the work carried out in \cite{ger2},
in which the sixth of our secondary results was given.
Hence,
we will state that result now,
together with some helpful details from the proof:

\begin{Proposition}[\cite{ger2}, 1.2] \label{ger2 1.2}
Let $c \in (1,3]$ and suppose that $m\!=\!m(n)\!=\!cn +~\!O(1)$ as $n \to \infty$.
Then 
\begin{displaymath}
\emph{add}(n,m)=\mu(c)n+O(1),
\end{displaymath} 
where $\mu(c)$, given explicitly in \cite{ger2}, 
satisfies $\mu(3)=0$, $\mu(c)>0$~$\forall c<3$ and $\mu(c) \to \infty$ as $c \to 1$. 
\end{Proposition}
\textbf{Sketch of Proof}
It is shown that there must exist a graph attaining add$(n,m)$ that belongs to a particular family for which
it happens to be relatively simple to obtain a lower bound for add.
To be more precise,
it is shown in Lemmas $3.5$ and $3.10$ of \cite{ger2} that,
if $n \geq 6$ and $8 \leq m \leq 3n-6$, 
there is a graph $G \in \mathcal{P}(n,m)$ with \\
\\
(i) $i(=i_{n})$ isolated vertices \\
and (ii) a plane embedding with $f_{k}(=f_{k,n})$ faces of size $k$ \\
such that \\
(iii) $G$ is connected apart from the isolated vertices \\
and (iv) add$(n,m) = |\textrm{add}(G)| \geq 
\frac{1}{4} \sum_{k=3}^{n}(k-3)(k+2)f_{k} +i(n- \frac{i+1}{2})$. 

A lower bound for add$(n,m)$ is then obtained using (iv),
while a close upper bound is demonstrated by constructing a subdivided triangulation with few addable edges
(following an idea from \cite{joh}).
$\phantom{qwerty}$ 
\setlength{\unitlength}{.25cm}
\begin{picture}(1,1)
\put(0,0){\line(1,0){1}}
\put(0,0){\line(0,1){1}}
\put(1,1){\line(-1,0){1}}
\put(1,1){\line(0,-1){1}}
\end{picture} \\
\\

A corollary which will be very useful in later sections is:

\begin{Corollary} \label{add41}
Let $m=m(n) \leq (1+o(1))n$.
Then 
\begin{displaymath}
\emph{add}(n,m)=\omega(n).
\end{displaymath}
\end{Corollary}
\textbf{Proof}
This follows from the monotonicity of add(n,m).
\phantom{qwerty} 
\setlength{\unitlength}{.25cm}
\begin{picture}(1,1)
\put(0,0){\line(1,0){1}}
\put(0,0){\line(0,1){1}}
\put(1,1){\line(-1,0){1}}
\put(1,1){\line(0,-1){1}}
\end{picture} 
\begin{displaymath}
\end{displaymath}

We will now see how to modify the proof of Proposition~\ref{ger2 1.2}
to obtain the lower bound for the second of our four main results,
which improves on Lemma~\ref{add33} for the case when $n-m$ is small.
Simultaneously, we will be provided with the lower bounds for the third and fourth of our secondary results.
The idea of the proof is due to Colin McDiarmid.

\begin{Theorem} \label{add3}
Let $\lambda$ be a (not necessarily positive) fixed constant and 
let $m=m(n)$ be such that $m \leq n + \lambda n^{1/2} + o(n^{1/2})$. 
Then
\begin{displaymath}
\emph{add}(n,m) \geq (1+o(1))
\left\{ \begin{array}{ll}
\left( 2 - \lambda \right) n^{3/2} & \textrm{if } \lambda \leq 1 \\
\frac{1}{\lambda} n^{3/2} & \textrm{if } \lambda \geq 1 \\
\end{array} \right. 
\end{displaymath}
\end{Theorem}
\textbf{Proof} 
Note that given any $\alpha > \lambda$, 
we have add$(n,m) \geq$ add$(n,n + \lfloor \alpha n^{1/2} \rfloor )$
for~all~sufficiently large $n$,
by the monotonicity of add.
Thus, since this holds with $\alpha$ arbitrarily close to $\lambda$ and since 
the function
$f(\alpha) = \left\{ \begin{array}{ll}
\left( 2 - \alpha \right) & \textrm{if } \alpha \leq 1 \\
\frac{1}{\alpha} & \textrm{if } \alpha \geq 1 \\
\end{array} \right.$
is continuous in $\alpha$,
it suffices to show add$(n,n + \lfloor \alpha n^{1/2} \rfloor) \geq (1+o(1)) f(\alpha)n^{3/2}$.
Hence, without loss of generality, we may assume that
$m=n + \lfloor \lambda n^{1/2} \rfloor$ $\forall n$
(in fact, for this case we shall be able to show add$(n,m) \geq f(\lambda)n^{3/2} + O(n)$). \\

Clearly, it suffices for us to consider the case when $n$ is large enough that 
$n \geq 6$ and $8 \leq n+ \lfloor \lambda n^{1/2} \rfloor \leq 3n-6$.
Thus, using the proof of Proposition~\ref{ger2 1.2},
we know there is a graph $G \in \mathcal{P}(n,n +\lfloor \lambda n^{1/2} \rfloor)$ with \\
(i) $i(=i_{n})$ isolated vertices \\
and (ii) a plane embedding with $f_{k}(=f_{k,n})$ faces of size $k$ \\
such that \\
(iii) $G$ is connected apart from the isolated vertices \\
and (iv) add$(n,n+\lfloor \lambda n^{1/2} \rfloor) = |\textrm{add}(G)| \geq 
\frac{1}{4} \sum_{k=3}^{n}(k-3)(k+2)f_{k} +i(n- \frac{i+1}{2})$. \\

Note that, by (ii), we have
\begin{equation}
\sum_{k=3}^{n}kf_{k}=2(n+\lfloor \lambda n^{1/2} \rfloor), \label{eq:2}
\end{equation}
and, by (iii) and Euler's formula, we have
\begin{eqnarray}
\sum_{k=3}^{n}f_{k}=\lfloor \lambda n^{1/2} \rfloor+i+2. \label{eq:3} \\ \nonumber
\end{eqnarray} 

Using (iv) and the fact that $i \! \leq \! n$,
we have add$(n,n+\lfloor \lambda n^{1/2} \rfloor) \! \geq \!
i \left( n- \frac{i+1}{2} \right) \! \geq \! i \left( \frac{n-1}{2} \right) > \frac{in}{3}$.
Let $x= \max \{6, 6 - 3 \lambda \}$.
Then, for those values of $n$ for which
$i = i_{n} \geq xn^{1/2}$,
we have 
\begin{eqnarray}
\textrm{add}(n,n+\lfloor \lambda n^{1/2} \rfloor) 
& > & \frac{xn^{3/2}}{3} \label{eq:pre1} \\
& = & \max \{2, 2- \lambda \}n^{3/2} \nonumber \\
& \geq & f(\lambda)n^{3/2}. \nonumber \\ \nonumber 
\end{eqnarray}

We shall now also obtain a lower bound on add$(n,n + \lambda n^{1/2})$ for those values of $n$ for which 
$i = i_{n} < xn^{1/2}$. 
By (iv), it suffices to minimize 
$h = h(i,k_{3},k_{4}, \ldots, k_{n}) = \frac{1}{4} \sum_{k=3}^{n}(k-3)(k+2)f_{k} +i(n- \frac{i+1}{2})$ 
over all remaining choices of non-negative integers $i,k_{3},k_{4}, \ldots, k_{n}$,
subject to constraints (\ref{eq:2}) and~(\ref{eq:3}).

Note that the formula $(k-3)(k+2)$ yields a convex function of $k$.
Thus, as observed in the proof of Theorem 1.2 of \cite{ger2},
the minimum value of $h$ under constraints (\ref{eq:2}) and (\ref{eq:3})
is attained with either a single non-zero value $f_{l}$
or with two `adjacent' non-zero values $f_{l}$ and $f_{l+1}$.
Thus (with $f_{l+1}=0$ if appropriate), 
by (\ref{eq:2}) we have 
$2(n+\lfloor \lambda n^{1/2} \rfloor)=lf_{l}+(l+1)f_{l+1} < (l+~1)(f_{l}+f_{l+1}) = 
(l+1)(\lfloor \lambda n^{1/2} \rfloor+i+2)$ by (\ref{eq:3}).
Hence, $l-3 > \frac{2(n+\lfloor \lambda n^{1/2} \rfloor)}{\lfloor \lambda n^{1/2} \rfloor+i+2}-4$.
Therefore, 
\begin{eqnarray}
& & \textrm{add}(n,n+\lfloor \lambda n^{1/2} \rfloor) \nonumber \\
& \geq & \frac{1}{4} \left( (l-3)(l+2)f_{l} + (l-2)(l+3)f_{l+1} \right) + i \left( n - \frac{i+1}{2} \right) \nonumber \\ 
& > & \left(\frac{(n+\lfloor \lambda n^{1/2} \rfloor)}{2(\lfloor \lambda n^{1/2} \rfloor+i+2)}-1\right) 
\left( (l+2)f_{l} + (l+3)f_{l+1} \right) + i \left( n - \frac{i+1}{2} \right) \nonumber \\
& > & \left(\frac{(n+\lfloor \lambda n^{1/2} \rfloor)}{2(\lfloor \lambda n^{1/2} \rfloor+i+2)}-1\right) 
\left( lf_{l} + (l+1)f_{l+1} \right) + i \left( n - \frac{i+1}{2} \right) \nonumber \\
& \stackrel{\textrm{\small{by (\ref{eq:2})}}}> & 
\frac{(n+\lfloor \lambda n^{1/2} \rfloor)^{2}}{(\lfloor \lambda n^{1/2} \rfloor+i+2)}
-2(n+\lfloor \lambda n^{1/2} \rfloor)  + i \left( n - \frac{i+1}{2} \right) \label{eq:pre2} \\
& = & \frac{n^{2}}{\lambda n^{1/2} + i + 2} + in + O(n) \phantom{w} \textrm{ if } i < xn^{1/2} \nonumber \\
& \geq & \left\{ \begin{array}{ll}
\left( 2 - \lambda \right) n^{3/2} + O(n) & \textrm{if } \lambda < 1 \\
\frac{1}{\lambda} n^{3/2} + O(n) & \textrm{if } \lambda \geq 1, \\
\end{array} \right. \nonumber \\
& & \textrm{since, by differentiation, } i = \max \{ 0, (1- \lambda)n^{1/2} -2 \}
\nonumber \\
& & \textrm{minimizes } 
\frac{n^{2}}{\lambda n^{1/2} + i + 2} + in
\textrm{ in the region } i \geq 0. \nonumber \\ \nonumber 
\end{eqnarray} 

Thus, regardless of whether or not $i_{n} \geq xn^{1/2}$, 
we have 
\begin{eqnarray*}
\textrm{add}(n,n+\lfloor \lambda n^{1/2} \rfloor) 
& \geq & 
\left\{ \begin{array}{ll}
\left( 2 - \lambda \right) n^{3/2} + O(n) & \textrm{if } \lambda < 1 \\
\frac{1}{\lambda} n^{3/2} + O(n) & \textrm{if } \lambda \geq 1. 
\phantom{qwerty} 
\setlength{\unitlength}{.25cm}
\begin{picture}(1,1)
\put(0,0){\line(1,0){1}}
\put(0,0){\line(0,1){1}}
\put(1,1){\line(-1,0){1}}
\put(1,1){\line(0,-1){1}}
\end{picture} \\
\end{array} \right. 
\end{eqnarray*}

By the same method as in the proof of Theorem~\ref{add3}, 
we also obtain the lower bound for our third main result:

\begin{Theorem} \label{add54}
Let $d = d(n)>0$ be such that $d=\Omega(n^{1/2})$ and $\limsup \frac{d}{n}<2$. 
Then 
\begin{displaymath}
\emph{add}(n,n+d) = \Omega \left( \frac{n^{2}}{d} \right). 
\end{displaymath}
\end{Theorem}
\textbf{Proof}
Let $\epsilon \in \{ 0,1 \}$ be arbitrary.
Then for those values of $n$ for which $d \geq \epsilon n$,
the result is equivalent to showing that add$(n,n+d)=\Omega(n)$,
which follows immediately from Proposition~\ref{ger2 1.2}
(since $\limsup_{n \to \infty} \frac{d}{n} < 2$ and add is monotonic).
Hence, we may assume that $d \leq \epsilon n$~$\forall n$.

We now follow the proof of Theorem~\ref{add3},
with $\lfloor \lambda n^{1/2} \rfloor$ replaced by $d$.
As before, we find (by (\ref{eq:pre1}))
that if $i \geq xn^{1/2}$ 
(where $x$ is an arbitrary constant)
then add$(n,n+d) > \frac{xn^{3/2}}{3} = \Omega \left( n^{3/2} \right)$,
and (analogously to (\ref{eq:pre2}))
that if $i <~xn^{1/2}$ then add$(n,n+d) > \frac{(n+d)^{2}}{d+i+2} - 2(n+d) + i \left( n - \frac{i+1}{2} \right)
\geq \frac{(n+d)^{2}}{d+i+2} - 2(n+d)
=~\Omega \left( \frac{n^{2}}{d} \right)$,
since $\frac{d}{n} \leq \epsilon <1$
(and so $\frac{(n+d)^{2}}{d+i+1} > 2(n+d)$).
Thus, we obtain
add$(n,n+d) \geq \min \left \{ \Omega \left( n^{3/2} \right), \Omega \left( \frac{n^{2}}{d} \right) \right \} 
= \Omega \left( \frac{n^{2}}{d} \right)$.
\phantom{qwerty}
\begin{picture}(1,1)
\put(0,0){\line(1,0){1}}
\put(0,0){\line(0,1){1}}
\put(1,1){\line(-1,0){1}}
\put(1,1){\line(0,-1){1}}
\end{picture} 
\begin{displaymath}
\end{displaymath}

Similarly,
we may obtain the lower bound for our fifth secondary result:

\begin{Proposition} \label{add11}
Let $d = d(n)>0$ be such that $d=\omega(n^{1/2})$ and $d=o(n)$. 
Then 
\begin{displaymath}
\emph{add}(n,n+d) \geq \frac{n^{2}}{d} + O(n). 
\end{displaymath}
\end{Proposition}
\textbf{Proof}
We copy the proof of Theorem~\ref{add54}, 
except that for the case when~$i <~\!xn^{1/2}$ we now have 
add$(n,n+d) > \frac{(n+d)^{2}}{d+i+2} - 2(n+d) + i \left( n - \frac{i+1}{2} \right) \geq \frac{n^{2}}{d} + O(n)$.~
\begin{picture}(1,1)
\put(0,0){\line(1,0){1}}
\put(0,0){\line(0,1){1}}
\put(1,1){\line(-1,0){1}}
\put(1,1){\line(0,-1){1}}
\end{picture} 
\begin{displaymath}
\end{displaymath}

Finally for this part of the section,
we shall obtain the lower bound for our last secondary result
(which also gives us the lower bound for our final main result):

\begin{Proposition} \label{add712}
Let $n \geq 6$ and let $m \geq 8$.
Then 
\begin{displaymath}
\emph{add}(n,m) \geq \left \lceil \frac{3}{2} (3n-6-m) \right \rceil.
\end{displaymath}
\end{Proposition}
\textbf{Proof}
Since add$(n,m)$ is integral,
it clearly suffices to show that
add$(n,m) \geq \frac{3}{2} (3n-6-m)$.

As in the proof of Proposition~\ref{ger2 1.2},
since $n \geq 6$ and $m \geq 8$,
we have a graph $G \in \mathcal{P}(n,m)$ with \\
(i) $i(=i_{n})$ isolated vertices \\
and (ii) a plane embedding with $f_{k}(=f_{k,n})$ faces of size $k$ \\
such that \\
(iii) $G$ is connected apart from the isolated vertices \\
and (iv) add$(n,m) = |\textrm{add}(G)| \geq 
\frac{1}{4} \sum_{k=3}^{n}(k-3)(k+2)f_{k} +i(n- \frac{i+1}{2})$. \\

Thus,
\begin{eqnarray*}
\sum_{k \geq 3} (k-3) f_{k} & = & \sum_{k \geq 3} kf_{k} - 3 \sum_{k \geq 3} f_{k} \\
& = & 2m - 3(m-n+i+2) \phantom{ww} \textrm{ by (iii) and Euler's formula.} \\
& = & 3n-6-m-3i.
\end{eqnarray*}
Hence, 
$\frac{1}{4} \sum_{k=3}^{n}(k-3)(k+2)f_{k} \geq \frac{3}{2} \sum_{k \geq 3} (k-3) f_{k} = \frac{3}{2}(3n-6-m-3i)$,
and so add$(n,m) \geq \frac{3}{2}(3n-6-m-3i) + i \left( n - \frac{i+1}{2} \right) 
= \frac{3}{2}(3n-6-m) + i \left( n - \frac{i}{2} -5 \right)$.
Thus, we certainly have add$(n,m) \geq \frac{3}{2} (3n-6-m)$
for those values of $n$ for which $i=i_{n} \leq 2(n-5)$. \\

It now only remains to consider the values of $n$ for which $i=i_{n} > 2(n-5)$.
But recall add$(n,m) \geq \frac{1}{4} \sum_{k=3}^{n}(k-3)(k+2)f_{k} +i(n- \frac{i+1}{2}) \geq i(n- \frac{i+1}{2})$
and note that $i=2n-9$ minimizes $i \left( n - \frac{i+1}{2} \right)$
in the region $i \in \{ 2n-9,n \}$.
Thus, for those values of $n$ for which $i=i_{n} > 2(n-5)$,
we must have add$(n,m) \geq (2n-9) \left( n - \frac{2n-8}{2} \right) = 8n-36 \geq \frac{3}{2} (3n-6-m)$,
since $n \geq 6$ and $m \geq 8$.
\phantom{ww}
\begin{picture}(1,1)
\put(0,0){\line(1,0){1}}
\put(0,0){\line(0,1){1}}
\put(1,1){\line(-1,0){1}}
\put(1,1){\line(0,-1){1}}
\end{picture} 
\begin{displaymath}
\end{displaymath}

We shall now spend the remainder of this section providing upper bounds 
to show that our lower bounds are actually all tight.

We start with the upper bound for the first of our main results:

\begin{Theorem} \label{add55} 
Let $d=d(n)>0$ be such that $d= \Omega(n^{1/2})$.
Then 
\begin{displaymath}
\emph{add}(n,n-d) =O(dn).
\end{displaymath}
\end{Theorem}
\textbf{Proof}
Clearly, add$(n,m) \leq \left( ^{n}_{2} \right) = O(n^{2})$.
Thus, it suffices to prove the result for the case when $d \leq \epsilon n$ $\forall n$,
where $\epsilon >0$ is an arbitrary (small) constant.

By the monotonicity of add, it suffices to construct a planar graph, $G$,
with $e(G) \leq~n -d$
and add$(G)$ sufficiently small.
We shall use the same construction as in the proof of Proposition~\ref{ger2 1.2} (\cite{ger2}, 1.2),
which itself follows an idea from \cite{joh}.

Let $x > 0$ be an arbitrary constant 
and consider the following triangulation on $\lfloor xd \rfloor$ vertices:

\begin{figure}[ht]
\setlength{\unitlength}{1cm}
\begin{picture}(2,5.25)
\put(4,0){\line(1,0){4}}
\put(6,0.5){\line(0,1){2.15}}
\put(6,2.8){\line(0,1){0.25}}
\put(6,3.2){\line(0,1){0.25}}
\put(6,3.55){\line(0,1){0.25}}
\put(6,3.9){\line(0,1){1.1}}
\put(4,0){\line(4,1){2}}
\put(8,0){\line(-4,1){2}}
\put(4,0){\line(4,3){2}}
\put(8,0){\line(-4,3){2}}
\put(4,0){\line(4,5){2}}
\put(8,0){\line(-4,5){2}}
\put(4,0){\line(1,2){2}}
\put(8,0){\line(-1,2){2}}
\put(4,0){\line(2,5){2}}
\put(8,0){\line(-2,5){2}}

\put(4,0){\circle*{0.1}}
\put(8,0){\circle*{0.1}}
\put(6,0.5){\circle*{0.1}}
\put(6,2.5){\circle*{0.1}}
\put(6,1.5){\circle*{0.1}}
\put(6,4){\circle*{0.1}}
\put(6,5){\circle*{0.1}}
\end{picture} 
\caption{Our triangulation.}
\end{figure}
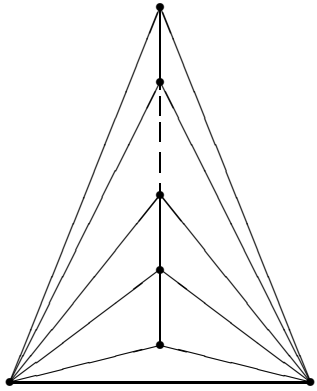

From now on, we shall refer to the vertices in this triangulation as `core vertices',
and to the two vertices of degree~ $>4$ as `base vertices'.
We shall use the term `spine line' to denote an edge between two non-base core vertices or between the two base vertices.

Let $i=2x+1$ and note that we may assume that $n - \lfloor xd \rfloor - \lfloor id \rfloor >~0$,
since $x$ may be taken to be arbitrarily small and $d \leq \epsilon n$.
Let us insert $n - \lfloor xd \rfloor - \lfloor id \rfloor$ new vertices 
as evenly as possible on the $\lfloor xd \rfloor - 2$ spine lines
(i.e. so that between 
$\left \lfloor \frac{n - \lfloor xd \rfloor - \lfloor id \rfloor}
{\lfloor xd \rfloor - 2} \right \rfloor$
and 
$\left \lceil \frac{n - \lfloor xd \rfloor - \lfloor id \rfloor}
{\lfloor xd \rfloor - 2} \right \rceil$
new vertices are inserted on each spine line).
Then our new graph is a subdivision of our triangulation, which was $3$-connected.
Thus, by a theorem of Whitney \cite{whi}, 
our new graph has the following \textit{unique} embedding in the plane:

\begin{figure}[ht]
\setlength{\unitlength}{1cm}
\begin{picture}(2,5.5)
\put(4,0){\line(1,0){4}}
\put(6,0.5){\line(0,1){2.15}}
\put(6,2.8){\line(0,1){0.25}}
\put(6,3.2){\line(0,1){0.25}}
\put(6,3.55){\line(0,1){0.25}}
\put(6,3.9){\line(0,1){1.1}}
\put(4,0){\line(4,1){2}}
\put(8,0){\line(-4,1){2}}
\put(4,0){\line(4,3){2}}
\put(8,0){\line(-4,3){2}}
\put(4,0){\line(4,5){2}}
\put(8,0){\line(-4,5){2}}
\put(4,0){\line(1,2){2}}
\put(8,0){\line(-1,2){2}}
\put(4,0){\line(2,5){2}}
\put(8,0){\line(-2,5){2}}

\put(4,0){\circle*{0.1}}
\put(8,0){\circle*{0.1}}
\put(6,0.5){\circle*{0.1}}
\put(6,2.5){\circle*{0.1}}
\put(6,1.5){\circle*{0.1}}
\put(6,4){\circle*{0.1}}
\put(6,5){\circle*{0.1}}

\put(5,0){\circle*{0.1}}
\put(7,0){\circle*{0.1}} 
\put(6,0.8333){\circle*{0.1}}
\put(6,4.6667){\circle*{0.1}}
\put(6,1.1667){\circle*{0.1}}
\put(6,2){\circle*{0.1}}
\put(6,4.3333){\circle*{0.1}}

\end{picture} 

\caption{The unique embedding in the plane of our new graph.}
\end{figure}
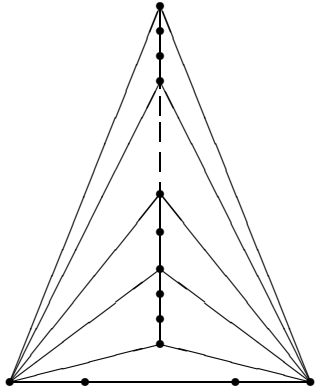

Let us also include $\lfloor id \rfloor$ isolated vertices in our graph, which we shall call $G$.
Then the total number of vertices in $G$ is
$\lfloor xd \rfloor + (n - \lfloor xd \rfloor  - \lfloor id \rfloor) + \lfloor id \rfloor = n$
and the total number of edges is 
$(3 \lfloor xd \rfloor - 6) + (n - \lfloor xd \rfloor - \lfloor id \rfloor)
= n + 2\lfloor xd \rfloor - \lfloor id \rfloor - 6$.

Note that 
$n-d > n + 2 \lfloor xd \rfloor - \lfloor (2x+1)d \rfloor -6 = n + 2 \lfloor xd \rfloor - \lfloor id \rfloor - 6$,
so add$(n,n-d) \leq \textrm{add}(G)$.

Let us now consider add$(G)$.
The only addable non-edges are those between:
(i) two isolated vertices;
(ii) an isolated vertex and another vertex;
(iii) a non-isolated new vertex and another new vertex on the same spine line;
(iv) a non-base core vertex and a new vertex from one of two spine lines;
(v) a non-base core vertex and at most two other non-base core vertices;
or (vi) a base vertex and another vertex.

Thus,
\begin{eqnarray*}
\textrm{add}(n,n-d) 
& \leq & \textrm{add}(G) \\
& \leq & \left( ^{\lfloor id \rfloor}_{\phantom{w}2} \right) \\
& & + \lfloor id \rfloor (n - \lfloor id \rfloor) \\
& & +\frac{1}{2} \left( \Big( n - \lfloor id \rfloor - \lfloor xd \rfloor \Big)
\left( \left \lceil 
\frac{n - \lfloor id \rfloor - \lfloor xd \rfloor}{ \lfloor xd \rfloor - 2} 
\right \rceil \right) \right) \\
& & + \lfloor xd \rfloor 
\left( 2 \left \lceil 
\frac{n - \lfloor id \rfloor - \lfloor xd \rfloor}{ \lfloor xd \rfloor - 2} 
\right \rceil \right) \\
& & + \lfloor xd \rfloor \\
& & + 2n \\
& = & O \left( dn \right), \textrm{ since } d = \Omega \left( n^{1/2} \right).
\phantom{qwerty} 
\setlength{\unitlength}{.25cm}
\begin{picture}(1,1)
\put(0,0){\line(1,0){1}}
\put(0,0){\line(0,1){1}}
\put(1,1){\line(-1,0){1}}
\put(1,1){\line(0,-1){1}}
\end{picture}
\end{eqnarray*}
\begin{displaymath}
\end{displaymath}

By the same method, 
we may also obtain the upper bounds for the first two of our secondary results:

\begin{Proposition} \label{add61}
Let $A<1$ be a fixed constant and let $m=m(n)$ be such that $m \geq An+o(n)$ for all large $n$.
Then 
\begin{displaymath}
\emph{add}(n,m) \leq (1+o(1)) \frac{(1-A)(1+A)}{2} n^{2}. 
\end{displaymath}
\end{Proposition}
\textbf{Proof}
As in the proof of Theorem~\ref{add55}, with $d=(1-A)n$, we obtain
\begin{eqnarray*}
\textrm{add}(n,An) 
& \leq & \left( ^{\lfloor i(1-A)n \rfloor}_{\phantom{www}2} \right) \\
& & + \lfloor i(1-A)n \rfloor (n - \lfloor i(1-A)n \rfloor) \\
& & +\frac{1}{2} \Big( n - \lfloor i(1-A)n \rfloor - \lfloor x(1-A)n \rfloor \Big) \\
& & \phantom{w} \cdot
\left( \left \lceil 
\frac{n - \lfloor i(1-A)n \rfloor - \lfloor x(1-A)n \rfloor}{ \lfloor x(1-A)n \rfloor - 2} 
\right \rceil \right) \\
& & + \lfloor x(1-A)n \rfloor 
\left( 2 \left \lceil 
\frac{n - \lfloor i(1-A)n \rfloor - \lfloor x(1-A)n \rfloor}{ \lfloor x(1-A)n \rfloor - 2} 
\right \rceil \right) \\
& & + \lfloor x(1-A)n \rfloor \\
& & + 2n \\
& = & (1+o(1)) \left( \frac{i^{2}(1-A)^{2}n^{2}}{2} + i(1-A)n(n-i(1-A)n) \right) \\
& = & (1+o(1)) i(1-A)n^{2} \left( 1- \frac{i(1-A)}{2} \right). 
\end{eqnarray*}
Since this holds with $i$ arbitrarily close to $1$,
we must have 
\begin{eqnarray*}
\textrm{add}(n,An) & \leq & (1+o(1)) (1-A)n^{2} \left( 1- \frac{1-A}{2} \right) \\
& = & (1+o(1)) \frac{(1-A)(1+A)}{2}n^{2}.
\phantom{qwerty} 
\setlength{\unitlength}{.25cm}
\begin{picture}(1,1)
\put(0,0){\line(1,0){1}}
\put(0,0){\line(0,1){1}}
\put(1,1){\line(-1,0){1}}
\put(1,1){\line(0,-1){1}}
\end{picture}
\end{eqnarray*} 
\phantom{p}

\begin{Proposition} \label{add52} 
Let $d=d(n)>0$ be such that $d= \omega(n^{1/2})$ and $d=o(n)$.
Then 
\begin{displaymath}
\emph{add}(n,n-d) \leq (1+o(1))dn.
\end{displaymath}
\end{Proposition}
\textbf{Proof}
By the same method as in Theorem~\ref{add55}, we obtain
\begin{eqnarray*}
\textrm{add}(n,n-d) 
& \leq & \left( ^{\lfloor id \rfloor}_{\phantom{w}2} \right) \\
& & + \lfloor id \rfloor (n - \lfloor id \rfloor) \\
& & +\frac{1}{2} \left( \Big( n - \lfloor id \rfloor - \lfloor xd \rfloor \Big)
\left( \left \lceil 
\frac{n - \lfloor id \rfloor - \lfloor xd \rfloor}{\lfloor xd \rfloor - 2} 
\right \rceil \right) \right) \\
& & + \lfloor xd \rfloor 
\left( 2 \left \lceil 
\frac{n - \lfloor id \rfloor - \lfloor xd \rfloor}{ \lfloor xd \rfloor - 2} 
\right \rceil \right) \\
& & + \lfloor xd \rfloor \\
& & + 2n \\
& = & (1+o(1)) idn. 
\end{eqnarray*}
Since this holds with $i$ arbitrarily close to $1$,
we must have add$(n,n-d) \leq (1+o(1))dn$.
\phantom{qwerty}
\setlength{\unitlength}{.25cm}
\begin{picture}(1,1)
\put(0,0){\line(1,0){1}}
\put(0,0){\line(0,1){1}}
\put(1,1){\line(-1,0){1}}
\put(1,1){\line(0,-1){1}}
\end{picture}
\begin{displaymath}
\end{displaymath} 
\phantom{p}

We may also use the same method to obtain the upper bound for the second of our main results,
and simultaneously the third and fourth of our secondary results:

\begin{Theorem} \label{add21}
Let $\lambda$ be a (not necessarily positive) fixed constant and let $m=m(n)$ be such that 
$m \geq n+\lambda n^{1/2} + o(n^{1/2})$. 
Then
\begin{displaymath}
\emph{add}(n,m) \leq (1+o(1))
\left\{ \begin{array}{ll}
\left( 2 - \lambda \right) n^{3/2} & \textrm{if } \lambda \leq 1 \\
\frac{1}{ \lambda} n^{3/2} & \textrm{if } \lambda \geq 1 \\
\end{array} \right. 
\end{displaymath}
\end{Theorem}
\textbf{Proof}
By the same argument as in the proof of Theorem~\ref{add3}, 
we may without loss of generality assume that
$m=n + \lfloor \lambda n^{1/2} \rfloor$ $\forall n$.
Again, for this case we shall be able to show that 
$\textrm{add}(n,m) \leq 
\left\{ \begin{array}{ll}
\left( 2 - \lambda \right) n^{3/2} + O(n) & \textrm{if } \lambda \leq 1 \\
\frac{1}{ \lambda} n^{3/2} + O(n) & \textrm{if } \lambda \geq 1. \\
\end{array} \right.$

By the monotonicity of add, it suffices to construct a planar graph, $G$,
with $e(G) \leq~n +~ \lambda n^{1/2}$
and add$(G)$ sufficiently small.

Let $x= \max \{ \frac{1}{2}, \frac{\lambda}{2} \}$ and let $i=2x- \lambda$.
We take the triangulation of Theorem~\ref{add55} with $\lfloor xn^{1/2} \rfloor$ `core vertices'
and insert $n - \lfloor xn^{1/2} \rfloor - \lfloor in^{1/2} \rfloor$ 
new vertices as evenly as possible on the `spine lines'.
We also include $\lfloor in^{1/2} \rfloor$ isolated vertices in our graph.

The total number of vertices in our new graph, $G$, is $n$ 
and the number of edges is 
$(3 \lfloor x n^{1/2} \rfloor - 6) + (n - \lfloor x n^{1/2} \rfloor - \lfloor i n^{1/2} \rfloor)
= n + 2\lfloor x n^{1/2} \rfloor -~\lfloor i n^{1/2} \rfloor -~6$.
Note that 
$\lfloor \lambda n^{1/2} \rfloor = \lfloor (2x-i) n^{1/2} \rfloor > 
2 \lfloor x n^{1/2} \rfloor - \lfloor in^{1/2} \rfloor - 6$, 
so
add$(n,m) \leq~\textrm{add}(G)$.

Similarly to with the proof of Theorem~\ref{add55}, we obtain
\begin{eqnarray*}
\textrm{add}(G) & \leq & \left( ^{\lfloor in^{1/2} \rfloor}_{\phantom{ww}2} \right) \\
& & + \lfloor in^{1/2} \rfloor (n - \lfloor in^{1/2} \rfloor) \\
& & +\frac{1}{2} \left( \left( n - \lfloor in^{1/2} \rfloor - \lfloor xn^{1/2} \rfloor \right)
\left( \left \lceil 
\frac{n - \lfloor in^{1/2} \rfloor - \lfloor xn^{1/2} \rfloor}{ \lfloor xn^{1/2} \rfloor - 2} 
\right \rceil \right) \right) \\
& & + \lfloor xn^{1/2} \rfloor 
\left( 2 \left \lceil 
\frac{n - \lfloor in^{1/2} \rfloor - \lfloor xn^{1/2} \rfloor}{ \lfloor xn^{1/2} \rfloor - 2} 
\right \rceil \right) \\
& & + \lfloor xn^{1/2} \rfloor \\
& & + 2n \\
& = & \left( i + \frac{1}{2x} \right) n^{3/2} + O(n) \\
& = & 
\left\{ \begin{array}{ll}
\left( 2 - \lambda \right) n^{3/2} + O(n) & \textrm{if } \lambda \leq 1 \\
\frac{1}{ \lambda} n^{3/2} + O(n) & \textrm{if } \lambda \geq 1. 
\phantom{qwerty}
\setlength{\unitlength}{.25cm}
\begin{picture}(1,1)
\put(0,0){\line(1,0){1}}
\put(0,0){\line(0,1){1}}
\put(1,1){\line(-1,0){1}}
\put(1,1){\line(0,-1){1}}
\end{picture} \\
\end{array} \right. 
\end{eqnarray*} \\
\\

Similarly, we may use the same method to obtain the upper bound for the third of our main results:

\begin{Theorem} \label{add56}
Let $d=d(n)>0$ be such that $d=O(n)$.
Then
\begin{displaymath}
\emph{add}(n,n + d) = O \left( \frac{n^{2}}{d} \right).
\end{displaymath}
\end{Theorem}
\textbf{Proof}
We take the triangulation of Theorem~\ref{add55} with $\lfloor \frac{d}{2} \rfloor$ `core vertices' 
and insert $n - \lfloor \frac{d}{2} \rfloor$ new vertices as evenly as possible on the `spine lines'.

The total number of vertices in our new graph, $G$, is $n$
and the total number of edges is 
$\left( 3 \lfloor \frac{d}{2} \rfloor - 6\right) + \left(n - \lfloor \frac{d}{2} \rfloor \right) 
= n + 2 \lfloor \frac{d}{2} \rfloor - 6$.

Similarly to with the proof of Theorem~\ref{add55}, we have
\begin{eqnarray*}
\textrm{add}(n,n+d) & \leq & \textrm{add}(n,n+ 2 \left \lfloor d/2 \right \rfloor - 6) \\
& \leq & \textrm{add}(G) \\
& \leq & \frac{1}{2} \left( \left( n - \left \lfloor \frac{d}{2} \right \rfloor \right) 
\left( \left \lceil \frac{n - \left \lfloor \frac{d}{2} \right \rfloor}
{ \left \lfloor \frac{d}{2} \right \rfloor - 2} \right \rceil \right) \right) \\
& & + \left \lfloor \frac{d}{2} \right \rfloor 
\left( 2 \left \lceil \frac{n - \left \lfloor \frac{d}{2} \right \rfloor}
{\left \lfloor \frac{d}{2} \right \rfloor - 2} \right \rceil \right)
+ \left \lfloor \frac{d}{2} \right \rfloor + 2n \\
& = & O \left( \frac{n^{2}}{d} \right).
\phantom{qwerty}
\setlength{\unitlength}{.25cm}
\begin{picture}(1,1)
\put(0,0){\line(1,0){1}}
\put(0,0){\line(0,1){1}}
\put(1,1){\line(-1,0){1}}
\put(1,1){\line(0,-1){1}}
\end{picture}
\end{eqnarray*}
\begin{displaymath}
\end{displaymath} 
\phantom{p}

The same calculations also provide us with the upper bound for our fifth secondary result:

\begin{Proposition} \label{add22}
Let $d=d(n)>0$ be such that $d=o(n)$.
Then
\begin{displaymath}
\emph{add}(n,n +d) \leq \frac{n^{2}}{d}+O(n). 
\end{displaymath}
\end{Proposition}

Finally,
we may obtain the upper bound for our last secondary result
(which combines with Theorem~\ref{add56} to also give us the upper bound for our final main result):

\begin{Proposition} \label{add713}
Let $m \geq 2n-3$.
Then 
\begin{displaymath}
\emph{add}(n,m) \leq \left \lceil \frac{3}{2} (3n-6-m) \right \rceil.
\end{displaymath}
\end{Proposition}
\textbf{Proof}
We take the triangulation of Theorem~\ref{add55} with 
$\left \lfloor \frac{m-n}{2} \right \rfloor + 3$ `core vertices'
and insert $\left \lceil \frac{3n-m}{2} \right \rceil - 3$ new vertices
as evenly as possible on the `spine lines'.
If $m-n$ is odd,
we also insert a new edge from one new vertex to one base vertex.

The total number of vertices in our new graph, $G$,
is $n$ (since $m-n$ and $3n-m$ have the same parity)
and the total number of edges is
$3 \left( \lfloor \frac{m-n}{2} \rfloor + 3 \right) - 6
+ \left \lceil \frac{3n-m}{2} \right \rceil - 3 + \mathbf{1} \{m-n \textrm{ odd} \}
= 3 \lfloor \frac{m-n}{2} \rfloor + \left \lceil \frac{3n-m}{2} \right \rceil + \mathbf{1} \{m-n \textrm{ odd} \} = m$.

Note that the number of spine lines is 
$\left \lfloor \frac{m-n}{2} \right \rfloor + 1 \geq \left \lceil \frac{3n-m}{2} \right \rceil - 3$,
since $m \geq 2n-3$.
Thus, we will have at most one new vertex on each spine line,
and so $G$ is as shown in Figure~\ref{lastadd}.
\begin{figure} [ht]
\setlength{\unitlength}{0.9cm}

\begin{picture}(2,4.4)(-0.6,0.5)
\put(4,0){\line(1,0){4}}
\put(6,0.5){\line(0,1){2.15}}
\put(6,2.8){\line(0,1){0.25}}
\put(6,3.2){\line(0,1){0.25}}
\put(6,3.55){\line(0,1){0.25}}
\put(6,3.9){\line(0,1){1.1}}
\put(4,0){\line(4,1){2}}
\put(8,0){\line(-4,1){2}}
\put(4,0){\line(4,3){2}}
\put(8,0){\line(-4,3){2}}
\put(4,0){\line(4,5){2}}
\put(8,0){\line(-4,5){2}}
\put(4,0){\line(1,2){2}}
\put(8,0){\line(-1,2){2}}
\put(4,0){\line(2,5){2}}
\put(8,0){\line(-2,5){2}}

\put(4,0){\circle*{0.1}}
\put(8,0){\circle*{0.1}}
\put(6,0.5){\circle*{0.1}}
\put(6,2.5){\circle*{0.1}}
\put(6,1.5){\circle*{0.1}}
\put(6,4){\circle*{0.1}}
\put(6,5){\circle*{0.1}}

\put(6,0){\circle*{0.1}}
\put(6,1){\circle*{0.1}}
\put(6,4.5){\circle*{0.1}}

\put(4,0){\line(2,1){2}}

\end{picture} 

\caption{The unique embedding in the plane of our new graph, $G$.} \label{lastadd}
\end{figure}
It is then clear that
\begin{eqnarray*}
|\textrm{add}(G)| & = & 3 \times (\textrm{number of new vertices}) - \mathbf{1} \{m-n \textrm{ odd} \} \\
& = & 3 \left \lceil \frac{3n-m}{2} \right \rceil - 9 - \mathbf{1} \{m-n \textrm{ odd} \} \\
& = & \left \lceil \frac{9n-3m}{2} -9 \right \rceil \\
& = & \left \lceil \frac{3}{2} (3n-6-m) \right \rceil.
\phantom{qwerty}
\setlength{\unitlength}{.25cm}
\begin{picture}(1,1)
\put(0,0){\line(1,0){1}}
\put(0,0){\line(0,1){1}}
\put(1,1){\line(-1,0){1}}
\put(1,1){\line(0,-1){1}}
\end{picture}
\end{eqnarray*}

\newpage
\section{Components I: Lower Bounds} \label{cptlow}

We now come to our first main section, 
where we shall start to look at lower bounds for
$\mathbf{P}[P_{n,m} \textrm{ will have a component isomorphic to } H]$.
We will see (in Theorem~\ref{gen3}) that,
for any connected planar graph $H$,
the probability that
$P_{n,m} \textrm{ will have a component isomorphic to } H$
is bounded away from $0$ for sufficiently large $n$
if $\frac{m}{n}$ is bounded below by $b>1$ and above by $B<3$.
We shall then see (in Theorem~\ref{gen4}) that, in fact, 
the lower bound need only be $c>0$ if $H$ has at most one cycle. 
Finally,
we will discover
(in Corollary~\ref{tree6} via Theorem~\ref{tree4})
that the probability actually converges to $1$ if $H$ is a tree
and $\frac{m}{n} \in [c,1+o(1)]$~as~$n \to \infty$.

The proofs of our results will be based on counting:
we will construct graphs with components isomorphic to $H$ from other graphs in $\mathcal{P}(n,m)$
by deleting and inserting `suitable' edges in carefully chosen ways,
and we will then show that there isn't too much double-counting.
The properties shown in Sections~\ref{pen} and~\ref{add} will play a crucial role.

\begin{figure} [ht]
\setlength{\unitlength}{1cm}
\begin{picture}(20,4.75)(1.7,0)

\put(5.075,0.425){\circle{0.85}}
\put(6.825,0.425){\circle{0.85}}

\put(5.075,1.95){\vector(0,-1){1}}
\put(5.075,3.4){\vector(0,-1){1}}
\put(6.825,3.4){\vector(0,-1){2.45}}

\put(4.725,0.3){\textbf{T\ref{gen3}}}
\put(6.475,0.3){\textbf{T\ref{gen4}}}
\put(4.725,2.05){\textbf{C\ref{gen1}}}
\put(6.6,3.5){T\ref{pen5}}
\put(4.8,3.5){L\ref{gen2}}

\put(9.4,0.45){\circle{0.8}}
\put(7.85,2.05){\textbf{T\ref{tree4}}}
\put(8.85,3.5){Table 1}
\put(10.3,2.05){C\ref{add41}}
\put(9.05,0.325){\textbf{C\ref{tree6}}}

\put(7,3.4){\vector(1,-1){1}}
\put(9.4,3.4){\vector(-1,-1){1}}
\put(8.2,1.95){\vector(1,-1){1}}
\put(10.6,1.95){\vector(-1,-1){1}}
\end{picture}

\caption{The structure of Section~\ref{cptlow}.}
\end{figure}
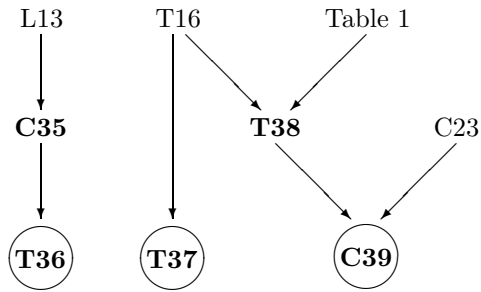

\begin{displaymath}
\end{displaymath} 

We start by noting a simple corollary to Lemma~\ref{gen2}:

\begin{Corollary} \label{gen1}
Let $b>1$ and $B<3$ be fixed constants and let $m(n) \! \in \! [bn,Bn]$~$\! \forall n$.
Then there exist constants $N(b,B)$ and $\delta (b,B) >0$ such that
\begin{displaymath}
\mathbf{P}[P_{n,m} \textrm{ will have $< \delta n$ edge-disjoint appearances of } K_{4}] < 
e^{- \delta n} \textrm{ } \forall n \geq N.
\end{displaymath}
\end{Corollary}
\textbf{Proof}
Any appearances of $K_{4}$ must be edge-disjoint, by $2$-edge-connectivity
(in fact, they must be vertex-disjoint).
Hence, the result follows from Lemma~\ref{gen2}.~
\begin{picture}(1,1)
\put(0,0){\line(1,0){1}}
\put(0,0){\line(0,1){1}}
\put(1,1){\line(-1,0){1}}
\put(1,1){\line(0,-1){1}}
\end{picture} \\

We will now use Corollary~\ref{gen1} to prove the first of our aforementioned theorems.
In fact,
we shall actually show a slightly stronger version than that advertised,
involving several components at once:

\begin{Theorem} \label{gen3}
Let $b>1$ and $B<3$ be fixed constants and let $m(n) \in [bn,Bn]$ for all large $n$.
Let $H_{1},H_{2}, \ldots, H_{l}$ 
be (fixed) connected planar graphs and let $t$ be a fixed constant.
Then there exist constants $\epsilon > 0$ and $N$ such that
\begin{eqnarray*}
& & \mathbf{P} 
\Big[ \bigcap_{i=1}^{l} 
(
\textrm{$P_{n,m}$
will have $\geq t$ components with} \\
& & \phantom{wwwww}\textrm{an order-preserving isomorphism to $H_{i}$} ) ] 
\geq \epsilon \textrm{ } \forall n \geq N.  
\end{eqnarray*} 
\end{Theorem}
\textbf{Sketch of Proof}
By symmetry,
it suffices to prove the result without the order-preserving condition.
The proof is then by induction on $l$.
We shall suppose that the result is true for $l=j$, but false for $l=j+1$.
Thus, using Corollary~\ref{gen1},
there must be a decent proportion of graphs in $\mathcal{P}(n,m)$ that have
(a) $\geq t$ components isomorphic to $H_{i}$ $\forall i \leq j$, 
(b) $<t$ components isomorphic to $H_{j+1}$,
and (c) many edge-disjoint appearances of $K_{4}$.

For each such graph, 
we can delete edges from some of these appearances of~$K_{4}$ to create isolated vertices,
on which we may then construct $t$ components isomorphic to $H_{j+1}$.
By inserting extra edges in appropriate places elsewhere,
we may hence obtain graphs in $\mathcal{P}(n,m)$.
The fact that the original graphs contained 
few components isomorphic to $H_{j+1}$
can then be used to show that there isn't too much double-counting,
and so we find that we have actually constructed a decent proportion of \textit{distinct} graphs in $\mathcal{P}(n,m)$.

By carefully selecting where we delete/insert edges,
we may ensure that these new graphs still have $\geq t$ components 
isomorphic to $H_{i}$~$\forall i \leq j$,
and so we are done. \\
\\
\textbf{Full Proof} 
As mentioned, 
it suffices to prove the result without the order-preserving condition,
since given any collection of $tl$ components such that exactly $t$ are isomorphic to $H_{i}$~$\forall i$,
the probability that they are all order-preserving is
$\prod_{i=1}^{l} \left( \left( \frac{| \textrm{\footnotesize{Aut}}(H_{i}) |}{|H_{i}|!} \right)^{t} \right)$.

We shall now prove the simplified version of the result by induction on $l$.
Suppose it is true for $l=j$,
i.e.~$\exists \epsilon_{j} >0$ and $\exists N_{j}$ such that
\begin{displaymath}
\mathbf{P}
\left[ \bigcap_{i=1}^{j} 
\left(P_{n,m} \textrm{ will have $\geq t$ components isomorphic to } H_{i} \right) \right]
\geq \epsilon_{j} \textrm{ } \forall n \geq N_{j}.
\end{displaymath} 
Let $\mathcal{L}_{n,r}$ denote the set of graphs in $\mathcal{P}(n,m)$ that have $\geq t$ components
isomorphic to $H_{i}$ $\forall i \leq r$.
Then $|\mathcal{L}_{n,j}| \geq \epsilon_{j} |\mathcal{P}(n,m)|$ $\forall n \geq N_{j}$. 

We have $\frac{m}{n} \in [b,B]$ $\forall$ large $n$.
Thus, by Corollary~\ref{gen1}, 
there are constants $\delta = \delta(b,B) >0$ and $N^{\prime}(b,B)$ such that, for all $n \geq N^{\prime}$,
\begin{displaymath}
\mathbf{P}[P_{n,m} \textrm{ will have at least $\delta n$ edge-disjoint appearances of } K_{4}] 
\geq \left(1- \frac{\epsilon_{j}}{3} \right).
\end{displaymath}
Let $\mathcal{I}_{n}$ denote the set of graphs in $\mathcal{P}(n,m)$ 
that have $\geq \delta n$ edge-disjoint appearances of $K_{4}$.
Then $|\mathcal{I}_{n} \cap \mathcal{L}_{n,j}| \geq 
\frac{2 \epsilon_{j}}{3} |\mathcal{P}(n,m)|$ $\forall n \geq N = \max \{N^{\prime}, N_{j} \}$. 

Consider an $n \geq N$ and suppose that 
$|\mathcal{L}_{n,j+1}| \leq \frac{\epsilon_{j}}{3} |\mathcal{P}(n,m)|$
(if not, then we are done).
Let $\mathcal{G}_{n,j}$ denote the set of graphs in $\mathcal{L}_{n,j}$ 
with (i) $<t$ components isomorphic to $H_{j+1}$ 
and (ii) at least $\delta n$ edge-disjoint appearances of~$K_{4}$.
Then, under our assumption,
we have $|\mathcal{G}_{n,j}| \geq \frac{\epsilon_{j}}{3} |\mathcal{P}(n,m)|$.
We shall use $\mathcal{G}_{n,j}$ to construct graphs in $\mathcal{L}_{n,j}$
with $\geq t$ components isomorphic to $H_{j+1}$. \\

Consider a graph $G \in \mathcal{G}_{n,j}$.
Then $G$ contains a subgraph $F$ consisting of $jt$ components such that,
for all $i \leq j$,
exactly $t$ of these components isomorphic to $H_{i}$.
Also, $G$ has a set of at least $\delta n$ edge-disjoint appearances of $K_{4}$.
Note that at least $\delta n - t \sum_{i=1}^{j} \left \lfloor \frac{e(H_{i})}{6} \right \rfloor$
of these edge-disjoint appearances of $K_{4}$ must lie in $G \setminus F$.
Let $s= t \sum_{i=1}^{j} \left \lfloor \frac{e(H_{i})}{6} \right \rfloor$,
let $H=H_{j+1}$,
and let $k=|H|$.
We may assume that $n$ is large enough that $\delta n - s \geq tk$.
Thus, we may choose $tk$ of the edge-disjoint appearances of $K_{4}$ in $G \setminus F$
$\left( \textrm{at least} \left( ^{\lceil \delta n \rceil -s} _{\phantom{wq}tk} \right) \textrm{ choices} \right)$,
and for each of these chosen appearances we may choose a `special' vertex in the $K_{4}$ that is not the 
root~$\left( 3^{tk} \textrm{ choices} \right)$.
Let us then delete all $3tk$ edges that are incident to the `special' vertices and 
insert edges between these $tk$ newly isolated vertices in such a way that they now form $t$ components
isomorphic to~$H$ 
$\left( \textrm{at least } \left( ^{\phantom{qi}tk} _{k, \ldots, k} \right) \frac{1}{t!} \textrm{ choices} \right)$.
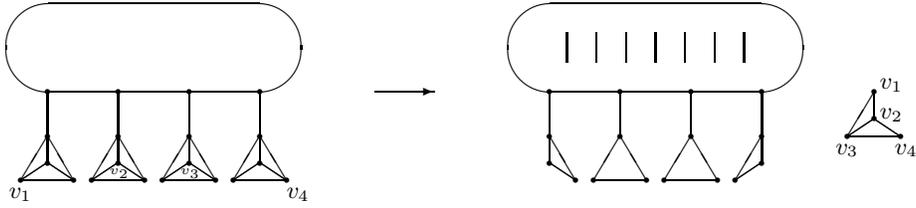
\begin{figure} [ht]
\setlength{\unitlength}{0.785cm}
\begin{picture}(20,3.85)(-0.25,-0.35)

\put(0,0){\line(1,0){0.9}}
\put(0,0){\line(3,2){0.45}}
\put(0,0){\line(3,5){0.45}}
\put(0.45,0.3){\line(0,1){0.45}}
\put(0.9,0){\line(-3,2){0.45}}
\put(0.9,0){\line(-3,5){0.45}}
\put(0,0){\circle*{0.1}}
\put(0.9,0){\circle*{0.1}}
\put(0.45,0.3){\circle*{0.1}}
\put(0.45,0.75){\circle*{0.1}}

\put(1.2,0){\line(1,0){0.9}}
\put(1.2,0){\line(3,2){0.45}}
\put(1.2,0){\line(3,5){0.45}}
\put(1.65,0.3){\line(0,1){0.45}}
\put(2.1,0){\line(-3,2){0.45}}
\put(2.1,0){\line(-3,5){0.45}}
\put(1.2,0){\circle*{0.1}}
\put(2.1,0){\circle*{0.1}}
\put(1.65,0.3){\circle*{0.1}}
\put(1.65,0.75){\circle*{0.1}}

\put(2.4,0){\line(1,0){0.9}}
\put(2.4,0){\line(3,2){0.45}}
\put(2.4,0){\line(3,5){0.45}}
\put(2.85,0.3){\line(0,1){0.45}}
\put(3.3,0){\line(-3,2){0.45}}
\put(3.3,0){\line(-3,5){0.45}}
\put(2.4,0){\circle*{0.1}}
\put(3.3,0){\circle*{0.1}}
\put(2.85,0.3){\circle*{0.1}}
\put(2.85,0.75){\circle*{0.1}}

\put(3.6,0){\line(1,0){0.9}}
\put(3.6,0){\line(3,2){0.45}}
\put(3.6,0){\line(3,5){0.45}}
\put(4.05,0.3){\line(0,1){0.45}}
\put(4.5,0){\line(-3,2){0.45}}
\put(4.5,0){\line(-3,5){0.45}}
\put(3.6,0){\circle*{0.1}}
\put(4.5,0){\circle*{0.1}}
\put(4.05,0.3){\circle*{0.1}}
\put(4.05,0.75){\circle*{0.1}}

\put(2.25,2.25){\oval(5,1.5)}

\put(0.45,0.75){\line(0,1){0.75}}
\put(0.45,1.5){\circle*{0.1}}
\put(1.65,0.75){\line(0,1){0.75}}
\put(1.65,1.5){\circle*{0.1}}
\put(2.85,0.75){\line(0,1){0.75}}
\put(2.85,1.5){\circle*{0.1}}
\put(4.05,0.75){\line(0,1){0.75}}
\put(4.05,1.5){\circle*{0.1}}

\put(6,1.5){\vector(1,0){1}}

\put(8.95,0.3){\line(0,1){0.45}}
\put(9.4,0){\line(-3,2){0.45}}
\put(9.4,0){\line(-3,5){0.45}}
\put(9.4,0){\circle*{0.1}}
\put(8.95,0.3){\circle*{0.1}}
\put(8.95,0.75){\circle*{0.1}}

\put(9.7,0){\line(1,0){0.9}}
\put(9.7,0){\line(3,5){0.45}}
\put(10.6,0){\line(-3,5){0.45}}
\put(9.7,0){\circle*{0.1}}
\put(10.6,0){\circle*{0.1}}
\put(10.15,0.75){\circle*{0.1}}

\put(10.9,0){\line(1,0){0.9}}
\put(10.9,0){\line(3,5){0.45}}
\put(11.8,0){\line(-3,5){0.45}}
\put(10.9,0){\circle*{0.1}}
\put(11.8,0){\circle*{0.1}}
\put(11.35,0.75){\circle*{0.1}}

\put(12.1,0){\line(3,2){0.45}}
\put(12.1,0){\line(3,5){0.45}}
\put(12.55,0.3){\line(0,1){0.45}}
\put(12.1,0){\circle*{0.1}}
\put(12.55,0.3){\circle*{0.1}}
\put(12.55,0.75){\circle*{0.1}}

\put(10.75,2.25){\oval(5,1.5)}

\put(8.95,0.75){\line(0,1){0.75}}
\put(8.95,1.5){\circle*{0.1}}
\put(10.15,0.75){\line(0,1){0.75}}
\put(10.15,1.5){\circle*{0.1}}
\put(11.35,0.75){\line(0,1){0.75}}
\put(11.35,1.5){\circle*{0.1}}
\put(12.55,0.75){\line(0,1){0.75}}
\put(12.55,1.5){\circle*{0.1}}

\put(14,0.75){\line(1,0){0.9}}
\put(14,0.75){\line(3,2){0.45}}
\put(14,0.75){\line(3,5){0.45}}
\put(14.45,1.05){\line(0,1){0.45}}
\put(14.9,0.75){\line(-3,2){0.45}}
\put(14,0.75){\circle*{0.1}}
\put(14.9,0.75){\circle*{0.1}}
\put(14.45,1.05){\circle*{0.1}}
\put(14.45,1.5){\circle*{0.1}}

\put(9.25,2){\line(0,1){0.5}}
\put(9.75,2){\line(0,1){0.5}}
\put(10.25,2){\line(0,1){0.5}}
\put(10.75,2){\line(0,1){0.5}}
\put(11.25,2){\line(0,1){0.5}}
\put(11.75,2){\line(0,1){0.5}}
\put(12.25,2){\line(0,1){0.5}}

\put(13.8,0.45){\small{$v_{3}$}}
\put(14.8,0.45){\small{$v_{4}$}}
\put(14.55,1.05){\small{$v_{2}$}}
\put(14.55,1.55){\small{$v_{1}$}}

\put(-0.2,-0.3){\small{$v_{1}$}}
\put(1.5,0.1){\tiny{$v_{2}$}}
\put(2.7,0.1){\tiny{$v_{3}$}}
\put(4.5,-0.3){\small{$v_{4}$}}

\end{picture}
\caption{Constructing a component isomorphic to $H$.}

\end{figure}

To maintain the correct number of edges, we should insert $t(3k-e(H))$ extra ones somewhere into the graph, 
making sure that we maintain planarity.
We will do this in such a way that we do not interfere with our new components,
with the chosen appearances of $K_{4}$
(which are now appearances of $K_{3}$),
or with $F$.
Thus, the part of the graph where we wish to insert edges contains $n-4tk-|F|$ vertices and $m - 7tk - e(F)$ edges.
We know that there exists a triangulation on these vertices containing these edges,
and clearly inserting an edge from this triangulation would not violate planarity.
Thus, we have at least 
$\left( ^{3(n-4tk-|F|)-6-(m -7tk-e(F))} _{\phantom{wwwwwwww} t(3k-l)} \right)$
choices for where to add the edges, where~$l=~e(H)$.

Therefore, we have at least
\begin{eqnarray*}
&& |\mathcal{G}_{n,j}|
\left( ^{\lceil \delta n \rceil -s} _{\phantom{wq}tk} \right)
3^{tk}
\left( ^{\phantom{qi}tk} _{k, \ldots, k} \right)
\frac{1}{t!}
\left( \! \! \begin{array}{c} 3(n-4tk-|F|)-6-(m -7tk-e(F)) \\ t(3k-l) \end{array} \! \! \right) \\
& = & |\mathcal{G}_{n,j}|
\Theta \left( n^{t(4k-l)} \right)
\phantom{w}
\textrm{ (recalling that $\frac{m}{n}$ is bounded away from $3$)}
\end{eqnarray*}
ways to build (not necessarily distinct) graphs in $\mathcal{L}_{n,j}$ 
that have at least $t$ components isomorphic to $H$. \\

We will now consider the amount of double-counting: \\
Each of our constructed graphs will contain at most $t(4k-l+2)-1$ components isomorphic to $H$ 
(since there were at most $t-1$ already;
we have deliberately added $t$; and 
we may have created at most one extra one each time we cut a `special' vertex away from its $K_{4}$ or 
added an edge in the rest of the graph).
Hence, we have at most $\left( ^{t(4k-l+2)-1} _{\phantom{wwww}t} \right)$ possibilities 
for which were our $tk$ `special' vertices.
Since appearances of $K_{3}$ must be vertex-disjoint, by $2$-edge-connectedness, 
we have at most $\frac{n}{3}$ of them and hence
at most 
$\left( \frac{n}{3} \right)^{tk}$ 
possibilities
for where the `special' vertices were originally.
There are then at most 
$\left(^{m -tl-4tk}_{\phantom{w}t(3k-l)}\right)$ 
possibilities for which edges were added
in the rest of the graph 
(i.e. away from the constructed components isomorphic to $H$ and these appearances of $K_{3}$). 
Thus, the amount of double-counting is at most 
$\left( ^{t(4k-l+2)-1} _{\phantom{wwww}t} \right)
\left( \frac{n}{3} \right)^{tk}
\left(^{m -tl-4tk}_{\phantom{w}t(3k-l)}\right) 
= \Theta \left( n^{t(4k-l)} \right)$,
recalling that $m = \Theta (n)$. \\

Hence, under our assumption that
$|\mathcal{G}_{n,j}| \geq \frac{\epsilon_{j} |\mathcal{P}(n,m)|}{3}$,
we find that the number of graphs in $\mathcal{L}_{n,j}$ 
that have $\geq t$ components isomorphic to $H$ is at least~$\Theta (|\mathcal{P}(n,m)|)$. 

But the set of graphs in $\mathcal{L}_{n,j}$ that have $\geq t$ components isomorphic to $H$ 
is exactly $\mathcal{L}_{n,j+1}$.
Thus, by assuming that 
$|\mathcal{L}_{n,j+1}| \leq \frac{\epsilon_{j}}{3} |\mathcal{P}(n,m)|$,
we have proved $\exists \zeta >0$ such that 
$\frac{|\mathcal{L}_{n,j+1}|}{|\mathcal{P}(n,m)|} \geq \zeta$.
Therefore, 
$\frac{|\mathcal{L}_{n,j+1}|}{|\mathcal{P}(n,m)|} \geq \epsilon_{j+1} =\min \{ \zeta, \frac{\epsilon_{j}}{3} \} >0$.

Hence, we are done, by induction.
\phantom{qwerty}
\begin{picture}(1,1)
\put(0,0){\line(1,0){1}}
\put(0,0){\line(0,1){1}}
\put(1,1){\line(-1,0){1}}
\put(1,1){\line(0,-1){1}}
\end{picture}

Note that in the previous proof,
we could have constructed components isomorphic to $H$ directly from appearances of $H$.
We chose to instead build the components from isolated vertices cut from appearances of $K_{4}$,
as this technique generalises more easily to our next proof,
as we shall now explain.

Recall that when we cut the isolated vertices from the appearances of $K_{4}$,
this involved deleting three edges for each isolated vertex that we created,
which crucially meant that we had enough edges to play with when we wanted to turn these isolated vertices
into components isomorphic to $H$.
Notice, though, that the proof was only made possible by the fact that we had lots of appearances of~$K_{4}$
to choose from,
which was why we needed to restrict $\frac{m}{n}$ to the region~$[b,B]$, where $b>1$ and $B<3$.

However, if $e(H) \leq |H|$ then we would have enough edges to play with 
even if we only deleted one edge for each isolated vertex that we created.
Thus, we may replace the role of the appearances of $K_{4}$ by pendant edges,
which we know are plentiful even for small values of $\frac{m}{n}$,
by Theorem~\ref{pen5}.
Hence, we may obtain:

\begin{Theorem} \label{gen4}
Let $c>0$ and $B<3$ be fixed constants and let $m(n) \in [cn,Bn]$ for all large $n$.
Let $H_{1},H_{2}, \ldots, H_{l}$ 
be (fixed) connected planar graphs with at most one cycle each and let $t$ be a fixed constant.
Then there exist constants $\epsilon > 0$ and $N$ such that
\begin{eqnarray*}
& & \mathbf{P} 
\Big[ \bigcap_{i=1}^{l} 
(
\textrm{$P_{n,m}$
will have $\geq t$ components with} \\
& & \phantom{wwwww}\textrm{an order-preserving isomorphism to $H_{i}$} ) ] 
\geq \epsilon \textrm{ } \forall n \geq N.
\end{eqnarray*} 
\end{Theorem}
\textbf{Proof}
As with Theorem~\ref{gen3},
it suffices to prove the result without the order-preserving condition,
and again the proof is by induction on $l$.
Suppose the result is true for $l=j$,
but false for $l=j+1$.
Then, similarly to with the proof of Theorem~\ref{gen3}, 
we have a set $\mathcal{G}_{n,j}$ of at least $\frac{\epsilon_{j}}{3} |\mathcal{P}(n,m)|$ graphs with
(i) $\geq~t$ components isomorphic to $H_{i}$ $\forall i \leq j$,
(ii) $<t$ components isomorphic to $H_{j+1}$ 
and (iii) at least $\delta n$ pendant edges (using Theorem~\ref{pen5}).

Given a graph $G \in \mathcal{G}_{n,j}$,
we may delete $t|H_{j+1}|$ of the pendant edges
(taking care not to interfere with our components isomorphic to $H_{i}$ for $i \leq j$)
and use the resulting isolated vertices to construct $t$ components isomorphic to $H_{j+1}$.
If $H_{j+1}$ is a tree, then we should also add $t$ edges in suitable places somewhere in the rest of the graph.

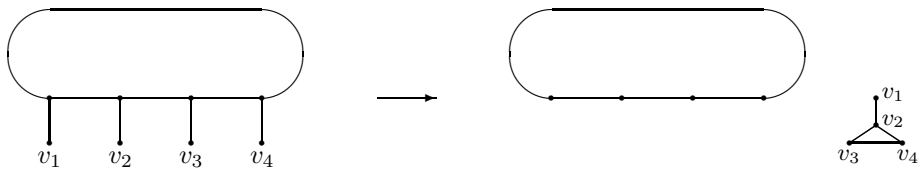
\begin{figure} [ht]
\setlength{\unitlength}{0.785cm}
\begin{picture}(20,2.75)(-0.25,0.5)

\put(0.45,0.75){\circle*{0.1}}
\put(1.65,0.75){\circle*{0.1}}
\put(2.85,0.75){\circle*{0.1}}
\put(4.05,0.75){\circle*{0.1}}

\put(2.25,2.25){\oval(5,1.5)}

\put(0.45,0.75){\line(0,1){0.75}}
\put(0.45,1.5){\circle*{0.1}}
\put(1.65,0.75){\line(0,1){0.75}}
\put(1.65,1.5){\circle*{0.1}}
\put(2.85,0.75){\line(0,1){0.75}}
\put(2.85,1.5){\circle*{0.1}}
\put(4.05,0.75){\line(0,1){0.75}}
\put(4.05,1.5){\circle*{0.1}}

\put(6,1.5){\vector(1,0){1}}

\put(10.75,2.25){\oval(5,1.5)}

\put(14,0.75){\line(1,0){0.9}}
\put(14,0.75){\line(3,2){0.45}}
\put(14.45,1.05){\line(0,1){0.45}}
\put(14.9,0.75){\line(-3,2){0.45}}
\put(14,0.75){\circle*{0.1}}
\put(14.9,0.75){\circle*{0.1}}
\put(14.45,1.05){\circle*{0.1}}
\put(14.45,1.5){\circle*{0.1}}

\put(13.8,0.45){\small{$v_{3}$}}
\put(14.8,0.45){\small{$v_{4}$}}
\put(14.55,1.05){\small{$v_{2}$}}
\put(14.55,1.5){\small{$v_{1}$}}

\put(8.95,1.5){\circle*{0.1}}
\put(10.15,1.5){\circle*{0.1}}
\put(11.35,1.5){\circle*{0.1}}
\put(12.55,1.5){\circle*{0.1}}

\put(0.25,0.4){$v_{1}$}
\put(1.45,0.4){$v_{2}$}
\put(2.65,0.4){$v_{3}$}
\put(3.85,0.4){$v_{4}$}

\end{picture}
\caption{Constructing a component isomorphic to $H_{j+1}$.}

\end{figure}

By similar counting arguments to those used in the proof of Theorem~\ref{gen3},
we achieve our result. 
\phantom{qwerty}
\begin{picture}(1,1)
\put(0,0){\line(1,0){1}}
\put(0,0){\line(0,1){1}}
\put(1,1){\line(-1,0){1}}
\put(1,1){\line(0,-1){1}}
\end{picture} \\
\\
\\

We shall now finish this section by looking specifically at the case when $H$ is a tree.
We already know from Theorem~\ref{gen4} that 
the probability that $P_{n,m}$ will contain a component isomorphic to a given fixed tree 
is certainly bounded away from $0$ for large $n$ 
if $\frac{m}{n}$ is bounded below by $c > 0$ and above by $B < 3$,
We will now see (in Corollary~\ref{tree6}) that the limiting probability is, in fact, $1$ 
if $\frac{m}{n} \in [c,1+o(1)]$ as $n \to \infty$.
Note that this result can actually be shown by exactly the same proof as for Theorem~\ref{gen4},
using the additional ingredient that add$(n,m) = \omega(n)$ if $\frac{m}{n} \leq 1 + o(1)$
(from Corollary~\ref{add41}).
However,
we shall instead aim to give (in Theorem~\ref{tree4}) 
a more detailed account of the number of components isomorphic to $H$,
which will include Corollary~\ref{tree6}.
The proof will again involve deleting pendant edges and inserting edges elsewhere,
but the calculations will now be more complicated:

\begin{Theorem} \label{tree4}
Let $H$ be a (fixed) tree, let $c>0$ be a fixed constant
and let $m=m(n) \in [cn, (1+o(1))n]$.
Then $\exists \lambda (H,c) > 0$ such that
\begin{eqnarray*}
& & \mathbf{P}\Big[
\textrm{$P_{n,m}$
will have \footnotesize{$< \left \lceil \frac{\lambda \emph{\tiny{add}}(n,m)}{n} \right \rceil$}} 
\textrm{ components with} \\
& & \phantom{ww}\textrm{an order-preserving isomorphism to $H$}\Big] 
< e^{- \left \lceil \frac{\lambda \emph{\tiny{add}}(n,m)}{n} \right \rceil} \textrm{ for all large $n$}.
\end{eqnarray*}
\end{Theorem}
\textbf{Sketch of Proof}
We suppose that the result is false.
Thus, using Theorem~\ref{pen5},
there must be a decent proportion of graphs in $\mathcal{P}(n,m)$ that have
(i) `few' components with an order-preserving isomorphism to $H$
and (ii) many pendant edges.
For each such graph, 
we can delete some of these pendant edges and
use the resulting isolated vertices to construct components with an order-preserving isomorphism to $H$.
By inserting extra edges elsewhere,
we may thus construct lots of graphs in $\mathcal{P}(n,m)$.
The fact that the original graphs contained 
few components with order-preserving isomorphism to $H$
can then be used to show that there is not much double-counting,
and so we find that we have actually constructed more than $|\mathcal{P}(n,m)|$ distinct graphs in $\mathcal{P}(n,m)$,
which is a contradiction. \\
\\
\textbf{Full Proof} 
By Theorem~\ref{pen5}, there exist constants $\alpha >0$ and $n_{0}$ such that 
\begin{displaymath}
\mathbf{P}[P_{n,m} \textrm{ will have } < \alpha n \textrm{ pendant edges}] < e^{- \alpha n} 
\textrm{ } \forall n \geq n_{0}.
\end{displaymath}
Let $\lambda$ be a small positive constant (whose value we shall choose later),
let $t=~t(n)= \left \lceil \frac{\lambda \textrm{\scriptsize{add}}(n,m)}{n} \right \rceil$
and suppose $\exists n \geq n_{0}$ such that 
\begin{eqnarray*}
\mathbf{P}[P_{n,m} 
\textrm{ will have $<t$ 
components with an order-preserving isomorphism to } H] \\
\geq e^{-t}.
\end{eqnarray*}
Then there is a set $\mathcal{G}_{n}$ of at least a proportion 
$e^{-t} - e^{- \alpha n}$ of the graphs in 
$\mathcal{P}(n,m)$ 
with (i) $<t$ components with an order-preserving isomorphism to $H$ and (ii) at least $\alpha n$ pendant edges.

Let $k=|H|$.
Since add$(n,m) \leq \left( ^{n}_{2} \right)$,
we may assume that $n$ is large enough (and $\lambda$ small enough)
that $\alpha n \geq tk$ and 
$e^{-t} - e^{- \alpha n} 
\geq \frac{1}{2}e^{-t}$. 

\phantom{p}

To build graphs with $\geq t$ components with an order-preserving isomorphism to $H$,
one can start with a graph $G \in \mathcal{G}_{n}$ ($|\mathcal{G}_{n}|$ choices), 
delete $tk$ of the pendant edges 
$\left( \textrm{ at least } \left(^{\lceil \alpha n \rceil} _{\phantom{ii}tk} \right) \textrm{ choices } \right)$,
and insert edges between $tk$ of the newly-isolated vertices 
(choosing one from each pendant edge) 
in such a way that they now form $t$ components, 
each with an order-preserving isomorphism to $H$
$\left( \textrm{at least } \left( ^{\phantom{qi}tk} _{k, \ldots, k} \right) \frac{1}{t!} \textrm{ choices} \right)$.
We should then add $t$ edges somewhere in the rest of the graph 
(i.e.~away from our newly constructed components)
to maintain the correct number of edges overall
(at least 
$\prod_{i=0}^{t-1} \textrm{add}(n-tk,m-tk+i) \geq 
(\textrm{add}(n-tk,m-tk+t-1))^{t}$
choices).

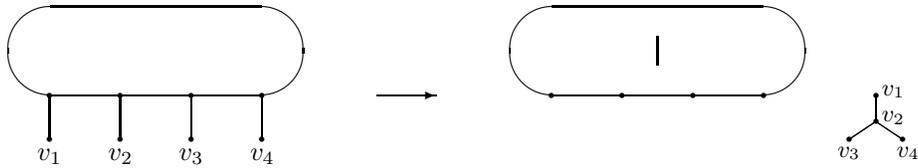
\begin{figure} [ht]
\setlength{\unitlength}{0.785cm}
\begin{picture}(20,3)(-0.25,0.5)

\put(0.45,0.75){\circle*{0.1}}
\put(1.65,0.75){\circle*{0.1}}
\put(2.85,0.75){\circle*{0.1}}
\put(4.05,0.75){\circle*{0.1}}

\put(2.25,2.25){\oval(5,1.5)}

\put(0.45,0.75){\line(0,1){0.75}}
\put(0.45,1.5){\circle*{0.1}}
\put(1.65,0.75){\line(0,1){0.75}}
\put(1.65,1.5){\circle*{0.1}}
\put(2.85,0.75){\line(0,1){0.75}}
\put(2.85,1.5){\circle*{0.1}}
\put(4.05,0.75){\line(0,1){0.75}}
\put(4.05,1.5){\circle*{0.1}}

\put(6,1.5){\vector(1,0){1}}

\put(10.75,2.25){\oval(5,1.5)}

\put(14,0.75){\line(3,2){0.45}}
\put(14.45,1.05){\line(0,1){0.45}}
\put(14.9,0.75){\line(-3,2){0.45}}
\put(14,0.75){\circle*{0.1}}
\put(14.9,0.75){\circle*{0.1}}
\put(14.45,1.05){\circle*{0.1}}
\put(14.45,1.5){\circle*{0.1}}

\put(13.8,0.45){\small{$v_{3}$}}
\put(14.8,0.45){\small{$v_{4}$}}
\put(14.55,1.05){\small{$v_{2}$}}
\put(14.55,1.5){\small{$v_{1}$}}

\put(8.95,1.5){\circle*{0.1}}
\put(10.15,1.5){\circle*{0.1}}
\put(11.35,1.5){\circle*{0.1}}
\put(12.55,1.5){\circle*{0.1}}

\put(0.25,0.4){$v_{1}$}
\put(1.45,0.4){$v_{2}$}
\put(2.65,0.4){$v_{3}$}
\put(3.85,0.4){$v_{4}$}

\put(10.75,2){\line(0,1){0.5}}

\end{picture}
\caption{Constructing a component isomorphic to $H$.}

\end{figure}

Let us now consider the amount of double-counting: \\
Each of our constructed graphs will contain at most $t(k+3)-1$ components with an order-preserving isomorphism to $H$
(since there were at most $t-1$ already in $G$; 
we have deliberately added~$t$;
and we may have created at most one extra one each time 
we deleted a pendant edge or added an edge in the rest of the graph),
so we have at most 
$\left(^{t(k+3)-1}_{\phantom{wwp}t}\right) \leq \frac{1}{t!} (t(k+3))^{t}$ 
possibilities for which are our created components.
We then have at most $n^{tk}$ 
possibilities for where the vertices in our created components were attached originally 
and at most $\left(^{m}_{t}\right) \leq (3n)^{t}$ possibilities for which edges was added.

Hence, putting everything together,
we find that the number of \textit{distinct} graphs in $\mathcal{P}(n,m)$ that have $\geq t$ components 
with an order-preserving isomorphism to $H$ is at least 
\begin{eqnarray*}
& & \frac{ \left( ^{\phantom{wwqq} \lceil \alpha n \rceil}_{k, \ldots, k, \phantom{p} \lceil \alpha n \rceil -tk} \right)}
{n^{tk}}
\left( \frac{ (\textrm{add}(n-tk,m-tk+t-1))^{t}}{(3n)^{t}} \right)
\frac{|\mathcal{G}_{n}|}{(t(k+3))^{t}} \\
& \geq & 
\left( \frac{\alpha^{k}}{2^{k}k!3(k+3)} \cdot \frac{\textrm{add}(n-tk,m-tk+t-1)}{tn} \right)^{t} |\mathcal{G}_{n}|, 
\phantom{ww} \textrm{since we may} \\
& & \textrm{assume $\lambda$ is sufficiently small and $n$ sufficiently large 
that $\alpha n - tk \geq \frac{\alpha n}{2}$}.
\end{eqnarray*} 
\phantom{p}

The main thrust of the proof is now over,
and we are left with the fiddly task of evaluating
add$(n-tk,m-tk+t-1)$.
We shall see that it is at least $C\textrm{add}(n,m)$
(for an appropriate constant $C>0$,
independent of $\lambda$),
which is roughly $\frac{C}{\lambda}tn$.
Our desired contradiction will then follow by taking $\lambda$ to be sufficiently small. \\

Let $i=n-tk$ and let $s=t-1$.
Then
\begin{eqnarray}
\textrm{add}(n-tk,m-tk+t-1) 
& = & \textrm{add}(i,i+(m-n+s)). \label{eq:a1} 
\end{eqnarray}
Since add$(i,i+l)$ is quite sensitive to changes in $l$,
we shall need to investigate the value of $m-n+s$ in detail. 

Recall that $s=\left \lceil \frac{\lambda \textrm{\scriptsize{add}}(n,m)}{n} \right \rceil - 1$.
From the upper bounds of Section~\ref{add}, 
we know there exist constants $B>0$ and $N$ such that $\forall n \geq N$ we have
\begin{displaymath}
\textrm{add}(n,m) \leq
\left\{ \begin{array}{ll}
Bn(n-m) & \textrm{if $m-n<0$ and $n^{1/2} \leq n-m \leq n$} \\
Bn^{3/2} & \textrm{if $|m-n| \leq n^{1/2}$} \\
B \frac{n^{2}}{m-n} & \textrm{if $m-n>0$ and $n^{1/2} \leq m-n \leq \frac{n}{2}$}. 
\end{array} \right. 
\end{displaymath} 

Thus, since we may assume that $\lambda$ is sufficiently small that $\lambda B \leq 1/2$, 
$\forall n \geq N$ we have 
\begin{eqnarray}
s < \frac{\lambda \textrm{add}(n,m) }{n} \leq
\left\{ \begin{array}{ll}
\frac{n-m}{2} & \textrm{if $m-n<0$ and $n^{1/2} \leq |m-n| \leq n$} \\
\frac{n^{1/2}}{2} & \textrm{if $|m-n| \leq n^{1/2}$} \\
\frac{n}{2(m-n)} & \textrm{if $m-n>0$ and $n^{1/2} \leq |m-n| \leq \frac{n}{2}$}. \label{eq:a2}
\end{array} \right\} 
\end{eqnarray} 

Hence
(combining (\ref{eq:a2}) with the fact that we may assume that $\lambda$ is sufficiently small that $i \geq \frac{n}{2}$),
for all $n \geq N$ we have 
\begin{eqnarray}
\left. \begin{array}{lll}
m-n<0 \textrm{ and } n^{1/2} \leq |m-n| \leq n 
& \Rightarrow &
m-n+s < 0 \textrm{ and } \\
& & \frac{i^{1/2}}{2} \leq \frac{n^{1/2}}{2} \leq |m-n+s| \\
& & \leq n \leq 2i \\
|m-n| \leq n^{1/2} 
& \Rightarrow &
|m-n+s| \leq \frac{3n^{1/2}}{2} \leq 3i^{1/2} \\
m-n>0 \textrm{ and } n^{1/2} \leq |m-n| \leq \frac{n}{2} 
& \Rightarrow &
m-n+s>0 \textrm{ and } \\
& & i^{1/2} \leq n^{1/2} \leq |m-n+s| \\
& & \leq \frac{n}{2} + \frac{n^{1/2}}{2} \leq i + i^{1/2}. \label{eq:a3}
\end{array} \right \}
\end{eqnarray} 
\phantom{p}

Recall that we are interested in add$(i,i+(m-n+s))$,
where $i=n-tk$.
We know from the lower bounds of Section~\ref{add} 
that there exist constants $b>0$ and~$N_{2}$ such that $\forall i \geq N_{2}$
we have
\begin{eqnarray*}
\textrm{add}(i,i\!+\!(m\!-\!n\!+\!s)) 
& \geq &
\left\{ \begin{array}{ll}
bi|m-n+s| & \textrm{if $m-n+s<0$} \\
& \textrm{and $\frac{i^{1/2}}{2} \leq |m-n+s| \leq 2i$} \\
bi^{3/2} & \textrm{if $|m-n+s| \leq 3i^{1/2}$} \\
b \frac{i^{2}}{|m-n+s|} & \textrm{if $m-n+s>0$} \\
& \textrm{and $i^{1/2} \leq |m-n+s| \leq i + i^{1/2}$}. 
\end{array} \right. 
\end{eqnarray*}
\phantom{p}

Thus, combining this with $(\ref{eq:a1})$ and $(\ref{eq:a3})$
(and the fact that $i \geq \frac{n}{2}$),
for all large~$n$ we have
\begin{eqnarray*}
\textrm{add}(n-tk,m-tk+s) 
& \geq &
\left\{ \begin{array}{ll}
\frac{bn}{2}|m-n+s| & \textrm{if $m-n<0$} \\
& \textrm{and $n^{1/2} \leq |m-n| \leq n$} \\
\frac{bn^{3/2}}{2} & \textrm{if $|m-n| \leq n^{1/2}$} \\
\frac{b}{2} \frac{n^{2}}{|m-n+s|} 
& \textrm{if $m-n>0$} \\
& \textrm{and $n^{1/2} \leq |m-n| \leq \frac{n}{2}$}
\end{array} \right. \\
& \stackrel{\textrm{\small{by (\ref{eq:a2})}}}{\geq} &
\left\{ \begin{array}{ll}
\frac{bn}{4}|m-n| & \textrm{if $m-n<0$} \\
& \textrm{and $n^{1/2} \leq |m-n| \leq n$} \\
\frac{bn^{3/2}}{2} & \textrm{if $|m-n| \leq n^{1/2}$} \\
\frac{b}{3} \frac{n^{2}}{|m-n|} 
& \textrm{if $m-n>0$} \\
& \textrm{and $n^{1/2} \leq |m-n| \leq \frac{n}{2}$} \\
\end{array} \right. \\
& \geq & C \textrm{add}(n,m) 
\textrm{ if } m \leq \frac{3n}{2}, 
\textrm{ by Section~\ref{add}} \\
& & \textrm{(for some } C>0 \textrm{ and sufficiently large } n). 
\end{eqnarray*} 
\phantom{p}

Hence, continuing from where we left off,
the number of graphs in $\mathcal{P}(n,m)$ that have $\geq t$ components 
with an order-preserving isomorphism to $H$ is at least 
$
\left( \frac{\alpha^{k}}{2^{k}k!3(k+3)} \cdot \frac{C\textrm{\scriptsize{add}}(n,m)}{tn} \right)^{t} |\mathcal{G}_{n}| 
= \left( (1+o(1)) \frac{\alpha^{k}C}{2^{k}k!3(k+3)\lambda} \right)^{t} |\mathcal{G}_{n}|.
$ 
But this is more than $|\mathcal{P}(n,m)|$ for large $n$,
if $\lambda$ is sufficiently small, 
since we recall our assumption that
$|\mathcal{G}_{n}| \geq \frac{1}{2} e^{-t} |\mathcal{P}(n,m)|$. 

Thus, by proof by contradiction, it must be that
\begin{eqnarray*}
& & \mathbf{P}\Big[
\textrm{$P_{n,m}$
will have \footnotesize{$< \left \lceil \frac{\lambda \textrm{\tiny{add}}(n,m)}{n} \right \rceil$}} 
\textrm{ components with} \\
& & \phantom{wwwwww}\textrm{an order-preserving isomorphism to $H$}\Big] \\
& & \phantom{wwwwwwwwwwwwwww}
< e^{- \left \lceil \frac{\lambda \textrm{\tiny{add}}(n,m)}{n} \right \rceil} \textrm{ for all large $n$}.
\phantom{qwerty}
\setlength{\unitlength}{.25cm}
\begin{picture}(1,1)
\put(0,0){\line(1,0){1}}
\put(0,0){\line(0,1){1}}
\put(1,1){\line(-1,0){1}}
\put(1,1){\line(0,-1){1}}
\end{picture}
\end{eqnarray*} \\

Our aforementioned corollary now follows easily:

\begin{Corollary} \label{tree6}
Let $H$ be a (fixed) tree, let $t$ and $c>0$ be fixed constants and let $m=m(n) \in [cn, (1+o(1))n]$ as $n \to \infty$.
Then
\begin{eqnarray*}
\mathbf{P}[P_{n,m} \textrm{ will have } \geq t \textrm{ components with an order-preserving isomorphism to } H] \\
\to 1 \textrm{ as } n \to \infty.
\end{eqnarray*}
\end{Corollary}
\textbf{Proof}
By Corollary~\ref{add41},
we know add$(n,m)=\omega(n)$.
Thus, the result follows from Theorem~\ref{tree4}. 
$\phantom{qwerty}$\begin{picture}(1,1)
\put(0,0){\line(1,0){1}}
\put(0,0){\line(0,1){1}}
\put(1,1){\line(-1,0){1}}
\put(1,1){\line(0,-1){1}}
\end{picture}

\newpage
\section{Connectivity \& \boldmath{$\kappa(P_{n,m})$}} \label{conkappa}

In this section, 
we shall look indirectly at the chances of $P_{n,m}$ containing specific components 
by investigating the probability that it will be connected
(in which case it clearly won't contain any component of order $<n$).
This is also an interesting topic in its own right.

We already know (from Corollary~\ref{tree6}) that the probability that $P_{n,m}$ will be connected
converges to $0$ if $m \leq (1+o(1))n$,
and (from Theorem~\ref{gen4}) 
that it is bounded away from $1$ for sufficiently large $n$ if $\limsup_{n \to \infty} \frac{m}{n} < 3$.
Conversely, we will now see (in Theorem~\ref{conn4}) 
that the probability is bounded away from $0$
if~$\liminf_{n \to \infty} \frac{m}{n} >1$,
and (in Corollary~\ref{conn5}) that it converges to $1$ if $\frac{m}{n} \to 3$.
Note that we shall then have a complete description of
$\mathbf{P}[P_{n,m} \textrm{ will be connected}]$,
in terms of exactly when it is bounded away from $0$ and $1$.

The proofs will be in two parts, Lemmas~\ref{conn1} and~\ref{conn2},
dealing with the cases when $\frac{m}{n}$ is or isn't, respectively, bounded away from $3$.
As a tool for obtaining~Lemma~\ref{conn1},
we shall first work towards a result (Lemma~\ref{cpt12}) on $\kappa(P_{n,m})$,
the number of components in $P_{n,m}$.
Although not one of the main objectives of this thesis,
it is interesting to note that we may obtain quite a lot of information on $\kappa(P_{n,m})$,
and so we will explore this topic in more detail during the second half of this section
(in particular with Propositions~\ref{cpt2031} and~\ref{cpt31}).
We will then collect up all our results on $\kappa(P_{n,m})$,
including some lower bounds derived from~Section~\ref{cptlow},
on page~\pageref{cptsum}.

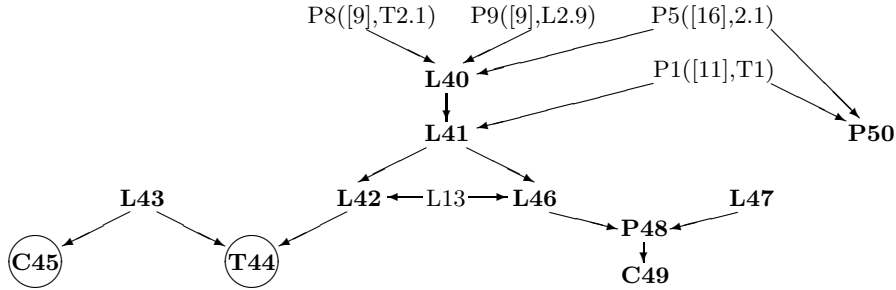
\begin{figure} [ht]
\setlength{\unitlength}{0.9cm}
\begin{picture}(20,3.5)(0,0.6)

\put(5.525,3.925){\vector(2,-1){1}}
\put(7.95,3.925){\vector(-2,-1){1}}
\put(6.7,2.975){\vector(0,-1){0.4}}
\put(9.75,3.925){\vector(-4,-1){2.6}}
\put(11.5,3.925){\vector(1,-1){1.3}}
\put(9.75,3.125){\vector(-4,-1){2.6}}
\put(11.5,3.125){\vector(2,-1){1.1}}

\put(4.65,4.025){\small{P\ref{ger T2.1}(\cite{ger},T2.1)}}
\put(7.05,4.025){\small{P\ref{ger L2.9}(\cite{ger},L2.9)}}
\put(6.375,3.075){\small{\textbf{L\ref{cpt11}}}}
\put(9.725,4.025){\small{P\ref{mcd 2.1}(\cite{mcd},2.1)}}
\put(9.75,3.225){\small{P\ref{gim T1}(\cite{gim},T1)}}
\put(6.375,2.275){\small{\textbf{L\ref{cpt12}}}}
\put(12.625,2.275){\small{\textbf{P\ref{cpt31}}}}

\put(6.4,1.325){\small{L\ref{gen2}}}
\put(5.075,1.325){\small{\textbf{L\ref{conn1}}}}
\put(1.875,1.325){\small{\textbf{L\ref{conn2}}}}
\put(3.47,0.375){\small{\textbf{T\ref{conn4}}}}
\put(0.27,0.375){\small{\textbf{C\ref{conn5}}}}

\put(3.825,0.5){\circle{0.8}}
\put(0.625,0.5){\circle{0.8}}

\put(6.35,1.45){\vector(-1,0){0.525}}
\put(6.4,2.175){\vector(-2,-1){1}}
\put(5.225,1.225){\vector(-2,-1){1}}
\put(2.025,1.225){\vector(-2,-1){1}}
\put(2.425,1.225){\vector(2,-1){1}}

\put(7.675,1.325){\small{\textbf{L\ref{cpt201}}}}
\put(10.875,1.325){\small{\textbf{L\ref{cpt202}}}}
\put(9.275,0.875){\small{\textbf{P\ref{cpt2031}}}}
\put(9.275,0.175){\small{\textbf{C\ref{cpt203}}}}

\put(7.025,1.45){\vector(1,0){0.575}}
\put(7,2.175){\vector(2,-1){1}}
\put(8.225,1.225){\vector(4,-1){1}}
\put(11,1.225){\vector(-4,-1){1}}
\put(9.6125,0.775){\vector(0,-1){0.3}}

\end{picture}

\caption{The structure of Section~\ref{conkappa}.}
\end{figure}

We start by obtaining a lower bound for
$\left( \frac{|\mathcal{P}(n,m)|}{n!} \right)^{1/n}$,
which we will use in the proof of Lemma~\ref{cpt12} on $\kappa(P_{n,m})$.

\begin{Lemma} \label{cpt11}
Let $m=m(n) \in [n,3n-6]$ for all large $n$.
Then
\begin{displaymath}
\liminf_{n \to \infty} \left( \frac{|\mathcal{P}(n,m)|}{n!} \right)^{1/n} \geq \gamma (1) = e.
\end{displaymath}
\end{Lemma}
\textbf{Proof}
The case when $\frac{m}{n}$ is bounded away from $1$ follows fairly simply from 
uniform convergence:

Let $\eta >0$.
Then, by Proposition~\ref{ger L2.9} with $a=\frac{5}{4}$ (for example), 
$\exists n_{0}$ such that 
$
\left( \frac{|\mathcal{P}(n,m)|}{n!} \right)^{1/n} > \gamma \left( \frac{m}{n} \right) - \eta \phantom{q}
\textrm{ } \forall m \in \left[ \frac{5n}{4},3n-6 \right] \textrm{ } \forall n \geq n_{0}.
$
Recall (from Proposition~\ref{ger T2.1}) that $\gamma \left( \frac{m}{n} \right) \geq \gamma (1)$.
Thus, since $\eta$ was an arbitrary positive constant, we are done. 

We shall now consider the case $m \in [n,\frac{5n}{4}]$:

Let $m<3n-6$.
Given a graph, $G$, in $\mathcal{P}(n,m)$,
we may construct a graph in~$\mathcal{P}(n,m+1)$ by inserting an extra edge into $G$
in such a way that we maintain planarity.
Since there exists a (planar) triangulation of order $n$ that contains~$G$ as a subgraph,
we have at least $3n-6-m$ choices for the edge to insert.
Thus, taking any possible double-counting into account,
we have $|\mathcal{P}(n,m+~1)| \geq~\frac{3n-6-m}{m+1} |\mathcal{P}(n,m)|$.
Hence, if $m \leq \frac{3n-7}{2}$ then $|\mathcal{P}(n,m)| \leq |\mathcal{P}(n,m+1)|$.
Thus, since $\frac{5n}{4} \leq \frac{3n-7}{2}$ for sufficiently large $n$,
it suffices to consider the case when $m=n$.
But this holds by Proposition~\ref{ger T2.1}. 

Thus, the result holds for all $m \in [n,3n-6]$,
by considering the cases \mbox{$m \in \left[ \frac{5n}{4}, 3n-6 \right]$} and 
$m \in \left[ n, \frac{5n}{4} \right]$ separately.
\phantom{qwerty}
\setlength{\unitlength}{.25cm}
\begin{picture}(1,1)
\put(0,0){\line(1,0){1}}
\put(0,0){\line(0,1){1}}
\put(1,1){\line(-1,0){1}}
\put(1,1){\line(0,-1){1}}
\end{picture} \\
\\

We are now ready to obtain our aforementioned first result on the number of components in $P_{n,m}$,
which will turn out to be useful to us later when investigating
$\mathbf{P}[P_{n,m} \textrm{ will be connected}]$.
We follow the method of proof of Proposition~\ref{ger L2.6}~(\cite{ger}, Lemma 2.6),
which dealt with the case $m = \lfloor qn \rfloor$ for fixed $q \in [1,3)$.

\begin{Lemma} \label{cpt12}
Let $m=m(n) \in [n,3n-6]$ for all large $n$ and let the constant $c$ satisfy $c> \ln \gamma_{l} - 1$.
Then
\begin{displaymath}
\mathbf{P} \left[ \kappa (P_{n,m}) > \left \lceil \frac{cn}{\ln n} \right \rceil \right] = e^{- \Omega(n)}.
\end{displaymath}
\end{Lemma}
\textbf{Proof}
Let $k=k(n)=\lceil \frac{cn}{\ln n} \rceil$.
Then we have 
$|G \in \mathcal{P}(n,m): \kappa(G)>k| \leq |G \in~\mathcal{P}(n): \kappa(G)>k| \leq \frac{|\mathcal{P}(n)|}{k!}$, 
using Proposition~\ref{mcd 2.1}.
Hence,
it must be that
$
\left(\frac{|G \in \mathcal{P}(n,m):\kappa(G)>k|}{|\mathcal{P}(n,m)|} \right)^{\frac{1}{n}} \leq
\left( \frac{1}{k!} \frac{|\mathcal{P}(n)|}{|\mathcal{P}(n,m)|} \right)^{\frac{1}{n}}. 
$

As $n \to \infty$,
recall that
we have $\left( \frac{|\mathcal{P}(n)|}{n!} \right)^{\frac{1}{n}} \to \gamma_{l}$ (by Proposition~\ref{gim T1})
and $\liminf \left( \frac{|\mathcal{P}(n,m)|}{n!} \right)^{\frac{1}{n}} \!\geq\! e$ (by Lemma~\ref{cpt11}).
Thus, 
$\limsup_{n \to \infty} \left( \frac{|\mathcal{P}(n)|}{|\mathcal{P}(n,m)|} \right)^{\frac{1}{n}} 
\!\leq\! e^{-1} \gamma_{l}$. 

By Stirling's formula,
$k! \sim \sqrt{2 \pi k} k^{k} e^{-k}$ as $n$ (and hence $k$) $\to \infty$.
Thus,
\begin{eqnarray*}
(k!)^{\frac{1}{n}} & \sim & k^{\frac{k}{n}} e^{-\frac{k}{n}} \phantom{www} \textrm{ as } n \to \infty \\
& \sim & k^{\frac{k}{n}}, \phantom{www} \textrm{ since } \frac{k}{n} \to 0 \textrm{ as } n \to \infty \\
& \sim & \left( \frac{cn}{\ln n} \right) ^{\frac{c}{\ln n}}, \phantom{www} \textrm{ by definition of } k \\
& \sim & n^{\frac{c}{\ln n}}, \phantom{www} \textrm{ since } 
\left( \frac{c}{\ln n} \right) ^{\frac{c}{\ln n}} \to 1 \textrm { as } \frac{c}{\ln n} \to 0 \\
& = & e^{\frac{c}{\ln n } \ln n} \\
& = & e^{c}.
\end{eqnarray*}

Therefore,
$
\limsup_{n \to \infty} \left( \frac{1}{k!} \frac{|\mathcal{P}(n)|}{|\mathcal{P}(n,m)|} \right)^{\frac{1}{n}} \leq 
e^{-(1+c)} \gamma_{l},
$
and the theorem follows from the fact that $c > \ln \gamma_{l} - 1$.
$\phantom{qwerty}$
\setlength{\unitlength}{.25cm}
\begin{picture}(1,1)
\put(0,0){\line(1,0){1}}
\put(0,0){\line(0,1){1}}
\put(1,1){\line(-1,0){1}}
\put(1,1){\line(0,-1){1}}
\end{picture} \\
\\

We will now use Lemma~\ref{cpt12} to work towards Theorem~\ref{conn4},
where we shall show that 
$\mathbf{P}[P_{n,m} \textrm{ will be connected}]$
is bounded away from $0$ if $\liminf_{n \to \infty} \frac{m}{n} > 1$.
At first,
we shall restrict ourselves to the case when we also have $\limsup_{n \to \infty} \frac{m}{n} < 3$,
so that we can use the topic of appearances.

\begin{Lemma} \label{conn1}
Let $b>1$ and $B<3$ be fixed constants and let $m(n) \in [bn,Bn]$~$\forall n$.
Then $\exists c(b,B)>0$ such that 
\begin{displaymath}
\mathbf{P}[P_{n,m} \textrm{ will be connected}] \geq c \textrm{ } \forall n.
\end{displaymath}
\end{Lemma}
\textbf{Sketch of Proof}
By Lemmas~\ref{gen2} and~\ref{cpt12}, 
it suffices to consider the set of graphs in $\mathcal{P}(n,m)$ with `few' components
and with `many' appearances of a suitable $H$.
Given one of these graphs, 
we may construct another graph in $\mathcal{P}(n,m)$ 
with one less component by deleting a non-cut-edge from an appearance of $H$
(to create an appearance of $H-f$, for a suitable $f$) and inserting an edge between two components.
By cascading this result downwards,
we find that the proportion of graphs in $\mathcal{P}(n,m)$ with exactly one component must be quite decent.
The crucial ingredient in this counting argument is that we have $\Omega (n)$ appearances of $H$ in our given graph
and only $O(n)$ appearances of $H-f$ in our constructed graph. \\
\\
\textbf{Full Proof}
Let $H$ be a $2$-edge-connected planar graph
(for example, we could take $H$ to be $K_{4}$)
and let $f \in E(H)$.
By Lemma~\ref{gen2}, $\exists \alpha = \alpha (b,B) > 0$ such that
$\mathbf{P}[f_{H}(P_{n,m}) \leq \alpha n] = e^{-\Omega(n)}.$
Let us define $\mathcal{G}_{n}$ to denote the set of graphs in $\mathcal{P}(n,m)$
such that $G \in G_{n}$ iff $f_{H}(G) \geq \alpha n$.
Then we have 
$|\mathcal{G}_{n}| > \left( 1 - e^{-\Omega(n)} \right) |\mathcal{P}(n,m)|$.

Let $\mathcal{H}_{n}$ denote the set of graphs in $\mathcal{P}(n,m)$
with less than $\frac{\alpha n}{6}$ components.
Then, by Lemma~\ref{cpt12},
we also have 
$|\mathcal{H}_{n}| > \left( 1 - e^{-\Omega(n)} \right) |\mathcal{P}(n,m)|$.

Let $\mathcal{L}_{n}$ denote the set of graphs in $\mathcal{P}(n,m)$
such that $G \in \mathcal{L}_{n}$ iff $f_{H}(G) \geq \frac{\alpha n}{2} + 3 \kappa (G)$
and $\kappa (G) \leq \frac{\alpha n}{6}$.
Note $\mathcal{L}_{n} \supset \mathcal{G}_{n} \cap \mathcal{H}_{n}$,
so $|\mathcal{L}_{n}| > \left( 1 \!-\! e^{-\Omega(n)} \right) |\mathcal{P}(n,m)|$.
Thus, $\exists N = N(b,B)$ such that
$|\mathcal{L}_{n}| > \frac{1}{2} |\mathcal{P}(n,m)| \textrm{ } \forall n \geq N$.

\begin{figure} [ht]
\setlength{\unitlength}{1cm}
\begin{picture}(20,2.5)(0,0.25)

\put(0.55,0){\line(1,0){0.9}}
\put(0.55,0){\line(3,2){0.45}}
\put(0.55,0){\line(3,5){0.45}}
\put(1,0.3){\line(0,1){0.45}}
\put(1.45,0){\line(-3,2){0.45}}
\put(1.45,0){\line(-3,5){0.45}}
\put(0.55,0){\circle*{0.1}}
\put(1.45,0){\circle*{0.1}}
\put(1,0.3){\circle*{0.1}}
\put(1,0.75){\circle*{0.1}}

\put(1,2){\oval(2,1)}
\put(3.5,2){\oval(2,1)}

\put(1,0.75){\line(0,1){0.75}}
\put(1,1.5){\circle*{0.1}}

\put(5.5,2){\vector(1,0){1}}

\put(8.05,0){\line(1,0){0.9}}
\put(8.05,0){\line(3,2){0.45}}
\put(8.05,0){\line(3,5){0.45}}
\put(8.5,0.3){\line(0,1){0.45}}
\put(8.95,0){\line(-3,5){0.45}}
\put(8.05,0){\circle*{0.1}}
\put(8.95,0){\circle*{0.1}}
\put(8.5,0.3){\circle*{0.1}}
\put(8.5,0.75){\circle*{0.1}}

\put(8.5,2){\oval(2,1)}
\put(11,2){\oval(2,1)}

\put(9.5,2){\line(1,0){0.5}}
\put(9.5,2){\circle*{0.1}}
\put(10,2){\circle*{0.1}}

\put(8.5,0.75){\line(0,1){0.75}}
\put(8.5,1.5){\circle*{0.1}}

\end{picture}

\caption{Constructing a graph in $\mathcal{L}_{n,k}$ from a graph in $\mathcal{L}_{n,k+1}$.} \label{conn1fig}
\end{figure}
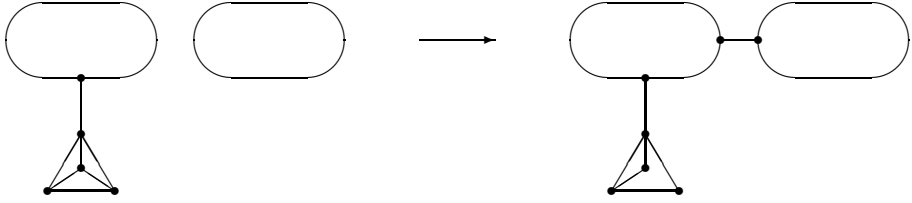

Let $n \geq N$ and 
let $\mathcal{L}_{n,k}$ denote the set of graphs in $\mathcal{L}_{n}$ with exactly $k$ components.
If $2 \leq k+1 \leq \frac{\alpha n}{6}$, we may construct a graph in $\mathcal{L}_{n,k}$
from a graph in $\mathcal{L}_{n,k+1}$ by the following method (see Figure~\ref{conn1fig}):
choose a graph $G \in \mathcal{L}_{n,k+1}$;
choose an appearance of $H$ in $G$ and delete the edge corresponding to $f$ in this appearance
(we have at least $\frac{\alpha n}{2} + 3(k+1) \geq \frac{\alpha n}{2}$ choices for this edge, 
since clearly all appearances of $H$ are disjoint,
by $2$-edge-connectivity);
and insert an edge between two vertices in different components, 
making sure that we don't interfere with a vertex that was in our chosen appearance of $H$
(we have $a_{k}$, say, choices for this edge).

It is clear that the number of components in our new graph, $G^{\prime}$, will be~$k$.
Also, we have $f_{H}(G^{\prime}) \geq \frac{\alpha n}{2} + 3(k+1) - 3 = \frac{\alpha n}{2} - 3k$
(since we had $f_{H}(G) \geq \frac{\alpha n}{2} + 3(k+1)$), 
we deliberately interfered with one appearance in the set $Z$,
and we may have interfered with at most two more when we inserted an edge).
Thus, our new graph is indeed in $\mathcal{L}_{n,k}$.

Note that the number of possible edges between disjoint sets $X$ and $Y$ is~$|X||Y|$ and that
if $|X| \leq |Y|$ then $|X||Y| > (|X|-~1)(|Y|+~1)$,
so it follows that the number of choices for the edge to insert is minimized when 
the $n-|H|$ vertices that are not in our chosen appearance of $H$ are $k$ isolated vertices and 
one component of $n-|H|-k$ vertices.
Thus, \mbox{$a_{k} \geq \left( ^{k} _{2} \right) + k(n-|H|-k)$.}
Hence, we have created at least 
$|\mathcal{L}_{n,k+1}| \frac{\alpha n}{2} \left( \left( ^{k} _{2} \right) + k(n-|H|-k) \right)$ 
graphs in~$\mathcal{L}_{n,k}$.

Given one of our created graphs in $\mathcal{L}_{n,k}$, 
there are at most $6n$ possibilities for~where the deleted edge was originally, 
since it must have been in what is now an appearance of $H - f$
(and clearly $f_{H-f} \left( G^{\prime} \right) \leq 2m$,
since each edge in~$G^{\prime}$ can only be the cut-edge for at most $2$ appearances).
There are then $b_{k}$, say, possibilities for which was the inserted edge.
Since this edge must be a cut-edge that doesn't interfere with the appearance of $H - f$,
we have $b_{k} \leq n-|H|-k$ 
(since the number of edges in a spanning forest of a graph with $k$ components and $n-|H|$ vertices is $n-|H|-k$).
Hence, we have built each graph at most $6n(n-|H|-k)$ times.

Clearly, 
$
\frac{|\mathcal{L}_{n,k+1}| \frac{1}{2} \alpha n \left( \left( ^{k} _{2} \right) + k(n-|H|-k) \right)} {6n(n-|H|-k)} 
\geq \frac{1}{12} |\mathcal{L}_{n,k+1}| \alpha k.
$
Thus, for $2 \!\leq\! k+1 \!\leq\! \frac{\alpha n}{6}$, we have
$|\mathcal{L}_{n,k+1}| \leq \frac{12|\mathcal{L}_{n,k}|}{\alpha k}$.

Let $p_{k} = \frac{|\mathcal{L}_{n,k+1}|}{|\mathcal{L}_{n}|}$
and let $p = p_{0} = \frac{|\mathcal{L}_{n,1}|}{|\mathcal{L}_{n}|}$.
Then (remembering that $|\mathcal{L}_{n,k+1}| = 0$ for $k+1 > \frac{\alpha n}{6}$),
we have $p_{k} \leq \frac{p(\frac{12}{\alpha})^{k}}{k!} \forall k$.
Since $\sum_{k \geq 0} p_{k} = 1$, we must have $\sum_{k \geq 0} \frac{p(\frac{12}{\alpha})^{k}}{k!} \geq 1$
and so $p \geq \left( \sum_{k \geq 0} \frac{(\frac{12}{\alpha})^{k}}{k!} \right)^{-1} = e^{-\frac{12}{\alpha}}$.

Since $\mathcal{L}_{n} > \frac{1}{2} |\mathcal{P}(n,m)|$ and 
$\mathcal{L}_{n,1} \subset \mathcal{P}_{c} (n,m)$,
we have 
$
\frac{|\mathcal{P}_{c}(n,m)|} {|\mathcal{P}(n,m)|} \geq 
\frac{|\mathcal{L}_{n,1}|}{2|\mathcal{L}_{n}|}
= \frac{p}{2} \geq \frac{1}{2}e^{-\frac{12}{\alpha}}.
$
Thus, 
$\mathbf{P}[P_{n,m} \textrm{ will be connected}] 
\geq \frac{1}{2}e^{-\frac{12}{\alpha (b,B)}} > 0$
~$\forall n \geq N(b,B)$.

Clearly, $|\mathcal{P}_{c}(n,m)| > 1$ $\forall n$.
Thus, for all $n$ we have
\begin{displaymath}
\mathbf{P}[P_{n,m} \textrm{ will be connected}] \!\geq\! c(b,B) \!=\! 
\min 
\!\left\{ 
\frac{e^{\frac{-12}{\alpha (b,B)}}}{2}, \min_{n<N(b,B)} \frac{1}{|\mathcal{P}(n,m)|} 
\right\}\! > 0.\!~
\setlength{\unitlength}{.25cm}
\begin{picture}(1,1)
\put(0,0){\line(1,0){1}}
\put(0,0){\line(0,1){1}}
\put(1,1){\line(-1,0){1}}
\put(1,1){\line(0,-1){1}}
\end{picture}
\end{displaymath} \\

The proof of Lemma~\ref{conn1} may seem unnecessarily complicated,
but this is due to the fact that we only know
$\mathbf{P}[f_{H}(P_{n,m})>\beta n] \to 1$
rather than that
$\mathbf{P}[f_{H}(P_{n,m})>\beta n~|~\kappa(P_{n,m})=k+1] \to1$
(in fact,
the second result is trivially false if $k+1=n$,
for example).
Thus,
if we let $\mathcal{P}(n,m,r)$ denote~$\{ G \in \mathcal{P}(n,m) : \kappa (G) = r \}$,
then although we could have used our construction to obtain a lower bound for
$\frac{ | \{ G \in \mathcal{P}(n,m,k+1) : f_{H}(G)>\beta n \} |}
{|\mathcal{P}(n,m,k)|}$,
this would not have given us a lower bound for 
$\frac{|\mathcal{P}(n,m,k+1)|}{|\mathcal{P}(n,m,k)|}$.
Hence,
we had to instead aim to obtain a lower bound for
$\frac{ | \{ G \in \mathcal{P}(n,m,k+1) : f_{H}(G) \textrm{ \scriptsize{is `large'}} \} |}
{ | \{ G \in \mathcal{P}(n,m,k) : f_{H}(G) \textrm{ \scriptsize{is `large'}} \} |}$,
and it is then clear that our construction forces us to define `large' to be something like
$\frac{\alpha n}{2} + 3 \kappa (G)$,
rather than $\alpha n$,
which in turn forced us to use Lemma~\ref{cpt12}
to make sure that 
$\frac{ | \{ G \in \mathcal{P}(n,m) : f_{H}(G) \textrm{ \scriptsize{is `large'}} \} |}
{|\mathcal{P}(n,m)|}$
is bounded away from $0$. \\
\\
\\

We will now see an analogous result to Lemma~\ref{conn1} for when $\frac{m}{n}$ is close to $3$.
The main difference in the construction is that
instead of removing edges from suitable appearances, 
we shall now remove edges from suitable triangles.
Also,
we will work with the entire set of disconnected graphs at once,
rather than conditioning on the exact number of components.

\begin{Lemma} \label{conn2}
Let $C \in \left( \frac{41}{14},3 \right]$ and let $m=m(n) \in [Cn-o(n), 3n-6]$.
Then
\begin{displaymath}
\limsup _{n \to \infty} \mathbf{P}[P_{n,m} \textrm{ will \emph{not} be connected}] 
\leq \frac{15(3-C)}{\frac{2}{7} - 12 + 4C}.
\end{displaymath}
\end{Lemma}
\textbf{Proof}
Let $G \in \mathcal{P}(n,m)$ and let us consider how many triangles in $G$ contain at least one vertex with degree $\leq 6$.
We shall call such triangles `good' triangles.

First, note that (assuming $n \geq 3$)
$G$ may be extended to a triangulation by inserting $3n-6-m \leq (3-C)n + o(n)$ `phantom' edges.
Let $d_{i}$ denote the number of vertices of degree $i$ in such a triangulation.
Then $7 \sum_{i \geq 7} d_{i} \leq~\sum_{i \geq 1} i d_{i} = 2(3n-6)$.
Thus, $\sum_{i \geq 7} d_{i} < \frac{6n}{7}$ and so $\sum_{i \leq 6} d_{i} > \frac{n}{7}$.

For $n>3$, each of these vertices of small degree will be in at least three faces of the triangulation,
all of which will be good triangles.
This counts each good triangle at most three times
(once for each vertex),
so our triangulation must have at least 
$3 \cdot \frac{n}{7} \cdot \frac{1}{3} = \frac{n}{7}$ good triangles that are \textit{faces}.
Each of our phantom edges is in exactly two faces of the triangulation,
so our original graph $G$ must contain at least $\frac{n}{7} - 2(3-C)n + o(n)$ of our good triangles
(note that these triangles will still be `good', 
since the degrees of the vertices will be at most what they were in the triangulation).

We will now consider how many cut-edges a graph in $\mathcal{P}(n,m)$ may have.
If we delete all $c$ cut-edges,
then the remaining graph will consist of $b$, say, blocks,
each of which is either $2$-edge-connected or is an isolated vertex.
Note that the graph formed by condensing each block to a single node and re-inserting the cut-edges must be acyclic,
so $c \leq b-1$.
Label the blocks $1,2, \ldots b$ and let $n_{i}$ denote the number of vertices in block $i$.
Then the number of edges in block $i$ is at most $3n_{i}-6$ if $n_{i} \geq 3$ and is $0=3n_{i}-3$ otherwise
(since $n_{i} < 3$ implies that $n_{i} = 1$).
Thus, $m \leq \sum_{i=1}^{b} (3n_{i}-3) + c = 3n-3b+c < 3n-2c$,
and so $c < \frac{3n-m}{2} < \frac{(3-C)n}{2} + o(n)$.

We now come to the main part of the proof.
Let $\mathcal{G}_{n}$ denote the set of graphs in $\mathcal{P}(n,m)$ that are not connected,
and choose a graph $G \in \mathcal{G}_{n}$.
Choose a good triangle in $G$ (at least $\frac{n}{7} - 2(3-C)n + o(n)$ choices)
and delete an edge that is opposite a vertex with degree $\leq 6$.
Then insert an edge between two vertices in different components 
(we have $a$, say, choices for this edge).

\begin{figure} [ht]
\setlength{\unitlength}{1cm}
\begin{picture}(20,0.8)(0,0.3)

\put(0.7,0.25){\line(1,0){0.6}}
\put(0.7,0.25){\line(3,5){0.3}}
\put(1.3,0.25){\line(-3,5){0.3}}
\put(0.7,0.25){\circle*{0.1}}
\put(1.3,0.25){\circle*{0.1}}
\put(1,0.75){\circle*{0.1}}
\put(1.15,0.7){\scriptsize{$v$}}
\put(0.6,-0.3){\scriptsize{$d(v) \leq 6$}}

\put(1,0.5){\oval(2,1)}
\put(3.5,0.5){\oval(2,1)}

\put(5.5,0.5){\vector(1,0){1}}

\put(8.2,0.25){\line(3,5){0.3}}
\put(8.8,0.25){\line(-3,5){0.3}}
\put(8.2,0.25){\circle*{0.1}}
\put(8.8,0.25){\circle*{0.1}}
\put(8.5,0.75){\circle*{0.1}}
\put(8.65,0.7){\scriptsize{$v$}}

\put(8.5,0.5){\oval(2,1)}
\put(11,0.5){\oval(2,1)}

\put(9.5,0.5){\line(1,0){0.5}}
\put(9.5,0.5){\circle*{0.1}}
\put(10,0.5){\circle*{0.1}}

\end{picture}

\caption{Constructing our new graph.}
\end{figure}
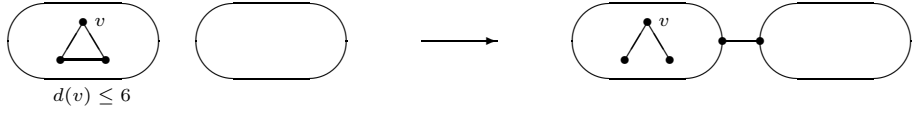

As mentioned in the previous proof, 
the number of possible edges between disjoint sets $X$ and $Y$ is $|X||Y|$ and
if $|X| \leq |Y|$ then $|X||Y| > (|X|-~1)(|Y|+~1)$,
so it follows that the number of choices for the edge to insert is minimized when 
we have one isolated vertex and 
one component of $n-1$ vertices.
Thus, $a \geq n-1$ and so we have created at least 
$|\mathcal{G}_{n}| \left( \frac{n}{7} - 2(3-C)n + o(n) \right) (n-1)$ 
graphs in~$\mathcal{P}(n,m)$.

Given one of our created graphs, 
there are at most $\frac{(3-C)n}{2} + o(n)$ possibilities for which edge was inserted,
since it must be a cut-edge.
There are then at most~$\left( ^{6}_{2} \right)n$ possibilities for where the deleted edge was originally,
since it must have been between two neighbours of a vertex with degree $\leq 6$
(we have at most $n$ possibilities for this vertex and then
at most $\left( ^{6}_{2} \right)$ possibilities for its neighbours).
Hence, we have built each graph at most $\frac{15(3-C)n^{2}}{2} + o(n^{2})$ times.

Thus,
$|\mathcal{P}(n,m)| \geq
\frac{(\frac{1}{7} - 2(3-C))n^{2} + o \left( n^{2} \right)}{ \frac{15(3-C)n^{2}}{2} + o \left( n^{2} \right)}
|\mathcal{G}_{n}|$
and so
\begin{eqnarray*}
\frac{|\mathcal{G}_{n}|}{|\mathcal{P}(n,m)|} & \leq & \frac{\frac{15(3-C)n^{2}}{2} + o(n^{2})}
{\left( \frac{1}{7}-2(3-C) \right) n^{2} + o(n^{2})} 
\phantom{w} \textrm{ since }
\frac{1}{7} \!-\! 2(3\!-\!C) \!>\! 0 \textrm{ for } C \!>\! \frac{41}{14} \\
& \to & \frac{15(3-C)}{\frac{2}{7} - 12 + 4C} \phantom{w} \textrm{ as } n \to \infty.
\phantom{qwerty}
\setlength{\unitlength}{.25cm}
\begin{picture}(1,1)
\put(0,0){\line(1,0){1}}
\put(0,0){\line(0,1){1}}
\put(1,1){\line(-1,0){1}}
\put(1,1){\line(0,-1){1}}
\end{picture}
\end{eqnarray*} \\

We may now combine Lemmas~\ref{conn1} and~\ref{conn2}
to obtain our first main result:

\begin{Theorem} \label{conn4}
Let $b>1$ and let $m=m(n) \in [bn, 3n-6]$. Then $\exists c(b) > 0$ such that 
\begin{displaymath}
\mathbf{P}[P_{n,m} \textrm{ will be connected}] > c \textrm{  } \forall n.
\end{displaymath}
\end{Theorem}
\textbf{Proof}
Clearly, $\exists C<3$ such that $\frac{15(3-C)}{\frac{2}{7} - 12 + 4C} < 1$.
Thus, by considering separately the values of $n$ for which $m(n) \in [bn,Cn]$
and the values for which $m(n) \in~[Cn,3n-6]$,
the result follows from Lemmas~\ref{conn1} and~\ref{conn2}.
$\phantom{qwerty}
\setlength{\unitlength}{.25cm}
\begin{picture}(1,1)
\put(0,0){\line(1,0){1}}
\put(0,0){\line(0,1){1}}
\put(1,1){\line(-1,0){1}}
\put(1,1){\line(0,-1){1}}
\end{picture}$
\begin{displaymath}
\end{displaymath}

As a Corollary to Lemma~\ref{conn2}, by taking $C=3$, we also obtain the following result:

\begin{Corollary} \label{conn5}
Let $m=m(n)=3n-o(n)$.
Then
\begin{displaymath}
\mathbf{P}[P_{n,m} \textrm{ will be connected}] \to 1 \textrm{ as } n \to \infty.
\end{displaymath} \\
\end{Corollary}

As mentioned,
we shall now take a break from the main themes of Part I
to observe that the proofs of this section
can be used to obtain bounds on $\kappa(P_{n,m})$.
We start by observing what was implicitly shown about $\kappa (P_{n,m})$ in the proof of Lemma~\ref{conn1}:

\begin{Lemma} \label{cpt201}
Let $m=m(n)$ satisfy 
$1 < \liminf_{n \to \infty} \frac{m}{n} \leq \limsup_{n \to \infty} \frac{m}{n} < 3$.
Then there exists a constant $K$ such that
\begin{displaymath}
\mathbf{E}[\kappa (P_{n,m})] < K \textrm{ } \forall n.
\end{displaymath}
\end{Lemma}
\textbf{Proof}
Clearly,
it suffices to prove that the result holds for all sufficiently large~$n$.
We define $\mathcal{L}_{n}$ as in the proof of Lemma~\ref{conn1} and recall that,
by Lemmas~\ref{gen2} and~\ref{cpt12},
we have $|\mathcal{L}_{n}| > \left( 1 - e^{- \Omega(n)} \right) |\mathcal{P}(n,m)|$.
Thus,
\begin{eqnarray*}
\mathbf{E}[\kappa (P_{n,m})] 
& = & \mathbf{P}[P_{n,m} \in \mathcal{L}_{n}] \cdot \mathbf{E}[\kappa (P_{n,m}) | P_{n,m} \in \mathcal{L}_{n}] \\
& & + \mathbf{P}[P_{n,m} \notin \mathcal{L}_{n}] \cdot \mathbf{E}[\kappa (P_{n,m}) | P_{n,m} \notin \mathcal{L}_{n}] \\
& \leq & \mathbf{E}[\kappa (P_{n,m}) | P_{n,m} \in \mathcal{L}_{n}]
+ e^{- \Omega(n)} n \\
& = & \mathbf{E}[\kappa (P_{n,m}) | P_{n,m} \in \mathcal{L}_{n}]
+ e^{- \Omega(n)}.
\end{eqnarray*}

It now only remains to show that 
$\mathbf{E}[\kappa (P_{n,m}) | P_{n,m} \in \mathcal{L}_{n}]$
is bounded by a constant.
But recall we showed in the proof of Lemma~\ref{conn1} that 
$|\mathcal{L}_{n,k+1}| \!\leq~\!\frac{12 |\mathcal{L}_{n,k}|}{ \alpha k}$~$\forall k$,
for a suitable constant $\alpha$.
Thus, we have
$\mathbf{E}[\kappa (P_{n,m}) | P_{n,m} \in \mathcal{L}_{n}] \leq
\mathbf{E}[X]$,
where $X \sim~\textrm{Poi} \left( \frac{12}{\alpha} \right)$.
Since $\mathbf{E}[X] = \frac{12}{\alpha}$,
we are done.
$\phantom{qwerty}
\setlength{\unitlength}{.25cm}
\begin{picture}(1,1)
\put(0,0){\line(1,0){1}}
\put(0,0){\line(0,1){1}}
\put(1,1){\line(-1,0){1}}
\put(1,1){\line(0,-1){1}}
\end{picture}$
\\
\\

Similarly,
we may obtain a result for when $\frac{m}{n}$ is close to $3$,
by analysing the proof of Lemma~\ref{conn2}:

\begin{Lemma} \label{cpt202}
Let $C \in (\frac{41}{14},3]$ and let $m=m(n) \in [Cn-o(n), 3n-6]$.
Then there exists a constant $K$ such that
\begin{displaymath}
\mathbf{E}[\kappa (P_{n,m})] < K \textrm{ } \forall n.
\end{displaymath}
\end{Lemma}
\textbf{Proof}
Again,
it suffices to prove that the result holds for all sufficiently large~$n$.
Let $\mathcal{G}_{n,k}$ denote the set of graphs in $\mathcal{P}(n,m)$ with exactly $k$ components.
We follow the proof of Lemma~\ref{conn2} to construct graphs in $\mathcal{G}_{n,k}$ from graphs in~$\mathcal{G}_{n,k+1}$.
The details are as before, 
except that the number of ways to insert an edge between two vertices in different components is now
$\left( ^{k}_{2} \right)+ k(n-k)$.
Thus, we obtain 
\begin{eqnarray*}
\frac{|\mathcal{G}_{n,k+1}|}{|\mathcal{G}_{n,k}|} 
& \leq &  \frac{\frac{15(3-C)n^{2}}{2} + o(n^{2})}
{\left( \frac{1}{7}-2(3-C) \right)k n^{2} + o(n^{2})} \\
& \to & \frac{15(3-C)}{\left( \frac{2}{7} - 12 + 4C \right)k} \phantom{w} \textrm{ as } n \to \infty.
\end{eqnarray*}
Therefore, for sufficiently large $n$, we have
$\mathbf{E}[\kappa (P_{n,m})] \leq \mathbf{E}[X]$,
where $X \!\sim~\!\textrm{Poi} (\lambda)$ for any $\lambda > \frac{15(3-C)}{\frac{2}{7} - 12 + 4C}$.
But $\mathbf{E}[X] = \lambda$, and so we are done.
$\phantom{qwerty}
\setlength{\unitlength}{.25cm}
\begin{picture}(1,1)
\put(0,0){\line(1,0){1}}
\put(0,0){\line(0,1){1}}
\put(1,1){\line(-1,0){1}}
\put(1,1){\line(0,-1){1}}
\end{picture}$
\\
\\

We may combine Lemmas~\ref{cpt201} and~\ref{cpt202}:

\begin{Proposition} \label{cpt2031}
Let $m=m(n)$ satisfy $\liminf_{n \to \infty} \frac{m}{n} > 1$.
Then there exists a constant $K$ such that
\begin{displaymath}
\mathbf{E}[\kappa (P_{n,m})] < K \textrm{ } \forall n.
\end{displaymath} 
\end{Proposition}

By Markov's inequality,
we also obtain the following result,
which improves on Lemma~\ref{cpt12} 
(in terms of decreasing the bound on $\kappa(P_{n,m})$)
for the case when $\liminf_{n \to \infty} \frac{m}{n} > 1$:

\begin{Corollary} \label{cpt203}
Let $m=m(n)$ satisfy $\liminf_{n \to \infty} \frac{m}{n} > 1$
and let $g(n) \to \infty$ as~$n \to \infty$.
Then
\begin{displaymath}
\mathbf{P}[\kappa (P_{n,m}) < g(n)] \to 1 \textrm{ as } n \to \infty.
\end{displaymath} \\
\end{Corollary}

We shall shortly give a summary of all our results on $\kappa(P_{n,m})$,
but first let us obtain one final bound,
for when $m<n$
(this will in fact be useful later on,
in Section~\ref{cycsub}):

\begin{Proposition} \label{cpt31}
Let $m=m(n) \in [qn,n-1]$ for some fixed $q \in (0,1)$
and let the constant $c$ satisfy $c> \ln \left( \frac{\gamma _{l}}{q} \right) + e^{-1}$.
Then
\begin{displaymath}
\mathbf{P} \left[ \kappa (P_{n,m}) > \left \lceil \frac{cn}{\ln n} + n-m \right \rceil \right] = e^{- \Omega(n)}.
\end{displaymath}
\end{Proposition}
\textbf{Proof}
Note that
$\mathbf{P} \left[ \kappa (P_{n,m}) > \left \lceil \frac{cn}{\ln n} \!+\! n\!-\!m \right \rceil \right]
= \frac { \left | \left \{ G \in \mathcal{P}(n,m):
\kappa(G)> \left \lceil \frac{cn}{\ln n} + n-m \right \rceil \right \} \right|}{|\mathcal{P}(n,m)|}$,
by definition.
We shall obtain bounds for both parts of this fraction by following
the proof of Proposition~\ref{ger3 6.6} (\cite{ger3}, 6.6),
which deals with the case when \mbox{$m= n-(\beta + o(1))(n/ \ln n)$} for fixed $\beta >0$.

Let $x=x(n)=n-m$ and let $k=k(n)= \left \lceil \frac{cn}{\ln n} +x \right \rceil$.
Then
\begin{displaymath}
| \{ G \in \mathcal{P}(n,m):\kappa(G)>k \}| \leq | \{ G \in \mathcal{P}(n):\kappa(G)>k \}| \leq 
\frac{|\mathcal{P}(n)|}{k!},
\end{displaymath}
using Proposition~\ref{mcd 2.1}.

Let $\mathcal{F}(i,j)$ denote the set of forests with $i$ vertices and $j$ edges.
Then
\begin{eqnarray*}
|\mathcal{P}(n,m)| & \geq & |\mathcal{F}(n,m)| \\
& \geq & |\mathcal{F}(m+1,m)|, \phantom{w} \textrm{ since } n \geq m+1 \phantom{ww} \\
& & \textrm{(consider adding $n-(m+1)$ isolated vertices)} \\
& \geq & (m+1)^{m-1}, \phantom{www} \textrm{by Cayley's Theorem}. 
\end{eqnarray*}

Thus, 
\begin{displaymath}
\mathbf{P} [\kappa (P_{n,m}) > k] \leq \frac{|\mathcal{P}(n)|}{k!(m+1)^{m-1}}. 
\end{displaymath}

Note that 
\begin{eqnarray}
\ln \left( (m+1)^{m-1} \right) & = & (m-1) \ln (m+1) \nonumber \\
& = & m \ln (m+1) + o(n) \nonumber \\
& = & m \ln n + m \ln \left( \frac{m+1}{n} \right) + o(n) \nonumber \\
& \geq & m \ln n + m \ln q + o(n) \nonumber \\
& \geq & m \ln n + n \ln q + o(n), \textrm{ since } q<1\label{eq:c1} 
\end{eqnarray}

Also, note that
\begin{eqnarray} 
k \ln k - x \ln n & = & (k-x) \ln n + k \ln k - k \ln n \nonumber \\
& = & (k-x) \ln n - n \frac{ \ln \left( \frac{n}{k} \right) }{\frac{n}{k}} \nonumber \\
& \geq & (k-x) \ln n - ne^{-1}, \phantom{www}
\textrm{ since } \frac{\ln y}{y} \leq e^{-1} \textrm{ } \forall y \nonumber \\
& = & \frac{cn}{\ln n} \ln n - ne^{-1} + o(n) \textrm{ by definition of } k \nonumber \\
& = & \left(c- e^{-1} \right) n + o(n). \label{eq:c2} 
\nonumber 
\end{eqnarray}

Thus, we have
\begin{eqnarray*}
& & \limsup_{n \to \infty} \frac{|\mathcal{P}(n)|}{k!(m+1)^{m-1}} \\
& \leq & \frac{(\gamma_{l})^{n} n!}{ \left( \frac{k}{e} \right)^{k} e^{m \ln n + n \ln q + o(n)}}, 
\textrm{$\phantom{w}$ by Proposition~\ref{gim T1}, Stirling's formula and (\ref{eq:c1})} \\
& = & \exp (n \ln \gamma_{l} + n \ln n -n +k - k \ln k - m \ln n - n \ln q + o(n)), \\
& & \textrm{using Stirling's formula for $n!$} \\
& = & \exp (n \ln \gamma_{l} -n +k - k \ln k + x \ln n - n \ln q + o(n)), \phantom{w} 
\textrm{ by definition of } x \\
& \leq & \exp (n \ln \gamma_{l} -n +k -\left(c-e^{-1}\right)n - n \ln q + o(n)), 
\phantom{www} \textrm{ by } (\ref{eq:c2}) \\
& \leq & \exp \left(n \left(\ln \left( \frac{\gamma_{l}}{q} \right) -c + e^{-1} + o(1) \right) \right),
\phantom{www} \textrm{ since } k \leq n +o(n) \\
& = & e^{- \Omega (n)} 
\textrm{, $\phantom{www}$ since $c> \ln \left( \frac{\gamma_{l}}{q} \right) + e^{-1}$}.
\phantom{qwerty}
\setlength{\unitlength}{.25cm}
\begin{picture}(1,1)
\put(0,0){\line(1,0){1}}
\put(0,0){\line(0,1){1}}
\put(1,1){\line(-1,0){1}}
\put(1,1){\line(0,-1){1}}
\end{picture}
\end{eqnarray*}

\newpage
\label{cptsum}
\subsection*{Summary of Results on the Number of Components}

$\exists \lambda >0$ such that

\begin{displaymath}
\kappa(P_{n,m}) \geq
\left\{ \begin{array}{ll}
n-m & \forall m \phantom{www} \textrm{(trivial observation)} \\
\lambda n^{1/2} \phantom{w} \textrm{a.a.s.} &
\textrm{if $|m-n|=O(n^{1/2})$ \phantom{ww} (Theorems~\ref{tree4} and~\ref{add3}) }\\
\lambda \frac{n}{d} \phantom{w} \textrm{a.a.s.} &
\textrm{if $d\!=\!m\!-\!n\!>\!0$ is such that $d=\Omega(n^{1/2})$ and $O(n)$ } \\
& \textrm{(Theorems~\ref{tree4} and~\ref{add54}) }
\end{array} \right. 
\end{displaymath} 
and $\exists c$ such that
\begin{displaymath}
\kappa(P_{n,m}) \leq
\left\{ \begin{array}{ll}
n-m + \frac{cn}{\ln n} \phantom{w} \textrm{a.a.s.} 
\phantom{w} \textrm{ if } 0<q \leq \frac{m}{n} <1 &
\textrm{(Proposition~\ref{cpt31})} \\
\frac{cn}{\ln n} \phantom{w} \textrm{a.a.s.} 
\phantom{wwwqqqqqi} \textrm{if } \frac{m}{n} \geq 1 &
\textrm{(Lemma~\ref{cpt12})} \\
g(n) \phantom{w} \textrm{a.a.s.} 
\phantom{wwwwqqii} \textrm{if } \liminf_{n \to \infty} \frac{m}{n} > 1 &
\textrm{(Corollary~\ref{cpt203})} \\
(\textrm{where $g$ is any function with $g(n) \to \infty$}) 
\end{array} \right. 
\end{displaymath} \\

These results may be represented as below,
where $\epsilon$ is an arbitrarily small positive constant: \\

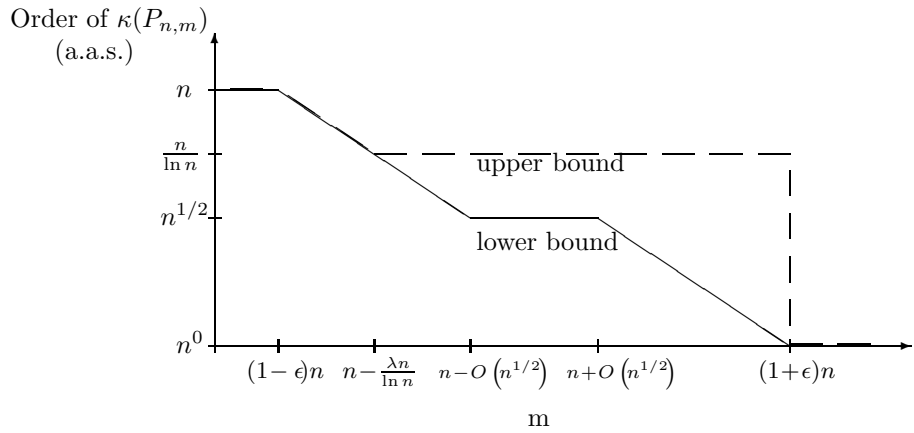
\begin{figure} [ht]
\setlength{\unitlength}{0.85cm}
\begin{picture}(2,6.4)(-1.2,0)
\put(-1.2,6.25){\textrm{Order of $\kappa(P_{n,m})$}} 
\put(-0.5,5.75){\textrm{(a.a.s.)}}
\put(1.4,1.15){$n^{0}$}
\put(1.15,3.15){$n^{1/2}$}
\put(1.15,4.15){$\frac{n}{\ln n}$}
\put(1.4,5.15){$n$}
\put(2.5,0.75){\small{$(1\!-\epsilon\!)n$}}
\put(4,0.75){\small{$n\!-\!\frac{\lambda n}{\ln n}$}}
\put(5.5,0.75){\scriptsize{$n\!-\!O\left(\!n^{1/2}\!\right)$}}
\put(7.5,0.75){\scriptsize{$n\!+\!O\left(\!n^{1/2}\!\right)$}}
\put(10.5,0.75){\small{$(1\!+\!\epsilon)n$}}
\put(2,1.15){\vector(0,1){5}}
\put(1.9,1.25){\vector(1,0){11}}
\put(3,1.15){\line(0,1){0.2}}
\put(4.5,1.15){\line(0,1){0.2}}
\put(6,1.15){\line(0,1){0.2}}
\put(8,1.15){\line(0,1){0.2}}
\put(11,1.15){\line(0,1){0.2}}
\put(1.9,3.25){\line(1,0){0.2}}
\put(1.9,4.25){\line(1,0){0.2}}
\put(1.9,5.25){\line(1,0){1.1}}
\put(3,5.25){\line(3,-2){3}}
\put(6,3.25){\line(1,0){2}}
\put(8,3.25){\line(3,-2){3}}
\put(4.5,4.25){\line(1,0){0.5}}
\put(5.25,4.25){\line(1,0){0.5}}
\put(6,4.25){\line(1,0){0.5}}
\put(6.75,4.25){\line(1,0){0.5}}
\put(7.5,4.25){\line(1,0){0.5}}
\put(8.25,4.25){\line(1,0){0.5}}
\put(9,4.25){\line(1,0){0.5}}
\put(9.75,4.25){\line(1,0){0.5}}
\put(10.5,4.25){\line(1,0){0.5}}
\put(11,1.27){\line(1,0){0.5}}
\put(11.75,1.27){\line(1,0){0.5}}
\put(2.25,5.27){\line(1,0){0.5}}
\put(3.15,5.17){\line(3,-2){0.5}}
\put(3.9,4.67){\line(3,-2){0.5}}

\put(11,4){\line(0,-1){0.4}}
\put(11,3.325){\line(0,-1){0.4}}
\put(11,2.65){\line(0,-1){0.4}}
\put(11,1.975){\line(0,-1){0.4}}

\put(6.1,4){\textrm{upper bound}}
\put(6.1,2.75){\textrm{lower bound}}
\put(6.9,0){m}
\end{picture} 
\caption{Summary of results on the number of components.}
\end{figure}

\newpage
\section{Components II: Upper Bounds} \label{cyccpt}

We now return to looking at 
$\mathbf{P}[P_{n,m} \textrm{ will have a component isomorphic to } H]$. 
Note that we already have a full description
(in terms of knowing exactly when the probability is or isn't bounded away from $0$ or $1$)
for when $H$ is a tree,
by combining the lower bounds of Section~\ref{cptlow} 
with the upper bounds implied by our results on 
$\mathbf{P}[P_{n,m} \textrm{ will be connected}]$.
In this section,
we will complete matters by obtaining further upper bounds for when $H$ is not acyclic.

First, we will look at the case when $H$ is `multicyclic',
i.e.~when it has more edges than vertices.
We have already seen (in Theorems~\ref{gen3} and~\ref{conn4}) 
that the probability of $P_{n,m}$ containing $H$ as a component is bounded away from both~$0$ and $1$
for large~$n$
if $\frac{m}{n}$ is bounded below by $b > 1$ and above by $B < 3$,
and (from Corollary~\ref{conn5}) that the probability converges to $0$ if $\frac{m}{n} \to 3$.
We shall now see (in Theorem~\ref{cyc41}) that the limiting probability is also $0$ if $\frac{m}{n} \leq 1+o(1)$. 

We shall then look at the case when $H$ is unicyclic.
We have already seen (in Corollary~\ref{conn5}) 
that the probability of $P_{n,m}$ containing $H$ as a component tends to $0$ if $\frac{m}{n} \to 3$,
that (in Theorem~\ref{gen4}) it is bounded away from $0$ for large~$n$
if $\frac{m}{n}$ is bounded away from both $0$ and $3$,
and (in Theorem~\ref{conn4}) 
that it is bounded away from $1$ if $\frac{m}{n}$ is bounded below by~$b>1$.
In this section,
we shall complete the picture by seeing (in Lemma~\ref{cyc52}) that, in fact, 
the probability is \textit{always} bounded away from $1$ for large $n$. 

We shall finish this section (in Theorem~\ref{cyc53}) with an extension of Lemma~\ref{cyc52}
in which we will show that, for any fixed $k$,
the probability that 
$P_{n,m}$ will contain \textit{any} unicyclic or multicyclic components of order $\leq k$ 
is also always bounded away from $1$ for large $n$.

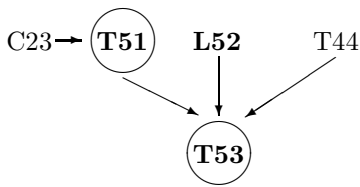
\begin{figure} [ht]
\setlength{\unitlength}{1cm}
\begin{picture}(20,1.5)(0.3,0.5)

\put(5.5,1.9){\circle{0.8}}
\put(6.775,0.4){\circle{0.8}}

\put(4.6,1.9){\vector(1,0){0.35}}
\put(5.5,1.4){\vector(2,-1){1}}
\put(6.775,1.675){\vector(0,-1){0.775}}
\put(8.325,1.675){\vector(-3,-2){1.2}}

\put(3.95,1.775){C\ref{add41}}
\put(5.15,1.775){\textbf{T\ref{cyc41}}}
\put(6.425,1.775){\textbf{L\ref{cyc52}}}
\put(8.025,1.775){T\ref{conn4}}
\put(6.425,0.275){\textbf{T\ref{cyc53}}}

\end{picture}

\caption{The structure of Section~\ref{cyccpt}.}
\end{figure}

We start with our aforementioned result for multicyclic components:

\begin{Theorem} \label{cyc41}
Let $m \leq (1+o(1))n$ and let $H$ be a (fixed) multicyclic connected planar graph.
Then
\begin{displaymath}
\mathbf{P}[P_{n,m} \textrm{ will have a component isomorphic to } H] \to 0 \textrm{ as } n \to \infty.
\end{displaymath}
\end{Theorem}
\textbf{Proof}
Let $\mathcal{G}_{n}$ denote the set of graphs in $\mathcal{P}(n,m)$ with a component isomorphic to $H$.
For each graph $G \in \mathcal{G}_{n}$,
let us delete $2$ edges from a component $H^{\prime} (=H^{\prime}_{G})$ isomorphic to $H$
in such a way that we do not disconnect the component.
Let us then insert one edge between a vertex in the remaining component and a vertex elsewhere in the graph.
We have $|H|(n-|H|)$ ways to do this, and planarity is maintained.
Let us then also insert one other edge into the graph,
without violating planarity.
We have at least $(\textrm{add}(n,m))=\omega(n)$ 
choices for where to place this second edge,
by Corollary~\ref{add41}.
Thus, we can construct $|\mathcal{G}_{n}| \omega \left( n^{2} \right)$ 
(not necessarily distinct) graphs in $\mathcal{P}(n,m)$.

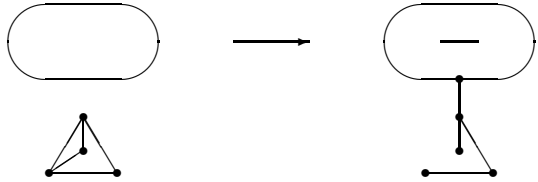
\begin{figure} [ht]
\setlength{\unitlength}{1cm}
\begin{picture}(20,2.75)

\put(3.05,0){\line(1,0){0.9}}
\put(3.05,0){\line(3,2){0.45}}
\put(3.05,0){\line(3,5){0.45}}
\put(3.5,0.3){\line(0,1){0.45}}
\put(3.95,0){\line(-3,5){0.45}}
\put(3.05,0){\circle*{0.1}}
\put(3.95,0){\circle*{0.1}}
\put(3.5,0.3){\circle*{0.1}}
\put(3.5,0.75){\circle*{0.1}}

\put(3.5,1.75){\oval(2,1)}

\put(5.5,1.75){\vector(1,0){1}}

\put(8.05,0){\line(1,0){0.9}}
\put(8.5,0.3){\line(0,1){0.45}}
\put(8.95,0){\line(-3,5){0.45}}
\put(8.05,0){\circle*{0.1}}
\put(8.95,0){\circle*{0.1}}
\put(8.5,0.3){\circle*{0.1}}
\put(8.5,0.75){\circle*{0.1}}

\put(8.5,1.75){\oval(2,1)}

\put(8.25,1.75){\line(1,0){0.5}}

\put(8.5,0.75){\line(0,1){0.5}}
\put(8.5,1.25){\circle*{0.1}}

\end{picture}

\caption{Redistributing edges from our multicyclic component.}
\end{figure}

Given one of our constructed graphs, there are $m=O(n)$
possibilities for the edge that was inserted last.
There are then at most $m-1=O(n)$ possibilities for the other edge that was inserted.
Since one of the two vertices incident with this edge must belong to $V(H^{\prime})$,
we then have at most two possibilities for~$V(H^{\prime})$
and then at most $\left( ^{\left( ^{|H|}_{\phantom{q}2} \right)}_{\phantom{qq}2} \right) = O(1)$ 
possibilities for $E(H^{\prime})$.
Thus, we have built each graph at most $O \left( n^{2} \right)$ times,
and so 
$\frac{|\mathcal{G}_{n}|}{|\mathcal{P}(n,m)|} = \frac{O \left( n^{2} \right)}{\omega \left( n^{2} \right)} \to 0$~
as~$n \to~\infty.\!$~
\setlength{\unitlength}{.25cm}
\begin{picture}(1,1)
\put(0,0){\line(1,0){1}}
\put(0,0){\line(0,1){1}}
\put(1,1){\line(-1,0){1}}
\put(1,1){\line(0,-1){1}}
\end{picture}
\phantom{p}

We shall now look at unicyclic components.
The basic argument will be the same as that of Theorem~\ref{cyc41},
i.e.~we will start with $\mathcal{G}_{n}$,
the set of graphs in~$\mathcal{P}(n,m)$ with a component isomorphic to $H$,
and redistribute edges from such components to construct other graphs in $\mathcal{P}(n,m)$.
This time we will only be able to delete one edge from each component,
and so we will only be able to show
$\frac{\mathcal{G}_{n}|}{|\mathcal{P}(n,m)|} = O(1)$,
rather than $o(1)$.
Hence,
it will be crucial to keep track of the constants involved in the calculations,
so that we can try to show~$\limsup_{n \to \infty} \frac{\mathcal{G}_{n}|}{|\mathcal{P}(n,m)|} < 1$.
Unfortunately,
it turns out that the constants are actually only small enough for certain $H$,
such as when $H$ is a cycle.
However,
it is simple to relate the probability that a given component is isomorphic to~$H$ 
to the probability that it is a cycle,
and so we may deduce that the result actually holds for any unicyclic $H$.

A further complication is that we will actually prove a stronger version of Lemma~\ref{cyc52}
than that advertised at the start of this section,
involving many unicyclic graphs at once,
in preparation for Theorem~\ref{cyc53}. \\

\begin{Lemma} \label{cyc52}
Let $k$ be a fixed constant.
Then, given any $m(n)$,
\begin{displaymath}
\limsup_{n \to \infty}
\mathbf{P}[P_{n,m} \textrm{ will contain any unicyclic components of order } \leq k] < 1.
\end{displaymath}
\end{Lemma}
\textbf{Sketch of Proof}
Let $C_{i}$ be used to denote a component that is a cycle of order~$i$ 
and let $H_{i}$ be used to denote all other unicyclic components of order $i$. 

Step $1$: 
By an induction hypothesis, 
we may assume that there is a decent proportion of graphs with no unicyclic components of order $<l$.

Step $2$:
If we have a $C_{l}$,
we may remove an edge from it and then insert an edge
between the remainder of the cycle and the rest of the graph,
analogously to the proof of Theorem~\ref{cyc41}.
We find that we may create so many graphs by doing this
that the proportion of graphs \textit{without} a $C_{l}$ (and still with no unicyclic components of order $<l$)
must be quite decent.

Step $3$:
We then show that the proportion of graphs with lots of $H_{l}$'s and no~$C_{l}$'s is small.
Thus, using the result of Step $2$,
there must be a decent proportion of graphs with few $H_{l}$'s, no $C_{l}$'s and no unicyclic components of order $<l$.

Step $4$:
We replace the few $H_{l}$'s with $C_{l}$'s to obtain a decent proportion of graphs 
with no $H_{l}$'s and no unicyclic components of order $<l$.

Step $5$:
We then repeat the argument of Step $2$ to find that the proportion of graphs without a $C_{l}$,
and still with no $H_{l}$'s and no unicyclic components of order $<l$,
must be quite decent.
Thus, we have a decent proportion of graphs with no unicyclic components of order $\leq l$,
and so we are done by induction. \\
\\
\textbf{Full Proof}
For $\mathcal{A}_{n} \subset \mathcal{P}(n,m)$,
we shall use $\mathcal{A}_{n}^{c}$ to denote the set of graphs in $\mathcal{P}(n,m)$ that are not in $\mathcal{A}_{n}$. \\
\\
\underline{Step $1$} 

We will prove the result by induction on $k$.
It is trivial for $k=2$.
Let us assume that the result is true for $k=l-1$.

Let $\mathcal{L}_{n}$ denote the set of graphs in $\mathcal{P}(n,m)$ with no unicyclic component of order $<l$.
Then, by our induction hypothesis,
$\exists \epsilon (l) > 0$ such that
$|\mathcal{L}_{n}| \!\geq~\!\epsilon |\mathcal{P}(n,m)|$ for all sufficiently large $n$. \\
\\
\underline{Step $2$} 

Let $\mathcal{G}_{n}$ denote the set of graphs in $\mathcal{L}_{n}$ with a component isomorphic to $C_{l}$,
the cycle of order $l$.
For each graph $G \in \mathcal{G}_{n}$,
let us delete an edge from a component $C^{\prime} (=C^{\prime}_{G})$ isomorphic to $C_{l}$,
leaving a spanning path. 
We have $l$ choices for this edge.
Let us then insert an edge between a vertex in the spanning path and a vertex elsewhere.
We have $l(n-l)$ ways to do this, planarity is maintained, and no unicyclic components of order $<l$ will be created,
since the only new component will have order $>l$.
Thus, we have constructed $|\mathcal{G}_{n}|l^{2}(n-l)$ graphs in $\mathcal{L}_{n}$.

Given one of our constructed graphs,
there are at most $n-1$ possibilities for which is the inserted edge,
since it must be a cut-edge.
There are then at most $2$ possibilities for $V(C^{\prime})$ (at most one possibility for each endpoint of the edge),
and then only one possibility for where the inserted edge was originally.
Thus, we have built each graph at most $2(n-1)$ times.

Therefore, $|\mathcal{G}_{n}| \leq \frac{2(n-1)}{l^{2}(n-l)} |\mathcal{L}_{n}|$.
Thus, for any fixed constant $\delta \in \left( 0, 1-\frac{2}{l^{2}} \right)$,
we have $|\mathcal{G}_{n}| < (1- \delta) |\mathcal{L}_{n}|$ for all sufficiently large $n$,
and so we have
$|\mathcal{L}_{n} \cap \mathcal{G}_{n}^{c}| > \delta |\mathcal{L}_{n}| \geq \delta \epsilon |\mathcal{P}(n,m)|$
for all sufficiently large $n$. \\
\\
\underline{Step $3$} 

Let $r = r(l) > \frac{4 \left( ^{ \left( ^{l}_{2} \right)} _{\phantom{q}l} \right)}{\delta \epsilon (l-1)!}$
and let $\mathcal{J}_{n}$ denote the set of graphs in $\mathcal{P}(n,m)$ with $>r$ unicyclic components of order $l$.

Consider the set $\mathcal{J}_{n} \cap \mathcal{G}_{n}^{c}$,
i.e.~the set of graphs in $\mathcal{P}(n,m)$ with $>r$ unicyclic components of order $l$,
but no components isomorphic to $C_{l}$.
For each graph $J \in \mathcal{J}_{n} \cap \mathcal{G}_{n}^{c}$,
delete a unicyclic component of order $l$ (we have $>r$ choices for this)
and replace it with a component isomorphic to $C_{l}$ (we have $\frac{(l-1)!}{2}$ choices for this).
We have constructed at least $|\mathcal{J}_{n} \cap \mathcal{G}_{n}^{c}| \frac{r(l-1)!}{2}$
(not necessarily distinct) graphs in $\mathcal{P}(n,m)$.

Given one of our constructed graphs, there are at most 
$\left( ^{ \left( ^{l}_{2} \right)} _{\phantom{q}l} \right)$ 
possibilities for what the component isomorphic to $C_{l}$ was originally,
since it must have been a component with $l$ edges.
Thus, we have built each graph at most $\left( ^{ \left( ^{l}_{2} \right)} _{\phantom{q}l} \right)$ times,
and so 
$
|\mathcal{J}_{n} \cap \mathcal{G}_{n}^{c}| \leq 
\frac{ 2 \left( ^{ \left( ^{l}_{2} \right)} _{\phantom{q}l} \right) |\mathcal{P}(n,m)| } {r(l-1)!} 
< \frac{\delta \epsilon}{2} |\mathcal{P}(n,m)|, \textrm{ since }
r > \frac{4 \left( ^{ \left( ^{l}_{2} \right)} _{\phantom{q}l} \right)}{\delta \epsilon (l-1)!}.
$

Since we also know (from Step $2$) that 
$|\mathcal{L}_{n} \cap \mathcal{G}_{n}^{c}| > \delta \epsilon |\mathcal{P}(n,m)|$ for all sufficiently large~$n$,
we must have 
$|\mathcal{L}_{n} \cap \mathcal{J}_{n}^{c} \cap \mathcal{G}_{n}^{c}| > \frac{\delta \epsilon}{2} |\mathcal{P}(n,m)|$.
Thus, 
$|\mathcal{L}_{n} \cap\mathcal{J}_{n}^{c}| > \frac{\delta \epsilon}{2} |\mathcal{P}(n,m)|$ 
for all sufficiently large $n$. \\
\\
\underline{Step $4$} 

Let $\mathcal{H}_{n}$ denote the set of graphs in $\mathcal{P}(n,m)$
with a unicyclic component of order $l$ that is not isomorphic to $C_{l}$,
and recall that $\mathcal{L}_{n} \cap \mathcal{J}_{n}^{c}$ is the set of graphs in $\mathcal{P}(n,m)$
with no unicyclic component of order $<l$ and with $\leq r$ unicyclic components of order $l$.
For each graph $L \in \mathcal{L}_{n} \cap \mathcal{J}_{n}^{c}$,
delete all the unicyclic components of order $l$ and replace them with components isomorphic to $C_{l}$.
Thus, we have constructed at least $|\mathcal{L}_{n} \cap \mathcal{J}_{n}^{c}|$ 
(not necessarily distinct) graphs in $\mathcal{L}_{n} \cap \mathcal{H}_{n}^{c}$.

Suppose we are given one of our constructed graphs 
and suppose it has exactly $s$ components that are isomorphic to $C_{l}$
(note that $s \leq r$).
The $ls$ vertices in these components must all have been in unicyclic components originally,
so there are at most
$\left( ^{ \left( ^{ls}_{\phantom{l}2} \right)} _{\phantom{q}ls} \right) \leq (ls)^{2ls} \leq (lr)^{2lr}$ 
possibilities for what the original graph was.

Thus, 
\begin{eqnarray*}
|\mathcal{L}_{n} \cap \mathcal{H}_{n}^{c}| 
& \geq & \frac{|\mathcal{L}_{n} \cap \mathcal{J}_{n}^{c}|} {(lr)^{2lr}} \\
& > & \frac{\delta \epsilon |\mathcal{P}(n,m)|} {2(lr)^{2lr}} 
\textrm{ for all sufficiently large $n$, by Step $3$} \\
& = & \lambda |\mathcal{P}(n,m)|, \textrm{ where } \lambda = \frac{\delta \epsilon} {2(lr)^{2lr}}.
\end{eqnarray*} 
\\
\underline{Step $5$} 

By applying the same argument as in Step $2$,
but to $\mathcal{L}_{n} \cap \mathcal{H}_{n}^{c}$ rather than~$\mathcal{L}_{n}$,
we may obtain 
$|\mathcal{L}_{n} \cap \mathcal{H}_{n}^{c} \cap \mathcal{G}_{n}^{c}| >
\delta^{\prime} |\mathcal{L}_{n} \cap \mathcal{H}_{n}^{c}| \geq
\delta^{\prime} \lambda |\mathcal{P}(n,m)|$,
where $\delta^{\prime}$ is a suitable positive constant.
But $\mathcal{L}_{n} \cap \mathcal{H}_{n}^{c} \cap \mathcal{G}_{n}^{c}$ is the set of graphs in $\mathcal{P}(n,m)$
with no unicyclic components of order $\leq l$.
Thus, the induction hypothesis holds for~$k=l$, and so we are done.
$\phantom{qwerty}$
\setlength{\unitlength}{.25cm}
\begin{picture}(1,1)
\put(0,0){\line(1,0){1}}
\put(0,0){\line(0,1){1}}
\put(1,1){\line(-1,0){1}}
\put(1,1){\line(0,-1){1}}
\end{picture}
\begin{displaymath}
\end{displaymath} \\

By using Theorems~\ref{conn4} and~\ref{cyc41},
we may extend Lemma~\ref{cyc52} to cover all \textit{non-acyclic} components of order $\leq k$:

\begin{Theorem} \label{cyc53}
Let $k$ be a fixed constant.
Then, given any $m(n)$,
\begin{displaymath}
\limsup_{n \to \infty}
\mathbf{P}[P_{n,m} \textrm{ will contain any non-acyclic components of order } \leq k] < 1.
\end{displaymath}
\end{Theorem}
\textbf{Proof}
We shall suppose,
hoping for a contradiction,
that there exists a function $m(n)$ such that
$
\limsup_{n \to \infty} \mathbf{P}[P_{n,m} \textrm{ will contain any non-acyclic components of}$
$\textrm{order} \leq k] = 1.
$
Thus,
there exists an infinite subsequence $n_{1},n_{2}, \ldots$
such that 
$
\mathbf{P}[P_{n_{i},m(n_{i})} \textrm{ will contain any non-acyclic components of order} \!\leq\! k] 
\!\to\! 1
\textrm{ as } i \!\to~\!\!\infty.
$

Clearly,
by combining Theorem~\ref{cyc41} with Lemma~\ref{cyc52},
we can't have $m(n_{i}) \leq (1+o(1))n_{i}$ as $i \to \infty$.
Thus,
$\exists \epsilon > 0$ such that $m(n_{i}) \geq (1+\epsilon)n_{i}$
for arbitrarily many values of $i$.
But then this contradicts Theorem~\ref{conn4},
where we showed that the probability of $P_{n,m}$ being connected is bounded away from $0$ for such $n.\!$~
\setlength{\unitlength}{.25cm}
\begin{picture}(1,1)
\put(0,0){\line(1,0){1}}
\put(0,0){\line(0,1){1}}
\put(1,1){\line(-1,0){1}}
\put(1,1){\line(0,-1){1}}
\end{picture}
\begin{displaymath}
\end{displaymath} \\

\newpage
\section{General Subgraphs \& Acyclic Subgraphs} \label{sub}

We have now finished looking at when $P_{n,m}$ will contain given components.
In the remainder of Part~\ref{I},
we shall instead look at the probability that $P_{n,m}$ will contain a given 
(not necessarily connected) \textit{subgraph}.

In this section, we will see (in Theorem~\ref{sub4}) that the limiting probability is~$1$ 
if $\liminf_{n \to \infty} \frac{m}{n} > 1$.
We shall then see (in Theorem~\ref{sub6})
that this result can be extended to $\liminf_{n \to \infty} \frac{m}{n} > 0$ if the subgraph is acyclic.

The proof of Theorem~\ref{sub4} will be in two parts, Lemmas~\ref{sub3} and~\ref{sub2},
dealing with the cases when $\frac{m}{n}$ is or isn't, respectively, bounded away from $3$.

We will first see (in Lemma~\ref{sub3}) that the part of the result dealing with the case when 
$\limsup_{n \to \infty} \frac{m}{n} < 3$ 
may be deduced easily from our knowledge of `appearances'.
We shall then prove (in Lemma~\ref{sub2}) the result for when $\frac{m}{n}$ is close to $3$ by modifying
appearance proofs to deal instead with the useful concept of `$6$-appearances'
(which we will define shortly).
This will involve first proving (in Lemmas~\ref{sub1} and~\ref{sub7})
some basic properties concerning $6$-appearances and triangulations.

By combining Theorem~\ref{sub4} with our knowledge of acyclic components
(from Section~\ref{cptlow}),
we shall then achieve (in Theorem~\ref{sub6}) our aforementioned extension for acyclic subgraphs.

\begin{figure} [ht]
\setlength{\unitlength}{1cm}
\begin{picture}(20,4.85)

\put(3.175,4.85){L\ref{gen2}}
\put(3.125,4.125){\textbf{L\ref{gen12}}}
\put(4.4,4.85){\textbf{L\ref{sub1}}}
\put(5.25,4.85){\textbf{L\ref{sub7}}}
\put(7.1,4.85){P\ref{ger L2.9}}
\put(6.3,4.85){P\ref{ger T2.1}}
\put(8.275,4.85){C\ref{tree6}}
\put(3.125,3.4){\textbf{L\ref{sub3}}}
\put(5.675,3.4){\textbf{L\ref{sub2}}}
\put(8.225,3.4){\textbf{L\ref{sub5}}}
\put(4.4,1.675){\textbf{T\ref{sub4}}}
\put(5.95,-0.325){\textbf{T\ref{sub6}}}

\put(4.75,1.8){\circle{0.8}}
\put(6.3,-0.2){\circle{0.8}}

\put(3.475,4.75){\vector(0,-1){0.3}}
\put(3.475,4.025){\vector(0,-1){0.3}}
\put(4.75,4.75){\vector(1,-1){1}}
\put(6.475,4.75){\vector(-1,-3){0.333}}
\put(5.575,4.75){\vector(1,-3){0.333}}
\put(7.3,4.75){\vector(-1,-1){1}}
\put(8.575,4.75){\vector(0,-1){1}}
\put(3.475,3.3){\vector(1,-1){1}}
\put(6.025,3.3){\vector(-1,-1){1}}
\put(8.575,3.3){\vector(-2,-3){2}}
\put(5.025,1.3){\vector(1,-1){1}}
\end{picture}

\caption{The structure of Section~\ref{sub}.}
\end{figure}
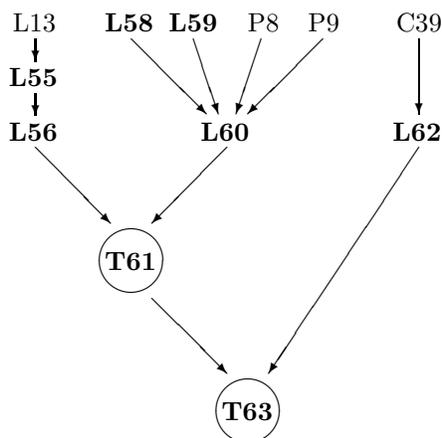

We start with the case when $\limsup_{n \to \infty} \frac{m}{n} < 3$,
which will involve recalling the concept of `appearances' from Definition~\ref{defapps}.
As usual,
we will actually aim to show a slightly stronger result than that advertised at the beginning of this section,
so we shall find it useful to first state vertex-disjoint versions of both Definition~\ref{defapps}
and Lemma~\ref{gen2}:

\begin{Definition}
Given a connected planar graph $H$
and another planar graph~$G$,
let \emph{\boldmath{${f_{H}^{\prime}}(G)$}}
denote the maximum size of a set of vertex-disjoint appearances of~$H$ in $G$.
\end{Definition}

\begin{Lemma} \label{gen12}
Let $H$ be a fixed connected planar graph on the vertices $\{1,2, \ldots ,h\}$,
let $b>1$ and $B<3$ be given constants,
and let $m=m(n) \in [bn,Bn]$ $\forall n$.
Then there exist constants $N(H,b,B)$ and $\beta(H,b,B) >0$ such that
\begin{displaymath}
\mathbf{P}[f_{H}^{\prime}(P_{n,m}) \leq \beta n] < e^{- \beta n} \textrm{ } \forall n \geq N.
\end{displaymath}
\end{Lemma}
\textbf{Proof}
For given $G$, let $\mathcal{W}$ be the collection of all sets $W \subset V(G)$ such that $H$ appears at $W$ in $G$.
As noted in \cite{ger}, each set $W \in \mathcal{W}$ meets at most $|H|-1$ other sets in $\mathcal{W}$
(since there are at most $|H|-1$ cut-edges in $H$
and each of these can have at most one `orientation' that provides an appearance of $H$),
and so we must be able to find a collection of at least 
$\frac{\mathcal{W}}{|H|}$ pairwise vertex disjoint sets in $\mathcal{W}$.
Thus, the result follows from Lemma~\ref{gen2} by taking $\beta = \frac{\alpha}{|H|}$.
$\phantom{qwerty}
\setlength{\unitlength}{.25cm}
\begin{picture}(1,1)
\put(0,0){\line(1,0){1}}
\put(0,0){\line(0,1){1}}
\put(1,1){\line(-1,0){1}}
\put(1,1){\line(0,-1){1}}
\end{picture}$ 
\\
\\

Our aforementioned result now follows fairly easily:

\begin{Lemma} \label{sub3}
Let $H$ be a (fixed) planar graph and let $m=m(n)$ satisfy \\
$1 < \liminf_{n \to \infty} \frac{m}{n} \leq \limsup_{n \to \infty} \frac{m}{n} < 3$.
Then $\exists \alpha > 0$ and $\exists N$ such that
\begin{eqnarray*}
& & \mathbf{P}\Big[\textrm{$P_{n,m}$
will \emph{not} have a set of $\geq \alpha n$ vertex-disjoint} \\
& & \phantom{wwwwwwwwi}\textrm{induced order-preserving copies of $H$}\Big] 
< e^{- \alpha n} 
\textrm{ } \forall n \geq N.
\end{eqnarray*}
\end{Lemma}
\textbf{Proof}
Let the components of $H$ be $H_{1}, H_{2}, \ldots, H_{k}$
and let us choose vertices
$\{ v_{1}, v_{2}, , \ldots, v_{k} \} \subset V(H)$
such that $v_{i} \in V(H_{i})$ $\forall i$.

Without loss of generality, 
we may assume that
$V(H)= \{1,2, \ldots, |H| \}$.
Let us define $H^{\prime}$ to be the graph with vertex set
$V(H^{\prime}) = \{1,2, \ldots, |H|, |H|+1 \}$
and edge set
$E(H^{\prime}) = E(H) \cup \{ (v_{1},|H|+1), (v_{2},|H|+1), (v_{3},|H|+1) \}$.
\begin{figure} [ht]
\setlength{\unitlength}{1cm}
\begin{picture}(20,4.5)(-0.2,-0.5)

\put(4.9,-0.5){\circle{1}}
\put(4.9,2){\circle{1}}
\put(4.9,3.5){\circle{1}}

\put(5.4,-0.5){\line(2,5){1}}
\put(5.4,2){\line(1,0){1}}
\put(5.4,3.5){\line(2,-3){1}}

\put(5.4,-0.5){\circle*{0.1}}
\put(5.4,2){\circle*{0.1}}
\put(5.4,3.5){\circle*{0.1}}
\put(6.4,2){\circle*{0.1}}

\put(4.7,3.4){$H_{1}$}
\put(4.7,1.9){$H_{2}$}
\put(4.7,-0.6){$H_{k}$}
\put(6.5,2.1){\scriptsize{$|H|+1$}}
\put(5.5,3.6){\footnotesize{$v_{1}$}}
\put(5.5,2.1){\footnotesize{$v_{2}$}}
\put(5.5,-0.6){\footnotesize{$v_{k}$}}

\put(4.9,0.25){\circle*{0.05}}
\put(4.9,0.5){\circle*{0.05}}
\put(4.9,0.75){\circle*{0.05}}
\put(4.9,1){\circle*{0.05}}
\put(4.9,1.25){\circle*{0.05}}
\end{picture}

\caption{The graph $H^{\prime}$.}
\end{figure}
\\ Then $H^{\prime}$ is a planar graph containing an induced order-preserving copy of $H$.
Thus, it suffices to prove the result for $H^{\prime}$.
But this follows from Lemma~\ref{gen12}.
~$\setlength{\unitlength}{.25cm}
\begin{picture}(1,1)
\put(0,0){\line(1,0){1}}
\put(0,0){\line(0,1){1}}
\put(1,1){\line(-1,0){1}}
\put(1,1){\line(0,-1){1}}
\end{picture}$
\\
\\
\\

\begin{figure} [ht]
\setlength{\unitlength}{1cm}
\begin{picture}(20,2.8)(-0.3,-0.7)
\put(4.2,0.7){\oval(2,2.8)}
\put(7.2,0.7){\oval(1,1.4)}
\put(7.1,0.6){$W$}
\put(4.8,1.3){$v_{3}$}
\put(4.8,0.6){$v_{1}$}
\put(4.8,0){$v_{2}$}
\put(6.8,1){\small{$r_{1}$}}
\put(7.75,0.85){\small{$r_{3}$}}
\put(6.8,0.3){\small{$r_{2}$}}
\put(5.2,1.3){\circle*{0.1}}
\put(5.2,0.1){\circle*{0.1}}
\put(5.2,0.7){\circle*{0.1}}
\put(6.7,0.4){\circle*{0.1}}
\put(7.7,0.7){\circle*{0.1}}
\put(6.7,1){\circle*{0.1}}
\put(5.2,1.3){\line(5,-1){1.5}}
\put(5.2,0.7){\line(5,1){1.5}}
\put(5.2,0.7){\line(5,-1){1.5}}
\put(5.2,0.1){\line(5,1){1.5}}
\put(7.7,1.2){\oval(1,1)[r]}
\put(7.7,0.2){\oval(1,1)[r]}
\put(5.2,1.3){\line(5,1){2}}
\put(5.2,0.1){\line(5,-1){2}}
\put(7.2,1.7){\line(1,0){0.5}}
\put(7.2,-0.3){\line(1,0){0.5}}
\end{picture} 

\caption{A $6$-appearance at $W$.} \label{6app}
\end{figure}
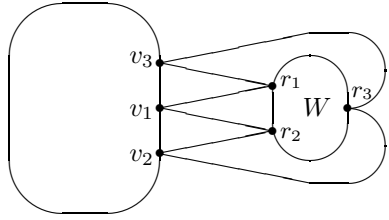

We shall now start working towards the case when $\frac{m}{n}$ is close to $3$.
As mentioned, the proof will involve the concept of `$6$-appearances':

\begin{Definition}
We say that a graph $H$ with vertex set $\{ 1,2, \ldots, h \}$ \textbf{\emph{\boldmath{$6$}-appears}} 
at $W \subset V(G)$ if
(a) the increasing bijection from $\{ 1,2, \ldots, h \}$ to $W$
gives an isomorphism between $H$ and the induced subgraph $G[W]$ of $G$;
and (b) there are exactly six edges in $G$ between $W$ and the rest of $G$,
and these are of the form 
$E_{W} =~\{ r_{1}v_{1}, v_{1}r_{2}, r_{2}v_{2}, v_{2}r_{3}, r_{3}v_{3}, v_{3}r_{1} \}$,
where $\{ r_{1}, r_{2}, r_{3} \} \subset W$,
$\{ v_{1}, v_{2}, v_{3} \} \subset~V(G)~\setminus~W$,
and $\{ v_{2}, v_{3} \} \subset \Gamma (v_{1})$ (see Figure~\ref{6app}).
We shall call $E(G[W]) \cup E_{W}$ the \textbf{\emph{total edge set}} of the $6$-appearance. 
\end{Definition}

\phantom{p}

The following result about $6$-appearances shall be useful:

\begin{Lemma} \label{sub1}
The total edge set of a $6$-appearance of a $2$-vertex-connected graph of order $|T|$
will intersect (i.e. have an edge in common with)
the total edge set of at most 
$2 \left( ^{|T|+3} _{\phantom{qw}3} \right)$ other $6$-appearances of connected graphs of order $|T|$.
\end{Lemma}
\textbf{Proof}
Suppose we have a $6$-appearance of a $2$-vertex-connected graph of order $|T|$ at $W \subset V(G)$, 
as in Figure~\ref{6app}.
Note that the vertices $\{ v_{1}, v_{2}, v_{3} \}$ form a $3$-vertex-cut.
Suppose that $G$ also contains another $6$-appearance, at $W_{2}$, of a connected graph of order $|T|$
and let $\{ w_{1}, w_{2}, w_{3} \}$ 
denote the $3$-vertex-cut in $V(G) \setminus W_{2}$ that is associated with $W_{2}$. \\

(a) Suppose $\{ u_{1}, u_{2}, u_{3} \} \subset V(G) \setminus W$.
If $\{ v_{1}, v_{2}, v_{3} \} = \{ u_{1}, u_{2}, u_{3} \}$,
then either $W_{2} = W$ or $W_{2} \subset V(G) \setminus (W \cup \{ v_{1}, v_{2}, v_{3} \})$,
in which case the total edge set of $W_{2}$ would not meet the total edge set of $W$.

If $\{ v_{1}, v_{2}, v_{3} \} \neq \{ u_{1}, u_{2}, u_{3} \}$,
then without loss of generality $v_{1} \notin \{ u_{1}, u_{2}, u_{3} \}$ 
and so 
$W \cup v_{1}$ will all be in one component of the graph $G \setminus \{ u_{1},u_{2},u_{3} \}$.
Thus, either $(W \cup v_{1}) \subset W_{2}$,
in which case $|W_{2}|>|T|$,
or $W_{2} \subset V(G) \setminus~(W \cup~\{ v_{1}, v_{2}, v_{3} \})$
(since $v_{2}$ and $v_{3}$ must each either be in the same component of $G \setminus \{ u_{1},u_{2},u_{3} \}$
as $W \cup v_{1}$ or must belong to $\{ u_{1},u_{2},u_{3} \}$,
in which case the total edge set of $W_{2}$ would not meet the total edge set of $W$. \\

(b) Suppose $\{ u_{1}, u_{2}, u_{3} \} \subset (W \cup \{ v_{1}, v_{2}, v_{3} \})$.
Then there are at most $\left( ^{|T|+3} _{\phantom{qw}3} \right)$ choices for $\{ u_{1}, u_{2}, u_{3} \}$.
Note that once we have fixed on a particular choice of $\{ u_{1},u_{2},u_{3} \}$
then there can only be at most $2$ choices for $W_{2}$,
since if $Z_{1}$, $Z_{2}$ and $Z_{3}$ were all possibilities for $W_{2}$
then the graph obtained from $G$ by condensing each of $Z_{1}$, $Z_{2}$ and $Z_{3}$
down to a single point
would contain $K_{3,3}$
(since each of these points would be adjacent to each of the $u_{i}$'s),
and this would contradict planarity.
Hence,
if $\{ u_{1}, u_{2}, u_{3} \} \subset (W \cup \{ v_{1}, v_{2}, v_{3} \})$
then we have at most~$2 \left( ^{|T|+3} _{\phantom{qw}3} \right)$ choices for $W_{2}$ in total. \\

(c) If neither (a) nor (b) holds,
$\exists u_{i} \in W$ and $\exists u_{j} \in V(G) \setminus (W \cup~\{ v_{1}, v_{2}, v_{3} \})$.
By the definition of a $6$-appearance, we must have $\{ u_{2}, u_{3} \} \subset \Gamma(u_{1})$.
Since $u_{i}$ and $u_{j}$ are not adjacent, 
we must have $u_{1} \in \{ v_{1}, v_{2}, v_{3} \}$ and $\{ u_{i}, u_{j} \} = \{ u_{2}, u_{3} \}$. 

Let $C(W \setminus u_{i})$ denote the component containing $W \setminus u_{i}$ after the $3$-vertex-cut
defined by $\{ u_{1}, u_{2}, u_{3} \}$
(this is well-defined,
since the $2$-vertex-connectivity of $W$ implies that $W \setminus u_{i}$ is connected).
By the definition of a $6$-appearance,
$\exists x \in W_{2}$ such that $x \in \Gamma (u_{i}) \cap \Gamma (u_{j})$.
Since $u_{i} \in W$ and $u_{j} \in V(G) \setminus (W \cup \{ v_{1}, v_{2}, v_{3} \} )$,
we must have $x \in  \{ v_{1}, v_{2}, v_{3} \}$.
But then $x$ is certainly adjacent to a vertex in~$W \setminus u_{i}$, 
and so $x \in C(W \setminus u_{i})$.

Since $x \in W_{2}$
and $x \in C(W \setminus u_{i})$, 
we must have
$C(W \setminus u_{i}) = W_{2}$.
But $W \setminus u_{i}$ will still be connected to $\{ v_{1}, v_{2}, v_{3} \} \setminus u_{1}$ 
in $G \setminus \{ u_{1},u_{2},u_{3} \}$,
so $|C(W \setminus u_{i})| \geq |W \setminus u_{i}| + | \{v_{1},v_{2},v_{3} \} \setminus u_{1} | = |T| + 1 > |W_{2}|.$
\phantom{qwerty}
\setlength{\unitlength}{.25cm}
\begin{picture}(1,1)
\put(0,0){\line(1,0){1}}
\put(0,0){\line(0,1){1}}
\put(1,1){\line(-1,0){1}}
\put(1,1){\line(0,-1){1}}
\end{picture}
\begin{displaymath}
\end{displaymath}

Before proceeding with Lemma~\ref{sub2},
we also need to note the following useful result:

\begin{Lemma} \label{sub7}
Let $H$ be a (fixed) planar graph.
Then there exists a triangulation~$T$ with $|T| \geq \max \{ |H|+1, 4 \}$
such that $T$ contains an induced order-preserving copy of $H$. 
\end{Lemma}
\textbf{Proof}
Let $T_{1}$ be a triangulation with
$|T_{1}| \geq \max \{ |H|+1, 4 \}$
such that $L \subset T_{1}$, where $L$ is an order-preserving copy of $H$.
For each edge $e \in E(T_{1}) \setminus E(L)$ such that $e$ is between two vertices in $V(L)$,
let us subdivide $e$ by inserting a new vertex on it.

Let the new graph (which is not necessarily a triangulation)
be called $T_{2}$ and note that $L$ is an \textit{induced} subgraph of $T_{2}$.
Since each face of $T_{1}$ was a triangle, 
each face of $T_{2}$ will be one of the following, 
where $\bullet$ denotes a new vertex:

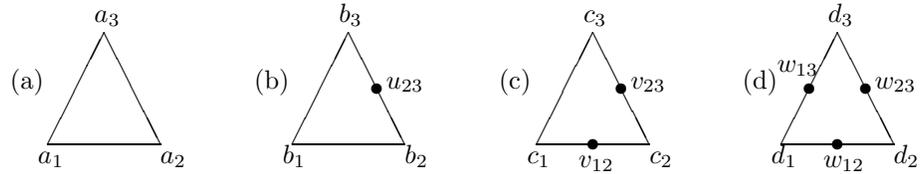
\begin{figure} [ht]
\begin{picture}(90,6.5)
\put(0,3){(a)}
\put(2,0){\line(1,0){6}}
\put(2,0){\line(1,2){3}}
\put(5,6){\line(1,-2){3}}
\put(1.5,-1){$a_{1}$}
\put(8,-1){$a_{2}$}
\put(4.5,6.5){$a_{3}$}
\put(13,3){(b)}
\put(19.5,3){\circle*{0.5}}
\put(20,3){$u_{23}$}
\put(15,0){\line(1,0){6}}
\put(15,0){\line(1,2){3}}
\put(18,6){\line(1,-2){3}}
\put(14.5,-1){$b_{1}$}
\put(21,-1){$b_{2}$}
\put(17.5,6.5){$b_{3}$}
\put(26,3){(c)}
\put(32.5,3){\circle*{0.5}}
\put(33,3){$v_{23}$}
\put(31,0){\circle*{0.5}}
\put(30.25,-1.2){$v_{12}$}
\put(28,0){\line(1,0){6}}
\put(28,0){\line(1,2){3}}
\put(31,6){\line(1,-2){3}}
\put(27.5,-1){$c_{1}$}
\put(34,-1){$c_{2}$}
\put(30.5,6.5){$c_{3}$}
\put(38.95,3){(d)}
\put(42.5,3){\circle*{0.5}}
\put(44,0){\circle*{0.5}}
\put(45.5,3){\circle*{0.5}}
\put(41,0){\line(1,0){6}}
\put(41,0){\line(1,2){3}}
\put(44,6){\line(1,-2){3}}
\put(40.5,-1){$d_{1}$}
\put(47,-1){$d_{2}$}
\put(43.5,6.5){$d_{3}$}
\put(43.25,-1.2){$w_{12}$}
\put(40.8,3.9){$w_{13}$}
\put(46,3){$w_{23}$}
\end{picture} 

\caption{The possibilities for faces of $T_{2}$.} 
\end{figure}

We can create a new \textit{triangulation} $T \supset T_{2}$ by replacing the faces of type (b),~(c)~and~(d)
with (b$^{\prime}$), (c$^{\prime}$) and (d$^{\prime}$):

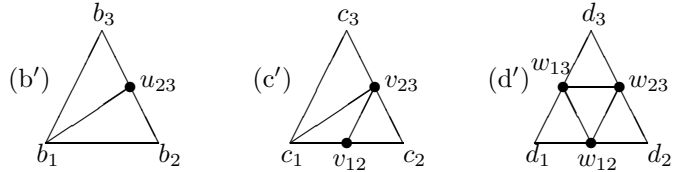
\begin{figure} [ht]
\begin{picture}(90,6.5)
\put(6.5,3){(b$^{\prime}$)}
\put(13,3){\circle*{0.5}}
\put(8.5,0){\line(1,0){6}}
\put(8.5,0){\line(1,2){3}}
\put(11.5,6){\line(1,-2){3}}
\put(8.5,0){\line(3,2){4.5}}
\put(8,-1){$b_{1}$}
\put(14.5,-1){$b_{2}$}
\put(11,6.5){$b_{3}$}
\put(13.5,3){$u_{23}$}
\put(19.5,3){(c$^{\prime}$)}
\put(26,3){\circle*{0.5}}
\put(24.5,0){\circle*{0.5}}
\put(21.5,0){\line(1,0){6}}
\put(21.5,0){\line(1,2){3}}
\put(24.5,6){\line(1,-2){3}}
\put(21.5,0){\line(3,2){4.5}}
\put(24.5,0){\line(1,2){1.5}}
\put(26.5,3){$v_{23}$}
\put(23.75,-1.2){$v_{12}$}
\put(21,-1){$c_{1}$}
\put(27.5,-1){$c_{2}$}
\put(24,6.5){$c_{3}$}
\put(32,3){(d$^{\prime}$)}
\put(36,3){\circle*{0.5}}
\put(37.5,0){\circle*{0.5}}
\put(39,3){\circle*{0.5}}
\put(34.5,0){\line(1,0){6}}
\put(34.5,0){\line(1,2){3}}
\put(37.5,6){\line(1,-2){3}}
\put(37.5,0){\line(1,2){1.5}}
\put(37.5,0){\line(-1,2){1.5}}
\put(36,3){\line(1,0){3}}
\put(34,-1){$d_{1}$}
\put(40.5,-1){$d_{2}$}
\put(37,6.5){$d_{3}$}
\put(36.75,-1.2){$w_{12}$}
\put(34.3,3.9){$w_{13}$}
\put(39.5,3){$w_{23}$}
\end{picture} 

\caption{Our new faces.}
\end{figure}

Note that we haven't inserted the same edge in more than one face.
This is because there is only one triangle in $T_{1}$ containing both the edge $b_{2}b_{3}$ and the vertex $b_{1}$
(and hence only one face, since $|T_{1}| \geq 4$)
and only one triangle (and hence one face) containing both the edges $c_{2}c_{3}$ and $c_{1}c_{2}$.
Hence, the edges $b_{1}u_{23}$ and $v_{12}v_{23}$ are both inserted in only one face of $T_{2}$, 
and since all inserted edges are of this form we are okay.

Recall that $T_{2} \subset T$ and that $L$ is an \textit{induced} subgraph of $T_{2}$.
Note that all edges in $E(T) \setminus E(T_{2})$ are incident with at least one new vertex,
and so none can be between two vertices in $L$.
Thus, $L$ is an \textit{induced} subgraph of $T$.
$\phantom{qwerty}
\setlength{\unitlength}{.25cm}
\begin{picture}(1,1)
\put(0,0){\line(1,0){1}}
\put(0,0){\line(0,1){1}}
\put(1,1){\line(-1,0){1}}
\put(1,1){\line(0,-1){1}}
\end{picture}$
\\
\\
\\

We now come to our aforementioned lemma for the case when $\frac{m}{n}$ is close to~$3$.
The proof is based on that of Proposition~\ref{ger T3.1},
which is itself based on the proof of Theorem $4.1$ in \cite{mcd}.

\begin{Lemma} \label{sub2}
Let $H$ be a (fixed) planar graph.
Then $\exists A(H) < 3$ such that if
$m=m(n) \in [An,3n-6]$ for all large $n$,
then $\exists  \alpha >0$ and $\exists N$ such that
\begin{eqnarray}
& & \mathbf{P}\Big[\textrm{$P_{n,m}$
\emph{won't} have a set of $> \alpha n$ vertex-disjoint} \nonumber \\
& & \phantom{wwwww}\textrm{induced order-preserving copies of $H$}\Big]  
< e^{- \alpha n} 
\textrm{ } \forall n \geq N. \label{eq:dag}
\end{eqnarray}
\end{Lemma}
\textbf{Sketch of Proof} 
We choose a specific $\alpha$ and suppose that $(\ref{eq:dag})$ is false for~$n=k$,
where $k$ is suitably large.
We let $q= \frac{m(k)}{k}$
and, by using the results of~\cite{ger},
we find that (for a given $\epsilon > 0$)
\begin{eqnarray}
(1-\epsilon)^{n}n!(\gamma (q))^{n} \leq |\mathcal{P}(n, \lfloor qn \rfloor)| \leq 
(1+ \epsilon)^{n}n!(\gamma (q))^{n} \phantom{ww}\forall n \geq k.
\label{eq:*} 
\end{eqnarray}
Using the assumption that $(\ref{eq:dag})$ is false for $n=k$,
together with the left-hand inequality in $(\ref{eq:*})$,
it follows that there are many graphs in $\mathcal{P}(k,qk)$
without a set of $> \alpha k$ vertex-disjoint induced order-preserving copies of $H$.

We attach $6$-appearances of $T$ and $T^{\prime}$
(which both contain induced order-preserving copies of $H$ and have slightly more and slightly less, respectively,
than the appropriate ratio of edges to vertices)
to these graphs to construct a great many graphs in $\mathcal{P}((1+\delta)k, \lfloor q(1+ \delta)k \rfloor)$,
for some $\delta > 0$.

Lemma~\ref{sub1} is used to show that the original graphs didn't have many 
$6$-appearances of $T$ and $T^{\prime}$,
and this is used to show that there is not much double-counting of our created graphs.
Thus, we obtain a contradiction to the right-hand inequality in $(\ref{eq:*})$ for $n=(1+\delta)k$. \\
\\
\textbf{Full Proof}
By Lemma~\ref{sub7},
we know that there exists a triangulation~$T$ with $|T| \geq \max \{ |H|+1, 4 \}$
such that $T$ contains an induced order-preserving copy of~$H$. 
Let $\beta = e^{2} \gamma_{l}^{|T|} \left( 4 \left(^{|T|+3}_{\phantom{qw}3} \right) +2 \right) |T|!$
and let $\alpha$ be a fixed constant in $\left( 0, \frac{1}{\beta} \right)$.
Thus, we have $\alpha \beta <1$
and so $\exists \epsilon \in \left( 0, \frac{1}{3} \right)$
such that $(\alpha \beta)^{\alpha} = 1 - 3 \epsilon$.

Let $A = \max \left \{ \frac{3|T|-1}{|T|}, \frac{11}{4} \right \}$
(this value is chosen for ease with later parts of the proof)
and let $N_{1}$ be such that $m(n) \in [An,3n-6]$~$\forall n \geq N_{1}$.
Since~$A>1$, 
it follows from Proposition~\ref{ger L2.9}
(choosing $a \in (1,A)$)
that $\exists N_{2}$ such that for all $n \geq N_{2}$ and all
$s \in [\lfloor An \rfloor, 3n-6]$ we have
$\left|\left(\frac{|\mathcal{P}(n,s)|}{n!} \right)^{1/n} - 
\gamma \left(\frac{s}{n} \right)\right| < \frac{\epsilon e}{2}$.
Also, note from the uniform continuity of
Proposition~\ref{ger T2.1}
that $\exists N_{3}$ such that
$\left|\gamma \left(\frac{\lfloor rn \rfloor}{n} \right)- \gamma (r)\right| <
\frac{\epsilon e}{2}$~$\forall r \in [1,3]$~$\forall n \geq N_{3}$.

Suppose $(\ref{eq:dag})$ doesn't hold for some $k \!>\! N \!=\! \max \{\! N_{1}, N_{2}, N_{3} \!\}$
and let~$q =~\!\!\frac{m(k)}{k}$.
Then, by the previous paragraph, we have both
$\left|\left(\frac{|\mathcal{P}(n, \lfloor qn \rfloor)|}{n!} \right)^{1/n} - 
\gamma \left( \frac {\lfloor qn \rfloor}{n} \right)\right| < \frac{\epsilon e}{2}$~$\forall n \geq k$
and
$\left|\gamma \left(\frac{\lfloor qn \rfloor}{n} \right)- \gamma (q)\right| < 
\frac{\epsilon e}{2} $~$\forall n \geq k$.
We recall from Proposition~\ref{ger T2.1} that 
$\gamma (q) \geq e$,
so putting everything together we must have 
$\left|\left(\frac{|\mathcal{P}(n, \lfloor qn \rfloor)|}{n!} \right)^{1/n} \!-\! \gamma (q)\right| < 
\epsilon e \leq \epsilon \gamma (q)$~$\forall n \geq k$.
Thus, we have (\ref{eq:*}).
\\

Let $\mathcal{G}_{k}$ denote the set of graphs in 
$\mathcal{P}(k,m(k)) = \mathcal{P}(k, \lfloor qk \rfloor)$
such that $G \in~\mathcal{G}_{k}$ iff 
$G$ does not have a set of $> \alpha k$ vertex-disjoint induced order-preserving copies of $H$.
Then $|\mathcal{G}_{k}| \geq e^{- \alpha k} f(k, \lfloor qk \rfloor) 
\geq e^{- \alpha k} (1 - \epsilon)^{k} (\gamma (q))^{k} k!$,
using $(\ref{eq:*})$ and our assumption that $(\ref{eq:dag})$ does not hold for $k$.

Recall that $T$ is a triangulation with $|T| \geq |H|+1$
such that $T$ contains an induced order-preserving copy of $H$
and let $T^{\prime}$ be the result of deleting one edge from $T$ that does not interfere with this copy of $H$
(this is possible, since $|T| \geq |H|+1$).
Starting with graphs in $\mathcal{G}_{k}$,
we shall construct graphs in $\mathcal{P}((1+\delta)k, \lfloor q(1+~\delta)k \rfloor)$,
where $\delta = \frac{ \lceil \alpha k \rceil |T| }{k}$,
by attaching $k_{1}$ $6$-appearances of $T$ and $k_{2}$ $6$-appearances of $T^{\prime}$.
In order that our constructed graphs are in 
$\mathcal{P}((1+\delta)k, \lfloor q(1+\delta)k \rfloor)$,
we shall need to achieve the correct balance of $k_{1}$ and~$k_{2}$: \\
\\

Let us define
$k_{1}$ to be $\lfloor q(k+ \lceil \alpha k \rceil |T|) \rfloor - \lfloor qk \rfloor 
- \lceil \alpha k \rceil (3|T|-1)$
and $k_{2}$ to be $\lfloor qk \rfloor +\lceil \alpha k \rceil 3|T|
-\lfloor q(k+\lceil \alpha k \rceil |T|) \rfloor$.
Then $k_{1}$ and $k_{2}$ are integers such that $k_{1} + k_{2} = \lceil \alpha k \rceil$.

Recall $q \geq A \geq \frac{3|T|-1}{|T|}$.
Thus,
$q \lceil \alpha k \rceil |T| \geq~(3|T|-~1) \lceil \alpha k \rceil$
and so, since the right-hand-side is an integer,
we have 
$\lfloor q \lceil \alpha k \rceil |T| \rfloor \geq \lceil \alpha k \rceil (3|T|-~1)$.
Thus,
\begin{eqnarray*}
k_{1} & \geq & \lfloor qk \rfloor + \lfloor q \lceil \alpha k \rceil |T| \rfloor - \lfloor qk \rfloor 
- \lceil \alpha k \rceil (3|T|-1) \\
& = & \lfloor q \lceil \alpha k \rceil |T| \rfloor 
- \lceil \alpha k \rceil (3|T|-1) \\
& \geq & 0.
\end{eqnarray*}
Also, note that
\begin{eqnarray*}
k_{2} & \geq & \lfloor qk \rfloor + \lceil \alpha k \rceil 3|T|
- \lfloor qk+ 3 \lceil \alpha k \rceil |T| \rfloor \\
& = & \lfloor qk \rfloor + \lceil \alpha k \rceil 3|T|
- \lfloor qk \rfloor + 3 \lceil \alpha k \rceil |T| \\
& = & 0.
\end{eqnarray*}
Thus, $k_{1}$ and $k_{2}$ are \textit{positive} integers that sum to $\lceil \alpha k \rceil$.

Finally, note that the total edge set of a $6$-appearance of $T$ will have size~$3|T|$
and the total edge set of a $6$-appearance of $T^{\prime}$ will have size $3|T|-1$,
and note that
\begin{eqnarray*}
\lfloor qk \rfloor + k_{1} 3 |T| + k_{2} (3|T|-1) 
& = & \lfloor qk \rfloor + k_{1} 3 |T| + ( \lceil \alpha k \rceil - k_{1}) (3|T|-1) \\
& = & \lfloor qk \rfloor + \lceil \alpha k \rceil (3|T|-1) + k_{1} \\
& = & \lfloor q(k+ \lceil \alpha k \rceil |T|) \rfloor \\
& = & \lfloor q(1+\delta)k \rfloor. \\
\end{eqnarray*} 

We shall now construct graphs in $\mathcal{P}((1+\delta)k, \lfloor q(1+\delta)k \rfloor)$:

Choose $\delta k$ special vertices 
(we have $\left( ^{(1+\delta)k} _{\phantom{qq} \delta k} \right)$ choices for these)
and partition them into $\lceil \alpha k \rceil$ unordered blocks of size $|T|$
(we have 
$\left( ^{\phantom{qqw}\delta k}_{|T|, \ldots, |T|} \right) \frac{1}{\left \lceil \alpha k \right \rceil !}$ 
choices for this).
Divide the blocks into two sets of sizes $k_{1}$ and $k_{2}$.
On each of the first $k_{1}$ blocks,
we put a copy of $T$ such that the increasing bijection from $\{ 1,2, \ldots, |T| \}$ to the block 
is an isomorphism between $T$ and this copy.
We do the same for the set of $k_{2}$ blocks, except with $T^{\prime}$ instead of $T$.

On the remaining (i.e.~non-special) vertices,
choose a planar graph $L$ of size $\lfloor qk \rfloor = m$
that does not have a set of 
$> \alpha k$ vertex-disjoint induced order-preserving copies of $H$.
We have 
$|\mathcal{G}_{k}| \geq e^{- \alpha k} (1 - \epsilon)^{k} (\gamma (q))^{k} k!$ choices for this.
Note that (assuming $k \geq 3$) $L$ may be extended to a triangulation by inserting 
$3k-6-m > (3-A)k-6$ `phantom' edges.
This triangulation will contain $2k-4$ triangles that are faces.
Each of our phantom edges is in exactly two faces of this triangulation,
so when we remove these phantom edges we are left with an embedding of $L$ which contains 
$\geq 2k-4-2((3-A)k-6) \geq k(2A-4)$
triangles that are faces.
We may assume that $k$ is large enough
that $k(2A-4) - \lceil \alpha k \rceil > k$,
since $A \geq \frac{11}{4}$ and $\alpha < \frac{1}{\beta} < \frac{1}{2}$.

We may attach our copies of $T$ and $T^{\prime}$ inside $\lceil \alpha k \rceil$ of these triangles
in such a way that we create $6$-appearances of $T$ and $T^{\prime}$.
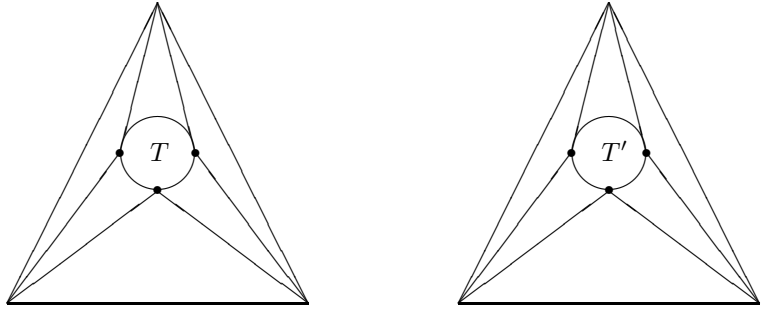
\begin{figure} [ht]
\setlength{\unitlength}{1cm}
\begin{picture}(90,4.25)

\put(1,0){\line(1,0){4}}
\put(1,0){\line(1,2){2}}
\put(5,0){\line(-1,2){2}}
\put(3,2){\circle{1}}
\put(1,0){\line(4,3){2}}
\put(1,0){\line(3,4){1.5}}
\put(5,0){\line(-4,3){2}}
\put(5,0){\line(-3,4){1.5}}
\put(3,4){\line(-1,-4){0.5}}
\put(3,4){\line(1,-4){0.5}}
\put(3,1.5){\circle*{0.1}}
\put(2.5,2){\circle*{0.1}}
\put(3.5,2){\circle*{0.1}}
\put(2.9,1.9){$T$}

\put(7,0){\line(1,0){4}}
\put(7,0){\line(1,2){2}}
\put(11,0){\line(-1,2){2}}
\put(9,2){\circle{1}}
\put(7,0){\line(4,3){2}}
\put(7,0){\line(3,4){1.5}}
\put(11,0){\line(-4,3){2}}
\put(11,0){\line(-3,4){1.5}}
\put(9,4){\line(-1,-4){0.5}}
\put(9,4){\line(1,-4){0.5}}
\put(9,1.5){\circle*{0.1}}
\put(8.5,2){\circle*{0.1}}
\put(9.5,2){\circle*{0.1}}
\put(8.9,1.9){$T^{\prime}$}

\end{picture} 

\caption{Creating $6$-appearances of $T$ and $T^{\prime}$ inside facial triangles.}
\end{figure}
\\We have at least 
$\left( ^{k(2A-4)}_{\phantom{qq}\lceil \alpha k \rceil} \right)$
choices for these triangles,
and we then have $\left \lceil \alpha k \right \rceil !$
choices for which copies of $T$ and $T^{\prime}$ to attach within which triangles.

Thus, for each choice of special vertices and each choice of $L$, we may construct at least
\begin{eqnarray*}
\left( ^{\phantom{qqw}\delta k}_{|T|, \ldots, |T|} \right) \left( ^{k(2A-4)}_{\phantom{qq}\lceil \alpha k \rceil} \right)
& \geq & \left( ^{\phantom{qqw}\delta k}_{|T|, \ldots, |T|} \right)
\frac{\left( k(2A-4) - \lceil \alpha k \rceil \right)^{\lceil \alpha k \rceil}}{\lceil \alpha k \rceil !} \\
& \geq & \left( ^{\phantom{qqw}\delta k}_{|T|, \ldots, |T|} \right)
\frac{k^{\lceil \alpha k \rceil}}{\lceil \alpha k \rceil !} \\
& = & \frac{(\delta k)! k^{\lceil \alpha k \rceil}}{(|T|!)^{\lceil \alpha k \rceil} \lceil \alpha k \rceil !} \\
& \geq & \frac{(\delta k)!} {(|T|! \alpha)^{\lceil \alpha k \rceil}} 
\textrm{ for $k$ large enough that } 
\left \lceil \alpha k \right \rceil ! \!\leq\! \left( \alpha k \right)^{\left \lceil \alpha k \right \rceil} 
\end{eqnarray*}
graphs in $\mathcal{P}((1+\delta)k, \lfloor q(1+\delta)k \rfloor)$.
Hence,
we may construct at least \\
$\left( ^{(1+\delta)k}_{\phantom{qq} \delta k} \right) e^{- \alpha k} (1 -~\epsilon)^{k} (\gamma(q))^{k} k! 
\frac{(\delta k)!}{(|T|! \alpha)^{\lceil \alpha k \rceil}}$
(not necessarily distinct) graphs in $\mathcal{P}((1+\delta)k, \lfloor q(1+\delta)k \rfloor)$ in total. \\

We shall now consider the amount of double-counting:

Recall that $L$ did not have a set of $> \alpha k$ vertex-disjoint induced order-preserving copies of~$H$
and that $T$ and $T^{\prime}$ both contain induced order-preserving copies of $H$.
Note that $T$ and $T^{\prime}$ are both $2$-vertex-connected, since $|T| \geq 4$.
Thus, by Lemma~\ref{sub1}, it must be that $L$ contained at most 
$\left( 2 \left( ^{|T|+3} _{\phantom{qw}3} \right) +1 \right) \alpha k$  
$6$-appearances of $|T|$ and $T^{\prime}$.

When we deliberately attach a copy of $T$ or $T^{\prime}$,
the number of `accidental' \\ $6$-appearances of $T$ of $T^{\prime}$ that we create in the graph
will be at most $2 \left( ^{|T|+3} _{\phantom{qw}3} \right)$, again using Lemma~\ref{sub1}.
Thus, the number of $6$-appearances of $T$ or $T^{\prime}$ will increase by at most 
$2 \left( ^{|T|+3} _{\phantom{qw}3} \right) + 1$ each time we attach one of our blocks.
Thus, our created graph will have 
$\leq \left( 2 \left( ^{|T|+3} _{\phantom{qw}3} \right) +1 \right) \alpha k + 
\left( 2 \left( ^{|T|+3} _{\phantom{qw}3} \right) + 1 \right) \lceil \alpha k \rceil
\leq \left( 4 \left( ^{|T|+3} _{\phantom{qw}3} \right) + 2 \right) \lceil \alpha k \rceil$
$6$-appearances of $T$ or $T^{\prime}$.

Let $x=4\left( ^{|T|+3} _{\phantom{qw}3} \right) + 2$.
Then, given one of our constructed graphs, we have at most
$\left( ^{x \lceil \alpha k \rceil}_{\phantom{i}\lceil \alpha k \rceil} \right) \leq (xe)^{\lceil \alpha k \rceil}$
choices for which were the special vertices.
Once we have identified these, we know what $L$ was.
Thus, each graph is constructed at most $(xe)^{\lceil \alpha k \rceil}$ times. \\
\\

Therefore, 
we find that
the number of distinct graphs that we have created in
$\mathcal{P}((1+\delta)k, \lfloor q(1+\delta)k \rfloor)$
is at least 
\begin{eqnarray*}
& & \left( ^{(1+\delta)k}_{\phantom{qq} \delta k} \right) e^{- \alpha k} (1 - \epsilon)^{k} (\gamma(q))^{k} k! 
\frac{(\delta k)!}{(|T|! \alpha)^{\lceil \alpha k \rceil}} (xe)^{- \lceil \alpha k \rceil} \\
& \geq & ((1+\delta)k)! (\gamma(q))^{(1+\delta)k} (1-\epsilon)^{k} 
\left( e^{2} (\gamma(q))^{|T|}x|T|! \alpha \right) ^{- \lceil \alpha k \rceil} 
\textrm{ since } \delta k \!=\! T \lceil \alpha k \rceil \\
& \geq & ((1+\delta)k)! (\gamma(q))^{(1+\delta)k} (1-\epsilon)^{k} 
(\alpha \beta)^{- \lceil \alpha k \rceil} \\
& &
\textrm{by definition of } \beta \textrm{ (and Proposition~\ref{ger T2.1})} \\
& \geq & f((1+\delta)k) (1+\epsilon)^{-(1+\delta)k} (1-\epsilon)^{k} (1-3\epsilon)^{-k}
\phantom{w} \textrm{ by } (\ref{eq:*}) \textrm{ with } n=(1+\delta)k \\
& \geq & f((1+\delta)k)  \left( \frac{(1-\epsilon)}{(1-3\epsilon)(1+\epsilon)^{2}} \right) ^{k} \\
& & \textrm{since we may assume $k$ is large enough that } \delta<1 \\
& > & f((1+\delta)k) 
\phantom{ww} \textrm{ since } (1-3\epsilon)(1+\epsilon)^{2}=1-\epsilon -5\epsilon^{2} - 3\epsilon^{3}
\end{eqnarray*}
which is a contradiction.
$\phantom{qwerty}
\setlength{\unitlength}{.25cm}
\begin{picture}(1,1)
\put(0,0){\line(1,0){1}}
\put(0,0){\line(0,1){1}}
\put(1,1){\line(-1,0){1}}
\put(1,1){\line(0,-1){1}}
\end{picture}$
\begin{displaymath}
\end{displaymath} 

\phantom{p}

Combining Lemmas~\ref{sub3} and~\ref{sub2},
we obtain our main result:

\begin{Theorem} \label{sub4}
Let $H$ be a (fixed) planar graph and let $\liminf_{n \to \infty} \frac{m}{n} > 1$.
Then $\exists \alpha > 0$ and $\exists N$ such that
\begin{eqnarray*}
& & \mathbf{P}\Big[\textrm{$P_{n,m}$
will \emph{not} have a set of $> \alpha n$ vertex-disjoint} \\
& & \phantom{wwwwwwwwi}\textrm{induced order-preserving copies of $H$}\Big] 
< e^{- \alpha n} 
\textrm{ } \forall n \geq N. \\
\end{eqnarray*} 
\end{Theorem}
\phantom{p}

It follows from Theorem~\ref{sub4} that
the probability that $P_{n,m}$ will contain an induced order-preserving copy of any given planar graph $H$
converges to $1$ 
if~$\liminf_{n \to \infty}\frac{m}{n}>1$.
We shall soon see that this can, in fact, be extended to 
$\liminf_{n \to \infty} \frac{m}{n}>0$ if $H$ is acyclic.
First, we need to note the following result, which deals with small values of $\frac{m}{n}$:

\begin{Lemma} \label{sub5}
Let $H$ be a (fixed) acyclic graph,
let $s$ be a fixed constant, 
and let $m=m(n)$ satisfy
$0<\liminf_{n \to \infty}\frac{m}{n} \leq \limsup_{n \to \infty}\frac{m}{n} \leq (1+o(1))n$ 
as $n \to \infty$.
Then
\begin{eqnarray*}
& & \mathbf{P}\Big[\textrm{$P_{n,m}$ 
will have a set of $\geq s$ vertex-disjoint} \\
& & \phantom{wwwww}\textrm{induced order-preserving copies of $H$}\Big] 
\to 1
\textrm{ as } n \to \infty. 
\end{eqnarray*}
\end{Lemma}
\textbf{Proof}
Let the components of $H$ be $H_{1}, H_{2}, \ldots, H_{k}$
and let us choose vertices
$\{ v_{1}, v_{2}, , \ldots, v_{k} \} \subset V(H)$
such that $v_{i} \in V(H_{i})$ $\forall i$.

Without loss of generality, 
we may assume that
$V(H)= \{1,2, \ldots, |H| \}$.
Let us define
$H^{\prime}$ be the graph with vertex set
$V(H^{\prime}) = \{1,2, \ldots, |H|, |H|+1 \}$
and edge set
$E(H^{\prime}) = E(H) \cup \{ (v_{1},|H|+1), (v_{2},|H|+~1), (v_{3},|H|+~1) \}$.
\begin{figure} [ht]
\setlength{\unitlength}{1cm}
\begin{picture}(20,2.5)(0,0.9)

\put(5.025,0.75){\oval(0.75,0.75)[r]}
\put(5.025,2){\oval(0.75,0.75)[r]}
\put(5.025,3.25){\oval(0.75,0.75)[r]}

\put(5.4,0.75){\line(4,5){1}}
\put(5.4,2){\line(1,0){1}}
\put(5.4,3.25){\line(4,-5){1}}

\put(5.4,0.75){\circle*{0.1}}
\put(5.4,2){\circle*{0.1}}
\put(5.4,3.25){\circle*{0.1}}
\put(6.4,2){\circle*{0.1}}

\put(4.85,3.15){$H_{1}$}
\put(4.85,1.9){$H_{2}$}
\put(4.85,0.65){$H_{k}$}
\put(6.5,2.1){\scriptsize{$|H|+1$}}
\put(5.5,3.35){\footnotesize{$v_{1}$}}
\put(5.5,2.1){\footnotesize{$v_{2}$}}
\put(5.5,0.65){\footnotesize{$v_{k}$}}
\end{picture}
\caption{The graph $H^{\prime}$.}
\end{figure} 
\\Then $H^{\prime}$ is a tree containing an induced order-preserving copy of $H$.
Thus, it suffices to prove the result for $H^{\prime}$.
But this follows from Corollary~\ref{tree6}.
\phantom{qwerty}
\begin{picture}(1,1)
\put(0,0){\line(1,0){1}}
\put(0,0){\line(0,1){1}}
\put(1,1){\line(-1,0){1}}
\put(1,1){\line(0,-1){1}}
\end{picture}
\\
\\

Combining Theorem~\ref{sub4} and Lemma~\ref{sub5},
we obtain our aforementioned result for acyclic graphs:

\begin{Theorem} \label{sub6}
Let $H$ be a (fixed) acyclic graph,
let $s$ be a fixed constant,
and let~$\liminf_{n \to \infty} \frac{m}{n} > 0$.
Then 
\begin{eqnarray*}
& & \mathbf{P}\Big[\textrm{$P_{n,m}$
will have a set of $\geq s$ vertex-disjoint} \\
& & \phantom{wwwww}\textrm{induced order-preserving copies of $H$}\Big] 
\to 1
\textrm{ as } n \to \infty. 
\end{eqnarray*}
\end{Theorem}
\textbf{Proof}
The proof is analogous to that of Theorem~\ref{cyc53}.

We shall suppose,
hoping for a contradiction,
that there exists a function~$m(n)$ such that 
$
\mathbf{P}[P_{n,m} 
\textrm{ will have a set of $\geq s$ vertex-disjoint induced order-}$
$\textrm{preserving copies of }H] \not \to 1.
$
Thus, 
$\exists \delta > 0$
and there exists an infinite subsequence $n_{1},n_{2}, \ldots $
such that
$
\mathbf{P}[P_{n_{i},m(n_{i})} 
\textrm{ will have a set of $\geq s$ vertex-disjoint}$
$\textrm{induced order-preserving copies of }H] 
\leq 1-\delta \textrm{ } \forall i.
$

Clearly,
by Lemma~\ref{sub5},
we can't have $m(n_{i}) \leq (1+o(1))n_{i}$ as $i \to \infty$.
Thus, $\exists \epsilon > 0$ such that 
$m(n_{i}) \geq (1+\epsilon)n_{i}$
for arbitrary many values of $i$.
But we then obtain a contradiction to Theorem~\ref{sub4}.
$\phantom{qwerty}$
\setlength{\unitlength}{.25cm}
\begin{picture}(1,1)
\put(0,0){\line(1,0){1}}
\put(0,0){\line(0,1){1}}
\put(1,1){\line(-1,0){1}}
\put(1,1){\line(0,-1){1}}
\end{picture}

\newpage
\section{Unicyclic Subgraphs} \label{cycsub}

We continue to look at the probability that $P_{n,m}$ will contain given subgraphs.
In Section~\ref{sub},
we saw that for any given planar subgraph the limiting probability is $1$ if~$\liminf_{n \to \infty} \frac{m}{n}>1$,
and that for acyclic subgraphs the limiting probability is $1$ if $\liminf_{n \to \infty} \frac{m}{n} >0$.
We shall now start to investigate the limiting probability for non-acyclic subgraphs when $\frac{m}{n} \leq 1+o(1)$.

In this section,
our main focus will be on unicyclic subgraphs.
We already know 
(by combining Theorems~\ref{gen4} and~\ref{sub4})
that the associated probability is bounded away from $0$ for all large $n$
if $\liminf_{n \to \infty} \frac{m}{n} >0$,
so it only remains to discover exactly when the probability converges to $1$.

Recall that in Section~\ref{cyccpt}
we showed that the probability that $P_{n,m}$ will contain a non-acyclic \textit{component} of order $\leq k$ 
is always bounded away from~$1$ for large~$n$, 
regardless of $m$.
We will now see (in Theorem~\ref{cyc54}) that,
if $\limsup_{n \to \infty} \frac{m}{n} <1$,
then the probability that $P_{n,m}$ will contain a non-acyclic \textit{subgraph} of order~$\leq k$
is also bounded away from $1$. 
We shall then complete matters by seeing 
(in Lemma~\ref{unisub1}) that,
by contrast to our results with components,
if $\frac{m}{n} \to 1$ then the probability that $P_{n,m}$ 
will contain any given connected unicyclic subgraph converges to $1$. 
Both proofs will make use of bounds on the number of isolated vertices,
which are implicit from Sections~\ref{cptlow} and~\ref{conkappa}.

In Theorem~\ref{unisub3},
we shall extend Lemma~\ref{unisub1} to cover any (not necessarily connected) given subgraph with no multicyclic components.

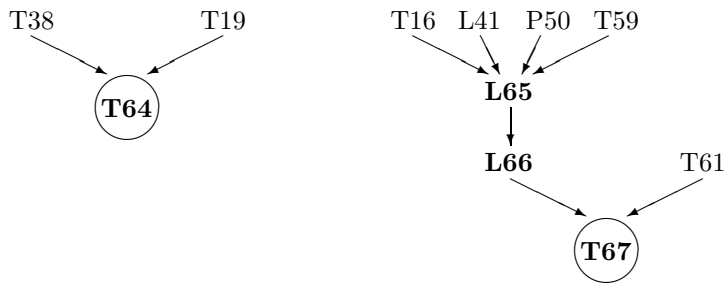
\begin{figure} [ht]
\setlength{\unitlength}{1cm}
\begin{picture}(20,3.5)(-0.1,2.9)

\put(1.225,5.85){T\ref{tree4}}
\put(3.7755,5.85){T\ref{add53}}
\put(9,5.85){T\ref{sub7}}
\put(7.2,5.85){L\ref{cpt12}}
\put(8.1,5.85){P\ref{cpt31}}
\put(2.45,4.675){\textbf{T\ref{cyc54}}}
\put(7.55,4.9){\textbf{L\ref{unisub1}}}
\put(7.55,3.95){\textbf{L\ref{unisub21}}}
\put(10.15,3.95){T\ref{sub4}}
\put(8.825,2.775){\textbf{T\ref{unisub3}}}
\put(6.3,5.85){T\ref{pen5}}

\put(2.8,4.8){\circle{0.8}}
\put(9.175,2.9){\circle{0.8}}

\put(1.525,5.75){\vector(2,-1){1}}
\put(4.075,5.75){\vector(-2,-1){1}}
\put(6.6,5.75){\vector(2,-1){1}}
\put(7.5,5.75){\vector(1,-2){0.25}}
\put(8.3,5.75){\vector(-1,-2){0.25}}
\put(9.2,5.75){\vector(-2,-1){1}}
\put(7.9,4.8){\vector(0,-1){0.5}}
\put(7.9,3.85){\vector(2,-1){1}}
\put(10.45,3.85){\vector(-2,-1){1}}
\end{picture}

\caption{The structure of Section~\ref{cycsub}.}
\end{figure}

\begin{displaymath}
\end{displaymath}

We shall go straight in with our result for when $\limsup_{n \to \infty} \frac{m}{n} < 1$:

\begin{Theorem} \label{cyc54}
Let $k$ be a fixed constant and let $m=m(n) \leq An$, where $A<1$.
Then
\begin{displaymath}
\limsup_{n \to \infty}
\mathbf{P}[P_{n,m} \textrm{ will have girth }\leq k] < 1.
\end{displaymath}
\end{Theorem}
\textbf{Sketch of Proof}
Let $C_{i}$ be used to denote cycles of order $i$ 
and let `$C_{i}$ component' be used to denote components that are cycles of order $i$
(thus, a $C_{i}$ component is a type of $C_{i}$). 

Step $1$: 
By an induction hypothesis, 
we may assume that there is a decent proportion of graphs with no cycles of order $<l$.

Step $2$:
If we have a $C_{l}$ \textit{component},
we may remove an edge from it and then insert an edge
between the remainder of the cycle and the rest of the graph.
We find that we may create so many graphs by doing this
that the proportion of graphs \textit{without} a $C_{l}$ component (and still with no cycles of order $<l$)
must be quite decent.

Step $3$: 
Using Theorem~\ref{tree4},
we find that there must be a decent proportion of graphs 
with no $C_{l}$ components, no cycles of order $<l$, and lots of isolated vertices.

Step $4$:
If we have a $C_{l}$ and lots of isolated vertices,
we may `transfer' the edges of the $C_{l}$ to the isolated vertices to construct a $C_{l}$ component.
If there were no $C_{l}$ components originally,
then the amount of double-counting is small,
and so we can use this idea to show that the proportion of graphs 
with lots of~$C_{l}$'s, lots of isolated vertices and no $C_{l}$ components is small.
Thus, using the result of Step $3$,
there must be a decent proportion of graphs with few $C_{l}$'s and no cycles of order $<l$.

Step $5$:
We then destroy the few $C_{l}$'s by deleting the edges from them and inserting edges between components of the graph. \\
\\
\textbf{Full Proof}
As in the sketch of the proof,
we shall let $C_{i}$ be used to denote cycles of order $i$ 
and $C_{i}$ component be used to denote components that are cycles of order $i$
(thus, a $C_{i}$ component is a type of $C_{i}$). 
For $\mathcal{A}_{n} \subset \mathcal{P}(n,m)$,
we shall use $\mathcal{A}_{n}^{c}$ to denote the set of graphs in $\mathcal{P}(n,m)$ that are not in $\mathcal{A}_{n}$. \\ 
\\
\underline{Step $1$}

We will prove the result by induction on $k$.
It is trivial for $k \leq 2$.
Let us assume that the result is true for $k=l-1$.

Let $\mathcal{L}_{n}$ denote the set of graphs in $\mathcal{P}(n,m)$ without a cycle of order $<l$.
Then, by our induction hypothesis,
$\exists \epsilon (l) > 0$ such that
$|\mathcal{L}_{n}| \geq~\epsilon |\mathcal{P}(n,m)|$ for all sufficiently large $n$. \\
\\
\underline{Step $2$}

Let $\mathcal{G}_{n}$ denote the set of graphs in $\mathcal{L}_{n}$ with a $C_{l}$ component.
For each graph~$G \in \mathcal{G}_{n}$,
let us delete an edge from a $C_{l}$ component, $C^{\prime} (=C^{\prime}_{G})$,
to leave a spanning path. 
We have $l$ choices for this edge.
Let us then insert an edge between a vertex in the spanning path and a vertex elsewhere.
\begin{figure} [ht]
\setlength{\unitlength}{1cm}
\begin{picture}(20,2.5)

\put(3.5,1.75){\oval(2,1)}
\put(3.5,0.375){\circle{0.75}}

\put(5.5,1.75){\vector(1,0){1}}

\put(8.5,0.75){\circle*{0.1}}

\put(8.5,1.75){\oval(2,1)}
\put(8.5,0.375){\oval(0.75,0.75)[t]}
\put(8.5,0.375){\oval(0.75,0.75)[br]}

\put(8.5,0.75){\line(0,1){0.5}}
\put(8.5,1.25){\circle*{0.1}}

\end{picture}

\caption{Destroying a $C_{l}$ component in Step $2$.}
\end{figure}
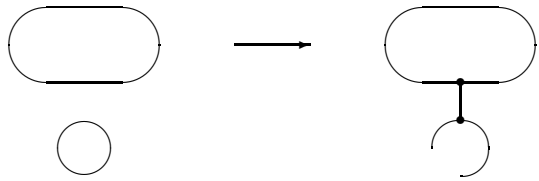
\\ We have $l(n-l)$ ways to do this, planarity is maintained and no cycles of order~$<l$ will be created.
Thus, we have constructed $|\mathcal{G}_{n}|l^{2}(n-l)$ graphs in $\mathcal{L}_{n}$.

Given one of our constructed graphs,
there are at most $n-1$ possibilities for which is the inserted edge,
since it must be a cut-edge.
There are then at most $2$ possibilities for $V(C^{\prime})$ (at most one possibility from each endpoint of the edge).
There is then only $1$ possibility for where the inserted edge was originally.
Thus, we have built each graph at most $2(n-1)$ times.

Therefore, $|\mathcal{G}_{n}| \leq \frac{2(n-1)}{l^{2}(n-l)} |\mathcal{L}_{n}|$.
Thus, for any fixed constant $\delta \in \left( 0, 1-\frac{2}{l^{2}} \right)$,
we have $|\mathcal{G}_{n}| < (1- \delta) |\mathcal{L}_{n}|$ for all sufficiently large $n$,
and so we have
$|\mathcal{L}_{n} \cap \mathcal{G}_{n}^{c}| > \delta |\mathcal{L}_{n}| \geq \delta \epsilon |\mathcal{P}(n,m)|$
for all sufficiently large $n$. \\
\\
\underline{Step $3$}

Let $c \in \left( 0, \frac{1}{2} \right)$ be a fixed constant and suppose $m \in [cn,An]$ $\forall n$.
Then, by Theorem~\ref{tree4} (with $|H|=1$) and Theorem~\ref{add53},
$\exists \beta > 0$ such that
\begin{displaymath}
\mathbf{P}[P_{n,m} \textrm{ will have $\geq \beta n$ isolated vertices}] \to 1 \textrm{ as } n \to \infty.
\end{displaymath}
Clearly, $P_{n,m}$ must also have at least $n-2cn$ isolated vertices if $m \leq cn$,
so by setting $\alpha = \min \{ \beta, 1-2c \}$
we find that whenever $m \leq An$ we have
\begin{displaymath}
\mathbf{P}[P_{n,m} \textrm{ will have $\geq \alpha n$ isolated vertices}] \to 1 \textrm{ as } n \to \infty.
\end{displaymath}

Let $\mathcal{I}_{n}$ denote the set of graphs in $\mathcal{P}(n,m)$ with $\geq \alpha n$ isolated vertices.
Then, using the result of Step $2$, 
we find that $\exists \epsilon^{\prime} \!>\! 0$ such that
$|\mathcal{G}_{n}^{c} \cap \mathcal{L}_{n} \cap~\!\mathcal{I}_{n}| \!\geq~\!\!\epsilon^{\prime} |\mathcal{P}(n,m)|$
for all sufficiently large $n$. \\
\\
\underline{Step $4$}

Let $r=r(l)> \frac{2l!}{\epsilon^{\prime}\alpha^{l}}$ 
and let $\mathcal{J}_{n}$ denote the set of graphs in $\mathcal{P}(n,m)$ 
whose maximal number of edge-disjoint $C_{l}$'s is $>r$.

Consider the set $\mathcal{J}_{n} \cap \mathcal{G}_{n}^{c} \cap \mathcal{I}_{n}$,
i.e. the set of graphs in $\mathcal{P}(n,m)$ with $\geq \alpha n$ isolated vertices,
$>r$ edge-disjoint $C_{l}$'s,
but no $C_{l}$ components.
For each graph $J \in \mathcal{J}_{n} \cap \mathcal{G}_{n}^{c} \cap \mathcal{I}_{n}$,
delete all $l$ edges from a $C_{l}$ 
(we have $>r$ choices for this)
and insert them between $l$ isolated vertices
(we have $\geq \left( ^{\alpha n}_{\phantom{q}l} \right)$ choices for these)
to form a $C_{l}$ component.
We can thus construct at least 
$|\mathcal{J}_{n} \cap \mathcal{G}_{n}^{c} \cap \mathcal{I}_{n}| r \left( ^{\alpha n}_{\phantom{q}l} \right)$
graphs in $\mathcal{P}(n,m)$.

Given one of our constructed graphs,
there are at most $l+1$ $C_{l}$ components,
since we have deliberately constructed one and we may have also created at most~$l$ when we deleted the edges
(since each vertex in the deleted $C_{l}$ might now be in a $C_{l}$ component).
Thus, there are at most $l+1$ possibilities for which is the deliberately constructed component
and hence at most $l+1$ possibilities for which are the inserted edges.
There are then at most 
$\left( ^{\phantom{i}n}_{\phantom{q}l} \right) \frac{l!}{|\textrm{\scriptsize{Aut}}(C_{l})|} = 
\left( ^{\phantom{i}n}_{\phantom{q}l} \right) \frac{(l-1)!}{2}$
possibilities for where these edges were originally.
Thus, we have built each graph at most
$(l+1) \left( ^{\phantom{i}n}_{\phantom{q}l} \right) \frac{(l-1)!}{2} < \frac{(l+1)n^{l}}{2l} < n^{l}$ times
and so 
\begin{eqnarray*}
|\mathcal{J}_{n} \cap \mathcal{G}_{n}^{c} \cap \mathcal{I}_{n}|
& < & \frac{n^{l} |\mathcal{P}(n,m)|}{ r \left( ^{\alpha n}_{\phantom{q}l} \right)} \\
& = & (1+o(1)) \frac{l!}{r \alpha ^{l}} |\mathcal{P}(n,m)| \\
& < & \frac{\epsilon^{\prime}}{2} |\mathcal{P}(n,m)| 
\textrm{ for sufficiently large $n$, since } r > \frac{2l!}{\epsilon^{\prime} \alpha^{l}}. 
\end{eqnarray*}

Since we also know that
$|\mathcal{G}_{n}^{c} \cap \mathcal{L}_{n} \cap \mathcal{I}_{n}| \geq \epsilon^{\prime} |\mathcal{P}(n,m)|$
for all sufficiently large~$n$,
from Step $3$,
we must have
$|\mathcal{J}_{n}^{c} \cap \mathcal{G}_{n}^{c} \cap \mathcal{L}_{n} \cap \mathcal{I}_{n}| > 
\frac{\epsilon^{\prime}}{2} |\mathcal{P}(n,m)|$.
Thus, $|\mathcal{J}_{n}^{c} \cap \mathcal{L}_{n}| > \frac{\epsilon^{\prime}}{2} |\mathcal{P}(n,m)|$
for all sufficiently large $n$. \\
\\
\underline{Step $5$}

Recall that $\mathcal{J}_{n}^{c} \cap \mathcal{L}_{n}$
is the set of graphs in $\mathcal{P}(n,m)$ without a copy of a cycle of order $<l$
and with $\leq r$ edge-disjoint $C_{l}$'s.
For $L \in \mathcal{J}_{n}^{c} \cap \mathcal{L}_{n}$,
let $S(L)$ denote a maximal set of edge-disjoint $C_{l}$'s (so $|S(L)| \leq r$)
and, for $s \leq r$, let $\mathcal{J}_{n,s}$ denote the set of graphs in $\mathcal{P}(n,m)$ with $|S(L)|=s$.

For each graph $L \in \mathcal{J}_{n,s} \cap \mathcal{L}_{n}$,
delete all $l$ edges from a $C_{l}$ that is in $S(L)$.
Note that the graphs will now have $|S|=s-1$,
by maximality of $S(L)$.

Clearly, we may insert an edge between any two vertices in different components
without introducing a copy of a cycle.
By the proof of Lemma~\ref{add33},
we have at least $(1+o(1)) \frac{(1-A)(1+A)}{2} n^{2}$ choices 
for where to insert an edge between two vertices in different components.
Doing this $l$ times,
we find that we may construct 
$\geq (1+o(1)) \frac{ \left( \frac{(1-A)(1+A)}{2} \right)^{l} n^{2l} }{l!} |\mathcal{J}_{n,s} \cap \mathcal{L}_{n}|$
graphs in $\mathcal{J}_{n,s-1} \cap \mathcal{L}_{n}$.

Given one of our created graphs,
there are $\leq \left( ^{m}_{\phantom{i}l} \right) \leq (An)^{l}$ possibilities for which edges were inserted 
and 
$\leq \left( ^{\phantom{i}n}_{\phantom{q}l} \right) \frac{l!}{|\textrm{\scriptsize{Aut}}(C_{l})|} < n^{l}$ 
possibilities for where they were originally.
Thus, we have built each graph at most $A^{l}n^{2l}$ times.

Thus, 
$|\mathcal{J}_{n,s-1} \cap \mathcal{L}_{n}| \geq 
(1+o(1)) \frac{ \left( \frac{(1-A)(1+A)}{2A} \right)^{l} }{l!} |\mathcal{J}_{n,s} \cap \mathcal{L}_{n}|$.

For $s \leq r$,
let $\mathcal{J}_{n, \leq s}$ denote the set of graphs in $\mathcal{P}(n,m)$ with $|S| \leq s$.
Then, with $z = \frac{ \left( \frac{(1-A)(1+A)}{2A} \right)^{l} }{l!}$, we have
$|\mathcal{J}_{n, \leq s-1} \cap \mathcal{L}_{n}| \geq 
(1+o(1)) \frac{z}{1+z} |\mathcal{J}_{n, \leq s} \cap~\!\mathcal{L}_{n}|$~$\forall s \!\leq~\!\!r$.
Thus, we may obtain 
\begin{eqnarray*}
|\mathcal{J}_{n, \leq 0} \cap \mathcal{L}_{n}|
& \geq & (1+o(1))^{r} \left( \frac{z}{1+z} \right)^{r} |\mathcal{J}_{n, \leq r} \cap \mathcal{L}_{n}| \\
& = & (1+o(1)) \left( \frac{z}{1+z} \right)^{r} |\mathcal{J}_{n}^{c} \cap \mathcal{L}_{n}| \\
& \geq & (1+o(1)) \left( \frac{z}{1+z} \right)^{r} \frac{\epsilon^{\prime}}{2} |\mathcal{P}(n,m)|
\phantom{ww} \textrm{ from Step $4$}.
\end{eqnarray*}

But $\mathcal{J}_{n, \leq 0} \cap \mathcal{L}_{n}$ is the set of graphs in $\mathcal{P}(n,m)$
without a cycle of order $\leq l$.
Thus, the induction hypothesis is true for $k=l$, and so we are done.~\phantom{qwerty}\begin{picture}(1,1)
\put(0,0){\line(1,0){1}}
\put(0,0){\line(0,1){1}}
\put(1,1){\line(-1,0){1}}
\put(1,1){\line(0,-1){1}}
\end{picture} 
\begin{displaymath}
\end{displaymath} 
\phantom{p}

We will now start to work towards Theorem~\ref{unisub3},
where we shall show that the probability that $P_{n,m}$ 
will contain any given subgraph with no multicyclic components converges to $1$
if $\frac{m}{n} \!\geq\! 1-o(1)$. 
We already know this for when~$\!\frac{m}{n} \!\geq~\!\!B\!>~\!\!1$,
by Theorem~\ref{sub4},
so it will suffice to now look at the case when $\frac{m}{n} \in [(1-o(1)), B]$.
We shall start by proving the result for a connected unicyclic graph,
showing that a.a.s.~$P_{n,m}$ will contain many vertex-disjoint appearances of any such graph: \\

\begin{Lemma} \label{unisub1}
Let $H$ be a (fixed) connected unicyclic graph,
let $t$ be a fixed constant,
and let $m=m(n) \in [(1-o(1))n, Bn]$, where $B<3$.
Then
\begin{displaymath}
\mathbf{P}[f^{\prime}_{H}(P_{n,m}) \geq t] \to 1 
\textrm{ as } n \to \infty.
\end{displaymath} 
\end{Lemma}
\textbf{Sketch of Proof}
By Theorem~\ref{pen5},
we may assume that we have lots of pendant edges.
We may delete $|H|$ of these edges,
and then use them with $|H|-1$ of the associated (now isolated) vertices
to convert another pendant edge into an appearance of $H$.
Note that we are also left with an extra isolated vertex.

By Lemma~\ref{cpt12} and Proposition~\ref{cpt31},
we may assume that there are not very many isolated vertices.
Thus, if our original graphs had few appearances of $H$,
then the amount of double-counting will be small,
and hence the size of our original set of graphs must have been small. \\
\\
\textbf{Full Proof}
Let $\mathcal{G}_{n}$ denote the set of graphs in $\mathcal{P}(n,m)$ 
with $f^{\prime}_{H} < t$,
and let $X$ denote the event that $f_{H}^{\prime}(P_{n,m}) < t$.

Recall, from Theorem~\ref{pen5},
$\exists \alpha >0$ such that
\begin{eqnarray}
\mathbf{P}[P_{n,m} 
\textrm{ will have $\geq \alpha n$ pendant edges}] \to 1 
\textrm{ as } n \to \infty. \label{eq:u*}
\end{eqnarray} 
Let $\mathcal{H}_{n}$ denote the set of graphs in $\mathcal{P}(n,m)$
with $\geq \alpha n$ pendant edges,
and let~$Y$ denote the event that $P_{n,m}$ will have $\geq \alpha n$ pendant edges.

Also, note that by Lemma~\ref{cpt12} and Proposition~\ref{cpt31},
$\exists c>0$ such that
\begin{eqnarray}
\mathbf{P} \left[ \kappa \left( P_{n,m} \right)
> 2 \max \left \{ \frac{cn}{\ln n}, n-m \right \} \right] = e^{- \Omega(n)}. \label{eq:udag}
\end{eqnarray} 
Let $x=x(n)= 2 \max \{ \frac{cn}{\ln n}, n-m \} = o(n)$,
let $\mathcal{I}_{n}$ denote the number of graphs in $\mathcal{P}(n,m)$ with $\leq x$ components,
and let $Z$ denote the event that $\kappa \left( P_{n,m} \right) \leq x$.

Note that 
$\mathbf{P}[X] \leq \mathbf{P}[X \cap Y \cap Z] + \mathbf{P} \left[\bar{Y}\right] + \mathbf{P}\left[\bar{Z}\right] 
\to \mathbf{P}[X \cap Y \cap Z] \textrm{ as } n \to \infty$, 
by (\ref{eq:u*}) and (\ref{eq:udag}).
Thus, it suffices to show that $\mathbf{P}[X \cap Y \cap Z] \to 0$ as $n \to \infty$. \\

Given a graph in $\mathcal{G}_{n} \cap \mathcal{H}_{n} \cap \mathcal{I}_{n}$,
let us choose $|H|+1$ pendant edges
(we have at least $\left( ^{\phantom{w} \alpha n}_{|H|+1} \right)$ choices for these).
We shall use these edges to create an appearance of $H$.

Out of our chosen $|H|+1$ pendant edges,
let the edge incident with the lowest labelled vertex of degree $1$ be called `special',
and let this associated lowest labelled vertex of degree $1$ be called the `root'.
Let us delete our $|H|$ non-special pendant edges to create at least $|H|$ isolated vertices.

We may choose $|H|-1$ of these newly isolated vertices in such a way
that no two were adjacent in the original graph
(i.e. we don't choose two vertices from the same pendant edge,
even if that is possible).
We may then use these chosen isolated vertices,
together with the root,
to construct an appearance of~$H$
(by inserting $|H|$ edges appropriately).

\begin{figure} [ht]
\setlength{\unitlength}{0.9cm}
\begin{picture}(20,3.25)(-0.75,0)

\put(-0.05,0.75){\circle*{0.1}}
\put(1.15,0.75){\circle*{0.1}}
\put(2.35,0.75){\circle*{0.1}}
\put(3.55,0.75){\circle*{0.1}}

\put(1.75,2.25){\oval(5,1.5)}

\put(-0.05,0.75){\line(0,1){0.75}}
\put(-0.05,1.5){\circle*{0.1}}
\put(1.15,0.75){\line(0,1){0.75}}
\put(1.15,1.5){\circle*{0.1}}
\put(2.35,0.75){\line(0,1){0.75}}
\put(2.35,1.5){\circle*{0.1}}
\put(3.55,0.75){\line(0,1){0.75}}
\put(3.55,1.5){\circle*{0.1}}

\put(5.5,1.5){\vector(1,0){1}}

\put(10.25,2.25){\oval(5,1.5)}

\put(8.45,1.5){\circle*{0.1}}
\put(9.65,1.5){\circle*{0.1}}
\put(10.85,1.5){\circle*{0.1}}
\put(12.05,1.5){\circle*{0.1}}

\put(-0.25,0.4){$v_{1}$}
\put(0.95,0.4){$v_{2}$}
\put(2.15,0.4){$v_{3}$}
\put(3.35,0.4){$v_{4}$}

\put(12.05,0.75){\circle*{0.1}}
\put(8.45,0.75){\circle*{0.1}}
\put(8.45,0.75){\line(0,1){0.75}}

\put(8.45,0.75){\line(3,-5){0.45}}
\put(8.45,0.75){\line(-3,-5){0.45}}
\put(8,0){\circle*{0.1}}
\put(8.9,0){\circle*{0.1}}
\put(7.8,-0.35){$v_{2}$}
\put(8.8,-0.35){$v_{3}$}
\put(8.55,0.75){$v_{1}$}
\put(11.85,0.4){$v_{4}$}
\put(8,0){\line(1,0){0.9}}

\end{picture}

\caption{Constructing an appearance of $H$.}
\end{figure}
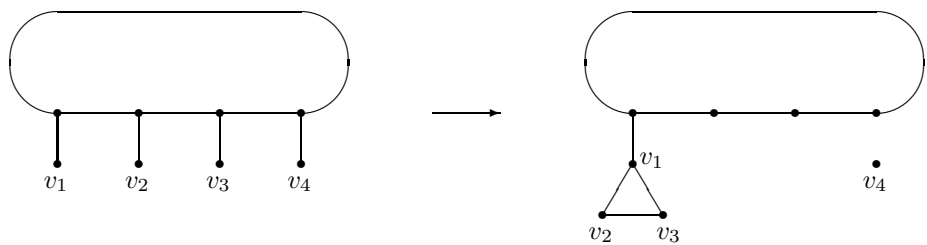

We shall now consider the amount of double-counting:

Suppose that when we deleted one of our pendant edges 
we `accidentally' created an appearance of $H$ at $W$.
Then the pendant edge must have been incident to a vertex in $W$.
Thus, if we simultaneously accidentally created appearances of $H$ at both $W_{1}$ and $W_{2}$,
then the pendant edge must have been incident to a vertex in $W_{1} \cap W_{2}$,
and so $W_{1} \cap W_{2} \neq \emptyset$.

As noted in Lemma~\ref{gen12},
an appearance of $H$ meets at most $|H|-1$ other appearances of $H$
(since there are at most $|H|-1$ cut-edges in $H$
and each of these can have at most one `orientation' that provides an appearance of $H$).
Thus 
(using the observation of the previous paragraph),
we can have created at most $|H|$ appearances of $H$ each time we deleted a pendant edge.
By the same argument,
our original graph must have satisfied $f_{H} < t|H|$
and when we deliberately created our appearance of $H$
we can have only increased the number of appearances of $H$ by at most $|H|$.
Therefore, 
in total we find that each of our constructed graphs will contain $<|H|(|H|+t+1)$ appearances of~$H$.

Thus, given one of our constructed graphs,
there are at most $|H|(|H|+~\!t+~\!1)$
possibilities for which is the constructed appearance of $H$.
We may then recover the original graph by deleting the $|H|$ edges from this appearance,
joining the $|H|-1$ non-root vertices back to the rest of the graph
(at most $n^{|H|-1}$ possibilities),
and joining the correct isolated vertex back to the rest of the graph
(at most $in$ possibilities,
where $i$ denotes the number of isolated vertices in the constructed graph).

We know that the number of components in the original graph was at most~$x$,
and each time we deleted an edge we can have only increased the number of components by at most $1$.
Thus, the number of components in our constructed graph is at most $x+|H|$,
and so $i \leq x+|H|$.

Thus, we have built each graph at most
$|H|(|H|+t+1)n^{|H|-1}(x+|H|)n$ times. \\

Therefore,
\begin{eqnarray*}
\mathbf{P}[X \cap Y \cap Z]
& = & \frac{|\mathcal{G}_{n} \cap \mathcal{H}_{n} \cap \mathcal{I}_{n}|}{|\mathcal{P}(n,m)|} \\
& \leq & \frac{|H|(|H|+t+1)n^{|H|}(x+|H|)}
{\left( ^{\phantom{w} \alpha n}_{|H|+1} \right) } \\
& = & \frac{o \left( n^{|H|+1} \right) }{\Theta \left( n^{|H|+1} \right) }, \textrm{ since } x=o(n)\\
& \to & 0 \textrm{ as } n \to \infty.
\phantom{qwerty}\begin{picture}(1,1)
\put(0,0){\line(1,0){1}}
\put(0,0){\line(0,1){1}}
\put(1,1){\line(-1,0){1}}
\put(1,1){\line(0,-1){1}}
\end{picture} 
\end{eqnarray*}
\\
\\

For $m \in [(1-o(1))n,Bn]$,
it follows from Lemma~\ref{unisub1} that a.a.s.~$P_{n,m}$ 
will have at least $t$ vertex-disjoint induced order-preserving copies of any connected unicyclic $H$.
We shall now extend this result to cover \textit{any} (not necessarily connected) subgraph without multicyclic components:
\\

\begin{Lemma} \label{unisub21}
Let $H$ be a (fixed) graph with $e(H_{i}) \leq |H_{i}|$ for all components $H_{i}$ of~$\!H$,
let $s$ be a fixed constant,
and let $m=m(n) \in [(1-o(1))n,Bn]$, where~$B<3$.
Then
\begin{eqnarray*}
& & \mathbf{P}\Big[\textrm{$P_{n,m}$ 
will have a set of $\geq s$ vertex-disjoint} \\
& & \phantom{wwwww}\textrm{induced order-preserving copies of $H$}\Big] 
\to 1 
\textrm{ as } n \to \infty.
\end{eqnarray*} 
\end{Lemma}
\textbf{Proof}
Since the probability of a copy of $H$ being order-preserving is at least~$\frac{1}{|H|!}$
and is independent of whether or not another vertex-disjoint copy of $H$ is also order-preserving,
it suffices to show 
\begin{displaymath}
\mathbf{P}[P_{n,m} 
\textrm{ will have a set of $\geq l$ vertex-disjoint induced copies of }H] \!\to\! 1 
\textrm{ as } n \!\to\! \infty
\end{displaymath}
for an arbitrary fixed constant $l$.

For every tree component $T$ of $H$,
let us introduce two new vertices $u_{T}$ and~$v_{T}$,
and three new edges $u_{T}v_{T}$, $u_{T}w_{T}$ and $v_{T}w_{T}$,
where $w_{T}$ is an arbitrary vertex in $V(T)$.
Then the resulting component $T^{\prime}$ is unicyclic and contains an induced copy of $T$.
\begin{figure} [ht]
\setlength{\unitlength}{1cm}
\begin{picture}(20,3.3)(-2,0)

\put(0,2){\circle*{0.1}}
\put(1,1){\circle*{0.1}}
\put(2,2){\circle*{0.1}}
\put(1,0){\circle*{0.1}}
\put(5,2){\circle*{0.1}}
\put(7,2){\circle*{0.1}}
\put(6,1){\circle*{0.1}}
\put(6,0){\circle*{0.1}}
\put(8,2.5){\circle*{0.1}}
\put(8,1.5){\circle*{0.1}}

\put(0,2){\line(1,-1){1}}
\put(2,2){\line(-1,-1){1}}
\put(1,1){\line(0,-1){1}}
\put(5,2){\line(1,-1){1}}
\put(7,2){\line(-1,-1){1}}
\put(6,1){\line(0,-1){1}}
\put(7,2){\line(2,1){1}}
\put(7,2){\line(2,-1){1}}
\put(8,2.5){\line(0,-1){1}}

\put(0.9,2.8){\Large{$T$}}
\put(5.9,2.8){\Large{$T^{\prime}$}}
\put(1.8,2.3){$w_{T}$}
\put(8.1,2.6){$u_{T}$}
\put(8.1,1.2){$v_{T}$}
\put(6.8,2.3){$w_{T}$}

\end{picture}

\caption{The graphs $T$ and $T^{\prime}$.}
\end{figure}
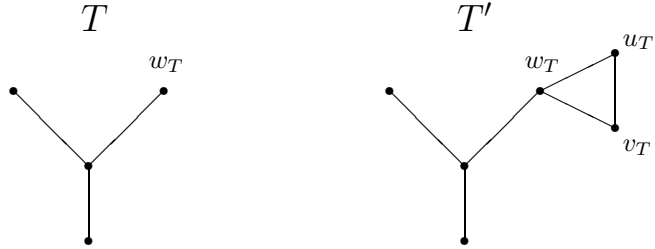
Thus, without loss of generality,
we may assume that $e(H_{i})=|H_{i}|$ for all components $H_{i}$ of $H$.

Let the components of $H$ be $H_{1},H_{2}, \ldots, H_{k}$,
for some $k$.
By Lemma~\ref{unisub1},
we have 
$\mathbf{P} \left[ f_{H_{i}}^{\prime}(P_{n,m}) \geq t~\forall i \right] \to 1$
as $n \to \infty$,
where $t$ may be chosen arbitrarily large.
Given a graph $G \in \mathcal{P}(n,m)$ with $f_{H_{i}}^{\prime}(G) \geq t$~$\forall i$,
let us select $l$ vertex-disjoint appearances of $H_{1}$.
Recall that $G$ contains a set of $t$ vertex-disjoint appearances of $H_{2}$
and notice that at most $l|H_{1}|$ of these 
can have a vertex in common with one of our selected appearances of $H_{1}$
(since the appearances of~$H_{2}$ are themselves vertex-disjoint)
and that at most $l$ can be attached by an edge to one of our selected appearances of $H_{1}$
(by definition of an appearance).
Hence,
if $t-l(|H_{1}|+1) \geq l$,
then we may select $l$ vertex-disjoint appearances of~$H_{2}$ in such a way
that the graph formed by these appearances and our selected appearances of $H_{1}$
will consist of $l$ vertex-disjoint induced copies of the graph with components $H_{1}$ and $H_{2}$.
Continuing in this manner,
we find that $G$ contains a set of $l$ vertex-disjoint induced copies of $H$,
and so we are done.
$\phantom{p}$
\begin{picture}(1,1)
\put(0,0){\line(1,0){1}}
\put(0,0){\line(0,1){1}}
\put(1,1){\line(-1,0){1}}
\put(1,1){\line(0,-1){1}}
\end{picture} 
\\
\\

Combining Lemma~\ref{unisub21} with Theorem~\ref{sub4},
we obtain our main result,
which holds for all $m \geq (1+o(1))n$:

\begin{Theorem} \label{unisub3}
Let $H$ be a (fixed) graph with $e(H_{i}) \leq |H_{i}|$ for all components $H_{i}$ of $H$,
let $s$ be a fixed constant,
and let $m=m(n) \in [(1-o(1))n, 3n-6]$.
Then 
\begin{eqnarray*}
& & \mathbf{P}\Big[\textrm{$P_{n,m}$ 
will have a set of $\geq s$ vertex-disjoint} \\
& & \phantom{wwwww}\textrm{induced order-preserving copies of $H$}\Big] 
\to 1
\textrm{ as } n \to \infty.
\end{eqnarray*} 
\end{Theorem}

\newpage
\section{Multicyclic Subgraphs} \label{multsub}

It now only remains to
look at the probability that $P_{n,m}$ will contain a given \textit{multicyclic} connected subgraph,
i.e.~a connected subgraph with more edges than vertices.
We already know from our work on appearances
(see Lemma~\ref{gen2}) that this probability converges to $1$ 
if $\liminf_{n \to \infty} \frac{m}{n} > 1$.
In this section,
we will show (in Corollary~\ref{msub3})
that the probability converges to $0$ if $\limsup_{n \to \infty} \frac{m}{n} < 1$,
and we shall also see some partial results
(in Theorems~\ref{msub2} and~\ref{msub1})
for the case~$\frac{m}{n} \to 1$.

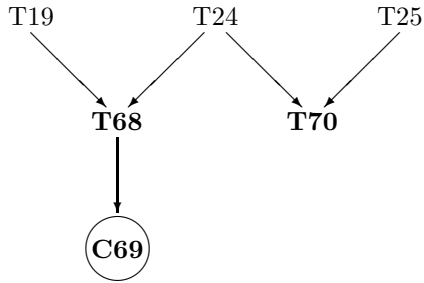
\begin{figure} [ht]
\setlength{\unitlength}{1cm}
\begin{picture}(20,3.6)(-2.55,0.925)

\put(0.7,4.2){T\ref{add53}}
\put(3.15,4.2){T\ref{add3}}
\put(5.6,4.2){T\ref{add54}}
\put(1.8,2.8){\textbf{T\ref{msub2}}}
\put(4.4,2.8){\textbf{T\ref{msub1}}}

\put(1,4.1){\vector(1,-1){1}}
\put(3.3,4.1){\vector(-1,-1){1}}
\put(3.6,4.1){\vector(1,-1){1}}
\put(5.9,4.1){\vector(-1,-1){1}}
\put(2.15,2.7){\vector(0,-1){1}}

\put(1.8,1.1){\textbf{C\ref{msub3}}}
\put(2.15,1.225){\circle{0.8}}
\end{picture}

\caption{The structure of Section~\ref{multsub}.}
\end{figure}

\phantom{p}

\begin{figure} [ht]
\setlength{\unitlength}{0.97cm}
\begin{picture}(20,2)(-0.25,0)

\put(1.8,0.4){\line(1,0){0.9}}
\put(1.8,0.4){\line(3,2){0.45}}
\put(1.8,0.4){\line(3,5){0.45}}
\put(2.25,0.7){\line(0,1){0.45}}
\put(2.7,0.4){\line(-3,2){0.45}}
\put(2.7,0.4){\line(-3,5){0.45}}
\put(1.8,0.4){\circle*{0.1}}
\put(2.7,0.4){\circle*{0.1}}
\put(2.25,0.7){\circle*{0.1}}
\put(2.25,1.15){\circle*{0.1}}

\put(2.25,0.75){\oval(5,1.5)}

\put(5.5,0.75){\vector(1,0){1}}

\put(8.25,0.4){\line(1,0){0.5}}
\put(8.25,0.75){\line(1,0){0.5}}
\put(8.25,1.1){\line(1,0){0.5}}
\put(10.75,0.4){\line(1,0){0.5}}
\put(10.75,0.75){\line(1,0){0.5}}
\put(10.75,1.1){\line(1,0){0.5}}

\put(9.3,0.4){\circle*{0.1}}
\put(10.2,0.4){\circle*{0.1}}
\put(9.75,0.7){\circle*{0.1}}
\put(9.75,1.15){\circle*{0.1}}

\put(9.75,0.75){\oval(5,1.5)}

\end{picture}

\caption{Redistributing the edges of $H$.} \label{msubfig}
\end{figure}
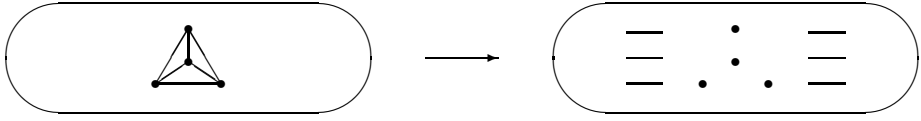

We start with a result for the case when either $\limsup_{n \to \infty} \frac{m}{n} < 1$
or $\frac{m}{n}$ converges to $1$ slowly from beneath:

\begin{Theorem} \label{msub2}
Let $H$ be a fixed multicyclic connected planar graph
and let $d=d(n)=n-m$ be such that $d>0$ and $d = \omega \left( n^{\frac{|H|}{e(H)}} \right)$.
Then 
\begin{displaymath}
\mathbf{P}[P_{n,m} \textrm{ will have a copy of }H] \to 0 \textrm{ as } n \to \infty.
\end{displaymath}
\end{Theorem}
\textbf{Proof}
Let $\mathcal{G}_{n}$ denote the set of graphs in $\mathcal{P}(n,m)$ with a copy of $H$.
For each graph $G \in \mathcal{G}_{n}$, let us delete all $e(H)$ edges from a copy of $H$ 
and then insert these edges back into the graph (see Figure~\ref{msubfig}).
By Theorems~\ref{add53} and~\ref{add3},
we have $\Omega \left( (dn)^{e(H)} \right)$ ways to do this, maintaining planarity.

There are then 
$\left(\left(^{\phantom{i}n}_{|H|}\right)\frac{|H|!}{|\textrm{\scriptsize{Aut}}(H)|} \right) = O\left( n^{|H|} \right)$ 
possibilities 
for where the copy of $H$ was originally 
and $\left(^{\phantom{n}m} _{e(H)} \right) = O\left( n^{e(H)} \right)$ 
possibilities for which edges were inserted, 
so we have built each graph $O\left( n^{(|H|+e(H))} \right)$ times.

Therefore, 
$\frac{|\mathcal{G}_{n}|}{|\mathcal{P}(n,m)|} \!=\! O \left(\! \frac{n^{(|H|+e(H))}}{(dn)^{e(H)}} \!\right)
\!=\! O \left(\! \frac{n^{|H|}}{d^{e(H)}} \!\right) \!\to\! 0$,
since $d \!=\! \omega \left(\! n^{\frac{|H|}{e(H)}} \!\right).\!$~
\setlength{\unitlength}{.25cm}
\begin{picture}(1,1)
\put(0,0){\line(1,0){1}}
\put(0,0){\line(0,1){1}}
\put(1,1){\line(-1,0){1}}
\put(1,1){\line(0,-1){1}}
\end{picture}
\\
\\

Our main result of this section follows immediately:

\begin{Corollary} \label{msub3}
Let $H$ be a fixed multicyclic connected planar graph
and let $m(n)$ be such that
$\limsup_{n \to \infty} \frac{m}{n} < 1$.
Then
\begin{displaymath}
\mathbf{P}[P_{n,m} \textrm{ will have a copy of }H] \to 0 \textrm{ as } n \to \infty. 
\end{displaymath} 
\end{Corollary}

\phantom{p}

\phantom{p}

By the same proof as Theorem~\ref{msub2}, 
using Theorem~\ref{add54} instead of Theorem~\ref{add53}, 
we may also obtain a stronger result for when $e(H)>2|H|$:

\begin{Theorem} \label{msub1}
Let $H$ be a fixed connected planar graph such that $e(H) > 2|H|$
and let $m=m(n)$ be such that $\max \{ 0,m-n \} = o \left( n^{1- \frac{|H|}{e(H)}} \right)$.
Then 
\begin{displaymath}
\mathbf{P}[P_{n,m} \textrm{ will have a copy of }H] \to 0 \textrm{ as } n \to \infty.
\end{displaymath} 
\end{Theorem}

\phantom{p}

\phantom{p}

Note that Theorem~\ref{msub1} includes the case when $\frac{m}{n}$ converges slowly to $1$ from above.
In particular,
it holds when $m=n$. \\
\\
\\

We now have a complete account of
$\mathbf{P}[P_{n,m} \textrm{ will have a copy of } H]$,
in terms of exactly when it is bounded away from $0$ and $1$,
except for the case when $H$ is multicyclic and $\frac{m}{n} \to 1$.
We shall now finish Part I with a discussion of this remaining question.

Recall that, in Theorems~\ref{msub2} and~\ref{msub1},
we showed that the probability that $P_{n,m}$ will contain a copy of any given multicyclic subgraph $H$
converges to $0$ even for some $m(n)$ such that $\frac{m}{n} \to 1$.
Our previous result on components and subgraphs have all produced thresholds with neat divisions,
so this may suggest that
$\mathbf{P}[P_{n,m} \textrm{ will contain a copy of }H] \to 0$
whenever $\frac{m}{n} \leq 1+o(1)$.

The proofs of Theorems~\ref{msub2} and~\ref{msub1}
used the work of Section~\ref{add},
where we defined add$(n,m)$ to be the \textit{minimum} value of $|\textrm{add}(G)|$
over all graphs~$G \!\in~\!\!\mathcal{P}(n,m)$.
We already know 
(from the second half of Section~\ref{add})
that these results cannot be improved,
but of course it might be possible to find some function $a(n,m)$ with
$a(n,m) = \omega(\textrm{add}(n,m))$
such that $|\textrm{add}(P_{n,m})| \geq a(n,m)$~\textit{a.a.s.}.
This could then be used to improve our current multicyclic subgraph results.
For example, if we could show that $|\textrm{add}(P_{n,m})|$
is usually of the order of $\frac{n^{2}}{\ln n}$ whenever~$\frac{m}{n} \leq 1+o(1)$,
then we could probably use the method of the proofs of Theorems~\ref{msub2} and~\ref{msub1}
to show that
$\mathbf{P}[P_{n,m} \textrm{ will contain a copy of }H] \to 0$
whenever $\frac{m}{n} \leq 1+o(1)$. 

Let us now consider for a moment how we might work towards finding such a function $a(n,m)$.
Note that if we did show that $|\textrm{add}(P_{n,m})| \geq a(n,m)$ a.a.s.,
then we could probably use the proof of Theorem~\ref{tree4} to show that
the number of components in $P_{n,m}$ that are isomorphic to any given tree is a.a.s.~of the order of $\frac{a(n,m)}{n}$.
Conversely,
we know from the proof of Lemma~\ref{add33} that 
$|\textrm{add}(G)| \geq~\left( \kappa(G)-1 \right) \left( n-\frac{\kappa(G)}{2} \right)$,
so if $\kappa(P_{n,m})$ is a.a.s.~of the order of $\frac{a(n,m)}{n}$
then we would indeed have $|\textrm{add}(P_{n,m})| \geq a(n,m)$~a.a.s..
Thus, if we are interested in the typical value of $|\textrm{add}(P_{n,m})|$,
then it seems that it would actually be equivalent for us to investigate the typical value of $\kappa(P_{n,m})$,
or indeed the typical number of components in $P_{n,m}$
that are isomorphic to any given tree
(such as an isolated vertex). 

It is possible, of course, that for some multicyclic connected $H$,
the probability that $P_{n,m} \textrm{ will contain a copy of }H$
is bounded away from $0$ for some functions $m(n)$ with $\frac{m}{n} \to 1$.
In Lemma~\ref{unisub1},
we showed that for \textit{unicyclic} connected $H$
we have 
$\mathbf{P}[P_{n,m} \textrm{ will have a copy of }H] \to 1$
if $\frac{m}{n} \geq 1-o(1)$.
The proof involved turning pendant edges into an appearance of $H$ plus one isolated vertex,
and the crucial ingredient was the fact that $\kappa(P_{n,m})$,
and hence the number of isolated vertices in $P_{n,m}$,
is typically $o(n)$.
If we tried to use the same method for a graph $H$ with $e(H) = |H|+1$,
we would be left with two isolated vertices, 
and hence the proof would only work if $m(n)$ is such that
$\kappa(P_{n,m})$ is typically $o \left( n^{1/2} \right)$.
In general,
for $H$ such that $e(H) = |H|+r$,
the proof would work if $m(n)$ is such that $\kappa(P_{n,m})$ is typically $o \left( n^{\frac{1}{r+1}} \right)$.
From page~\pageref{cptsum},
we can see that it is not impossible that $\kappa(P_{n,m})$
is small enough for certain $m$,
and thus we may in fact be able to use this method to show that
we sometimes have
$\mathbf{P}[P_{n,m} \textrm{ will have a copy of }H] \to 1$
for a connected multicyclic~$H$ even when~$\frac{m}{n} \to 1$. 

Note that, again,
the topic of
$\mathbf{P}[P_{n,m} \textrm{ will have a copy of }H]$
seems to be linked to that of $\kappa(P_{n,m})$.
Thus, to conclude,
in order to discover the behaviour of
$\mathbf{P}[P_{n,m} \textrm{ will have a copy of }H]$
for multicyclic $H$ when $\frac{m}{n}$ is close to $1$,
it seems that it will be necessary to obtain more precise results on $\kappa(P_{n,m})$,
or on the equivalent topics of add$(P_{n,m})$ 
or the number of isolated vertices in $P_{n,m}$.

\newpage
\part[Random Planar Graphs with Bounds on the Minimum~and Maximum Degrees]
{Random Planar Graphs with Bounds on the Minimum and Maximum Degrees} 
\label{II}
\section{Outline of Part~\ref{II}} \label{bintro}

In Part~\ref{I}, we saw how the typical properties of a random planar graph,
such as being connected or containing given subgraphs/components,
change depending on the number of edges or, equivalently, the average degree.
In Part~\ref{II}, we shall now instead look at how the minimum and maximum degrees influence these properties.

For functions $d_{1}(n)$, $d_{2}(n)$, $D_{1}(n)$ and $D_{2}(n)$,
we will use $\mathcal{P}(n,d_{1},d_{2},D_{1},D_{2})$ 
to denote the set of all labelled planar graphs on $\{ 1,2, \ldots, n \}$
with minimum degree between $d_{1}(n)$ and $d_{2}(n)$ inclusive
and maximum degree between $D_{1}(n)$ and $D_{2}(n)$ inclusive.
We shall use $P_{n,d_{1},d_{2},D_{1},D_{2}}$ to denote a graph taken uniformly at random from this class.

Note that the minimum degree of a planar graph must be at most $5$.
Thus, for example, $P_{n,1,5,3,\log n}$ would denote 
a random planar graph with minimum degree at least $1$ and maximum degree between $3$ and $\log n$,
while $P_{n,0,5,0,n-1}$ would simply denote a random planar graph
(with no bounds on the degrees at all).
For $D_{2}(n)<3$, we should note that our random \textit{planar} graph 
$P_{n,d_{1},d_{2},D_{1},D_{2}}$ is just the same as a \textit{general} random graph 
with the same degree constraints.
Thus, 
since general random graphs have already been extensively investigated,
we shall only bother to concern ourselves here with the case when $D_{2}(n) \geq 3$~$\forall n$. 

In fact, it will turn out that most of our results for $P_{n,d_{1},d_{2},D_{1},D_{2}}$
will follow just from consideration of the case 
when there is no upper bound on the minimum degree and no lower bound on the maximum degree.
Thus, for most of Part~\ref{II} we will only look at $P_{n,d_{1},5,0,D_{2}}$
(i.e.~a random planar graph with all degrees between $d_{1}(n)$ and $D_{2}(n)$),
before then extending our results in Section~\ref{general} 
to cover any functions $d_{1}(n)$, $d_{2}(n)$, $D_{1}(n)$ and $D_{2}(n)$ 
as long as $\limsup_{n \to \infty} D_{1}(n) < \infty$.
We shall not attempt to provide a description of what happens when 
$\limsup_{n \to \infty} D_{1}(n) = \infty$,
but we will note one partial result for this case in Section~\ref{unbounded}.

The structure of Part~\ref{II} shall be based on that of \cite{mcd},
where the graph $P_{n,0,5,0,n-1}$ was studied.
Hence, we will start in Section~\ref{bconn}
by establishing a lower bound for the probability that $P_{n,d_{1},5,0,D_{2}}$ will be connected.
In Section~\ref{growth}, we shall use this to show  
that there exists a non-zero finite `growth constant' for $|\mathcal{P}(n,d_{1},5,0,D_{2})|$,
and in Section~\ref{apps} we will use this second fact to show 
that~(a.a.s.)~$P_{n,d_{1},5,0,D_{2}}$ has many appearance-type copies of certain $H$.
In Section~\ref{cpts}, we shall then use this last result to
deduce a lower bound for the probability that $P_{n,d_{1},5,0,D_{2}}$ will contain a component isomorphic to $H$,
and in Section~\ref{subs}
we will prove that 
$P_{n,d_{1},5,0,D_{2}}$ has many vertex-disjoint induced copies of most, but not all, $H$.

During Section~\ref{subs},
we shall come across the topic of determining whether or not a given graph $H$ 
can \textit{ever} be a subgraph of a $4$-regular planar graph.
This issue turns out to be quite intricate,
and Section~\ref{4reg} (which is joint work with Louigi Addario-Berry) 
is devoted to giving a polynomial-time algorithm for this problem.

As already noted, in Section~\ref{general}
we shall then extend our results on $P_{n,d_{1},5,0,D_{2}}$ 
to also cover any functions $d_{2}(n)$ and $D_{1}(n)$ as long as $D_{1}(n)$ is bounded,
before then looking briefly at the case $\limsup_{n \to \infty} D_{1}(n) = \infty$ in Section~\ref{unbounded}.
A simplified summary of our main results 
is given on page~\pageref{sum},
although we may sometimes actually prove slightly stronger versions, as with Part~\ref{I}.

\newpage
\subsection*{Summary of Results} \label{sum}

Given functions $d_{1}(n)$, $d_{2}(n)$, $D_{1}(n)$ and $D_{2}(n)$, \\
let $H_{1}$ be a connected planar graph with 
$\overline{\lim} d_{1}(n) \leq \delta(H_{1}) \leq \Delta(H_{1}) \leq \underline{\lim} D_{2}(n)$ \\
and let $H_{2} \textrm{ be a connected planar graph with } 
\Delta(H_{2}) \leq \underline{\lim} D_{2}(n)$. \\
Let $\mathbf{P}_{c} = \mathbf{P}[P_{n,d_{1},d_{2},D_{1},D_{2}} \textrm{ will be connected}]$, \\
let $\mathbf{P} = 
\mathbf{P}[P_{n,d_{1},d_{2},D_{1},D_{2}} \textrm{ will have a component isomorphic to } H_{1}]$ \\
and let $\mathbf{P}_{s} = \mathbf{P}[P_{n,d_{1},d_{2},D_{1},D_{2}} \textrm{ will have a subgraph isomorphic to } H_{2}]$.
\\
\\
The following results hold for all choices of $d_{1}(n),d_{2}(n),D_{1}(n)$ and $D_{2}(n)$, 
subject to $\overline{\lim} D_{1}(n) < \infty$ (and $D_{2}(n) \geq 3$, as always), 
unless otherwise stated.
\\
\\
\\
Connectivity Results
$\left\{ \begin{array}{l}
d_{2}(n)=0~\forall n
\left\{ \begin{array}{l} 
\textrm{\textbf{Observation: }} \textbf{P}_{c}=0 \textrm{ } \forall n \geq 2.
\end{array} \right. \\
d_{2}(n)>0~\forall n
\left\{ \begin{array}{l}
\textrm{\textbf{Theorems~\ref{bounded311} \&~\ref{bounded791}: }}
\underline{\lim} \mathbf{P}_{c} >0. \\
\textrm{\textbf{Theorems~\ref{bounded404} \&~\ref{bounded791}: }}
\overline{\lim} \mathbf{P}_{c} <1. \\
\end{array} \right.
\end{array} \right.$ \\
\\
\\
Component Results
$\left\{ \begin{array}{l} \!\!
\textrm{\textbf{Theorems~\ref{bounded404},~\ref{bounded791},~\ref{bounded147} \&~\ref{bounded803}: }}
\underline{\lim} \mathbf{P} >0. \\
\!\!(d_{2}(n),|H_{1}|) = (0,1)~\forall n 
\left\{ \begin{array}{l} \!\!
\textrm{\textbf{Observation: }} \textbf{P}=1 \textrm{ } \forall n.
\end{array} \right. \\
\!\!\exists
\begin{picture}(0.2,1)
\put(-0.7,-0.2){\line(1,3){0.5}}
\end{picture}
n\!:\!(d_{2}(n),|H_{1}|) \!=\! (0,1) 
\left\{ \! \begin{array}{l} \!\!
\textrm{\textbf{Thms.~\ref{bounded311},~\ref{bounded791},~\ref{bounded147} \&~\ref{bounded807}: }} \\
\!\!\overline{\lim} \mathbf{P}<1. 
\end{array} \right.
\end{array} \right.$ \\
\\
\\
Subgraph Results
$\left\{ \begin{array}{l}
\!\!D_{2}(n) \!=\! \delta(H_{2})~\!\forall n \left\{ \begin{array}{l}
\!\!\!
\textrm{\small{\textbf{Thms.~\ref{bounded1002},~\ref{bounded791},~\ref{bounded147} \&~\ref{bounded806}: }}}
\underline{\lim} \mathbf{P}_{s} \!>\!0. \\
\!\!\!
\textrm{\small{\textbf{Thms.~\ref{bounded1002},~\ref{bounded791},~\ref{bounded147} \&~\ref{bounded806}: }}}
\overline{\lim} \mathbf{P}_{s} \!<\!1. \\
\end{array} \right. \\
\begin{array}{l}
\!\!\!\!D_{2}(n)\!>\!\delta(H_{2})~\forall n~\&\\
\!\!\!\!\exists
\begin{picture}(0.2,1)
\put(-0.7,-0.2){\line(1,3){0.5}}
\end{picture}
n\!:\!d_{1}(n)\!=\!D_{2}(n)\!=\!4 
\end{array}
\left\{ \begin{array}{l}
\!\!\textrm{\textbf{Thms.~\ref{bounded1001},~\ref{bounded791},~\ref{bounded147} 
\&~\ref{bounded802}: }} \\
\!\!\mathbf{P}_{s} \to 1. 
\end{array} \right. \\
\begin{array}{l}
\!\!\!\!D_{2}(n)\!>\!\delta(H_{2})~\forall n~\&\\
\!\!\!\!d_{1}(n)\!=\!D_{2}(n)\!=\!4~\forall n 
\end{array}
\left\{ \begin{array}{l}
\!\!\textrm{\textbf{Thms.~\ref{bounded1001},~\ref{bounded791} \&~\ref{bounded147}: }} \\
\!\!\mathbf{P}_{s} \!\to\! 1 \textrm{ if } \exists \textrm{ $4$-reg.~planar } H^{*} \!\supset\! H_{2}. \\ 
\!\!\textrm{\textbf{Observation: }} 
\mathbf{P}_{s}=0 \textrm{ otherwise.} \\
\end{array} \right. 
\end{array} \right. $ 
\label{sum2}

\newpage
\section{Connectivity} \label{bconn}

We will start by examining the probability that our random graph is connected.
Not only is this topic interesting in its own right,
but the results given here will also be important ingredients in later sections.

As mentioned in the introduction,
until Section~\ref{general} we shall restrict ourselves to the case when
we have no upper bound on the minimum degree and no lower bound on the maximum degree.
Thus, we will be looking at $P_{n,d_{1},5,0,D_{2}}$,
which is simply a random planar graph on $\{ 1,2, \ldots, n \}$
with all degrees between $d_{1}(n)$ and $D_{2}(n)$.

Recall that we must have $D_{2}(n) \geq 3$ for planarity to have any impact.
The main result of this section will be to show (in Theorem~\ref{bounded311}) that,
given any function $d_{1}(n)$ and any function $D_{2}(n)$ satisfying $D_{2}(n) \geq 3$~$\forall n$,
we have 
$\liminf_{n \to \infty} \mathbf{P}[P_{n,d_{1},5,0,D_{2}} \textrm{ will be connected}]>0$.
An upper bound for this probability shall be deduced later, in Section~\ref{cpts}, 
from results on components.

The proof of Theorem~\ref{bounded311} will copy that of Proposition~\ref{mcd 2.1}(\cite{mcd},~2.1),
but will be slightly more complicated,
as the upper bound on the maximum degree means that $\mathcal{P}(n,d_{1},5,0,D_{2})$ is not edge-addable
(i.e.~the class $\mathcal{P}(n,d_{1},5,0,D_{2})$ 
is not closed under the operation of inserting an edge between two components).
Hence, we shall first prove (in Lemma~\ref{bounded1}) a very helpful result on short cycles,
which will be extremely useful to us throughout Part~\ref{II}.

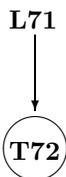
\begin{figure} [ht]
\setlength{\unitlength}{1cm}
\begin{picture}(20,3.2)(-5.6,0)

\put(0.4,0.4){\circle{0.8}}
\put(0.4,1.9){\vector(0,-1){1}}
\put(0.05,2){\textbf{L\ref{bounded1}}}
\put(0.05,0.25){\textbf{T\ref{bounded311}}}

\end{picture}

\caption{The structure of Section~\ref{bconn}.}
\end{figure}

\begin{displaymath}
\end{displaymath}

We shall now proceed with our aforementioned useful lemma:

\begin{Lemma} \label{bounded1}
Let $k < \frac{1}{15}$,
and let $S$ be a planar graph with at most $k|S|$ vertices of degree $\leq 2$.
Then $S$ must contain at least $\left( \frac{1-15k}{28} \right) |S|$ cycles of size $\leq 6$.
In particular, if $S$ has $\leq \frac{|S|}{43}$ vertices of degree $\leq 2$
then $S$ must contain $\geq \frac{|S|}{43}$ cycles of size $\leq 6$.
\end{Lemma}
\textbf{Proof}
Fix a planar embedding of $S$.
We shall first show that this embedding must have at least 
$\left( \frac{1-15k}{14} \right) |S|$ \textit{faces} of size $\leq 6$
(where, as usual, the `size' of a face denotes the number of edges in the associated facial boundary,
with an edge counted twice if it appears twice in the boundary),
and we will later deduce the lemma from this fact.

We shall argue by contradiction.
Let $f_{i}$ denote the number of faces of size $i$ and 
suppose that $\sum_{i \leq 6} f_{i} < \left( \frac{1-15k}{14} \right) |S|$.
We have
\begin{eqnarray*}
2e(S) & = & \sum_{i} i f_{i} \\
& \geq & 7 \sum_{i \geq 7} f_{i} \\
& > & 7 \left( \sum_{i} f_{i} - \left( \frac{1-15k}{14} \right) |S| \right), 
\textrm{ by our supposition} \\
& = & 7 \left( e(S) - |S| + \kappa(S) + 1- \left( \frac{1-15k}{14} \right) |S| \right),
\textrm{ by Euler's formula} \\
& > & 7 \left( e(S) - \left( \frac{15(1-k)}{14} \right) |S| \right). 
\end{eqnarray*}
Thus,
$\left( \frac{15(1-k)}{2} \right) |S| > 5e(S)$.
But $e(S) \geq \frac{3(1-k)|S|}{2}$,
since $S$ contains at least $(1-k)|S|$ vertices of degree $\geq 3$,
and so $5e(S) \geq \left( \frac{15(1-k)|S|}{2} \right)$.
Thus, we obtain our desired contradiction,
and so we must have at least $\left( \frac{1-15k}{14} \right) |S|$ faces of size~$\leq 6$.

Let us now consider these faces of size $\leq 6$.
Note that the boundary of a face of size $\leq 6$ must contain a cycle of size $\leq 6$
as a subgraph unless it is acyclic,
in which case it must be the entire graph $S$.
But if $S$ were acyclic,
then at least half of the vertices would have degree $\leq 2$
(since we would have $e(S) \leq |S|-1$),
and this would contradict the conditions of this lemma.
Thus, for each of our faces of size $\leq 6$,
the boundary must contain a cycle of size $\leq 6$ as a subgraph.

Each edge of $S$ can only be in at most two faces of the embedding,
and so each cycle can only be in at most two faces.
Thus, $S$ must contain at least $\left( \frac{1-15k}{28} \right) |S|$ \textit{distinct} cycles of size $\leq 6$.
$\phantom{qwerty}$ 
\setlength{\unitlength}{0.25cm}
\begin{picture}(1,1)
\put(0,0){\line(1,0){1}}
\put(0,0){\line(0,1){1}}
\put(1,1){\line(-1,0){1}}
\put(1,1){\line(0,-1){1}}
\end{picture} \\
\\

Note that Lemma~\ref{bounded1} does not hold for general graphs,
since it is known (see, for example, Corollary 11.2.3 of \cite{die})
that there exist graphs that have both arbitrarily large minimum degree and arbitrarily large girth. \\
\\
\\

We shall now use Lemma~\ref{bounded1} to obtain our aforementioned lower bound for
the probability that $P_{n,d_{1},5,0,D_{2}}$ will be connected
(we shall state the result in a slightly different form to that advertised earlier,
for ease with Section~\ref{growth}):

\begin{Theorem} \label{bounded311}
There exists a constant $c >0$ such that,
given any constants $r,d_{1},D_{2} \in \mathbf{N} \cup \{ 0 \}$ with $D_{2} \geq 3$
and $\mathcal{P}(r,d_{1},5,0,D_{2}) \neq \emptyset$,
\begin{displaymath}
\mathbf{P} [P_{r,d_{1},5,0,D_{2}} \textrm{ will be connected}] > c.
\end{displaymath}
\end{Theorem}
\textbf{Sketch of Proof} 
We shall choose any $r,d_{1},D_{2} \in \mathbf{N} \cup \{ 0 \}$ with $D_{2} \geq 3$ and 
show that there are many ways to construct a graph in $\mathcal{P}(r,d_{1},5,0,D_{2})$ with $k-1$ components
from a graph in $\mathcal{P}(r,d_{1},5,0,D_{2})$ with $k$ components,
by combining two components.
Our stated lower bound for the proportion of graphs with exactly one component will then follow
by `cascading' this result downwards.

If $D_{2} > 6$,
we shall see that we may obtain sufficiently many ways to combine components
simply by inserting edges between them that don't interfere with this upper bound on the maximum degree.

If $D_{2} \leq 6$,
we will sometimes also delete an edge from a small cycle
in order to maintain $\Delta \leq D_{2}$.
We shall use Lemma~\ref{bounded1} to show that
we have lots of choices for this small cycle,
and then the fact that it is small
(combined with the knowledge that $D_{2}<7$) 
will help us to bound the amount of double-counting. \\
\\
\textbf{Full Proof}
Choose any $r,d_{1},D_{2} \in \mathbf{N} \cup \{ 0 \}$ with $D_{2} \geq 3$.
We shall show that there exists a strictly positive constant $c$,
independent of $r,d_{1}$ and $D_{2}$,
such that
$\mathbf{P} [P_{r,d_{1},5,0,D_{2}} \textrm{ will be connected}] > c$.

Let $\mathcal{P}^{t}(r,d_{1},5,0,D_{2})$ denote the set of graphs in 
$\mathcal{P}(r,d_{1}, 5,0,D_{2})$ with exactly $t$ components.
For $k>1$, we shall construct graphs in $\mathcal{P}^{k-1}(r,d_{1},5,0,D_{2})$
from graphs in $\mathcal{P}^{k}(r,d_{1},5,0,D_{2})$.

Let $G \in \mathcal{P}^{k}(r,d_{1},5,0,D_{2})$
and let us denote the $k$ components of $G$ by $S_{1}, S_{2}, \ldots, S_{k}$,
where $|S_{i}|=n_{i}$ $\forall i$.
Without loss of generality,
we may assume that $S_{1}, S_{2}, \ldots, S_{k}$ are ordered so that 
$S_{i}$ contains $\geq \frac{n_{i}}{43}$ vertices of degree $<D_{2}$ iff $i \leq l$, 
for some fixed $l \in \{ 0,1, \ldots, k \}$.
Note that we must have $l=k$ if $D_{2} > 6$,
since (by planarity) $e(S_{i}) < 3n_{i}$
and so we can only have at most $\frac{6n_{i}}{7}$ vertices of degree~$\geq 7$.

For $1 \leq i < j \leq k$,
let us construct a new graph $G_{i,j} \in \mathcal{P}^{k-1}(r,d_{1},5,0,D_{2})$ as follows: \\
\\
Case (a): if $j \leq l$ (note that this is always the case if $D_{2} > 6$) \\
Insert an edge between a vertex in $S_{i}$ of degree $<D_{2}$
(we have at least $\frac{n_{i}}{43}$ choices for this)
and a vertex in $S_{j}$ of degree $<D_{2}$
(we have at least $\frac{n_{j}}{43}$ choices for this).

\begin{figure} [ht]
\setlength{\unitlength}{1cm}
\begin{picture}(20,2.75)(-1,-0.75)

\put(-0.5,0.5){\line(0,1){0.5}}
\put(1,0.5){\line(0,1){0.5}}
\put(2,0.5){\line(0,1){0.5}}
\put(3.5,0.5){\line(0,1){0.5}}
\put(6.5,0.5){\line(0,1){0.5}}
\put(8,0.5){\line(0,1){0.5}}
\put(9,0.5){\line(0,1){0.5}}
\put(10.5,0.5){\line(0,1){0.5}}

\put(0.25,1){\oval(1.5,1.2)[t]}
\put(2.75,1){\oval(1.5,1.2)[t]}
\put(7.25,1){\oval(1.5,1.2)[t]}
\put(9.75,1){\oval(1.5,1.2)[t]}

\put(0.25,0.5){\oval(1.5,1.2)[b]}
\put(2.75,0.5){\oval(1.5,1.2)[b]}
\put(7.25,0.5){\oval(1.5,1.2)[b]}
\put(9.75,0.5){\oval(1.5,1.2)[b]}

\put(4.5,0.75){\vector(1,0){1}}

\put(8,0.75){\line(1,0){1}}

\put(1,0.75){\circle*{0.1}}
\put(2,0.75){\circle*{0.1}}
\put(8,0.75){\circle*{0.1}}
\put(9,0.75){\circle*{0.1}}

\put(0.1,-0.75){$S_{i}$}
\put(2.6,-0.75){$S_{j}$}

\end{picture}

\caption{Constructing the graph $G_{i,j}$ in case (a).} \label{bconnfig}
\end{figure}
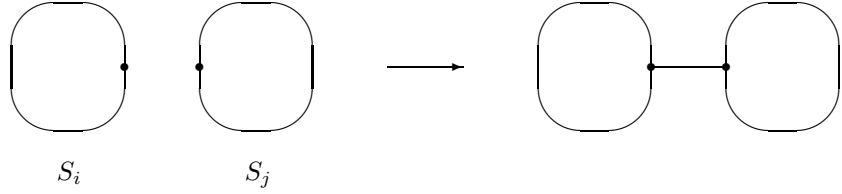

The constructed graph $G_{i,j}$ (see Figure~\ref{bconnfig}) is planar and has exactly $k-1$ components.
It is also clear that we have $d_{1} \leq \delta (G_{i,j}) \leq \Delta (G_{i,j}) \leq D_{2}$,
since we have not deleted any edges from the original graph
and have only inserted an edge between two vertices with degree $<D_{2}$.
Thus, $G_{i,j} \in \mathcal{P}^{k-1}(r,d_{1},5,0,D_{2})$. \\
\\
Case (b): if $j>l$ and $n_{i}>1$ (in which case $D_{2} \leq 6$) \\
If $j>l$,
then $S_{j}$ contains $<\frac{n_{j}}{43}$ vertices of degree $<D_{2}$.
Thus, by Lemma~\ref{bounded1},
$S_{j}$ must contain at least $\frac{n_{j}}{43}$ cycles of size $\leq 6$.
Delete an edge $uv$ in one of these cycles
(we have at least 
$\frac{3}{D_{2}+D_{2}^{2}+D_{2}^{3}+D_{2}^{4}}\frac{n_{j}}{43} \geq \frac{3}{6+6^{2}+6^{3}+6^{4}}\frac{n_{j}}{43}$ 
choices for this edge,
since each cycle must contain at least 
$3$ edges and each edge is in at most $(D_{2}-1)^{m-2}<D_{2}^{m-2}$ cycles of size~$m$),
insert an edge between $u$ and a vertex~$w \in S_{i}$
(we have $n_{i}$ choices for $w$),
delete an edge between $w$ and $x \in \Gamma(w)$
(we have at least one choice for $x$, since $n_{i}>1$),
and insert an edge between $x$ and $v$
(planarity is preserved, 
since we may draw $S_{j}$ so that the face containing $u$ and $v$ is on the outside,
and similarly we may draw $S_{i}$ so that the face containing $w$ and $x$ is on the outside).

\begin{figure} [ht]
\setlength{\unitlength}{1cm}
\begin{picture}(20,2.75)(-1,-0.75)

\put(-0.5,0.5){\line(0,1){0.5}}
\put(1,0.5){\line(0,1){0.5}}
\put(2,0.5){\line(0,1){0.5}}
\put(3.5,0.5){\line(0,1){0.5}}
\put(6.5,0.5){\line(0,1){0.5}}
\put(10.5,0.5){\line(0,1){0.5}}

\put(0.25,1){\oval(1.5,1.2)[t]}
\put(2.75,1){\oval(1.5,1.2)[t]}
\put(7.25,1){\oval(1.5,1.2)[t]}
\put(9.75,1){\oval(1.5,1.2)[t]}

\put(0.25,0.5){\oval(1.5,1.2)[b]}
\put(2.75,0.5){\oval(1.5,1.2)[b]}
\put(7.25,0.5){\oval(1.5,1.2)[b]}
\put(9.75,0.5){\oval(1.5,1.2)[b]}

\put(4.5,0.75){\vector(1,0){1}}

\put(8,0.5){\line(1,0){1}}
\put(8,1){\line(1,0){1}}

\put(1,0.5){\circle*{0.1}}
\put(2,0.5){\circle*{0.1}}
\put(8,0.5){\circle*{0.1}}
\put(9,0.5){\circle*{0.1}}
\put(1,1){\circle*{0.1}}
\put(2,1){\circle*{0.1}}
\put(8,1){\circle*{0.1}}
\put(9,1){\circle*{0.1}}

\put(2,0.75){\oval(1,0.5)[r]}
\put(9,0.75){\oval(1,0.5)[r]}

\put(1.75,1.1){$u$}
\put(1.75,0.2){$v$}
\put(1.1,1.1){$w$}
\put(1.1,0.2){$x$}
\put(8.75,1.1){$u$}
\put(8.75,0.2){$v$}
\put(8.1,1.1){$w$}
\put(8.1,0.2){$x$}

\put(0.1,-0.75){$S_{i}$}
\put(2.6,-0.75){$S_{j}$}

\end{picture}

\caption{Constructing the graph $G_{i,j}$ in case (b).}
\end{figure}
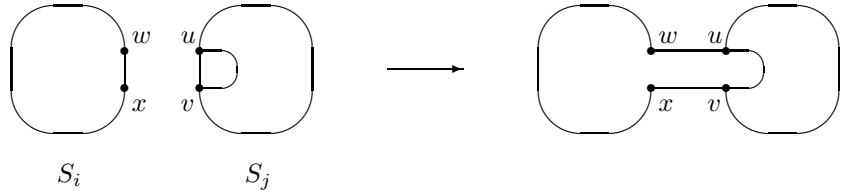

Since the deleted edge $uv$ was in a cycle,
it was not a cut-edge,
and so the vertex set $V(S_{j})$ is still connected.
The deleted edge $wx$ may have been a cut-edge in $S_{i}$,
but since we have also inserted edges from $w$ to $u \in V(S_{j})$ and from $x$ to $v \in V(S_{j})$
it must be that the vertex set $V(S_{i}) \cup V(S_{j})$ is now connected.
Thus, the constructed planar graph $G_{i,j}$ has exactly $k-1$ components.
By construction, the degrees of the vertices have not changed,
and so we have $d_{1} \leq \delta (G_{i,j}) \leq \Delta (G_{i,j}) \leq D_{2}$.
Thus, $G_{i,j} \in \mathcal{P}^{k-1}(r,d_{1},5,0,D_{2})$. \\
\\
Case (c): if $j>l$ and $n_{i}=1$ (in which case $D_{2} \leq 6$) \\
Delete any edge $uv$ in $S_{j}$
(we have at least $n_{j}$ choices for this,
since $S_{j}$ cannot be a forest if $j>l$)
and insert edges $uw$ and $vw$,
where $w$ is the unique vertex in $S_{i}$.

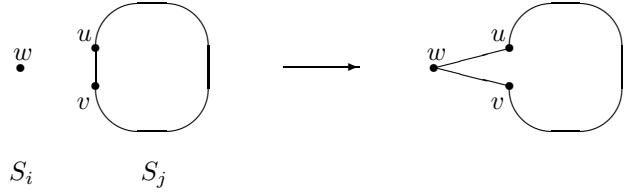
\begin{figure} [ht]
\setlength{\unitlength}{1cm}
\begin{picture}(20,2.3)(-1,-0.75)

\put(2,0.5){\line(0,1){0.5}}
\put(3.5,0.5){\line(0,1){0.5}}
\put(9,0.5){\line(0,1){0.5}}

\put(2.75,1){\oval(1.5,1.2)[t]}
\put(8.25,1){\oval(1.5,1.2)[t]}

\put(2.75,0.5){\oval(1.5,1.2)[b]}
\put(8.25,0.5){\oval(1.5,1.2)[b]}

\put(4.5,0.75){\vector(1,0){1}}

\put(6.5,0.75){\line(4,1){1}}
\put(6.5,0.75){\line(4,-1){1}}

\put(1,0.75){\circle*{0.1}}
\put(2,0.5){\circle*{0.1}}
\put(2,1){\circle*{0.1}}
\put(6.5,0.75){\circle*{0.1}}
\put(7.5,0.5){\circle*{0.1}}
\put(7.5,1){\circle*{0.1}}

\put(1.75,1.1){$u$}
\put(1.75,0.2){$v$}
\put(0.9,0.85){$w$}
\put(7.25,1.1){$u$}
\put(7.25,0.2){$v$}
\put(6.4,0.85){$w$}

\put(0.85,-0.75){$S_{i}$}
\put(2.6,-0.75){$S_{j}$}

\end{picture}

\caption{Constructing the graph $G_{i,j}$ in case (c).}
\end{figure}

The constructed graph $G_{i,j}$ is planar and has exactly $k-1$ components.
The degrees have not changed,
except that we now have deg$(w)\!=\!2$.
But~since~$D_{2} \!\geq~\!\!3$,
we still have $d_{1} \leq \delta(G_{i,j}) \leq \Delta(G_{i,j}) \leq D_{2}$.
Thus, $G_{i,j} \in~\mathcal{P}^{k-1}(r,d_{1},5,0,D_{2})$.
\\

Let $z= \frac{3}{43(6+6^{2}+6^{3}+6^{4})}
= \min \left\{ \left( \frac{1}{43} \right)^{2}, \frac{3}{43(6+6^{2}+6^{3}+6^{4})} \right\}$.
Then in all cases
we have at least $zn_{i}n_{j}$ choices when constructing the new graph.
Thus, from our initial graph $G$,
we have at least $\sum_{i<j} zn_{i}n_{j} = z \sum_{i<j} n_{i}n_{j}$
ways to construct a graph in $\mathcal{P}^{k-1}(r,d_{1},5,0,D_{2})$.
Note that if $x \leq y$ then $xy > (x-~1)(y+1)$,
so $\sum_{i<j} n_{i}n_{j}$
is at least what it would be if one component in $G$ had order $r-(k-1)$
and the other $k-1$ components were isolated vertices.
Thus,
$z \sum_{i<j} n_{i}n_{j} \geq z \left( \frac{1}{2}(k\!-\!1)(k\!-\!2) \!+\! (k\!-\!1)(r\!-\!k\!+\!1) \right) 
= (k\!-\!1) \left( r\!-\!\frac{k}{2} \right)z$.
Therefore, for $k>~\!1$,
we have at least 
$(k\!-\!1) \left( r\!-\!\frac{k}{2} \right) z|\mathcal{P}^{k}(r,d_{1},5,0,D_{2})|
\geq (k\!-\!1) \frac{r}{2} z|\mathcal{P}^{k}(r,d_{1},5,0,D_{2})|$
ways to construct a graph in $\mathcal{P}^{k-1}(r,d_{1},5,0,D_{2})$. \\

Given one of our constructed graphs,
there are at most $3$ possibilities for how the graph was obtained (case (a), (b) or (c)).

If case (a) was used (which must be so if $D_{2} > 6$),
then we can re-obtain the original graph simply by deleting the inserted edge,
for which there are at most $r-(k-1)<r$ possibilities,
since it must now be a cut-edge.
Thus, if case (a) was used, 
we have $<r$ possibilities for the original graph.

If case (b) was used, 
then we can re-obtain the original graph by locating the vertices $u,v,w$ and $x$,
deleting the two inserted edges ($uw$ and $vx$)
and re-inserting the two deleted edges ($uv$ and $wx$).
Note that we have at most $r$ possibilities for which vertex is $u$.
We know that $u$ and $v$ were originally on a cycle of size $\leq 6$,
and so $v$ is still at distance at most $5$ from $u$.
Since the graph has maximum degree at most $D_{2}$,
we therefore have at most $D_{2}^{2}+D_{2}^{3}+D_{2}^{4}+D_{2}^{5}$ possibilities for $v$.
Once we have located $u$ and $v$,
we then have at most $D_{2}$ possibilities for $w$
and at most $D_{2}$ possibilities for $x$,
since $w$ and $x$ are now neighbours of $u$ and $v$, respectively.
Thus, if case (b) was used,
we have at most
$D_{2}^{2}(D_{2}^{2}+D_{2}^{3}+D_{2}^{4}+D_{2}^{5})r \leq 36(6^{2}+6^{3}+6^{4}+6^{5})r$ 
possibilities for the original graph.

If case (c) was used,
then we can re-obtain the original graph by locating the vertices $u,v$ and $w$,
deleting the two inserted edges ($uw$ and $vw$)
and re-inserting the deleted edge ($uv$).
We have at most $r$ possibilities for which vertex is $w$,
and given $w$ we then know which edges to delete and insert,
as $v$ and $w$ are the only vertices adjacent to $u$.
Thus, if case (c) was used,
we have $\leq r$ possibilities for the original graph.

Therefore,
there are $<r$ possibilities for the original graph if $D_{2} > 6$,
since case (a) must have been used,
and 
$< r + 36(6^{2}+6^{3}+6^{4}+6^{5})r +r = 2r(1+18(6^{2}+6^{3}+6^{4}+6^{5}))$ 
possibilities for the original graph if $D_{2} \leq 6$,
since any of case (a), case (b) or case (c) may have been used. \\

Let 
$\alpha = \frac{z}{4(1+18(6^{2}+6^{3}+6^{4}+6^{5}))}
= \min \left\{ \frac{z}{2}, \frac{z}{4(1+18(6^{2}+6^{3}+6^{4}+6^{5}))} \right\}$.
Then we have shown that
we can construct at least
$\alpha(k-1) |\mathcal{P}^{k}(r,d_{1},5,0,D_{2})|$ 
\textit{distinct} graphs in $\mathcal{P}^{k-1}(r,d_{1},5,0,D_{2})|$,
and so we find that it must be that we have
$
|\mathcal{P}^{k-1}(r,d_{1},5,0,D_{2})| \geq \alpha(k-1) |\mathcal{P}^{k}(r,d_{1},5,0,D_{2})| \textrm{ } \forall k>1.
$

Let us define
$p_{k}$
to be
$\frac{|\mathcal{P}^{k+1}(r,d_{1},5,0,D_{2})|}{|\mathcal{P}(r,d_{1},5,0,D_{2})|}$
and let 
$p=p_{0}= \frac{|\mathcal{P}^{1}(r,d_{1},5,0,D_{2})|}{|\mathcal{P}(r,d_{1},5,0,D_{2})|} = 
\mathbf{P}[P_{r,d_{1},5,0,D_{2}} \textrm{ will be connected}]$.
From the previous paragraph,
we know that 
$|\mathcal{P}^{k+1}(r,d_{1},5,0,D_{2})| \!\leq\! \frac{|\mathcal{P}^{k}(r,d_{1},5,0,D_{2})|}{\alpha k}$~
$\!\forall k\!>\!0$,
and so $p_{k} \!\leq\! \frac{p}{\alpha^{k}k!}$~$\forall k \!\geq\! 0$.~
We~must have $\sum_{k \geq 0} p_{k} \!=\! 1$,
so $\sum_{k \geq 0} \frac{p}{\alpha^{k}k!} \!\geq\! 1$
and hence 
$p \!\geq\! \left( \sum_{k \geq 0} \frac{\left( \frac{1}{\alpha} \right)^{k}}{k!} \right)^{-1} 
\!=~\!e^{-\frac{1}{\alpha}}$.~
\setlength{\unitlength}{0.25cm}
\begin{picture}(1,1)
\put(0,0){\line(1,0){1}}
\put(0,0){\line(0,1){1}}
\put(1,1){\line(-1,0){1}}
\put(1,1){\line(0,-1){1}}
\end{picture} \\
\\

\newpage
\section{Growth Constants} \label{growth}

We shall now look at the topic of `growth constants',
which will play a vital role in the proofs of Section~\ref{apps}.
We already know from Section~\ref{previous} that there exists a finite growth constant $\gamma_{l}>0$
such that $\left( \frac{|\mathcal{P}(n,0,5,0,n-1)|}{n!} \right)^{1/n} \to \gamma_{l}$ as $n \to \infty$.
In this section,
we shall use our connectivity bound from Theorem~\ref{bounded311} 
to also obtain (in Theorems~\ref{bounded104} and~\ref{bounded105})
growth constants for $\mathcal{P}(n,d_{1},5,0,D_{2})$
for the case when $d_{1}(n)$ is a constant and $D_{2}(n)$ is any monotonically non-decreasing function
(it will turn out that the result for this restricted case is all that will be required for our later sections). 

Clearly, we shall have to treat the case when $d_{1}(n)=D_{2}(n) \in \{ 3,5 \}$
slightly differently,
since $\mathcal{P}(n,d_{1},5,0,D_{2})$ will be empty if $n$ is odd.
Hence, we will deal first with the more standard case in Theorem~\ref{bounded104},
and then separately with this special case in Theorem~\ref{bounded105}.
In both of these cases,
we shall follow the proof of Theorem 3.3 of \cite{mcd},
which will require us to first state two useful lemmas (Proposition~\ref{vanL 11.6} and Corollary~\ref{bounded103}).

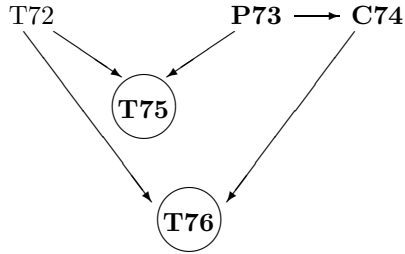
\begin{figure} [ht]
\setlength{\unitlength}{1cm}
\begin{picture}(20,4.5)(-2.6,0)

\put(0.8,3){T\ref{bounded311}}
\put(2.25,1.75){\textbf{T\ref{bounded104}}}
\put(2.85,0.25){\textbf{T\ref{bounded105}}}
\put(3.75,3){\textbf{P\ref{vanL 11.6}}}
\put(5.35,3){\textbf{C\ref{bounded103}}}

\put(2.6,1.9){\circle{0.8}}
\put(3.2,0.4){\circle{0.8}}

\put(1,2.9){\vector(3,-4){1.7}}
\put(1.4,2.9){\vector(3,-2){0.9}}
\put(4.6,3.1){\vector(1,0){0.6}}
\put(5.4,2.9){\vector(-3,-4){1.7}}
\put(3.8,2.9){\vector(-3,-2){0.9}}

\end{picture}

\caption{The structure of Section~\ref{growth}.}
\end{figure}

\begin{displaymath}
\end{displaymath}

As mentioned, the following result on supermultiplicative functions shall be very useful:

\begin{Proposition}[implicit in \cite{vanL}, 11.6] \label{vanL 11.6}
Let $f:\mathbf{N} \to \mathbf{R}^{+}$ be a function such that
$f(n)>0$ for all large $n$ and
$f(i+j) \geq f(i) \cdot f(j)$ $\forall i,j \in \mathbf{N}$.
Then $(f(n))^{1/n} \to~\sup_{n} \left( (f(n))^{1/n} \right)$ as $n \to \infty$. \\
\end{Proposition}

We should also note an even parity version,
which will be useful when $d_{1}(n)=~D_{2}(n) \in \{ 3,5 \}$:

\begin{Corollary} \label{bounded103}
Let $f:\mathbf{N} \to \mathbf{R}^{+}$ be a function such that
$f(2n)>0$ for all large $n$ and
$f(2(i+j)) \geq f(2i) \cdot f(2j)$ $\forall i,j \in \mathbf{N}$.
Then $(f(2n))^{1/2n} \to \sup_{n} \left( (f(2n))^{1/2n} \right)$ as $n \to \infty$. 
\end{Corollary}
\begin{displaymath}
\end{displaymath}

We may now use Proposition~\ref{vanL 11.6} to obtain our first growth constant result:

\begin{Theorem} \label{bounded104}
Let $d_{1} \in \{ 0,1, \ldots, 5 \}$ be a constant and 
let $D_{2}(n)$ be a monotonically non-decreasing integer-valued function that for all large $n$ satisfies
$D_{2}(n) \geq \max \{ d_{1}, 3 \}$ and $(d_{1},D_{2}(n)) \notin \{ (3,3), (5,5) \}$.
Then there exists a \emph{finite} constant $\gamma_{d_{1},D_{2}} >0$ such that
\begin{displaymath}
\left( \frac{|\mathcal{P}(n,d_{1},5,0,D_{2})|}{n!} \right)^{\frac{1}{n}} \to 
\gamma_{d_{1},D_{2}} 
\textrm{ as } n \to \infty.
\end{displaymath}
\end{Theorem}
\textbf{Proof}
We shall copy the method of proof of Theorem 3.3 of \cite{mcd}.
Let $c$ be the constant given by Theorem~\ref{bounded311} and
let $g(n,d_{1},D_{2}) = \frac{c^{2}|\mathcal{P}(n,d_{1},5,0,D_{2})|}{2 \cdot n!}$~$\forall n \in \mathbf{N}$.
We shall show that $g(n,d_{1},D_{2})$ satisfies the conditions of Proposition~\ref{vanL 11.6},
which we shall then use to deduce our result. \\

To show that $g(n,d_{1},D_{2})>0$ for all large $n$,
it suffices to show that $\mathcal{P}(n,d_{1},5,0,D_{2})$ is non-empty for all large $n$.
We shall now see that we only need to prove this for three values of $(d_{1},D_{2}(n))$: \\
\\
(i) If $d_{1}<3$,
then $\mathcal{P}(n,d_{1},5,0,D_{2}) \supset \mathcal{P}(n,2,5,0,2)$ for all large $n$, 
since $D_{2}(n) \geq 3$ for all large $n$. \\
(ii) If $d_{1} \in \{3,4\}$,
then $D_{2}(n)$ must be at least $4$ for all large $n$ 
(since $(d_{1},D_{2}(n))$ is not allowed to be $(3,3)$ for large $n$)
and so $\mathcal{P}(n,d_{1},5,0,D_{2}) \supset \mathcal{P}(n,4,5,0,4)$ for all large $n$. \\
(iii) If $d_{1}=5$, then $D_{2}(n)$ must be at least $6$ for all large $n$
(since $(d_{1},D_{2}(n))$ is not allowed to be $(5,5)$ for large $n$),
and so $\mathcal{P}(n,d_{1},5,0,D_{2}) \supset \mathcal{P}(n,5,5,0,6)$ for all large $n$. \\
Thus, it suffices just to show that 
$\mathcal{P}(n,2,5,0,2), \mathcal{P}(n,4,5,0,4)$ and $\mathcal{P}(n,5,5,0,6)$ 
are all non-empty for sufficiently large $n$.

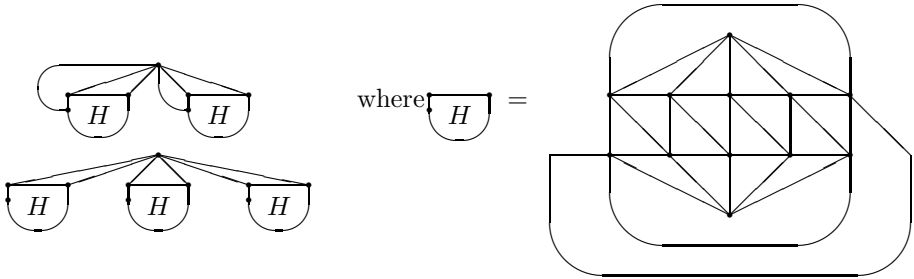
\begin{figure} [ht] 
\setlength{\unitlength}{0.8cm}
\begin{picture}(16,4.5)(-1,-3)

\put(1.5,0.5){\line(-3,-1){1.5}}
\put(1.5,0.5){\line(-1,-1){0.5}}
\put(1.5,0.5){\line(1,-1){0.5}}
\put(1.5,0.5){\line(3,-1){1.5}}
\put(0.5,-0.2){\oval(1,1)[b]}
\put(2.5,-0.2){\oval(1,1)[b]}
\put(0,0){\line(0,-1){0.25}}
\put(1,0){\line(0,-1){0.25}}
\put(0,-0.25){\circle*{0.1}}
\put(2,0){\line(0,-1){0.25}}
\put(3,0){\line(0,-1){0.25}}
\put(2,-0.25){\circle*{0.1}}
\put(0,0){\line(1,0){1}}
\put(2,0){\line(1,0){1}}
\put(1.5,0.5){\circle*{0.1}}
\put(0,0){\circle*{0.1}}
\put(1,0){\circle*{0.1}}
\put(2,0){\circle*{0.1}}
\put(3,0){\circle*{0.1}}
\put(0.31,-0.49){$H$}
\put(2.31,-0.49){$H$}
\put(0,-0.25){\line(-1,1){0.25}}
\put(0,0.125){\oval(1,0.75)[l]}
\put(1.5,0.5){\line(-1,0){1.5}}
\put(2,0.5){\oval(1,1.5)[bl]}
\put(1.5,-1){\line(-5,-1){2.5}}
\put(1.5,-1){\line(-3,-1){1.5}}
\put(1.5,-1){\line(-1,-1){0.5}}
\put(1.5,-1){\line(1,-1){0.5}}
\put(1.5,-1){\line(3,-1){1.5}}
\put(1.5,-1){\line(5,-1){2.5}}
\put(-0.5,-1.75){\oval(1,1)[b]}
\put(1.5,-1.75){\oval(1,1)[b]}
\put(3.5,-1.75){\oval(1,1)[b]}
\put(-1,-1.5){\line(0,-1){0.25}}
\put(0,-1.5){\line(0,-1){0.25}}
\put(1,-1.5){\line(0,-1){0.25}}
\put(2,-1.5){\line(0,-1){0.25}}
\put(3,-1.5){\line(0,-1){0.25}}
\put(4,-1.5){\line(0,-1){0.25}}
\put(-1,-1.75){\circle*{0.1}}
\put(1,-1.75){\circle*{0.1}}
\put(3,-1.75){\circle*{0.1}}
\put(-1,-1.5){\line(1,0){1}}
\put(1,-1.5){\line(1,0){1}}
\put(3,-1.5){\line(1,0){1}}
\put(1.5,-1){\circle*{0.1}}
\put(-1,-1.5){\circle*{0.1}}
\put(0,-1.5){\circle*{0.1}}
\put(1,-1.5){\circle*{0.1}}
\put(2,-1.5){\circle*{0.1}}
\put(3,-1.5){\circle*{0.1}}
\put(4,-1.5){\circle*{0.1}}
\put(-0.69,-1.99){$H$}
\put(1.31,-1.99){$H$}
\put(3.31,-1.99){$H$}

\put(4.8,-0.2){where}
\put(6,0){\line(1,0){1}}
\put(6,0){\line(0,-1){0.25}}
\put(7,0){\line(0,-1){0.25}}
\put(6.5,-0.25){\oval(1,1)[b]}
\put(6,0){\circle*{0.1}}
\put(7,0){\circle*{0.1}}
\put(6,-0.25){\circle*{0.1}}
\put(6.31,-0.49){$H$}
\put(7.3,-0.2){$=$}

\put(9,0){\line(1,0){4}}
\put(9,-1){\line(1,0){4}}

\put(9,0){\line(0,-1){1}}
\put(10,0){\line(0,-1){1}}
\put(11,1){\line(0,-1){3}}
\put(12,0){\line(0,-1){1}}
\put(13,0){\line(0,-1){1}}

\put(11,1){\line(-2,-1){2}}
\put(11,1){\line(-1,-1){1}}
\put(11,1){\line(1,-1){1}}
\put(11,1){\line(2,-1){2}}

\put(11,-2){\line(-2,1){2}}
\put(11,-2){\line(-1,1){1}}
\put(11,-2){\line(1,1){1}}
\put(11,-2){\line(2,1){2}}

\put(11,-1){\oval(4,3)[b]}
\put(11,0){\oval(4,3)[t]}
\put(11,-1){\oval(6,4)[b]}

\put(9,0){\line(1,-1){1}}
\put(10,0){\line(1,-1){1}}
\put(11,0){\line(1,-1){1}}
\put(12,0){\line(1,-1){1}}
\put(13,0){\line(1,-1){1}}

\put(9,0){\circle*{0.1}}
\put(10,0){\circle*{0.1}}
\put(11,0){\circle*{0.1}}
\put(12,0){\circle*{0.1}}
\put(13,0){\circle*{0.1}}

\put(9,-1){\circle*{0.1}}
\put(10,-1){\circle*{0.1}}
\put(11,-1){\circle*{0.1}}
\put(12,-1){\circle*{0.1}}
\put(13,-1){\circle*{0.1}}

\put(11,1){\circle*{0.1}}
\put(11,-2){\circle*{0.1}}

\put(8,-1){\line(1,0){1}}

\end{picture}
\caption{A graph in $\mathcal{P}(n,5,5,0,6)$.} \label{P56}
\end{figure}

Clearly $\mathcal{P}(n,2,5,0,2) \supset C_{n}$,
so certainly $|\mathcal{P}(n,2,5,0,2)| >0$ $\forall n \geq 3$.
It is shown in~\cite{sza} that there exist $4$-regular planar graphs of order $n$ $\forall n \geq 6$
apart from $n=7$,
so we also have $|\mathcal{P}(n,4,5,0,4)|>0$ $\forall n \geq 8$.
For $n$ of the form~$25x+37y$,
graphs in $\mathcal{P}(n,5,5,0,6)$  
can be constructed from $5$-regular planar graphs as shown in Figure~\ref{P56},
so (since $25$ and $37$ are co-prime) 
it follows that we have $|\mathcal{P}(n,5,5,0,6)|>0$ for all sufficiently large $n$ too.
Hence, we have $g(n,d_{1},D_{2})>0$ for all large $n$ for all allowable $(d_{1},D_{2}(n))$. \\

We shall now show that $g$ satisfies the supermultiplicative condition:

Let $i,j \in \mathbf{N}$ and
let us denote by 
$\mathcal{P}_{c}(i,d_{1},5,0,D_{2})$ and $\mathcal{P}_{c}(j,d_{1},5,0,D_{2})$
the set of connected graphs in $\mathcal{P}(i,d_{1},5,0,D_{2})$ and $\mathcal{P}(j,d_{1},5,0,D_{2})$,
respectively.
Then, by Theorem~\ref{bounded311},
we know that there exists a constant $c>0$ such that we have
$|\mathcal{P}_{c}(i,d_{1},5,0,D_{2})| \geq c |\mathcal{P}(i,d_{1},5,0,D_{2})|$
and 
$|\mathcal{P}_{c}(j,d_{1},5,0,D_{2})| \geq c |\mathcal{P}(j,d_{1},5,0,D_{2})|$.
We may form a graph in $\mathcal{P}(i+j,d_{1},5,0,D_{2})$ by choosing $i$ of the $i+j$ vertices 
$\left( \left(^{i+j}_{\phantom{q}j} \right) \textrm{ choices} \right)$,
placing a connected planar graph $G_{1}$ with $|G_{1}|=i$ and $d_{1} \leq \delta(G_{1}) \leq \Delta(G_{1}) \leq D_{2}(i)$
on the chosen vertices
($|\mathcal{P}_{c}(i,d_{1},5,0,D_{2})| \geq c |\mathcal{P}(i,d_{1},5,0,D_{2})|$ choices),
and then placing a connected planar graph $G_{2}$ with $|G_{2}|=j$ and 
$d_{1} \leq \delta(G_{2}) \leq \Delta(G_{2}) \leq D_{2}(j)$
on the remaining $j$ vertices
($|\mathcal{P}_{c}(j,d_{1},5,0,D_{2})| \geq c|\mathcal{P}(j,d_{1},5,0,D_{2})|$ choices).
If $i=j$,
then we need to divide by two to avoid double-counting.
Note that the constructed graph will have maximum degree at most $\max \{ D_{2}(i), D_{2}(j) \}$
and so will indeed be in $\mathcal{P}(i+j,d_{1},5,0,D_{2})$,
since $D_{2}$ is a monotonically non-decreasing function.
Thus,
\begin{eqnarray*}
|\mathcal{P}(i+j,d_{1},5,0,D_{2})| \geq 
\frac{c^{2}}{2} \left(^{i+j}_{\phantom{q}j} \right) 
|\mathcal{P}(i,d_{1},5,0,D_{2})| \cdot |\mathcal{P}(j,d_{1},5,0,D_{2})| \textrm{ } \forall i,j
\end{eqnarray*}
and, therefore,
\begin{eqnarray*}
g(i+j,d_{1},D_{2}) & = & \frac{c^{2} |\mathcal{P}(i+j,d_{1},5,0,D_{2})|}{2(i+j)!} \\
& \geq & \frac{ c^{4} \left(^{i+j}_{\phantom{q}j} \right) 
|\mathcal{P}(i,d_{1},5,0,D_{2})| |\mathcal{P}(j,d_{1},5,0,D_{2})|}
{4(i+j)!} \\
& = & \frac{ c^{2} |\mathcal{P}(i,d_{1},5,0,D_{2})|}{2 \cdot i!} 
\frac{c^{2} |\mathcal{P}(j,d_{1},5,0,D_{2})|}{2 \cdot j!} \\
& = & g(i,d_{1},D_{2}) \cdot g(j,d_{1},D_{2}). \\
\end{eqnarray*} 

Let $\gamma_{d_{1},D_{2}} = \sup_{n} \left( (g(n,d_{1},D_{2}))^{1/n} \right)$.
By Proposition~\ref{vanL 11.6}, 
it now only remains to show that $\gamma_{d_{1},D_{2}} \!<\! \infty$.
But clearly $\mathcal{P}(n,d_{1},5,0,D_{2}) \!\subset\! \mathcal{P}(n,0,5,0,n-~1)$,
so $\gamma_{d_{1},D_{2}} \leq \gamma_{l} < \infty$.
\phantom{qwerty}
\setlength{\unitlength}{0.25cm}
\begin{picture}(1,1)
\put(0,0){\line(1,0){1}}
\put(0,0){\line(0,1){1}}
\put(1,1){\line(-1,0){1}}
\put(1,1){\line(0,-1){1}}
\end{picture} \\
\\
\\

We may obtain an analogous result to Theorem~\ref{bounded104}
for the case when we have $D_{2}(n) = d_{1} \in \{ 3,5 \}$~$\forall n$
by using Corollary~\ref{bounded103} instead of Proposition~\ref{vanL 11.6}:

\begin{Theorem} \label{bounded105}
Let $D_{2} \in \{ 3,5 \}$ be a fixed constant.
Then there is a finite constant $\gamma_{D_{2},D_{2}}>0$ such that
\begin{displaymath}
\left( \frac{|\mathcal{P}(2n,D_{2},5,0,D_{2})|}{(2n)!} \right)^{\frac{1}{2n}} \to 
\gamma_{D_{2},D_{2}} \textrm{ as } n \to \infty.
\end{displaymath}
\end{Theorem}
\textbf{Proof}
The method of proof is exactly the same as that of Theorem~\ref{bounded104}.
Thus, it suffices to show that $\mathcal{P}(n,3,5,0,3)$ and $\mathcal{P}(n,5,5,0,5)$ 
are non-empty for sufficiently large even $n$.
But it has already been shown in \cite{sza} that there exist 
$3$-regular planar graphs of order $n$ for all even $n \geq 4$
and that there exist $5\textrm{-regular planar graphs}$ of order $n$ for all even $n \geq 12$ apart from $n=14$.
Thus, we are done.
$\phantom{qwerty}$ 
\setlength{\unitlength}{0.25cm}
\begin{picture}(1,1)
\put(0,0){\line(1,0){1}}
\put(0,0){\line(0,1){1}}
\put(1,1){\line(-1,0){1}}
\put(1,1){\line(0,-1){1}}
\end{picture}

\newpage
\section{Appearances} \label{apps}

We shall now look at appearances in $P_{n,d_{1},5,0,D_{2}}$.
Recall from Proposition~\ref{gim T4}~that, 
for any connected planar graph $H$, 
$P_{n,0,5,0,n-1}$ will a.a.s.~have at least linearly many appearances of $H$.
In this section,
we will prove two similar appearance-type results, 
Theorems~\ref{bounded11} and~\ref{bounded110},
for $P_{n,d_{1},5,0,D_{2}}$
(these final results will not require the monotonicity/constancy conditions that we imposed 
in the previous section,
but such restrictions will be used in lemmas along the way).
In Section~\ref{cpts},
we shall then aim to turn some of these appearances into components.

We will produce separate appearance-type results for the cases when we have
$d_{1}(n) < D_{2}(n)$~$\forall n$ (Theorem~\ref{bounded11})
and when $d_{1}(n) = D_{2}(n)$ $\forall n$ (Theorem~\ref{bounded110}). 
This is essentially because for the latter case it will be awkward 
to convert appearances into components in Section~\ref{cpts}
without violating our bound on $d_{1}(n)$,
and so we instead introduce the concept of `$2$-appearances' (see Definition~\ref{2appdefn}).

We will deal with the $d_{1}(n) < D_{2}(n)$ case first,
starting with a simple lemma (Lemma~\ref{bounded201}) concerning the number of intersections of appearances.
The main work will then be done in Lemma~\ref{bounded106},
where we shall copy a proof of \cite{mcd}
(using the growth constants of Section~\ref{growth})
to obtain an appearance result for the case when 
$d_{1}(n)$ is a constant and $D_{2}(n)$ is a monotonically non-decreasing function.
We shall then extend this to a more general result in Lemma~\ref{bounded781},
before finally achieving the full result of Theorem~\ref{bounded11}.

The $d_{1}(n)=D_{2}(n)$ case will follow a similar pattern.
First, we shall note two lemmas (Lemmas~\ref{bounded502} and Corollary~\ref{bounded108})
on the number of intersections of $2$-appearances,
before then using our growth constant results to obtain a weak $2$-appearance result (Lemma~\ref{bounded109}),
which we shall then extend to Theorem~\ref{bounded110}.

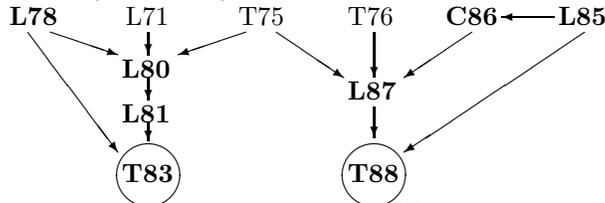
\begin{figure} [ht]
\setlength{\unitlength}{1cm}
\begin{picture}(20,1.5)(-1.8,2)

\put(2.1,1.6){\circle{0.8}}
\put(5.1,1.6){\circle{0.8}}

\put(0.5,3.5){\vector(3,-4){1.2}}
\put(0.8,3.5){\vector(3,-1){0.9}}
\put(2.1,3.5){\vector(0,-1){0.3}}
\put(3.4,3.5){\vector(-3,-1){0.9}}
\put(2.1,2.9){\vector(0,-1){0.3}}
\put(2.1,2.3){\vector(0,-1){0.25}}

\put(3.8,3.5){\vector(3,-2){0.9}}
\put(5.1,3.5){\vector(0,-1){0.6}}
\put(6.4,3.5){\vector(-3,-2){0.9}}
\put(5.1,2.5){\vector(0,-1){0.4}}
\put(7.9,3.5){\vector(-3,-2){2.4}}
\put(7.5,3.7){\vector(-1,0){0.7}}

\put(1.8,3.6){L\ref{bounded1}}
\put(3.3,3.6){T\ref{bounded104}}
\put(4.75,3.6){T\ref{bounded105}}

\put(0.25,3.6){\textbf{L\ref{bounded201}}}
\put(6.05,3.6){\textbf{C\ref{bounded108}}}
\put(7.55,3.6){\textbf{L\ref{bounded502}}}
\put(1.75,2.9){\textbf{L\ref{bounded106}}}
\put(4.75,2.6){\textbf{L\ref{bounded109}}}
\put(1.75,2.3){\textbf{L\ref{bounded781}}}
\put(4.75,1.5){\textbf{T\ref{bounded110}}}
\put(1.75,1.5){\textbf{T\ref{bounded11}}}

\end{picture}

\caption{The structure of Section~\ref{apps}.}
\end{figure}

We start with some definitions that we will find helpful:

\begin{Definition}
Suppose $H$ appears at $W$ in $G$ and let the unique edge in $G$ between $W$ and $V(G) \setminus W$
be $e_{W} = r_{W}v_{W}$, where $r_{W} \in W$ and $v_{W} \in V(G) \setminus~W$.
Let us call $e_{W}$ the \emph{\textbf{associated cut-edge}} of the appearance,
let us call $TV_{W} := W \cup \{ v_{W} \}$ the \emph{\textbf{total vertex set}} of the appearance
and let us call $TE_{W} := E(G[W]) \cup \{e_{W}\}$ the \emph{\textbf{total edge set}} of the appearance. 
\end{Definition}

\begin{figure} [ht]
\setlength{\unitlength}{1cm}
\begin{picture}(20,4.5)(-3.375,-0.75)

\put(0.75,0){\line(1,0){0.5}}
\put(0.75,1){\line(1,0){0.5}}
\put(0.75,2){\line(1,0){0.5}}
\put(0.75,3.5){\line(1,0){0.5}}
\put(4.25,0){\line(1,0){0.5}}
\put(4.25,1){\line(1,0){0.5}}

\put(1,1){\line(0,1){1}}
\put(4.5,1){\line(0,1){1}}

\put(0.75,0.5){\oval(0.5,1)[l]}
\put(0.75,2.75){\oval(2,1.5)[l]}
\put(4.25,0.5){\oval(0.5,1)[l]}
\put(1.25,0.5){\oval(0.5,1)[r]}
\put(1.25,2.75){\oval(2,1.5)[r]}
\put(4.75,0.5){\oval(0.5,1)[r]}

\put(1,1){\circle*{0.1}}
\put(1,2){\circle*{0.1}}
\put(4.5,1){\circle*{0.1}}
\put(4.5,2){\circle*{0.1}}

\put(0.6,2.7){$G \setminus W$}
\put(0.65,0.4){\small{$G[W]$}}
\put(4.15,0.4){\small{$G[W]$}}

\put(0.9,-0.75){\large{$G$}}
\put(3.8,-0.75){\large{$G[TV_{W}]$}}

\put(0.5,1.4){$e_{W}$}
\put(1.1,1.8){$v_{W}$}
\put(1.1,1.1){$r_{W}$}
\put(4,1.4){$e_{W}$}
\put(4.6,1.8){$v_{W}$}
\put(4.6,1.1){$r_{W}$}

\end{picture}

\caption{$\textrm{An appearance at $W$ in $G$ and its total vertex/edge set}$.}
\end{figure}
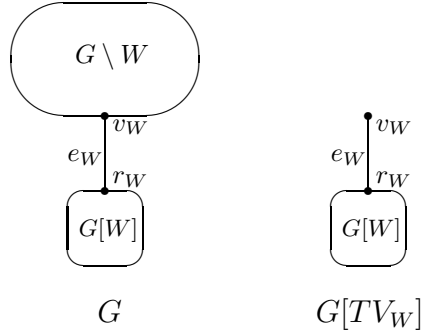

\phantom{q}

It is now easier to state the following useful result,
which is given implicitly in the proof of Theorem 4.1 of \cite{mcd}:

\begin{Lemma} \label{bounded201}
The total edge set of an appearance of a graph of order $|H|$ will intersect
(i.e. have an edge in common with)
the total edge set of at most $|H|$ other appearances of graphs of order $|H|$.
\end{Lemma}
\textbf{Proof}
Suppose we have an appearance of a graph of order $|H|$ at $W \subset V(G)$.
Let the associated cut-edge of the appearance
be $e_{W} = r_{W}v_{W}$,
where $r_{W} \in W$ and $v_{W} \in V(G) \setminus W$.
Suppose that $G$ also contains another appearance, at $W_{2}$,
of a graph of order $|H|$ 
and let the associated cut-edge of this appearance be $e_{W_{2}} = r_{W_{2}}v_{W_{2}}$,
where $r_{W_{2}} \in W_{2}$ and $v_{W_{2}} \in V(G) \setminus W_{2}$.

(a) Suppose $e_{W_{2}} \notin TE_{W}$.
Then $G[TV_{W}]$ is a connected subgraph in $G \setminus e_{W_{2}}$.
Note that $G[W_{2}]$ is a component in $G \setminus e_{W_{2}}$.
Thus, either $G[TV_{W}] \subset G[W_{2}]$,
in which case we would obtain the contradiction
$|W|+1 = |TV_{W}| \leq |W_{2}| = |W|$,
or $G[TV_{W}] \subset G \setminus W_{2}$,
in which case the total edge set of $W$ would not meet the total edge set of $W_{2}$.

(b) Suppose $e_{W_{2}} = e_{W}$.
Then, since $W_{2} \neq W$,
we must have $r_{W_{2}} = v_{W}$ and~$v_{W_{2}} = r_{W}$,
and so $W_{2}$ must be the component of $G \setminus e_{W}$ that contains $v_{W}$.

(c) Suppose $e_{W_{2}} \in E(G[W])$.
Since $e_{W_{2}}$ must be a cut-edge,
there are at most $|W|-1 = |H|-1$ possibilities for it.
Note that the `orientation' of $e_{W_{2}}$ must be chosen so that $r_{W} \in W_{2}$,
since otherwise we would have $W_{2} \subset W \setminus r_{W}$
and obtain the contradiction $|W \setminus r_{W}| \geq |W_{2}| = |W|$.
Hence, since knowing $e_{W_{2}}$ and its orientation determines $W_{2}$,
there can only be at most $|H|-1$ possibilities for $W_{2}$.

Thus, in total, there are at most 
$0+1+(|H|-1)=|H|$
possibilities for~$W_{2}$.~
\setlength{\unitlength}{0.25cm}
\begin{picture}(1,1)
\put(0,0){\line(1,0){1}}
\put(0,0){\line(0,1){1}}
\put(1,1){\line(-1,0){1}}
\put(1,1){\line(0,-1){1}}
\end{picture} \\
\\
\\

We shall shortly proceed with an appearance result for the case when we have
$d_{1}(n) < D_{2}(n)$~$\forall n$.
As already mentioned, this will later be used (in Section~\ref{cpts}) 
in proofs concerning the construction of components isomorphic to given $H$.
Hence, with this in mind, we will make one more definition
(which will help us later to delete certain edges without breaking our lower bound on the minimum degree):

\begin{Definition}
Given a connected planar graph $H$,
a function $d_{1}(n)$,
and a planar graph $G$,
let \emph{\boldmath{$f_{H}^{d_{1}}(G)$}} denote the number of appearances of $H$ in $G$ 
such that the associated cut-edge is between
two vertices with $\deg_{G} > d_{1}(|G|)$. 
\end{Definition}

\begin{displaymath}
\end{displaymath}

We are now finally ready to look at appearances in $P_{n,d_{1},5,0,D_{2}}$.
We shall start by assuming that $d_{1}(n)$ is constant and $D_{2}(n)$ is monotonically non-decreasing,
as in Section~\ref{growth},
but we will later (in Theorem~\ref{bounded781}) get rid of these conditions.
The statement of the result may seem complicated,
but basically it just asserts that for any `sensible' choice of $H$,
there will probably be lots of appearances of $H$ in $P_{n,d_{1},5,0,D_{2}}$
such that the associated cut-edge is between two vertices with degree $>d_{1}$.
Clearly, `sensible' entails that we must have $\delta(H) \geq d_{1}$,
$\Delta(H) \leq~\!D_{2}(n)$ and $\deg_{H}(1)+1 \leq D_{2}(n)$,
and as always we will require that $D_{2}(n) \geq 3$
(note that it follows from these conditions that we also have~$d_{1}<~\!D_{2}(n)$).
The proof is based on that of Theorem 4.1 of \cite{mcd}.

\begin{Lemma} \label{bounded106}
Let $H$ be a fixed connected planar graph on $\{ 1,2, \ldots, h \}$.
Then $\exists \alpha(h)>0$ such that,
given any constant $d_{1} \leq \delta(H)$
and any monotonically non-decreasing integer-valued function $D_{2}(n)$ satisfying
$\liminf_{n \to \infty}D_{2}(n) \geq \max \{ \Delta(H), \deg_{H}(1)+1, 3 \}$,
we have
\begin{displaymath}
\mathbf{P} [ f_{H}^{d_{1}} (P_{n,d_{1},5,0,D_{2}}) \leq \alpha n ] < e^{-\alpha n} \textrm{ for all sufficiently large }n.
\end{displaymath}
\end{Lemma}
\textbf{Sketch of Proof}
We choose a specific $\alpha$
and suppose that the result is false for~$n=k$,
where $k$ is suitably large.
Using Theorem~\ref{bounded104},
it then follows that there are many graphs $G \in \mathcal{P}(k,d_{1},5,0,D_{2})$ with $f_{H}^{d_{1}}(G) \leq \alpha k$.

From each such $G$,
we construct graphs in $\mathcal{P}((1+\delta)k,d_{1},5,0,D_{2})$, for a fixed~$\delta >0$.
If $G$ has lots of vertices with degree $<D_{2}(k)$,
then we do this simply by attaching appearances of $H$ to some of these vertices.
If $G$ has few vertices with degree~$<D_{2}(k)$,
then we attach appearances of $H$ to small cycles in $G$ 
and also delete appropriate edges.
By Lemma~\ref{bounded1},
we have lots of choices for these small cycles
and, since $G$ has few vertices with degree $<D_{2}(k)$,
we may assume that we don't interfere with any vertices of minimum degree.

The fact that the original graphs satisfied $f_{H}^{d_{1}} \leq \alpha k$,
together with Lemma~\ref{bounded201} and the knowledge that any deleted edges were in small cycles,
is then used to show that there is not much double-counting,
and so we find that we have constructed so many graphs in $\mathcal{P}((1+\delta)k,d_{1},5,0,D_{2})$
that we contradict Theorem~\ref{bounded104}. \\
\\
\textbf{Full Proof}
Let $p < \frac{1}{7(6^{2}+6^{3}+6^{4}+6^{5})}$,
let $\beta = \frac{344e^{2}(h+7)(6^{2}+6^{3}+6^{4}+6^{5})h! \left( \gamma_{l} \right)^{h}}{p}$,
and let $\alpha$ be a fixed constant in 
$\left( 0, \frac{1}{\beta} \right)$.
Then we have $\alpha \beta <1$ 
and so $\exists \epsilon \in \left( 0,\frac{1}{3} \right)$ such that $(\alpha \beta)^{\alpha} = 1-3 \epsilon$.

By Theorem~\ref{bounded104},
$\exists N$ such that
\begin{equation} \label{eq:gr}
(1-\epsilon)^{n}n! \left( \gamma_{d_{1},D_{2}} \right)^{n}
\leq |\mathcal{P}(n,d_{1},5,0,D_{2})|
\leq (1+\epsilon)^{n}n! \left( \gamma_{d_{1},D_{2}} \right)^{n} 
\textrm{ } \forall n \geq N.
\end{equation}
Suppose (aiming for a contradiction) 
that we can find a value
$k>N$ such that
$\mathbf{P} [ f_{H}^{d_{1}} (P_{k,d_{1},5,0,D_{2}}) \leq \alpha k ] \geq e^{-\alpha k}$,
and let $\mathcal{G}$ denote the set of graphs in $\mathcal{P}(k,d_{1},5,0,D_{2})$ such that
$G \in \mathcal{G}$ iff $f_{H}^{d_{1}} (G) \leq \alpha k$.
Then we must have 
$|\mathcal{G}| \geq e^{-\alpha k} |\mathcal{P}(k,d_{1},5,0,D_{2})| \geq 
e^{-\alpha k} (1-\epsilon)^{k}k! \left( \gamma_{d_{1},D_{2}} \right)^{k}$. 

Let $\delta = \frac{\lceil \alpha k \rceil h}{k}$.
We may assume that $k$ is sufficiently large that 
$\lceil \alpha k \rceil \leq 2 \alpha k$.
Thus, $\delta \leq 2 \alpha h < 1$ \label{delta}
(by our definition of $\alpha$).
This fact will be useful later. \\

We shall construct graphs in $\mathcal{P}((1+\delta)k,d_{1},5,0,D_{2})$:

Choose $\delta k$ special vertices 
(we have $\left(^{(1+\delta)k} _{\phantom{qq} \delta k} \right)$ choices for these)
and partition them into $\lceil \alpha k \rceil$ unordered blocks of size $h$
(we have $\left(^{\phantom{w} \delta k} _{h, \ldots, h} \right) \frac{1}{\lceil \alpha k \rceil !}$ choices for this).
On each of the blocks,
put a copy of $H$ such that the increasing bijection from $\{ 1,2, \ldots, h \}$ to the block
is an isomorphism between $H$ and this copy.
Note that we may assume that $k$ is large enough that $D_{2}(k) \geq \liminf_{n \to \infty}D_{2}(n)$,
and so the root, $r_{B}$, of a block
(i.e.~the lowest numbered vertex in it)
satisfies deg$(r_{B})<D_{2}(k)$,
by the conditions of the theorem. 
On the remaining $k$ vertices,
we place a planar graph $G$ with 
$d_{1} \leq \delta(G) \leq \Delta(G) \leq D_{2}(k)$
and $f_{H}^{d_{1}}(G) \leq \alpha k$
(we have at least $|\mathcal{G}|$ choices for this).

We shall continue our construction in one of two ways,
depending on the number of vertices of degree $D_{2}(k)$ in $G$: \\
\\
\\
Case (a): If $G$ has $\geq \frac{pk}{43}$ vertices of degree $<D_{2}(k)$
(note that this is certainly the case if $D_{2}(k) \geq 7$). \\
For each block $B$,
we choose a \textit{different} non-special vertex $v_{B} \in V(G)$ with deg$(v_{B})<D_{2}(k)$ 
(we have 
$\geq \left( ^{pk/43} _{\phantom{i} \lceil \alpha k \rceil} \right) \lceil \alpha k \rceil !$
choices for this,
since certainly $\alpha < \frac{p}{86}$
and we may assume that $k$ is large enough that 
$\lceil \alpha k \rceil \leq 2 \alpha k$),
and we insert the edge $r_{B}v_{B}$ from the root of the block to this vertex, 
creating an appearance of $H$ at $B$
(we should note for later use that $r_{B}$ and $v_{B}$ will now clearly both have degree $>d_{1}$).
\begin{figure} [ht]
\setlength{\unitlength}{1cm}
\begin{picture}(20,2.25)(-0.5,-0.25)

\put(0,0.5){\line(0,1){0.5}}
\put(1,0.5){\line(0,1){0.5}}
\put(2,0.5){\line(0,1){0.5}}
\put(4,0.5){\line(0,1){0.5}}
\put(7,0.5){\line(0,1){0.5}}
\put(8,0.5){\line(0,1){0.5}}
\put(9,0.5){\line(0,1){0.5}}
\put(11,0.5){\line(0,1){0.5}}

\put(8,0.75){\line(1,0){1}}

\put(5,0.75){\vector(1,0){1}}

\put(0.5,1){\oval(1,0.5)[t]}
\put(0.5,0.5){\oval(1,0.5)[b]}
\put(3,1){\oval(2,2)[t]}
\put(3,0.5){\oval(2,2)[b]}
\put(7.5,1){\oval(1,0.5)[t]}
\put(7.5,0.5){\oval(1,0.5)[b]}
\put(10,1){\oval(2,2)[t]}
\put(10,0.5){\oval(2,2)[b]}

\put(1,0.75){\circle*{0.1}}
\put(2,0.75){\circle*{0.1}}
\put(8,0.75){\circle*{0.1}}
\put(9,0.75){\circle*{0.1}}

\put(0.4,0.6){$B$}
\put(2.8,0.6){\Large{$G$}}
\put(1.05,0.55){\footnotesize{$r_{B}$}}
\put(1.6,0.55){\footnotesize{$v_{B}$}}

\put(7.4,0.6){$B$}
\put(9.8,0.6){\Large{$G$}}
\put(8.05,0.55){\footnotesize{$r_{B}$}}
\put(8.6,0.55){\footnotesize{$v_{B}$}}

\end{picture}

\caption{Creating an appearance of $H$ at $B$ in case (a).}
\end{figure}
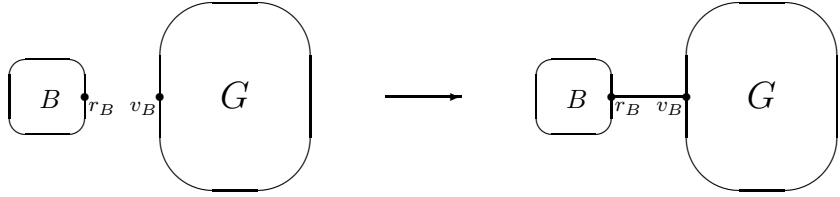
Note that we have not deleted any edges,
so we shall still have minimum degree at least $d_{1}$,
and we have only inserted edges between vertices of degree $<D_{2}(k)$,
so we still have maximum degree at most $D_{2}(k)$,
which is at most $D_{2}((1+ \delta)k)$ by monotonicity of $D_{2}$.
Thus, our new graph is indeed in $\mathcal{P}((1+\delta)k,d_{1},5,0,D_{2})$.

Hence, for each graph $G$ with $\geq \frac{pk}{43}$ vertices of degree $<D_{2}(k)$,
we find that we can construct at least
$
\left(^{(1+\delta)k} _{\phantom{qq} \delta k} \right)
\left(^{\phantom{p}\delta k}_{h \ldots h} \right) 
\cdot \frac{1}{\lceil \alpha k \rceil !} \cdot \left( ^{pk/43}_{\phantom{i} \lceil \alpha k \rceil} \right)
\lceil \alpha k \rceil !
$
different graphs in $\mathcal{P}((1+\delta)k,d_{1},5,0,D_{2})$. \\
\\
Case (b): If $G$ has $<\frac{pk}{43}$ vertices of degree $<D_{2}(k)$
(in which case $D_{2}(k)<7$). \\
Before describing the case (b) continuation of our construction,
it shall first be useful to investigate the number of short cycles in $G$: \\
If $G$ has $<\frac{pk}{43} < \frac{k}{43}$ vertices of degree $<D_{2}(k)$,
then (by Lemma~\ref{bounded1}) $G$ contains at least $\frac{k}{43}$ cycles of size at most $6$.
A vertex can only be in at most 
$(D_{2}(k))^{2}+(D_{2}(k))^{3}+(D_{2}(k))^{4}+(D_{2}(k))^{5} \leq 6^{2}+6^{3}+6^{4}+6^{5}$ cycles of size at most $6$,
so $G$ must have at most $\frac{pk(6^{2}+6^{3}+6^{4}+6^{5})}{43}$ cycles of size at most $6$
that contain a vertex of degree $<D_{2}(k)$.
In particular, $G$ must have at least $\frac{(1-(6^{2}+6^{3}+6^{4}+6^{5})p)k}{43}$ cycles of size at most $6$
that don't contain a vertex of degree~$d_{1}$, since $d_{1} \leq \delta(H) < \deg_{H}(1)+1 \leq D_{2}(k)$.
Since a vertex can only be in at most 
$6^{2}+6^{3}+6^{4}+6^{5}$ cycles of size at most $6$,
each cycle of size at most $6$ can only have a vertex in common with
at most $6(6^{2}+6^{3}+6^{4}+6^{5})$ other cycles of size at most $6$.
Thus, $G$ must have a set of at least 
$\frac{\left( \frac{1-(6^{2}+6^{3}+6^{4}+6^{5})p}{6(6^{2}+6^{3}+6^{4}+6^{5})} \right)k}{43} >~\frac{pk}{43}$ 
\textit{vertex-disjoint} cycles of size at most $6$ that don't contain a vertex of degree $d_{1}$
$\left( \textrm{using the fact that } 
p < \frac{1}{7(6^{2}+6^{3}+6^{4}+6^{5})} \right)$.
We shall call these cycles `special'.

Recall that we have $\lceil \alpha k \rceil$ blocks isomorphic to $H$.
For each block $B$,
choose a \textit{different} one of our `special' cycles
(we have $\geq \left( ^{pk/43}_{\phantom{i}\lceil \alpha k \rceil} \right) \lceil \alpha k \rceil !$
choices for this),
delete an edge $u_{B}v_{B}$ in the cycle 
and insert an edge $r_{B}v_{B}$ from the root of the block to a vertex $v_{B}$ that was incident to the deleted edge,
creating an appearance of $H$ at $B$.
\begin{figure} [ht]
\setlength{\unitlength}{1cm}
\begin{picture}(20,2.25)(-0.5,-0.25)

\put(0,0.5){\line(0,1){0.5}}
\put(1,0.5){\line(0,1){0.5}}
\put(2,0.5){\line(0,1){0.5}}
\put(4,0.5){\line(0,1){0.5}}
\put(7,0.5){\line(0,1){0.5}}
\put(8,0.5){\line(0,1){0.5}}
\put(11,0.5){\line(0,1){0.5}}

\put(8,0.75){\line(4,-1){1}}

\put(5,0.75){\vector(1,0){1}}

\put(0.5,1){\oval(1,0.5)[t]}
\put(0.5,0.5){\oval(1,0.5)[b]}
\put(3,1){\oval(2,2)[t]}
\put(3,0.5){\oval(2,2)[b]}
\put(7.5,1){\oval(1,0.5)[t]}
\put(7.5,0.5){\oval(1,0.5)[b]}
\put(10,1){\oval(2,2)[t]}
\put(10,0.5){\oval(2,2)[b]}

\put(2,0.75){\oval(1,0.5)[r]}
\put(9,0.75){\oval(1,0.5)[r]}

\put(1,0.75){\circle*{0.1}}
\put(2,0.5){\circle*{0.1}}
\put(2,1){\circle*{0.1}}
\put(8,0.75){\circle*{0.1}}
\put(9,0.5){\circle*{0.1}}
\put(9,1){\circle*{0.1}}

\put(0.4,0.6){$B$}
\put(3,0.6){\Large{$G$}}
\put(1.05,0.5){\footnotesize{$r_{B}$}}
\put(1.6,0.3){\footnotesize{$v_{B}$}}
\put(1.6,1.15){\footnotesize{$u_{B}$}}

\put(7.4,0.6){$B$}
\put(10,0.6){\Large{$G$}}
\put(8.05,0.5){\footnotesize{$r_{B}$}}
\put(8.6,0.35){\footnotesize{$v_{B}$}}
\put(8.6,1.15){\footnotesize{$u_{B}$}}
\end{picture}

\caption{Creating an appearance of $H$ at $B$ in case (b).}
\end{figure}
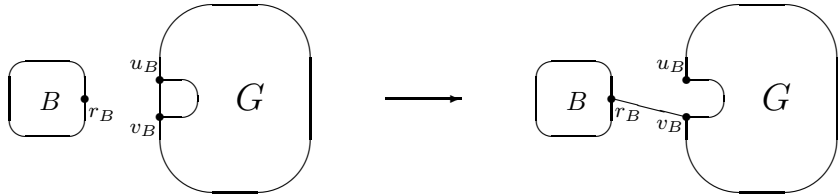
Note that the deleted edge was between two vertices of degree~$>d_{1}$,
so we still have minimum degree at least $d_{1}$
(we should also note for later use that $v_{B}$ will still have degree $>d_{1}$,
and that $r_{B}$ will now also have degree $>d_{1}$).
Recall that the root of each block has degree $<D_{2}(k)$,
so we still have maximum degree at most $D_{2}(k)$,
which is at most $D_{2}((1+~\delta)k)$ by monotonicity of $D_{2}$.
Thus, our constructed graph is indeed in $\mathcal{P}((1+\delta)k,d_{1},5,0,D_{2})$.

Thus, for each graph $G$ with $<\frac{pk}{43}$ vertices of degree $<D_{2}(k)$,
we find that we can construct at least
$
\left(^{(1+\delta)k} _{\phantom{qq} \delta k} \right)
\left(^{\phantom{p}\delta k}_{h \ldots h} \right) 
\cdot \frac{1}{\lceil \alpha k \rceil !} \cdot \left( ^{pk/43}_{\phantom{i} \lceil \alpha k \rceil} \right)
\lceil \alpha k \rceil !
$
different graphs in $\mathcal{P}((1+\delta)k,d_{1},5,0,D_{2})$. \\

We have shown that,
regardless of whether case (a) or case (b) is used,
for each $G$ we can construct at least 
\begin{eqnarray*}
& & \left(^{(1+\delta)k} _{\phantom{qq} \delta k} \right)
\left(^{\phantom{p}\delta k}_{h \ldots h} \right) 
\cdot \frac{1}{\lceil \alpha k \rceil !} \cdot \left( ^{pk/43}_{\phantom{i} \lceil \alpha k \rceil} \right)
\lceil \alpha k \rceil ! \\
& = & \frac{((1+\delta)k)!}{k!} \frac{1}{(h!)^{\lceil \alpha k \rceil}} 
\left(^{pk/43}_{\phantom{i} \lceil \alpha k \rceil} \right) \\
& \geq & \frac{((1+\delta)k)!}{k!} \frac{1}{(h!)^{\lceil \alpha k \rceil}}
\frac{(pk/43 - \lceil \alpha k \rceil +1)^{\lceil \alpha k \rceil}}{\lceil \alpha k \rceil !} \\
& \geq & \frac{((1+\delta)k)!}{k!} \frac{1}{(h!)^{\lceil \alpha k \rceil}} 
\left( \frac{pk}{86} \right)^{\lceil \alpha k \rceil} \frac{1}{\lceil \alpha k \rceil !} \\
& & \left( \textrm{ since } \alpha < \frac{p}{86}
\textrm{ and so } \frac{pk}{43} - \alpha k \geq \frac{pk}{86} \right) \\
& \geq & \frac{((1+\delta)k)!}{k!} \left( \frac{pk}{86h! \lceil \alpha k \rceil} \right)^{\lceil \alpha k \rceil} \\
& \geq & \frac{((1+\delta)k)!}{k!} \left( \frac{p}{172h! \alpha} \right)^{\lceil \alpha k \rceil} \\
& & \left( \textrm{since we may assume $k$ is large enough that } 
\lceil \alpha k \rceil \leq 2 \alpha k \right)
\end{eqnarray*}
different graphs in $\mathcal{P}((1+\delta)k,d_{1},5,0,D_{2})$.
Thus, recalling that we have at least
$e^{-\alpha k} (1-\epsilon)^{k}k! \left( \gamma_{d_{1},D_{2}} \right)^{k}$
choices for $G$,
we can in total construct at least
$
e^{-\alpha k} (1-\epsilon)^{k} ((1+\delta)k)! \left( \gamma_{d_{1},D_{2}} \right)^{k}
\left( \frac{p}{172h! \alpha} \right)^{\lceil \alpha k \rceil}
$
(not necessarily distinct) graphs in $\mathcal{P}((1+\delta)k,d_{1},5,0,D_{2})$. \\
\\

We are now at the half way point of our proof,
and it remains to investigate the amount of double-counting,
i.e.~how many times each of our constructed graphs will have been built. \\
\\

Given one of our constructed graphs, $G^{\prime}$,
there are at most $2$ possibilities for how the graph was obtained (case (a) or case(b)). \\

If case (a) was used, then we can re-obtain the original graph, $G$,
simply by deleting the $\lceil \alpha k \rceil$ appearances that were deliberately added.
Recall that these appearances were constructed in such a way that the associated cut-edges
are all between vertices with $\deg_{G^{\prime}}>d_{1}$.
Thus, in order to bound the amount of double-counting,
we only need to investigate $f_{H}^{d_{1}}(G^{\prime})$:

Suppose $W$ is an appearance of $H$ in $G^{\prime}$ such that the associated cut-edge $e_{W}$
is between two vertices of degree $>d_{1}$.
We shall consider how many possibilities there are for $W$:

(i) If we don't have $TE_{W} \subset E(G)$,
then the total edge set of $W$ must intersect the total edge set of one of our deliberately created appearances,
and so we have at most $(h+1) \lceil \alpha k \rceil$ possibilities for $W$
(including the possibility that $W$ is one of our deliberately created appearances),
by Lemma~\ref{bounded201}.

If $TE_{W} \subset E(G)$,
then $W$ must have been an appearance of $H$ in $G$:

(ii) If $W$ was an appearance of $H$ in $G$ such that
$e_{W}$ was already between two vertices of degree $>d_{1}$,
then there are at most $\lceil \alpha k \rceil$ possibilities for $W$,
by definition of $\mathcal{G}$.

(iii) If $W$ was an appearance of $H$ in $G$ such that 
$e_{W}$ was \textit{not} already between two vertices of degree $>d_{1}$,
then the unique vertex $v \in V(G^{\prime}) \setminus W$
incident to the root of $W$ must have had deg$(v) = d_{1}$ originally
and must have been chosen as $v_{B}$ by some block $B$.
Hence, we have at most $\lceil \alpha k \rceil$ possibilities for $v$
and thus at most $d_{1} \lceil \alpha k \rceil$ possibilities for $W$.

Thus, if case (a) was used, 
then $f_{H}^{d_{1}} (G^{\prime}) \leq (h+d_{1}+2) \lceil \alpha k \rceil$,
and so we have at most 
$\left( ^{(h+d_{1}+2) \lceil \alpha k \rceil}_{\phantom{www} \lceil \alpha k \rceil} \right)
\leq ((h+d_{1}+2)e)^{\lceil \alpha k \rceil}
\leq ((h+7)e)^{\lceil \alpha k \rceil}$
possibilities for~$G$. \\

If case (b) was used,
we can re-obtain the original graph, $G$,
by deleting the $\lceil \alpha k \rceil$ appearances that were deliberately added
and re-inserting the $\lceil \alpha k \rceil$ deleted edges.
Note that once we have identified the appearances that were deliberately added,
we have at most 
$\left( (D_{2}(k))^{2}+(D_{2}(k))^{3}+(D_{2}(k))^{4}+(D_{2}(k))^{5} \right)^{\lceil \alpha k \rceil}
\leq (6^{2}+ 6^{3}+ 6^{4}+ 6^{5})^{\lceil \alpha k \rceil}$
possibilities for the edges that were deleted,
since for each appearance we will automatically know one endpoint, $v$, of the corresponding deleted edge
and we know that the other endpoint, $u$,
will now be at most distance $5$ from $v$,
since $uv$ was originally part of a cycle of size $\leq 6$.
Hence, as with case (a), 
it now remains to examine how many possibilities there are for the $\lceil \alpha k \rceil$ appearances
that were deliberately added.

Suppose $W$ is an appearance of $H$ in $G^{\prime}$ such that the associated cut-edge~$e_{W}$
is between two vertices of degree $>d_{1}$.

(i) If we don't have $TE_{W} \subset E(G)$,
then we have at most $(h+1) \lceil \alpha k \rceil$ possibilities for $W$,
as with case (a).

(ii) If $TE_{W} \subset E(G)$ 
and $W$ was an appearance of $H$ in $G$,
then note that~$e_{W}$ must have already been between two vertices of degree $>d_{1}$,
since it is clear that we have $\deg_{G^{\prime}} \leq \deg_{G}$
for all vertices that were in $V(G)$.
Hence, there are at most $\lceil \alpha k \rceil$ possibilities for $W$,
by definition of $\mathcal{G}$.

(iii) If $TE_{W} \subset E(G)$ 
and $W$ was \textit{not} an appearance of $H$ in $G$,
then there must have originally been either another edge between $W$ and $V(G) \setminus W$ other than $e_{W}$,
or another edge between vertices in $W$.
This deleted edge must be of the form $u_{B}v_{B}$ for some block $B$,
and so $W$ must contain either $u_{B}$ or $v_{B}$ (or both).
However, if $v_{B} \in W$ then $r_{B}v_{B}$ would belong to the total edge set of~$W$,
which would contradict our assumption that $TE_{W} \subset E(G)$.
Thus, $u_{B} \in W$ and~$v_{B} \notin W$.
\begin{figure} [ht]
\setlength{\unitlength}{1cm}
\begin{picture}(20,3.5)(-5,0)

\put(0.75,0){\line(1,0){0.5}}
\put(0.75,1){\line(1,0){0.5}}
\put(0.75,2){\line(1,0){0.5}}
\put(0.75,3.5){\line(1,0){0.5}}

\put(1,1){\line(0,1){1}}

\put(0.75,0.5){\oval(0.5,1)[l]}
\put(0.75,2.75){\oval(2,1.5)[l]}
\put(1.25,0.5){\oval(0.5,1)[r]}
\put(1.25,2.75){\oval(2,1.5)[r]}

\put(1.5,0.5){\circle*{0.1}}
\put(1.2,2){\circle*{0.1}}

\put(0.8,0.4){$W$}

\put(0.5,1.4){$e_{W}$}
\put(1.3,1.8){$v_{B}$}
\put(1.6,0.4){$u_{B}$}

\end{picture}

\caption{The appearance $W$ in case (iii).}
\end{figure}
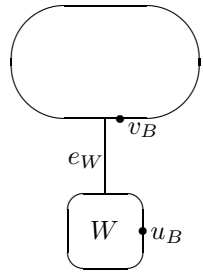
Recall that the deleted edge $u_{B}v_{B}$ was originally part of a cycle of size~$\leq 6$
and that no other edges from the same cycle were deleted.
Thus, there is still a $u_{B}-v_{B}$ path in $G^{\prime}$ consisting of the other edges in the cycle.
But, since $u_{B} \in W$, $v_{B} \notin W$ and $e_{W}$ is the unique edge between $W$ and $G^{\prime} \setminus W$,
it must be that $e_{W}$ belongs to this path,
i.e.~$e_{W}$ must have been one of the other (at most $5$) edges in the cycle.
Thus, we have at most $5 \lceil \alpha k \rceil$ possibilities for~$e_{W}$,
and hence for $W$
(since $e_{W}$ must be `oriented' so that $v_{B} \notin W$).

Thus, if case (b) was used, 
then $f_{H}^{d_{1}} (G^{\prime}) \leq (h+7) \lceil \alpha k \rceil$,
and so we have at most 
$\left( ^{(h+7) \lceil \alpha k \rceil}_{\phantom{ww} \lceil \alpha k \rceil} \right)
(6^{2}+6^{3}+6^{4}+6^{5})^{\lceil \alpha k \rceil}
\leq ((h+7)(6^{2}+~6^{3}+~6^{4}+~6^{5})e)^{\lceil \alpha k \rceil}$
possibilities for $G$. \\

We have shown each graph in $\mathcal{P}((1+\delta)k,d_{1},5,0,D_{2})$ is constructed at most
$((h+7)e)^{\lceil \alpha k \rceil} + ((h+7)(6^{2}+6^{3}+6^{4}+6^{5})e)^{\lceil \alpha k \rceil}
\leq 2 ((h+7)(6^{2}+6^{3}+6^{4}+6^{5})e)^{\lceil \alpha k \rceil}
\leq y^{\lceil \alpha k \rceil}$ times,
where $y$ denotes $2e(h+d_{1}+2)(6^{2}+6^{3}+6^{4}+6^{5})$.
Thus, the number of \textit{distinct} graphs 
that we have constructed in $\mathcal{P}((1+\delta)k,d_{1},5,0,D_{2})$ must be
\begin{eqnarray*}
& \geq & e^{-\alpha k} (1-\epsilon)^{k} ((1+\delta)k)! \left( \gamma_{d_{1},D_{2}} \right)^{k}
\left( \frac{p}{172h! \alpha y} \right)^{\lceil \alpha k \rceil} \\
& \geq & e^{-\alpha k} (1-\epsilon)^{k} ((1+\delta)k)! \left( \gamma_{d_{1},D_{2}} \right)^{(1+\delta)k}
\left( \gamma_{d_{1},D_{2}} \right)^{- \lceil \alpha k \rceil h}
\left( \frac{p}{172h! \alpha y} \right)^{\lceil \alpha k \rceil}, \\
& & \textrm{ since } \delta k = \lceil \alpha k \rceil h \\
& \geq & (1-\epsilon)^{k} ((1+\delta)k)! \left( \gamma_{d_{1},D_{2}} \right)^{(1+\delta)k}
\left( \frac{172 e h! \alpha y \left( \gamma_{d_{1},D_{2}} \right)^{h}}{p} \right)^{- \lceil \alpha k \rceil}, \\
& & \textrm{ since } e^{-\alpha k} \geq e^{- \lceil \alpha k \rceil} \\
& \geq & (1-\epsilon)^{k} ((1+\delta)k)! \left( \gamma_{d_{1},D_{2}} \right)^{(1+\delta)k}
(\alpha \beta)^{- \lceil \alpha k \rceil} \\
& \geq & (1-\epsilon)^{k} ((1+\delta)k)! \left( \gamma_{d_{1},D_{2}} \right)^{(1+\delta)k}
(\alpha \beta)^{- \alpha k} 
\textrm{, \phantom{www} since } \alpha \beta < 1 \\
& = & \left( \frac{1-\epsilon}{1-3\epsilon} \right)^{k} ((1+\delta)k)! \left( \gamma_{d_{1},D_{2}} \right)^{(1+\delta)k} 
\textrm{, \phantom{www} since } (\alpha \beta)^{\alpha} = 1-3 \epsilon \\
& \geq & \left( \frac{1-\epsilon}{1-3\epsilon} \right)^{k} 
\frac{|\mathcal{P}((1+\delta)k,d_{1},5,0,D_{2})|}{(1+\epsilon)^{(1+\delta)k}} 
\textrm{, \phantom{www} by (\ref{eq:gr})} \\
& > & \left( \frac{1-\epsilon}{(1-3\epsilon)(1+\epsilon)^{2}} \right)^{k} 
|\mathcal{P}((1+\delta)k,d_{1},5,0,D_{2})|, \phantom{w} \textrm{ since } \delta < 1 \textrm{ (page~\pageref{delta})} \\
& > & |\mathcal{P}((1+\delta)k,d_{1},5,0,D_{2})|, \phantom{www} \textrm{ since }
(1-3\epsilon)(1+\epsilon)^{2} = 1 - \epsilon - 5 \epsilon^{2} - 3 \epsilon^{3}.
\end{eqnarray*}
Hence, we have our desired contradiction.
$\phantom{qwerty}$ 
\setlength{\unitlength}{0.25cm}
\begin{picture}(1,1)
\put(0,0){\line(1,0){1}}
\put(0,0){\line(0,1){1}}
\put(1,1){\line(-1,0){1}}
\put(1,1){\line(0,-1){1}}
\end{picture}

As mentioned, we shall now see that we can actually drop the conditions that $d_{1}(n)$ is a constant
and $D_{2}(n)$ is monotonically non-decreasing:

\begin{Lemma} \label{bounded781}
Let $H$ be a fixed connected planar graph on $\{ 1,2, \ldots, h \}$.
Then $\exists \alpha (h) \!>\!0$ such that,
given any integer-valued functions $d_{1}(n)$ and $D_{2}(n)$ satisfying
$\limsup_{n \to \infty}\!d_{1}(n) \!\leq\! \delta(H)$ 
and $\liminf_{n \to \infty}\!D_{2}(n) \!\geq\! \max \{\Delta(H), \deg_{\!H}\!(1)+~\!1, 3 \},$
we have
\begin{displaymath}
\mathbf{P} [ f_{H}^{d_{1}} (P_{n,d_{1},5,0,D_{2}}) \leq \alpha n ] < e^{-\alpha n} 
\textrm{ for all sufficiently large } n.
\end{displaymath}
\end{Lemma}
\textbf{Proof}
Suppose we can find a graph $H$ and functions $d_{1}(n)$ and $D_{2}(n)$
that satisfy the conditions of this lemma,
but not the conclusion,
and let $\alpha = \alpha(h)$ be as given by Lemma~\ref{bounded106}.
Then there exist arbitrarily large `bad' $n$ for which 
$\mathbf{P} [ f_{H}^{d_{1}} (P_{n,d_{1},5,0,D_{2}}) \leq \alpha n ] \geq e^{-\alpha n}$.

Let $n_{1}$ be one of these bad $n$ 
and let us try to find a bad $n_{2}>n_{1}$ with $D_{2}(n_{2}) \geq D_{2}(n_{1})$.
Let us then try to find a bad $n_{3}>n_{2}$ with $D_{2}(n_{3}) \geq D_{2}(n_{2})$,
and so on.
We will either
(a) obtain an infinite sequence $n_{1},n_{2},n_{3} \ldots$
with $n_{1}<n_{2}<n_{3}< \ldots$
and $D_{2}(n_{1}) \leq D_{2}(n_{2}) \leq D_{2}(n_{3}) \leq \ldots$,
or (b) we will find a value $n_{k}$ such that all bad $n>n_{k}$ have $D_{2}(n) \leq D_{2}(n_{k})$. \\

Note that we must have $d_{1}(n) \in \{0,1,2,3,4,5\}$ $\forall n$.
Hence, in case (a) there must exist a constant $d$ such that infinitely many of our $n_{i}$ satisfy $d_{1}(n_{i})=d$
(we shall call these $n_{i}$ `special').
Let the function $D_{2}^{*}$ be defined by setting $D_{2}^{*}(n)=D_{2}(n_{1})$~$\forall n \leq n_{1}$
and $D_{2}^{*}(n) = D_{2}(n_{j})$~$\forall n \in \{ n_{j-1}+1,n_{j-1}+2, \ldots, n_{j} \}$~$\forall j>~1$.
Then $D_{2}^{*}$ is a monotonically non-decreasing integer-valued function satisfying
$\liminf_{n \to \infty} D_{2}^{*}(n) \geq \liminf_{n \to \infty} D_{2}(n) 
\geq \max \{ \Delta(H), \deg_{H}(1)+1, 3 \}$.
Hence, since $d \leq \limsup_{n \to \infty} d_{1}(n) \leq \delta(H)$,
by Lemma~\ref{bounded106}
it must be that we have
$\mathbf{P} [ f_{H}^{d} (P_{n,d,5,0,D_{2}^{*}}) \leq \alpha n ] < e^{-\alpha n}$
for all sufficiently large $n$.
But recall that our infinitely many `special' $n_{i}$ satisfy 
$(d, D_{2}^{*}(n_{i})) = (d_{1}(n_{i}), D_{2}(n_{i}))$,
and so
$\mathbf{P} [ f_{H}^{d} (P_{n_{i},d,5,0,D_{2}^{*}}) \leq \alpha n_{i} ] 
= \mathbf{P} [ f_{H}^{d_{1}} (P_{n_{i},d_{1},5,0,D_{2}}) \leq \alpha n_{i} ] \geq e^{-\alpha _{i}}$ 
for these $n_{i}$.
Thus, we obtain a contradiction. 

In case (b),
note that we have $d_{1}(n_{i}) \in \{ 0,1,2,3,4,5 \}$ $\forall i$
and that we also have
$D_{2}(n_{i}) \in \{ 3,4, \ldots, D_{2}(n_{k}) \}$~$\forall i \geq k$.
Hence, there must exist constants $d$ and $D$ such that infintely many of our $n_{i}$ satisfy
$(d_{1}(n_{i}),D_{2}(n_{i})) = (d,D)$.
But by Lemma~\ref{bounded106} we have
$\mathbf{P} [ f_{H}^{d} (P_{n,d,5,0,D}) \leq~\alpha n ] < e^{-\alpha n}$
for all large $n$,
and so we again obtain a contradiction.
$\phantom{qwerty}$ 
\setlength{\unitlength}{0.25cm}
\begin{picture}(1,1)
\put(0,0){\line(1,0){1}}
\put(0,0){\line(0,1){1}}
\put(1,1){\line(-1,0){1}}
\put(1,1){\line(0,-1){1}}
\end{picture} \\
\\
\\

We will state a stronger version of this last result after one more definition:

\begin{Definition} \label{maxdisjoint}
Given a connected planar graph $H$,
a function $d_{1}(n)$,
and a planar graph $G$,
let \emph{\boldmath{$\widehat{f_{H}^{d_{1}}}(G)$}}
denote the maximum size of a set of totally edge-disjoint appearances of $H$ in $G$
such that the associated cut-edges are all between vertices with $\deg_{G}>d_{1}(|G|)$. \\
\end{Definition}
\begin{displaymath}
\end{displaymath}

Finally, we may now obtain our main theorem,
which is the following stronger totally edge-disjoint version of Lemma~\ref{bounded781}:

\begin{Theorem} \label{bounded11}
Let $H$ be a fixed connected planar graph on $\{ 1,2, \ldots, h \}$.
Then $\exists \beta(H)\!>\!0$ such that,
given any integer-valued functions $d_{1}(n)$ and $D_{2}(n)$ satisfying
$\limsup_{n \to \infty} \!d_{1}(n) \!\leq\! \delta(H)$
and $\liminf_{n \to \infty} \!D_{2}(n) \!\geq\! \max \{ \Delta(H), \deg_{\!H}\!(1)+~\!1, 3 \},$
we have
\begin{displaymath}
\mathbf{P}[\widehat{f_{H}^{d_{1}}}(P_{n,d_{1},5,0,D_{2}}) \leq \beta n] < e^{- \beta n} 
\textrm{ for all sufficiently large } n.
\end{displaymath}
\end{Theorem}
\textbf{Proof}
The result follows from Lemmas~\ref{bounded201} and~\ref{bounded781} by taking 
$\beta = \frac{\alpha}{h+1}$.
$\phantom{www}$
$\setlength{\unitlength}{.25cm}
\begin{picture}(1,1)
\put(0,0){\line(1,0){1}}
\put(0,0){\line(0,1){1}}
\put(1,1){\line(-1,0){1}}
\put(1,1){\line(0,-1){1}}
\end{picture}$ \\
\\
\\

We will now look at the case when $d_{1}(n)=D_{2}(n)$ $\forall n$.
We shall work towards a result similar to Theorem~\ref{bounded11},
but this time (again, to help us in Section~\ref{cpts})
we will find it more convenient to look at the concept of `$2$-appearances':

\begin{Definition} \label{2appdefn}
Let $J$ be a connected graph on the vertices $\{ 1,2, \ldots, |J| \}$.
Given a graph $G$,
we say that $J$ \emph{\textbf{2-appears}} at $W \subset V(G)$ if 
(a) the increasing bijection from $\{ 1,2, \ldots, |J| \}$ to $W$ gives an isomorphism 
between $J$ and the induced subgraph $G[W]$ of $G$;
and (b) there are exactly two edges, $e_{1} = r_{1}v_{1}$ and $e_{2} = r_{2}v_{2}$,
in $G$ between $W \supset \{ r_{1},r_{2} \}$ and $V(G) \setminus W \supset \{ v_{1}, v_{2} \}$,
these edges are non-adjacent
(i.e.~$r_{1} \neq r_{2}$ and $v_{1} \neq v_{2}$),
and $v_{1}$ and $v_{2}$ are also non-adjacent.

Let us call $\{ e_{1},e_{2} \}$ the
\emph{\textbf{associated 2-edge-set}}
of the $2$-appearance, 
let us call
$TE_{W}^{2} := E(G[W]) \cup \{ e_{1},e_{2} \}$ the \emph{\textbf{total edge set}} of the $2$-appearance
and let us call
$TV_{W}^{2} := W \cup \{ v_{1}, v_{2} \}$
the \emph{\textbf{total vertex set}} of the $2$-appearance.

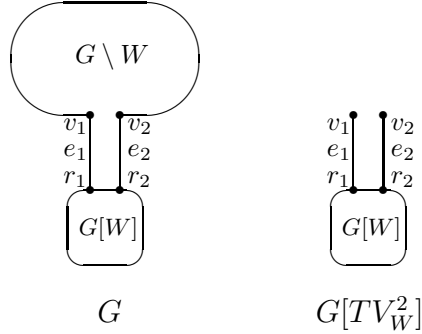
\begin{figure} [ht]
\setlength{\unitlength}{1cm}
\begin{picture}(20,4.5)(-3.375,-0.75)

\put(0.75,0){\line(1,0){0.5}}
\put(0.75,3.5){\line(1,0){0.5}}
\put(0.75,1){\line(1,0){0.5}}
\put(4.25,0){\line(1,0){0.5}}
\put(4.25,1){\line(1,0){0.5}}

\put(0.8,1){\line(0,1){1}}
\put(1.2,1){\line(0,1){1}}
\put(4.3,1){\line(0,1){1}}
\put(4.7,1){\line(0,1){1}}

\put(0.75,0.5){\oval(0.5,1)[l]}
\put(0.75,2.75){\oval(2,1.5)[l]}
\put(4.25,0.5){\oval(0.5,1)[l]}
\put(1.25,0.5){\oval(0.5,1)[r]}
\put(1.25,2.75){\oval(2,1.5)[r]}
\put(4.75,0.5){\oval(0.5,1)[r]}

\put(0.8,1){\circle*{0.1}}
\put(0.8,2){\circle*{0.1}}
\put(1.2,1){\circle*{0.1}}
\put(1.2,2){\circle*{0.1}}
\put(4.3,1){\circle*{0.1}}
\put(4.3,2){\circle*{0.1}}
\put(4.7,1){\circle*{0.1}}
\put(4.7,2){\circle*{0.1}}

\put(0.6,2.7){$G \setminus W$}
\put(0.65,0.4){\small{$G[W]$}}
\put(4.15,0.4){\small{$G[W]$}}

\put(0.9,-0.75){\large{$G$}}
\put(3.8,-0.75){\large{$G[TV^{2}_{W}]$}}

\put(1.3,1.45){$e_{2}$}
\put(1.3,1.8){$v_{2}$}
\put(1.3,1.1){$r_{2}$}
\put(4.8,1.45){$e_{2}$}
\put(4.8,1.8){$v_{2}$}
\put(4.8,1.1){$r_{2}$}

\put(0.45,1.45){$e_{1}$}
\put(0.45,1.8){$v_{1}$}
\put(0.45,1.1){$r_{1}$}
\put(3.95,1.45){$e_{1}$}
\put(3.95,1.8){$v_{1}$}
\put(3.95,1.1){$r_{1}$}

\end{picture}

\caption{$\textrm{A $2$-appearance at $W$ in $G$ and its total vertex/edge set.}$}
\end{figure}

\end{Definition}
\begin{displaymath}
\end{displaymath}

We will now prove results,
similar to Lemma~\ref{bounded201}, on the intersections of $2$-~appearances.
Since $d_{1}(n)=D_{2}(n)$ $\forall n$,
we must have $D_{2}(n) \leq 5$ $\forall n$ and so we can, in fact,
obtain a bound on the number of intersections of the total \textit{vertex} sets:

\begin{Lemma} \label{bounded502}
There exists $\Lambda \!=\! \Lambda(|J|)$ such that,
given any graph $G$ with $\Delta(G) \!\leq~\!\!5$,
the total vertex set of a $2$-appearance of a graph $J$ in $G$ will intersect
(i.e.~have a vertex in common with)
the total vertex set of at most $\Lambda$
other $2$-appearances of graphs of order~$|J|$.
\end{Lemma}
\textbf{Proof}
Suppose we have a $2$-appearance of a graph of order $|J|$ at $W \subset V(G)$,
and let $u \in TV_{W}^{2}$.
We shall show $u$ belongs to the total vertex set of at most
$\left( ^{|J|+2}_{\phantom{w}|J|} \right) \left( ^{5+5^{2}+ \ldots + 5^{|J|+1}}_{\phantom{wwww}|J|+1} \right)$
other $2$-appearances of graphs of order $|J|$.
Thus, since there are $|J|\!+\!2$ possibilities for $u$,
$TV_{W}^{2}$ only intersects the total vertex set~of at most
$(|J|\!+\!2) \left( ^{|J|+2}_{\phantom{w}|J|} \right) 
\left( ^{5+5^{2}+ \ldots + 5^{|J|+1}}_{\phantom{wwww}|J|+1} \right)$
other $2$-appearances of graphs of order~$|J|$.

Let us proceed with our argument by supposing that
$u$ also belongs to the total vertex set of a $2$-appearance of a graph of order $|J|$ at $W_{2}$.
Since $G[TV_{W_{2}}^{2}]$ is connected,
all the other $|J|+1$ vertices in $TV_{W_{2}}^{2}$
must be at most distance $|J|+1$ from $u$.
But, since $\Delta(G) \leq 5$,
there are at most $5+5^{2}+ \ldots + 5^{|J|+1}$ vertices in $V(G)$ that are at most distance $|J|+1$ from $u$.
Thus,
there are at most
$\left( ^{5+5^{2}+ \ldots + 5^{|J|+1}}_{\phantom{wwww}|J|+1} \right)$
possibilities for $TV_{W_{2}}^{2}$.
For each of these,
there are then at most 
$\left( ^{|J|+2}_{\phantom{w}|J|} \right)$
possibilities for $W_{2}$,
and hence in total there are at most
$\left( ^{|J|+2}_{\phantom{w}|J|} \right) \left( ^{5+5^{2}+ \ldots + 5^{|J|+1}}_{\phantom{wwww}|J|+1} \right)$
possibilities for $W_{2}$.
\phantom{qwerty}
\setlength{\unitlength}{0.25cm}
\begin{picture}(1,1)
\put(0,0){\line(1,0){1}}
\put(0,0){\line(0,1){1}}
\put(1,1){\line(-1,0){1}}
\put(1,1){\line(0,-1){1}}
\end{picture} \\
\\

It shall be useful for Section~\ref{cpts} to note that Lemma~\ref{bounded502}
immediately implies an analogous result for intersections of total edge sets:

\begin{Corollary} \label{bounded108}
There exists $\Lambda = \Lambda(|J|)$ such that,
given any graph $G$ with $\Delta(G) \leq 5$,
the total edge set of a $2$-appearance of a graph $J$ in $G$ will intersect
(i.e.~have an edge in common with)
the total edge set of at most $\Lambda$
other $2$-appearances of graphs of order $|J|$. 
\end{Corollary}

\phantom{p}

This second result actually holds even without the condition that $\Delta(G)$ is bounded,
but the details are rather fiddly.

We are now finally ready to obtain our first result on $2$-appearances in $P_{n,D_{2},5,0,D_{2}}$.
As with Lemma~\ref{bounded106},
we follow the method of Theorem 4.1 of \cite{mcd}:

\begin{Lemma} \label{bounded109}
Let $D_{2} \geq 3$ be a fixed constant, 
let $H$ be an $D_{2}$-regular connected planar graph on $\{ 1,2, \ldots,h \}$
and let $f \in E(H)$ be a non cut-edge.
Then $\exists \alpha >0$ and $\exists N$ such that
\begin{eqnarray*}
& & \mathbf{P}[P_{n,D_{2},5,0,D_{2}} \textrm{ will have } \leq \alpha n \textrm{ $2$-appearances of } H \setminus f] \\
& & \phantom{wwwwwwwwwwwwwwwww} < e^{- \alpha n} 
\left \{ \begin{array}{ll}
\forall n \geq N \textrm{ if } D_{2}=4 \\
\textrm{for all even } n \geq N \textrm{ if } D_{2} \in \{ 3,5 \}.
\end{array} \right.
\end{eqnarray*}
\end{Lemma}
\textbf{Sketch of Proof} 
In order to simplify parity matters,
we shall just prove the result for $D_{2}=4$,
but the $D_{2} \in \{ 3,5 \}$ cases will follow in a completely analogous way.

We choose a specific $\alpha$
and suppose that the result is false for $n=k$,
where $k$ is suitably large.
Using Theorem~\ref{bounded104},
it then follows that there are many graphs $G \in \mathcal{P}(k,D_{2},5,0,D_{2})$ with $\leq \alpha k$
$2$-appearances of $H \setminus f$.

From each such $G$,
we construct graphs in $\mathcal{P}((1+\delta)k,D_{2},5,0,D_{2})$, for a fixed $\delta >0$,
by replacing some edges in $G$ with $2$-appearances of $H \setminus f$.

The fact that the original graphs had few $2$-appearances of $H \setminus f$,
together with Corollary~\ref{bounded108},
is then used to show that there is not much double-counting,
and so we find that we have built so many graphs in $\mathcal{P}((1+\delta)k,D_{2},5,0,D_{2})$
that we contradict Theorem~\ref{bounded104}. \\
\\
\textbf{Full Proof}
As already mentioned in the sketch proof,
we shall just prove the $D_{2}=4$ case in order to avoid parity worries.
The $D_{2} \in \{ 3,5 \}$ cases will follow simply by substituting Theorem~\ref{bounded105} 
for every occurence of Theorem~\ref{bounded104}.

Let $\Lambda \!=\! \Lambda(h)$ be the constant given by Corollary~\ref{bounded108},
let $\beta \!=~\!\!\frac{ 8h! (\Lambda + 2) e^{2} \left( \gamma_{D_{2},D_{2}} \right)^{h} }{D_{2}}$
and let $\alpha$ be a fixed constant in 
$\left( 0, \min \left\{ \frac{D_{2}}{4}, \frac{1}{2h}, \frac{1}{\beta} \right\} \right)$.
Then we have $\alpha \beta <1$ 
and so $\exists \epsilon \in \left( 0,\frac{1}{3} \right)$ such that $(\alpha \beta)^{\alpha} = 1-3 \epsilon$.

By Theorem~\ref{bounded104},
$\exists N$ such that
\begin{equation} \label{eq:gr2}
(1-\epsilon)^{n}n! \left( \gamma_{D_{2},D_{2}} \right)^{n}
\leq |\mathcal{P}(n,D_{2},5,0,D_{2})|
\leq (1+\epsilon)^{n}n! \left( \gamma_{D_{2},D_{2}} \right)^{n} 
\textrm{ } \forall n \geq N.
\end{equation}
Suppose (aiming for a contradiction) that there exists a $k>N$ such that 
$
\mathbf{P} [P_{n,D_{2},5,0,D_{2}} \textrm{ will have } \leq \alpha k \textrm{ $2$-appearances of } H \setminus f] 
\geq e^{-\alpha k},
$
and let $\mathcal{G}$ denote the set of graphs in $\mathcal{P}(k,D_{2},5,0,D_{2})$ such that
$G \in \mathcal{G}$ iff 
$G \textrm{ has } \leq \alpha k$ $2$-appearances of $H \setminus f$.
Then we must have 
$|\mathcal{G}| \geq e^{-\alpha k} |\mathcal{P}(k,D_{2},5,0,D_{2})| \geq 
e^{-\alpha k} (1-\epsilon)^{k}k! \left( \gamma_{D_{2},D_{2}} \right)^{k}$. 

Let $\delta = \frac{\lceil \alpha k \rceil h}{k}$.
We may assume that $k$ is sufficiently large that 
$\lceil \alpha k \rceil \leq 2 \alpha k$.
Thus, $\delta \leq 2 \alpha h < 1$. \label{delta2}
This fact will be useful later. \\

We shall construct graphs in $\mathcal{P}((1+\delta)k,D_{2},5,0,D_{2})$:

Choose $\delta k$ special vertices 
(we have $\left(^{(1+\delta)k} _{\phantom{qq} \delta k} \right)$ choices for these)
and partition them into $\lceil \alpha k \rceil$ unordered blocks of size $h$
(we have $\left(^{\phantom{w} \delta k} _{h, \ldots, h} \right) \frac{1}{\lceil \alpha k \rceil !}$ choices for this).
On each of the blocks,
put a copy of $H \setminus f$ such that the increasing bijection from $\{ 1,2, \ldots, h \}$ to the block
is an isomorphism between $H \setminus f$ and this copy.
On the remaining $k$ vertices,
place an $D_{2}$-regular planar graph $G$ with at most
$\alpha k$ $2$-appearances of $H \setminus f$
(we have at least $|\mathcal{G}|$ choices for this).

Let $r_{B}$ and $s_{B}$ denote the two vertices in block $B$ with degree $D_{2}-1$.
For each of our $\lceil \alpha k \rceil$ blocks,
delete a different edge $u_{B}v_{B} \in E(G)$
(we have at least
$\left( ^{e(G)}_{\lceil \alpha k \rceil} \right)  \lceil \alpha k \rceil !
= \left( ^{kD_{2}/2}_{\phantom{`}\lceil \alpha k \rceil} \right) \lceil \alpha k \rceil !$
choices for this),
and insert edges $r_{B}v_{B}$ and $s_{B}u_{B}$
from the block to the vertices that were incident to the deleted edge
\begin{figure} [ht]
\setlength{\unitlength}{1cm}
\begin{picture}(20,2.5)(-0.5,-0.5)

\put(0,0.5){\line(0,1){0.5}}
\put(2,0.5){\line(0,1){0.5}}
\put(4,0.5){\line(0,1){0.5}}
\put(7,0.5){\line(0,1){0.5}}
\put(11,0.5){\line(0,1){0.5}}

\put(8,0.5){\line(1,0){1}}
\put(8,1){\line(1,0){1}}

\put(5,0.75){\vector(1,0){1}}

\put(0.5,1){\oval(1,0.5)[t]}
\put(0.5,0.5){\oval(1,0.5)[b]}
\put(3,1){\oval(2,2)[t]}
\put(3,0.5){\oval(2,2)[b]}
\put(7.5,1){\oval(1,0.5)[t]}
\put(7.5,0.5){\oval(1,0.5)[b]}
\put(10,1){\oval(2,2)[t]}
\put(10,0.5){\oval(2,2)[b]}

\put(1,0.75){\oval(0.5,0.5)[l]}
\put(8,0.75){\oval(0.5,0.5)[l]}

\put(1,0.5){\circle*{0.1}}
\put(1,1){\circle*{0.1}}
\put(2,0.5){\circle*{0.1}}
\put(2,1){\circle*{0.1}}
\put(8,0.5){\circle*{0.1}}
\put(8,1){\circle*{0.1}}
\put(9,0.5){\circle*{0.1}}
\put(9,1){\circle*{0.1}}

\put(0.3,0.6){$B$}
\put(2.8,0.6){\Large{$G$}}
\put(1,0.25){\footnotesize{$r_{B}$}}
\put(1,1.15){\footnotesize{$s_{B}$}}
\put(1.6,0.25){\footnotesize{$v_{B}$}}
\put(1.6,1.15){\footnotesize{$u_{B}$}}

\put(7.3,0.6){$B$}
\put(9.8,0.6){\Large{$G$}}
\put(8,0.25){\footnotesize{$r_{B}$}}
\put(8,1.15){\footnotesize{$s_{B}$}}
\put(8.6,0.25){\footnotesize{$v_{B}$}}
\put(8.6,1.15){\footnotesize{$u_{B}$}}

\end{picture}

\caption{Creating a $2$-appearance of $H \setminus f$ at $B$.} \label{6appfig}
\end{figure}
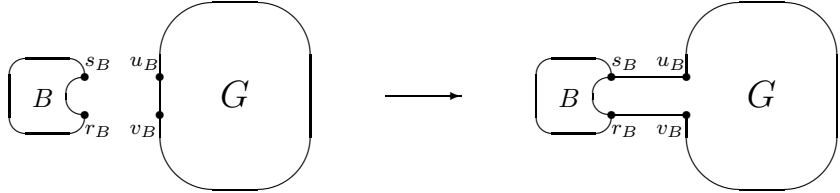 
(see Figure~\ref{6appfig}).
Planarity is maintained,
since if $H$ can be drawn with $f=rs$ on the outside,
then clearly~$H \setminus f$ can be drawn with $r$ and $s$ on the outside.

Thus, for each graph $G$
we can construct at least
\begin{eqnarray*}
& & \left(^{(1+\delta)k} _{\phantom{qq} \delta k} \right)
\left(^{\phantom{p}\delta k}_{h \ldots h} \right) 
\frac{1}{\lceil \alpha k \rceil !} \left( ^{kD_{2}/2}_{\phantom{`}\lceil \alpha k \rceil} \right)
\lceil \alpha k \rceil ! \\
& = & \frac{((1+\delta)k)!}{k!} \frac{1}{(h!)^{\lceil \alpha k \rceil}} 
\left( ^{kD_{2}/2}_{\phantom{`}\lceil \alpha k \rceil} \right) \\
& \geq & \frac{((1+\delta)k)!}{k!} \frac{1}{(h!)^{\lceil \alpha k \rceil}}
\frac{(kD_{2}/2 - \lceil \alpha k \rceil +1)^{\lceil \alpha k \rceil}}{\lceil \alpha k \rceil !} \\
& \geq & \frac{((1+\delta)k)!}{k!} \frac{1}{(h!)^{\lceil \alpha k \rceil}} 
\left( \frac{kD_{2}}{4} \right)^{\lceil \alpha k \rceil} \frac{1}{\lceil \alpha k \rceil !} \\
& & \left( \textrm{ since } \alpha < \frac{D_{2}}{4}
\textrm{ and so } kD_{2}/2 - \alpha k \geq \frac{kD_{2}}{4} \right) \\
& \geq & \frac{((1+\delta)k)!}{k!} \left( \frac{kD_{2}}{4h! \lceil \alpha k \rceil} \right)^{\lceil \alpha k \rceil} \\
& \geq & \frac{((1+\delta)k)!}{k!} \left( \frac{D_{2}}{8h! \alpha} \right)^{\lceil \alpha k \rceil} \\
& & \left( \textrm{ since we may assume $k$ is large enough that } 
\lceil \alpha k \rceil \leq 2 \alpha k \right)
\end{eqnarray*}
different graphs in $\mathcal{P}((1+\delta)k,D_{2},5,0,D_{2})$.

Therefore, recalling that we have at least
$e^{-\alpha k} (1-\epsilon)^{k}k! \left( \gamma_{D_{2},D_{2}} \right)^{k}$
choices for~$G$,
we can in total construct at least
$
e^{-\alpha k} (1\!-\!\epsilon)^{k} ((1\!+\!\delta)k)! \left( \gamma_{D_{2},D_{2}} \right)^{k}
\left( \frac{D_{2}}{8h! \alpha} \right)^{\lceil \alpha k \rceil}
$
(not necessarily distinct) graphs in $\mathcal{P}((1+\delta)k,D_{2},5,0,D_{2})$. \\

It now remains to investigate the amount of double-counting: 

Given one of our constructed graphs, $G^{\prime}$,
we can re-obtain the original graph,~$G$,
by deleting the $\lceil \alpha k \rceil$ $2$-appearances that were deliberately added
and re-inserting the $\lceil \alpha k \rceil$ deleted edges.
Notice that once we have identified the $2$-appearances that were deliberately added,
we will know where the deleted edges were.
Hence, it only remains to examine how many possibilities there are 
for the $\lceil \alpha k \rceil$ $2$-appearances that were deliberately added. 

Suppose there is a $2$-appearance of $H \setminus f$ at $W$ in $G^{\prime}$.
Let $\{ e_{1},e_{2} \}$ denote the associated $2$-edge-cut
and let $x$ and $y$ denote the vertices in $V(G^{\prime}) \setminus W$ 
that are incident to $e_{1}$ and $e_{2}$, respectively. 

(i) If there was a $2$-appearance of $H \setminus f$ at $W$ in $G$,
then there are at most~$\lceil \alpha k \rceil$ possibilities for $W$,
by definition of $\mathcal{G}$.

(ii) If the total vertex set of the $2$-appearance at $W$ intersects 
the total vertex set of one of our deliberately created $2 \textrm{-appearances}$,
then we have at most $(\Lambda + 1) \lceil \alpha k \rceil$ possibilities for $W$
(including the possibility that $W$ is one of our deliberately created $2$-appearances),
where $\Lambda$ is the constant given by Corollary~\ref{bounded108}.

(iii) Suppose (aiming for a contradiction) that neither (i) nor (ii) holds.
Then there was \textit{not} a $2$-appearance of $H \setminus f$ at $W$ in $G$,
but the total vertex set of the $2$-appearance at $W$ in $G^{\prime}$
does not intersect the total vertex set of any of our deliberately created $2$-appearances.
Note that we must then have $TE_{W}^{2} \subset E(G)$,
since it follows that $TV_{W}^{2} \subset V(G)$
and we know that we have not inserted any edges between vertices of $G$.
Thus, if there was not a $2$-appearance of $H \setminus f$ at~$W$ in $G$,
then there must have originally been either another edge between $W$ and $V(G) \setminus W$
other than $e_{1}$ and $e_{2}$,
or another edge between vertices in $W$,
or an edge between $x$ and $y$.
Note that in any of these three cases,
the deleted edge must have been incident to at least one vertex that is now in $TV_{W}^{2}$.
But the deleted edge must be of the form $u_{B}v_{B}$ for some block $B$,
and both $u_{B}$ and $v_{B}$ belong to the total vertex set of one of our deliberately created $2$-appearances.
Thus, this contradicts our assumption that $TV_{W}^{2}$ does not intersect
the total vertex set of any of our deliberately created $2$-appearances. \\

Thus, $G^{\prime}$ must have at most $(\Lambda + 2) \lceil \alpha k \rceil$
$2$-appearances of $H \setminus f$,
and so we have at most 
$\left( ^{(\Lambda + 2) \lceil \alpha k \rceil}_{\phantom{ww} \lceil \alpha k \rceil} \right)
\leq ((\Lambda + 2)e)^{\lceil \alpha k \rceil}$
possibilities for $G$. \\

We have now shown that each graph in $\mathcal{P}((1+\delta)k,D_{2},5,0,D_{2})$
may be constructed at most $((\Lambda + 2)e)^{\lceil \alpha k \rceil}$ times.
Thus, the number of \textit{distinct} graphs that we have constructed in $\mathcal{P}((1+\delta)k,D_{2},5,0,D_{2})$
must be at least \\
\begin{eqnarray*}
& & e^{-\alpha k} (1-\epsilon)^{k} ((1+\delta)k)! \left( \gamma_{D_{2},D_{2}} \right)^{k}
\left( \frac{D_{2}}{8h! \alpha (\Lambda +2)e} \right)^{\lceil \alpha k \rceil} \\
& \geq & e^{-\alpha k} (1\!-\!\epsilon)^{k} ((1\!+\!\delta)k)! \left( \gamma_{D_{2},D_{2}} \right)^{(1+\delta)k}
\left( \gamma_{D_{2},D_{2}} \right)^{- \lceil \alpha k \rceil h}
\left( \frac{D_{2}}{8h! \alpha (\Lambda +2)e} \right)^{\lceil \alpha k \rceil}\!\!, \\
& & \textrm{ since } \delta k = \lceil \alpha k \rceil h \\
& \geq & (1-\epsilon)^{k} ((1+\delta)k)! \left( \gamma_{D_{2},D_{2}} \right)^{(1+\delta)k}
\left( \frac{8h! \alpha (\Lambda +2) 
e^{2} \left( \gamma_{D_{2},D_{2}} \right)^{h}}{D_{2}} \right)^{- \lceil \alpha k \rceil}, \\
& & \textrm{ since } e^{- \alpha k} \geq e^{- \lceil \alpha k \rceil} \\
& = & (1-\epsilon)^{k} ((1+\delta)k)! \left( \gamma_{D_{2},D_{2}} \right)^{(1+\delta)k}
(\alpha \beta)^{- \lceil \alpha k \rceil} \\
& \geq & (1-\epsilon)^{k} ((1+\delta)k)! \left( \gamma_{D_{2},D_{2}} \right)^{(1+\delta)k}
(\alpha \beta)^{- \alpha k} 
\textrm{, \phantom{www} since } \alpha \beta < 1 \\
& = & \left( \frac{1-\epsilon}{1-3\epsilon} \right)^{k} ((1+\delta)k)! \left( \gamma_{D_{2},D_{2}} \right)^{(1+\delta)k} 
\textrm{, \phantom{www} since } (\alpha \beta)^{\alpha} = 1-3 \epsilon \\
& \geq & \left( \frac{1-\epsilon}{1-3\epsilon} \right)^{k} 
\frac{|\mathcal{P}((1+\delta)k,D_{2},5,0,D_{2})|}{(1+\epsilon)^{(1+\delta)k}} 
\textrm{, \phantom{www} by (\ref{eq:gr2})} \\
& > & \left( \frac{1-\epsilon}{(1-3\epsilon)(1+\epsilon)^{2}} \right)^{k} 
|\mathcal{P}((1+\delta)k,D_{2},5,0,D_{2})| 
\textrm{, \phantom{w} since } \delta < 1 \textrm{ (page~\pageref{delta2}}) \\
& > & |\mathcal{P}((1+\delta)k,D_{2},5,0,D_{2})|,
\phantom{ww} \textrm{ since }(1-3\epsilon)(1+\epsilon)^{2} = 1 - \epsilon - 5 \epsilon^{2} - 3 \epsilon^{3}.
\end{eqnarray*}
Hence, we have our desired contradiction.
$\phantom{qwerty}$ 
\setlength{\unitlength}{0.25cm}
\begin{picture}(1,1)
\put(0,0){\line(1,0){1}}
\put(0,0){\line(0,1){1}}
\put(1,1){\line(-1,0){1}}
\put(1,1){\line(0,-1){1}}
\end{picture} \\
\\

This time,
our main result is a totally \textit{vertex}-disjoint version:

\begin{Theorem} \label{bounded110}
Let $D_{2} \geq 3$ be a fixed constant, 
let $H$ be a $D_{2}$-regular connected planar graph on $\{1,2, \ldots,h \}$,
and let $f \in E(H)$ be a non cut-edge.
Then $\exists \beta$ and $\exists N$ such that
\begin{eqnarray*}
& & \mathbf{P}\Big[\textrm{$P_{n,D_{2},5,0,D_{2}}$ 
will not have a set of $\geq \beta n$ 
totally} \\
& & \phantom{wwwwwwq}\textrm{vertex-disjoint $2$-appearances of $H \setminus f$}\Big] \\
& & \phantom{wwwwwwwwwwwwwww}< e^{- \beta n} 
\left \{ \begin{array}{ll}
\forall n \geq N \textrm{ if } D_{2}=4 \\
\textrm{for all even } n \geq N \textrm{ if } D_{2} \in \{ 3,5 \}.
\end{array} \right.
\end{eqnarray*}
\end{Theorem}
\textbf{Proof}
This follows from Lemmas~\ref{bounded502} and~\ref{bounded109},
by taking $\beta = \frac{\alpha}{\Lambda+1}$.
$\phantom{qwerty}
\setlength{\unitlength}{.25cm}
\begin{picture}(1,1)
\put(0,0){\line(1,0){1}}
\put(0,0){\line(0,1){1}}
\put(1,1){\line(-1,0){1}}
\put(1,1){\line(0,-1){1}}
\end{picture}$

\newpage
\section{Components} \label{cpts}

We shall now use our appearance-type results from the previous section 
to investigate the probability of $P_{n,d_{1},5,0,D_{2}}$ having components isomorphic to given~$H$.
We already know from Section~\ref{bconn} that
(assuming $D_{2}(n) \geq 3$~$\forall n$, as always)
$\liminf_{n \to \infty} \mathbf{P}[P_{n,d_{1},5,0,D_{2}} \textrm{ will be connected}] >0$, 
so certainly we must have
$\limsup_{n \to \infty} \mathbf{P}[P_{n,d_{1},5,0,D_{2}} \textrm{ will have a component isomorphic to $H$}] <1$~$\forall H$.
In this section,
we will now see (in Theorem~\ref{bounded404}) that for all feasible $H$ we also have 
$\liminf_{n \to \infty} \mathbf{P}[P_{n,d_{1},5,0,D_{2}} \textrm{ will have a component isomorphic to $H$}] >0$.

As we are going to be using Theorems~\ref{bounded11} and~\ref{bounded110} from Section~\ref{apps},
we will start by dealing with the $d_{1}(n) < D_{2}(n)$ and $d_{1}(n)=D_{2}(n)$ cases separately
(in Lemmas~\ref{bounded7} and~\ref{bounded111}, respectively),
but we shall then combine these results in Theorem~\ref{bounded404}. \\

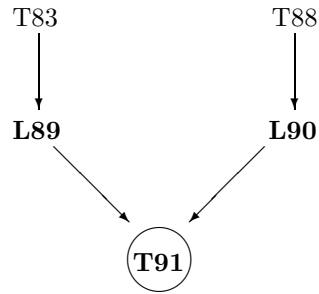
\begin{figure} [ht]
\setlength{\unitlength}{1cm}
\begin{picture}(20,4.2)(-3.4,0)

\put(2.6,0.4){\circle{0.8}}

\put(1.2,1.9){\vector(1,-1){1}}
\put(4,1.9){\vector(-1,-1){1}}
\put(1,3.4){\vector(0,-1){1}}
\put(4.4,3.4){\vector(0,-1){1}}

\put(0.65,3.5){T\ref{bounded11}}
\put(4.1,3.5){T\ref{bounded110}}
\put(0.65,2){\textbf{L\ref{bounded7}}}
\put(4.05,2){\textbf{L\ref{bounded111}}}
\put(2.25,0.25){\textbf{T\ref{bounded404}}}

\end{picture}

\caption{The structure of Section~\ref{cpts}.}
\end{figure}

\begin{displaymath}
\end{displaymath} 

We will start with the case when $d_{1}(n) < D_{2}(n)$~$\forall n$.
We shall prove a stronger result than that advertised at the beginning of this section,
as we are actually able to deal with several components at once:

\begin{Lemma} \label{bounded7}
Let $d_{1}(n)$ and $D_{2}(n)$ be any integer-valued functions 
that for all $n$ satisfy $D_{2}(n) \geq 3$ and $d_{1}(n)<D_{2}(n)$,
and let $t$ be a fixed constant.
Then, given any connected planar graphs $H_{1},H_{2}, \ldots, H_{k}$ with 
$\limsup_{n \to \infty}d_{1}(n) \leq \delta(H_{i}) \leq \Delta(H_{i}) \leq \liminf_{n \to \infty}D_{2}(n)$~$\forall i$,
we have
\begin{eqnarray*}
& & \liminf_{n \to \infty}
\mathbf{P}
\Big[\bigcap_{i \leq k} \left(\textrm{$P_{n,d_{1},5,0,D_{2}}$ 
will have $\geq t$ components}\right. \\
& & \textrm{\phantom{wwwwwwwwwww}with order-preserving isomorphisms to $H_{i}$} ) ] 
> 0.
\end{eqnarray*}
\end{Lemma}
\textbf{Sketch of Proof} 
The proof is by induction on $k$.
To prove the $k=l+1$ case,
we first consider the set of values of $n$ with $\delta (H_{l+1}) < D_{2}(n)$.
Without loss of generality,
$V(H_{l+1})= \{ 1,2, \ldots, l+1 \}$ and~$\deg_{H_{l+1}}(1)<D_{2}(n)$,
and so we can use Theorem~\ref{bounded11} to assume that $\widehat{f_{H_{l+1}}^{d_{1}}}$ is large.

We may then turn some of these appearances into components 
by deleting the associated cut-edges
(although we must take care not to interfere with our components isomorphic to 
$H_{1},H_{2}, \ldots, H_{l}$,
and also not to delete so many edges from the same vertex that we violate our bound on the minimum degree),
and we will see that the amount of double-counting is small 
unless there were already lots of components with order-preserving isomorphisms to $H_{l+1}$
originally. 

We shall then use a symmetry argument to deduce that the result also holds 
for the set of values of $n$ with $\delta(H_{l+1}) = D_{2}(n)$,
by relating the probability of having a component with an order-preserving isomorphism to $H_{l+1}$ 
to that of having a component with an order-preserving isomorphism to $H_{l+1} \setminus f$,
where $f$ is an arbitrary non cut-edge. \\
\\
\textbf{Full Proof}
Without loss of generality,
we may assume that $H_{1},H_{2}, \ldots, H_{k}$
are all distinct and that 
$e(H_{1}) \geq e(H_{2}) \geq \ldots \geq e(H_{k})$.
For brevity,
we will say that a component is `order-isomorphic' to $H_{i}$
to mean that it has an order-preserving isomorphism to $H_{i}$.
As mentioned in the sketch-proof,
we shall prove our result by using induction on $k$.
Since the statement is vacuous for $k=0$,
it suffices to suppose that it holds $\forall k \leq l$
and show that it must then also hold for $k=l+1$: \\

Let $\mathcal{N}^{<}$ denote the set of values of $n$ for which $\delta(H_{l+1})<D_{2}(n)$.
We shall start by proving that there exists a constant $c_{1}>0$ such that we have
$\mathbf{P}\left[ \bigcap_{i \leq l+1} \left(P_{n,d_{1},5,0,D_{2}} 
\textrm{ will have $\geq t$ components order-isomorphic to $H_{i}$} \right) \right] > c_{1}$
for all sufficiently large $n \in \mathcal{N}^{<}$.

Without loss of generality (by symmetry),
we may assume that we have
$V(H_{l+1}) = \{ 1,2, \ldots, |H_{l+1}| \}$
and that $\deg_{H_{l+1}}(1) = \delta (H_{l+1})$.
Thus, by Theorem~\ref{bounded11}, we know
$\exists \beta >0$ and $\exists N_{1}$ such that
$\mathbf{P}[\widehat{f_{H}^{d_{1}}}(P_{n,d_{1},5,0,D_{2}}) 
\leq \beta n] <~\!e^{- \beta n}$ for $\{ n \in \mathcal{N}^{<}:n \geq N_{1} \}$.

Let $\mathcal{A}(n,d_{1},D_{2})$ 
denote the set of graphs in 
$\mathcal{P}(n,d_{1},5,0,D_{2})$ 
that contain at least $t$ components order-isomorphic to $H_{i}$ $\forall i \leq l$,
and let $\mathcal{B}(n,d_{1},D_{2})$
denote the set of graphs in $\mathcal{A}(n,d_{1},D_{2})$ 
that also contain at least $t$ components order-isomorphic to~$H_{l+1}$.
By our induction hypothesis,
$\exists \epsilon >0$ and $\exists N_{2}$ such that
$|\mathcal{A}(n,d_{1},D_{2})| \geq \epsilon |\mathcal{P}(n,d_{1},5,0,D_{2})|$ $\forall n \geq N_{2}$.

Suppose there exists a value $n \in \mathcal{N}^{<}$ with $n \geq \max \{ N_{1},N_{2} \}$
such that 
$|\mathcal{B}(n,d_{1},D_{2})| < \frac{\epsilon |\mathcal{P}(n,d_{1},5,0,D_{2})|}{2}$
(if not, then we are done).
Let $\mathcal{G}(n,d_{1},D_{2})$ denote the collection of graphs in
$\mathcal{A}(n,d_{1},D_{2}) \setminus \mathcal{B}(n,d_{1},D_{2})$
that contain a set of at least $\beta n$ totally edge-disjoint appearances of $H_{l+1}$
such that the associated cut-edges are all between vertices of degree $>d_{1}(n)$.
Then
\begin{eqnarray*}
|\mathcal{G}(n,d_{1},D_{2})| 
& \geq & \epsilon |\mathcal{P}(n,d_{1},D_{2})| - \frac{ \epsilon |\mathcal{P}(n,d_{1},D_{2})| }{2}
- e^{- \beta n} |\mathcal{P}(n,d_{1},5,0,D_{2})| \\
& > & \frac{\epsilon |\mathcal{P}(n,d_{1},5,0,D_{2})|}{4},
\textrm{ if we assume $n$ is sufficiently large.} 
\end{eqnarray*}

Given a graph $G \in \mathcal{G}(n,d_{1},D_{2})$,
we may construct a graph in $\mathcal{B}(n,d_{1},D_{2})$ as follows: 

For each $i \leq l$,
we choose $t$ special components order-isomorphic to $H_{i}$.
Let $G^{\prime}$ denote the rest of the graph
(i.e.~away from these $tl$ special components).
We know that $G$ contains a set of $\lceil \beta n \rceil$
totally edge-disjoint appearances of $H_{l+1}$ 
such that the associated cut-edges are all between vertices of degree $>d_{1}(n)$.
At most $z:= \frac{ t \sum_{i \leq l}e(H_{i}) }{ e(H_{l+1})+1 }$
of these can be in our special components,
so (for sufficiently large $n$)
$G^{\prime}$ must contain a set, $S$,
of at least $\frac{\beta n}{2}$ 
totally edge-disjoint appearances of $H_{l+1}$ 
such that the associated cut-edges are all between vertices of degree $>d_{1}(n)$.

For $v \in V(G^{\prime})$,
let $a(v)$ denote the number of appearances of $H_{l+1}$ in $S$ 
such that $v$ is the unique vertex that is in the total vertex set of the appearance,
but not in the appearance itself.
Note that deg$(v)>d_{1}(n)$ if $a(v) \geq 1$.
Thus, for~$a(v) \geq 1$, we have
deg$(v)-d_{1}(n) \geq \max \{ 1, a(v)- d_{1}(n) \} \geq \frac{a(v)}{d_{1}(n)+1}$.
Therefore, $\min \{ a(v), \textrm{deg}(v)-d_{1}(n) \} \geq~\frac{a(v)}{d_{1}(n)+1}$~$\forall v$.

For each $v \in V(G^{\prime})$,
let us choose from $S$ exactly $\min \{ a(v), \textrm{deg}(v)-~\!d_{1}(n) \}$ appearances 
where $v$ is the unique vertex that is in the total vertex set of the appearance,
but not in the appearance itself.
Let $T$ denote the collection of all these chosen appearances.
Then, by the previous paragraph,
we have $|T| \geq~\sum_{v \in G} \frac{a(v)}{d_{1}(n)+1} =~\frac{|S|}{d_{1}(n)+1} \geq \frac{\beta n}{2(d_{1}(n)+1)}$. 

We may then complete our construction by choosing $t$ appearances from $T$
(at least $\left( ^{ \frac{\beta n}{2(d_{1}(n)+1)} } _{\phantom{www} t} \right)$ choices)
and simply deleting the $t$ associated cut-edges
(which are all distinct,
since the appearances are totally edge-disjoint).

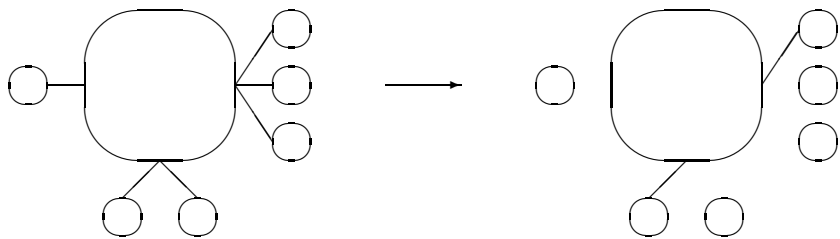
\begin{figure} [ht]
\setlength{\unitlength}{1cm}
\begin{picture}(20,3.25)(-0.5,0)

\put(2,2){\oval(2,2)}

\put(0.25,2){\oval(0.5,0.5)}
\put(1.5,0.25){\oval(0.5,0.5)}
\put(2.5,0.25){\oval(0.5,0.5)}
\put(3.75,1.25){\oval(0.5,0.5)}
\put(3.75,2){\oval(0.5,0.5)}
\put(3.75,2.75){\oval(0.5,0.5)}

\put(0.5,2){\line(1,0){0.5}}
\put(1.5,0.5){\line(1,1){0.5}}
\put(2.5,0.5){\line(-1,1){0.5}}
\put(3,2){\line(2,3){0.5}}
\put(3,2){\line(1,0){0.5}}
\put(3,2){\line(2,-3){0.5}}

\put(5,2){\vector(1,0){1}}

\put(9,2){\oval(2,2)}

\put(7.25,2){\oval(0.5,0.5)}
\put(8.5,0.25){\oval(0.5,0.5)}
\put(9.5,0.25){\oval(0.5,0.5)}
\put(10.75,1.25){\oval(0.5,0.5)}
\put(10.75,2){\oval(0.5,0.5)}
\put(10.75,2.75){\oval(0.5,0.5)}

\put(8.5,0.5){\line(1,1){0.5}}
\put(10,2){\line(2,3){0.5}}

\end{picture}

\caption{Using appearances to construct components isomorphic to $H_{l+1}$.}
\end{figure}

It remains to check that we have not violated our bound on the minimum degree:

Let $W_{1}$ and $W_{2}$ be two appearances in $T$
and suppose there exists a vertex~$u \in W_{1} \cap TV_{W_{2}}$.
Then we would have $TE_{W_{1}} \cap TE_{W_{2}} \neq \emptyset$,
which would be a contradiction,
since the appearances in $S$ (and hence $T$) are all totally edge-disjoint.
Thus, there exists no such vertex.
Therefore, the degree of a vertex $v$ can only have decreased by at most $1$ if it was in an appearance in $T$
(since it can't have been in the total vertex set of any other appearance in $T$)
and can only have decreased by at most deg$(v)-d_{1}(n)$
if it was not in an appearance in $T$
(since there are then at most deg$(v)-d_{1}(n)$ appearances in $T$
such that $v$ is in the total vertex set of $T$).
Thus, all vertices still have degree $\geq d_{1}(n)$,
and so our graph still has minimum degree at least $d_{1}(n)$. 

Thus, the constructed graphs are indeed in $\mathcal{B}(n,d_{1},D_{2})$.
Therefore, we have at least 
$\left( ^{\frac{\beta n}{2(d_{1}(n)+1)}}_{\phantom{www}t} \right)
|\mathcal{G}(n,d_{1},D_{2})|
\geq \left( ^{\frac{\beta n}{12}}_{\phantom{q}t} \right)
\frac{\epsilon |\mathcal{P}(n,d_{1},5,0,D_{2})|}{4}$
ways to construct a graph in~$\mathcal{B}(n,d_{1},D_{2})$.

Each time we deleted an edge,
we can have created at most $2$ components order-isomorphic to $H_{l+1}$.
Recall that each of our original graphs had $<t$ components order-isomorphic to $H_{l+1}$.
Thus, each of our constructed graphs will have at most $3t$ components order-isomorphic to $H_{l+1}$,
and so we will have at most $\left( ^{3t}_{\phantom{i}t} \right)$
possibilities for which are our deliberately created components.
For each of these deliberately created components,
we then have at most $n$ possibilities for which vertex in the rest of the graph was 
incident to the associated cut-edge.
Thus, each graph will have been constructed at most 
$\left( ^{3t}_{\phantom{i}t} \right) n^{t}$ times in total.

Therefore $\left( \textrm{under the assumption }
|\mathcal{B}(n,d_{1},D_{2})| < \frac{\epsilon |\mathcal{P}(n,d_{1},5,0,D_{2})|}{2} \right)$,
we have shown
$|\mathcal{B}(n,d_{1},D_{2})| \geq 
\frac{\left( ^{\frac{\beta n}{12}}_{\phantom{w}t} \right) \epsilon |\mathcal{P}(n,d_{1},5,0,D_{2})|}
{4 \left( ^{3t}_{\phantom{i}t} \right) n^{t}}$.
Thus, for a suitable $c_{1}>0$,
we have
$\mathbf{P}\big[ \bigcap_{i \leq l+1} (P_{n,d_{1},5,0,D_{2}}$ 
will have $\geq t$ components order-isomorphic to $H_{i}) \big] > c_{1}$
for all sufficiently large $n \in \mathcal{N}^{<}$. \\

Let $\mathcal{N}^{=}$ denote the set of values of $n$ for which $\delta(H_{l+1})=D_{2}(n)$.
We shall now prove $\exists c_{2}>0$ such that, 
for all sufficiently large $n \in \mathcal{N}^{=}$,
we have
$\mathbf{P}\big[ \bigcap_{i \leq l+1} (P_{n,d_{1},5,0,D_{2}}$ 
will have $\geq t$ components order-isomorphic to $H_{i}) \big] > c_{2}$.

If $\mathcal{N}^{=} \neq \emptyset$,
then we must have $\delta(H_{l+1}) \geq 3$,
since $D_{2}(n) \geq 3$.
Thus, there exists a non cut-edge $f \in E(H_{l+1})$.
Let $H_{l+1}^{\prime} = H_{l+1} \setminus f$
and note that $\delta(H_{l+1}^{\prime}) = D_{2}(n)-1 \geq d_{1}(n)$ $\forall n \in \mathcal{N}^{=}$.
Thus, we may use Theorem~\ref{bounded11} with~$H_{l+1}^{\prime}$
and then follow the same proof as with $\mathcal{N}^{<}$ to find $c_{2}>0$
such that 
\begin{eqnarray*}
\mathbf{P}\left[ \bigcap_{i \leq l}~(P_{n,d_{1},5,0,D_{2}} \right.
\textrm{ will have $\geq t$ components order-isomorphic to $H_{i}$} ) & \\
\cap~(P_{n,d_{1},5,0,D_{2}} 
\textrm{ will have $\geq t$ components order-isomorphic to $H_{l+1}^{\prime}$} )\biggr] & > & c_{2}
\end{eqnarray*}
for all sufficiently large $n \in \mathcal{N}^{=}$. 

Recall that $e(H_{i}) \!\geq\! e(H_{2}) \!\geq\! \ldots \!\geq\! e(H_{l+1})$.
Thus, since $e(H_{l+1}^{\prime}) \!=~\!\!e(H_{l+1}) \!-~\!\!1$,
it must be that $H_{1},H_{2}, \ldots, H_{l},H_{l+1},H_{l+1}^{\prime}$
are all \textit{distinct} graphs.
Hence, for $n \in \mathcal{N}^{=}$ and for any $x$ and $y$,
the number of graphs in $\mathcal{P}(n,d_{1},5,0,D_{2})$
that contain at least $t$ components order-isomorphic to $H_{i}$ $\forall i \leq l$,
exactly $x$ components order-isomorphic to~$H_{l+1}$
and exactly $y$ components order-isomorphic to $H_{l+1}^{\prime}$
is clearly exactly the same as the number of graphs in $\mathcal{P}(n,d_{1},5,0,D_{2})$
that contain at least $t$ components order-isomorphic to $H_{i}$ $\forall i \leq l$,
exactly $y$ components order-isomorphic to $H_{l+1}$
and exactly~$x$ components order-isomorphic to~$H_{l+1}^{\prime}$,
since swapping the components order-isomorphic to~$H_{l+1}$
and the components order-isomorphic to $H_{l+1}^{\prime}$ 
gives a bijection between such graphs.
Thus, for $n \in \mathcal{N}^{=}$, 
the number of graphs in $\mathcal{P}(n,d_{1},5,0,D_{2})$ that contain
at least $t$ components order-isomorphic to $H_{i}$ $\forall i \leq l+1$
must be exactly the same as the number of graphs in $\mathcal{P}(n,d_{1},5,0,D_{2})$ that contain
at least $t$ components order-isomorphic to $H_{i}$ $\forall i \leq l$
and at least $t$ components order-isomorphic to~$H_{l+1}^{\prime}$.

Therefore, to conclude,
we see that for all sufficiently large $n \in \mathcal{N}^{=}$ we have
$\mathbf{P}\left[ \bigcap_{i \leq l+1} \left(P_{n,d_{1},5,0,D_{2}} 
\textrm{ will have $\!\geq\! t$ components order-isomorphic to $H_{i}$} \right) \right] >~\!c_{2}$,
so 
$\mathbf{P}\left[ \bigcap_{i \leq l+1} \left(P_{n,d_{1},5,0,D_{2}} 
\textrm{ will have $\geq t$ components order-isomorphic to $H_{i}$} \right) \right] > \min\{c_{1},c_{2}\}$
for \textit{all} sufficiently large $n$.
\phantom{qwerty}
$\setlength{\unitlength}{.25cm}
\begin{picture}(1,1)
\put(0,0){\line(1,0){1}}
\put(0,0){\line(0,1){1}}
\put(1,1){\line(-1,0){1}}
\put(1,1){\line(0,-1){1}}
\end{picture}$ \\
\\

We shall now see an analogous result for when $d_{1}(n)=D_{2}(n)$~$\forall n$:

\begin{Lemma} \label{bounded111}
Let $D_{2} \geq 3$ and $t$ both be fixed.
Then, given any $D_{2}\textrm{-regular}$~connected planar graphs $H_{1},H_{2}, \ldots, H_{k}$,
there exist constants $\epsilon > 0$ and $N$ such~that
\begin{eqnarray*}
& &\mathbf{P}
\Big[ \bigcap_{i \leq k} ( P_{n,D_{2},5,0,D_{2}} 
\textrm{ will have $\geq t$ components} \\
& & \textrm{\phantom{wwwwww}with order-preserving isomorphisms to } H_{i} ) ] \\
& & \phantom{wwwwwwwwwwwwwwwwwww}> \epsilon
\left \{ \begin{array}{ll}
\forall n \geq N \textrm{ if } D_{2}=4 \\
\textrm{for all even } n \geq N \textrm{ if } D_{2} \in \{ 3,5 \}.
\end{array} \right.
\end{eqnarray*}
\end{Lemma}
\textbf{Sketch of Proof}
The proof is again by induction on $k$.
To prove the $k=l+1$ case,
we use Theorem~\ref{bounded110} to show that 
we may assume that there are lots of totally vertex-disjoint $2$-appearances of $H_{l+1} \setminus f$,
for some non cut-edge $f$,
and then use these $2$-appearances to create components isomorphic to $H_{l+1}$.
Again, the amount of double-counting is small unless there were already lots of components isomorphic to $H_{l+1}$
in the original graph. \\
\\
\textbf{Full Proof}
To simplify parity matters,
we will just consider the $D_{2}=4$ case. 

We shall prove the result by induction on $k$.
Since the statement is vacuous for $k=0$,
it suffices to suppose that it holds $\forall k \leq l$
and show that it must then also hold for $k=l+1$.

Without loss of generality,
we may assume that $V(H_{l+1}) = \{ 1,2, \ldots, |H_{l+1}| \}$.
Let $f \in E(H_{l+1})$ be an arbitrary non cut-edge.
Then, by Theorem~\ref{bounded110},
we know that there exists $\beta >0$ and there exists $N_{1}$ such that for all $n \geq N_{1}$ we have
$\mathbf{P}[P_{n,D_{2},5,0,D_{2}}$ does not have a set of 
$\geq \beta n$
totally vertex-disjoint $2$-appearances of $H_{l+1} \setminus f] < e^{- \beta n}$.

Let $\mathcal{A}(n,D_{2},D_{2})$ denote the set of graphs in $\mathcal{P}(n,D_{2},5,0,D_{2})$
that contain at least $t$ components with order-preserving isomorphisms to $H_{i}$ $\forall i \leq l$,
and let $\mathcal{B}(n,D_{2},D_{2})$ denote the set of graphs in $\mathcal{A}(n,D_{2},D_{2})$
that also contain at least $t$ components with order-preserving isomorphisms to $H_{l+1}$.
By our induction hypothesis,
$\exists \delta >0$ and $N_{2}$ such that
$|\mathcal{A}(n,D_{2},D_{2})| \!\geq~\!\!\delta |\mathcal{P}(n,D_{2},5,0,D_{2})|$~$\forall n \!\geq~\!\!N_{2}$.

Suppose $\exists n \geq \max \{ N_{1},N_{2} \}$ such that 
$|\mathcal{B}(n,D_{2},D_{2})| < \frac{\delta |\mathcal{P}(n,D_{2},5,0,D_{2})|}{2}$
(if not, then we are done).
Let $\mathcal{G}(n,D_{2},D_{2})$ denote the collection of graphs in 
$\mathcal{A}(n,D_{2},D_{2}) \setminus \mathcal{B}(n,D_{2},D_{2})$ that
contain a set of at least $\beta n$ totally vertex-disjoint $2$-appearances of $H_{l+1} \setminus f$.
Then we have 
\begin{eqnarray*}
|\mathcal{G}(n,D_{2},D_{2})| 
& \geq & \delta |\mathcal{P}(n,D_{2},D_{2})| \!-\! \frac{ \delta |\mathcal{P}(n,D_{2},D_{2})| }{2}
\!-\! e^{- \beta n} |\mathcal{P}(n,D_{2},5,0,D_{2})| \\
& > & \frac{\delta |\mathcal{P}(n,D_{2},5,0,D_{2})|}{4},
\textrm{ if we assume $n$ is sufficently large.}
\end{eqnarray*}

Given a graph $G \in \mathcal{G}(n,D_{2},D_{2})$,
we may construct a graph in $\mathcal{B}(n,D_{2},D_{2})$ as follows:
for each $i \leq l$,
we choose $t$ special components with order-preserving isomorphisms to $H_{i}$;
from the rest of the graph 
(i.e.~away from these $tl$ special components),
we choose $t$ totally vertex-disjoint $2$-appearances of $H_{l+1} \setminus f$
and denote these by $W_{1}, W_{2}, \ldots, W_{t}$
(we know that $G$ contains a set of at least $\beta n$ totally vertex-disjoint $2$-appearances of $H_{l+1}$,
and at most $z := \frac{ t \sum_{i \leq l}|H_{i}| }{ |H_{l+1}|+2 }$
of these can be in our special components,
so the number of choices that we have for our $2$-appearances is at least
$\left( ^{\beta n -z}_{\phantom{qw}t} \right) \geq \left( ^{\frac{\beta n}{2}}_{\phantom{q}t} \right)$,
for all sufficiently large~$n$);
for each $i$ we delete the $2$ edges of the form $v_{1}r_{1}$ and $v_{2}r_{2}$
for $\{ r_{1}, r_{2} \} \subset W_{i}$ and $ \{ v_{1}, v_{2} \} \subset V(G) \setminus W_{i}$
(these $2t$ edges are all distinct,
since the chosen $2$-appearances are totally vertex-disjoint);
and, finally,
for each $i$ we insert the two edges $v_{1}v_{2}$ and $e=r_{1}r_{2}$
(note that $v_{1}$ and $v_{2}$ were not originally adjacent,
by the definition of a $2$-appearance,
and that $r_{1}$ and $r_{2}$ were also not originally adjacent,
since they must be the two vertices of degree $D_{2}-1$ in $G[W_{i}]$,
by $D_{2}$-regularity of $G$,
and $G[W_{i}]$ is isomorphic to $H \setminus f$).

\begin{figure} [ht]
\setlength{\unitlength}{0.7cm}
\begin{picture}(20,2.3)(-4.95,0.8)

\put(0.75,0){\line(1,0){0.5}}
\put(0.75,3.5){\line(1,0){0.5}}

\put(0.8,1){\line(0,1){1}}
\put(1.2,1){\line(0,1){1}}

\put(0.75,0.5){\oval(0.5,1)[l]}
\put(0.75,2.75){\oval(2,1.5)[l]}
\put(1.25,0.5){\oval(0.5,1)[r]}
\put(1.25,2.75){\oval(2,1.5)[r]}

\put(0.8,1){\circle*{0.1}}
\put(0.8,2){\circle*{0.1}}
\put(1.2,1){\circle*{0.1}}
\put(1.2,2){\circle*{0.1}}

\put(1.27,1.72){\small{$v_{2}$}}
\put(1.27,1.1){\small{$r_{2}$}}

\put(0.33,1.72){\small{$v_{1}$}}
\put(0.33,1.1){\small{$r_{1}$}}

\put(3.25,1.75){\vector(1,0){1}}

\put(6.25,0){\line(1,0){0.5}}
\put(6.25,3.5){\line(1,0){0.5}}

\put(6.25,1){\line(1,0){0.5}}
\put(6.25,2){\line(1,0){0.5}}

\put(6.25,0.5){\oval(0.5,1)[l]}
\put(6.25,2.75){\oval(2,1.5)[l]}
\put(6.75,0.5){\oval(0.5,1)[r]}
\put(6.75,2.75){\oval(2,1.5)[r]}

\put(6.3,1){\circle*{0.1}}
\put(6.3,2){\circle*{0.1}}
\put(6.7,1){\circle*{0.1}}
\put(6.7,2){\circle*{0.1}}

\put(6.77,1.72){\small{$v_{2}$}}
\put(6.77,1.1){\small{$r_{2}$}}

\put(5.83,1.72){\small{$v_{1}$}}
\put(5.83,1.1){\small{$r_{1}$}}

\end{picture}

\caption{Using $2$-appearances to construct components isomorphic to $H_{l+1}$.} \label{bcptfig}
\end{figure}
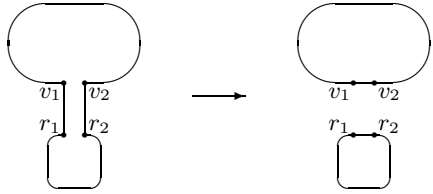

Note that the component,~$C$,
containing $r_{1}$ and $r_{2}$ will now have an order-preserving isomorphism to $H$,
since $C \setminus e$ has an order-preserving isomorphism to $H \setminus f$
and $C$ and $H$ are both $D_{2}$-regular.
Thus, our constructed graphs will now contain 
at least $t$ components with order-preserving isomorphisms to~$H_{i}$~$\forall i \leq~l+1$.
Note also that the degrees of all the vertices will remain unchanged and that planarity is preserved.
Thus, the constructed graphs are indeed in $\mathcal{P}(n,D_{2},5,0,D_{2})$.
Therefore, we have at least 
$\left( ^{\frac{\beta n}{2}}_{\phantom{q}t} \right)
|\mathcal{G}(n,D_{2},D_{2})|
\geq \left( ^{\frac{\beta n}{2}}_{\phantom{q}t} \right)
\frac{\delta |\mathcal{P}(n,D_{2},5,0,D_{2})|}{4}$
ways to construct a graph in~$\mathcal{B}(n,D_{2},D_{2})$.

Each time we performed the construction of Figure~\ref{bcptfig},
we can have created at most $2$ components with order-preserving isomorphisms to $H_{l+1}$.
Recall that each of our original graphs had $<t$ components with order-preserving isomorphisms to $H_{l+1}$.
Thus, each of our constructed graphs will have at most $3t$ components with order-preserving isomorphisms to $H_{l+1}$,
and so we will have at most $\left( ^{3t}_{\phantom{i}t} \right)$
possibilities for which are our deliberately created components.
We then know which edges were inserted into them and have at most 
$\left(^{3n}_{\phantom{i}t} \right)$,
by planarity,
possibilities for the edges that were inserted in the rest of the graph.
We then have $t!$ ways to pair up the edges inserted in the deliberately created components with
the edges inserted in the rest of the graph,
and hence $2^{t}t!$ possibilities for where the deleted edges were originally.
Thus, each graph will have been constructed at most 
$\left(^{3t}_{\phantom{i}t} \right) \left(^{3n}_{\phantom{i}t} \right) 2^{t}t!$ times.

Therefore
$\left( \textrm{under the assumption that }
|\mathcal{B}(n,D_{2},D_{2})| < \frac{\delta |\mathcal{P}(n,D_{2},5,0,D_{2})|}{2} \right)$,
we have shown 
$|\mathcal{B}(n,D_{2},D_{2})| \geq
\frac{\left( ^{\frac{\beta n}{2}}_{\phantom{q}t} \right) \delta|\mathcal{P}(n,D_{2},5,0,D_{2})|}
{4 \left(^{3t}_{\phantom{i}t} \right) \left(^{3n}_{\phantom{i}t} \right) 2^{t}t!}
= \Theta(1) |\mathcal{P}(n,D_{2},5,0,D_{2})|$.~$\setlength{\unitlength}{.25cm}
\begin{picture}(1,1)
\put(0,0){\line(1,0){1}}
\put(0,0){\line(0,1){1}}
\put(1,1){\line(-1,0){1}}
\put(1,1){\line(0,-1){1}}
\end{picture}$

\phantom{p}

We may now combine Lemmas~\ref{bounded7} and~\ref{bounded111} to obtain our full result:

\begin{Theorem} \label{bounded404}
Let $d_{1}(n)$ and $D_{2}(n)$ be any integer-valued functions, subject to $D_{2}(n) \geq 3$~$\forall n$
and $(d_{1}(n),D_{2}(n)) \notin \{ (3,3),(5,5) \}$ for odd $n$,
and let $t$ be a fixed constant.
Then, given any connected planar graphs 
$H_{1},H_{2}, \ldots, H_{k}$ with 
$\limsup_{n \to \infty} d_{1}(n) \leq \delta(H_{i}) \leq \Delta(H_{i}) 
\leq \liminf_{n \to \infty} D_{2}(n)$~$\forall i$,
we have
\begin{eqnarray*}
& & \liminf_{n \to \infty}
\mathbf{P}
\Big[ \bigcap_{i \leq k} (P_{n,d_{1},5,0,D_{2}} 
\textrm{ will have $\geq t$ components} \\
& & \textrm{\phantom{wwwwwwwwwww}with order-preserving isomorphisms to $H_{i}$} ) ] > 0.
\end{eqnarray*}
\end{Theorem}

\newpage
\section{Subgraphs} \label{subs}

We will now use the results of the previous two sections to investigate the probability of $P_{n,d_{1},5,0,D_{2}}$
having \textit{copies} of given (not necessarily connected) $H$.
As always, we shall assume throughout that $D_{2}(n) \geq 3$~$\forall n$.

Clearly, for those values of $n$ for which $D_{2}(n) < \Delta(H)$,
we must have $\mathbf{P}[P_{n,d_{1},5,0,D_{2}} \textrm{ will have a copy of } H] =0$.
For sufficiently large $n$,
it turns out that the only other time when we can have this
is if $d_{1}(n)=D_{2}(n)=4$ 
and $H$ happens to be a graph that can \textit{never} be a subgraph of a $4$-regular planar graph.
We will note a method for determining when this is so in Theorem~\ref{bounded991}.

Apart from the above exceptions,
we shall see that the matter of whether
$\mathbf{P}[P_{n,d_{1},5,0,D_{2}} \textrm{ will have a copy of } H]$
is bounded away from $0$ and/or $1$
actually depends only on whether $H$ has any $D_{2}(n)$-regular components.
For those values of $n$ for which this is the case,
it is already clear that the probability must be bounded away from $1$
(since we know from Theorem~\ref{bounded311} that
$\liminf_{n \to \infty} \mathbf{P}[P_{n,d_{1},5,0,D_{2}} \textrm{ will be connected}] >0$)
and in this section (Theorem~\ref{bounded1002})
we will deduce from Theorem~\ref{bounded404} that it is also bounded away from~$0$.
If there aren't arbitrarily large values of $n$ for which $H$ has $D_{2}(n)$-regular components,
we shall be able to use our appearance-type results of Section~\ref{apps} 
to see (in Theorem~\ref{bounded1001}) that 
$\mathbf{P}[P_{n,d_{1},5,0,D_{2}} \textrm{ will have a copy of } H] \to 1$
(again, with the exception of the cases given in the previous paragraph).

We will start by proceeding towards the
last result.
As mentioned,
the proof will use our appearance-type work,
and so we shall first split into $d_{1}(n) < D_{2}(n)$ and $d_{1}(n)=D_{2}(n)$ cases
(in Lemmas~\ref{bounded901} and~\ref{bounded904}, respectively),
before then combining these results in Theorem~\ref{bounded1001}.
We will then finish with the case when $H$ has $D_{2}(n)$-regular components,
in Theorem~\ref{bounded1002}.

\begin{figure} [ht]
\setlength{\unitlength}{1cm}
\begin{picture}(20,0.8)(-2.05,0.35)

\put(2.5,0){\circle{0.8}}
\put(6.4,0){\circle{0.8}}

\put(2.15,-0.15){\textbf{T\ref{bounded1001}}}
\put(6.05,-0.15){\textbf{T\ref{bounded1002}}}

\put(0.75,0.7){\textbf{L\ref{bounded901}}}
\put(2.15,0.7){\textbf{L\ref{bounded907}}}
\put(3.55,0.7){\textbf{L\ref{bounded904}}}

\put(5.5,0.7){T\ref{bounded311}}
\put(6.7,0.7){T\ref{bounded404}}

\put(0.75,1.3){T\ref{bounded110}}
\put(3.55,1.3){T\ref{bounded11}}

\put(3.1,0){\vector(1,0){2.7}}

\put(1.2,0.65){\vector(3,-1){0.9}}
\put(2.5,0.7){\vector(0,-1){0.3}}
\put(3.8,0.65){\vector(-3,-1){0.9}}

\put(5.8,0.68){\vector(2,-3){0.25}}
\put(7,0.68){\vector(-2,-3){0.25}}

\put(1.1,1.25){\vector(0,-1){0.3}}
\put(3.9,1.25){\vector(0,-1){0.3}}

\end{picture}

\caption{The structure of Section~\ref{subs}.}
\end{figure}
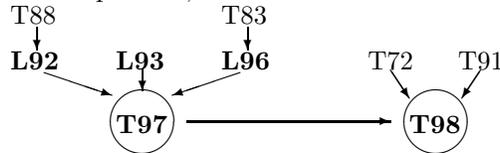

\begin{figure} [ht] 
\setlength{\unitlength}{0.85cm}
\begin{picture}(20,8.6)(-6.1,-5.8)

\put(-2,2){\line(1,0){2}}
\put(2,2){\line(1,0){2}}

\put(1,3){\line(-1,-1){1}}
\put(1,3){\line(1,-1){1}}

\put(1,1){\line(-1,1){1}}
\put(1,1){\line(1,1){1}}

\put(1,3){\line(0,-2){2}}

\put(0,2){\circle*{0.1}}
\put(2,2){\circle*{0.1}}
\put(1,3){\circle*{0.1}}
\put(1,1){\circle*{0.1}}

\put(-2,2){\circle*{0.1}}
\put(4,2){\circle*{0.1}}

\put(-2.1,2.2){$u_{i}$}
\put(3.9,2.2){$v_{i}$}

\put(-2,-0.5){\line(1,0){6}}

\put(1,0.5){\line(-1,-1){1}}
\put(1,0.5){\line(-1,-3){0.333333}}
\put(1,0.5){\line(1,-3){0.333333}}
\put(1,0.5){\line(1,-1){1}}

\put(1,-1.5){\line(-1,1){1}}
\put(1,-1.5){\line(-1,3){0.333333}}
\put(1,-1.5){\line(1,3){0.333333}}
\put(1,-1.5){\line(1,1){1}}

\put(0,-0.5){\circle*{0.1}}
\put(0.666667,-0.5){\circle*{0.1}}
\put(1.333333,-0.5){\circle*{0.1}}
\put(2,-0.5){\circle*{0.1}}
\put(1,0.5){\circle*{0.1}}
\put(1,-1.5){\circle*{0.1}}

\put(-2,-0.5){\circle*{0.1}}
\put(4,-0.5){\circle*{0.1}}

\put(-2.1,-0.3){$u_{i}$}
\put(3.9,-0.3){$v_{i}$}

\put(-1,-3.5){\line(1,0){4}}
\put(-1,-4.5){\line(1,0){4}}

\put(-1,-3.5){\line(0,-1){1}}
\put(0,-3.5){\line(0,-1){1}}
\put(1,-2.5){\line(0,-1){3}}
\put(2,-3.5){\line(0,-1){1}}
\put(3,-3.5){\line(0,-1){1}}

\put(1,-2.5){\line(-2,-1){2}}
\put(1,-2.5){\line(-1,-1){1}}
\put(1,-2.5){\line(1,-1){1}}
\put(1,-2.5){\line(2,-1){2}}

\put(1,-5.5){\line(-2,1){2}}
\put(1,-5.5){\line(-1,1){1}}
\put(1,-5.5){\line(1,1){1}}
\put(1,-5.5){\line(2,1){2}}

\put(1,-4.5){\oval(4,3)[b]}
\put(1,-3.5){\oval(4,3)[t]}

\put(-1,-3.5){\line(1,-1){1}}
\put(0,-3.5){\line(1,-1){1}}
\put(1,-3.5){\line(1,-1){1}}
\put(2,-3.5){\line(1,-1){1}}
\put(3,-3.5){\line(1,0){1}}

\put(-1,-3.5){\circle*{0.1}}
\put(0,-3.5){\circle*{0.1}}
\put(1,-3.5){\circle*{0.1}}
\put(2,-3.5){\circle*{0.1}}
\put(3,-3.5){\circle*{0.1}}

\put(-1,-4.5){\circle*{0.1}}
\put(0,-4.5){\circle*{0.1}}
\put(1,-4.5){\circle*{0.1}}
\put(2,-4.5){\circle*{0.1}}
\put(3,-4.5){\circle*{0.1}}

\put(1,-2.5){\circle*{0.1}}
\put(1,-5.5){\circle*{0.1}}

\put(-2,-4.5){\line(1,0){1}}

\put(-2.1,-4.3){$u_{i}$}
\put(3.9,-3.3){$v_{i}$}

\put(4,-3.5){\circle*{0.1}}
\put(-2,-4.5){\circle*{0.1}}

\end{picture}

\caption{Replacing the edge $f_{i}=u_{i}v_{i}$ (cases $D_{2}=3$, $D_{2}=4$ and $D_{2}=5$).} \label{edgefig}
\end{figure}
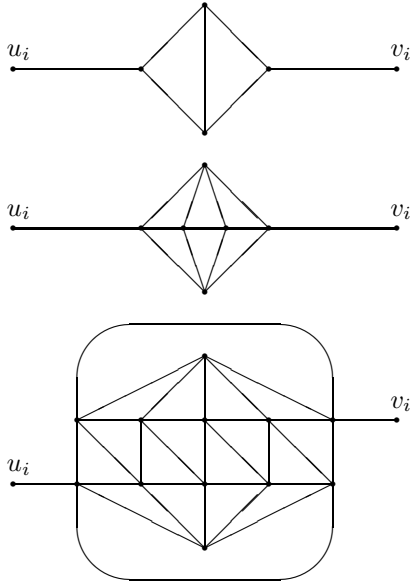

We start with the case when $d_{1}(n)=D_{2}(n)$~$\forall n$:

\begin{Lemma} \label{bounded901}
Let $H$ be a fixed planar graph with components $H_{1}, H_{2}, \ldots, H_{k}$, for some $k$,
and let $D_{2} \in \{3,4,5 \}$ be a fixed constant.
Suppose no component $H_{i} \textrm{ is $D_{2}$-regular}$,
but that for all $i$ there exists a $D_{2}$-regular planar graph $H_{i}^{*}$ that contains a copy of $H_{i}$.
Then $\exists \beta > 0$ and $\exists N$ such that
\begin{eqnarray*}
& & \mathbf{P} \Big[
\textrm{$P_{n,D_{2},5,0,D_{2}}$
will \emph{not} have a set of $\beta n$ vertex-disjoint} \\
& & \phantom{wwwwwwwwwwq}\textrm{induced order-preserving copies of $H$}\Big] \\
& & \phantom{wwwwwwwwwwwwwwww}< e^{- \beta n} 
\left \{ \begin{array}{ll}
\forall n \geq N \textrm{ if } D_{2}=4 \\
\textrm{for all even } n \geq N \textrm{ if } D_{2} \in \{ 3,5 \}.
\end{array} \right.
\end{eqnarray*}
\end{Lemma}
\textbf{Proof}
Without loss of generality,
we may assume that each $H_{i}^{*}$ is connected
(because each $H_{i}$ is connected).
Since no $H_{i}$ is $D_{2}$-regular,
it must be that each~$H_{i}^{*}$ contains an edge $f_{i}=u_{i}v_{i}$
such that $H_{i}^{*} \setminus f_{i}$ also contains a copy of $H_{i}$.
Without loss of generality,
we may assume that $f_{i}$ is not a cut-edge in $H_{i}^{*}$,
since we could replace $f_{i}$ with a copy of the appropriate graph from Figure~\ref{edgefig}
and in this way obtain a $D_{2}$-regular planar graph containing several non cut-edges
that don't interfere with our copy of $H_{i}$.
Thus,
we may assume that the graph formed from the $H_{i}^{*} \setminus f_{i}$'s by inserting the edges $u_{k}v_{1}$
and $u_{i}v_{i+1}$~$\forall i \leq k-1$
will be a \textit{connected} $D_{2}$-regular planar graph $H^{*}$ containing a copy of $H$.

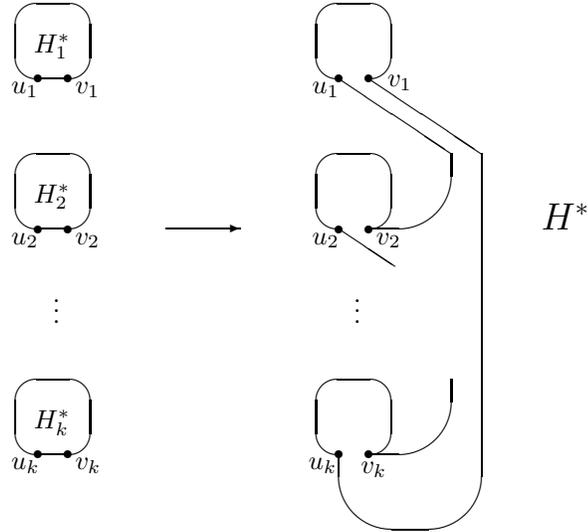
\begin{figure} [ht] 
\setlength{\unitlength}{1cm}
\begin{picture}(20,7.25)(-2.9,4)

\put(0.3,5){\line(1,0){0.4}}
\put(0.3,6){\line(1,0){0.4}}
\put(0.3,8){\line(1,0){0.4}}
\put(0.3,9){\line(1,0){0.4}}
\put(0.3,10){\line(1,0){0.4}}
\put(0.3,11){\line(1,0){0.4}}

\put(0.3,5.5){\oval(0.6,1)[l]}
\put(0.3,8.5){\oval(0.6,1)[l]}
\put(0.3,10.5){\oval(0.6,1)[l]}

\put(0.7,5.5){\oval(0.6,1)[r]}
\put(0.7,8.5){\oval(0.6,1)[r]}
\put(0.7,10.5){\oval(0.6,1)[r]}

\put(0.3,10){\circle*{0.1}}
\put(0.7,10){\circle*{0.1}}
\put(0.3,8){\circle*{0.1}}
\put(0.7,8){\circle*{0.1}}
\put(0.3,5){\circle*{0.1}}
\put(0.7,5){\circle*{0.1}}

\put(0.5,6.75){$\vdots$}

\put(0.25,10.35){$H_{1}^{*}$}
\put(0.25,8.35){$H_{2}^{*}$}
\put(0.25,5.35){$H_{k}^{*}$}

\put(-0.05,9.8){$u_{1}$}
\put(0.8,9.8){$v_{1}$}
\put(-0.05,7.8){$u_{2}$}
\put(0.8,7.8){$v_{2}$}
\put(-0.05,4.8){$u_{k}$}
\put(0.8,4.8){$v_{k}$}

\put(2,8){\vector(1,0){1}}

\put(7,8){\Large{$H^{*}$}}

\put(4.3,6){\line(1,0){0.4}}
\put(4.3,9){\line(1,0){0.4}}
\put(4.3,11){\line(1,0){0.4}}

\put(4.3,5.5){\oval(0.6,1)[l]}
\put(4.3,8.5){\oval(0.6,1)[l]}
\put(4.3,10.5){\oval(0.6,1)[l]}

\put(4.7,5.5){\oval(0.6,1)[r]}
\put(4.7,8.5){\oval(0.6,1)[r]}
\put(4.7,10.5){\oval(0.6,1)[r]}

\put(4.3,10){\circle*{0.1}}
\put(4.7,10){\circle*{0.1}}
\put(4.3,8){\circle*{0.1}}
\put(4.7,8){\circle*{0.1}}
\put(4.3,5){\circle*{0.1}}
\put(4.7,5){\circle*{0.1}}

\put(4.3,10){\line(3,-2){1.5}}
\put(4.7,10){\line(3,-2){1.5}}
\put(4.3,8){\line(3,-2){0.75}}

\put(4.7,9){\oval(2.2,2)[br]}
\put(4.7,6){\oval(2.2,2)[br]}

\put(6.2,9){\line(0,-1){4}}

\put(5.25,5){\oval(1.9,2)[b]}

\put(4.5,6.75){$\vdots$}

\put(3.95,9.8){$u_{1}$}
\put(4.95,9.9){$v_{1}$}
\put(3.95,7.8){$u_{2}$}
\put(4.8,7.8){$v_{2}$}
\put(3.9,4.8){$u_{k}$}
\put(4.6,4.75){$v_{k}$}

\end{picture}

\caption{Constructing $H^{*}$ from the $H_{i}^{*}$'s.}
\end{figure}

Without loss of generality,
we may assume that the copy of $H$ is order-preserving 
and also that it is induced 
(by again replacing appropriate edges with a copy of the relevant graph from Figure~\ref{edgefig}).
As before,
we may assume that $H^{*}$ contains a non cut-edge $f$ that doesn't interfere with this copy.
The result then follows from Theorem~\ref{bounded110}.
$\phantom{qwerty}
\setlength{\unitlength}{.25cm}
\begin{picture}(1,1)
\put(0,0){\line(1,0){1}}
\put(0,0){\line(0,1){1}}
\put(1,1){\line(-1,0){1}}
\put(1,1){\line(0,-1){1}}
\end{picture}$ \\
\\
\\

Lemma~\ref{bounded901} leaves us with the matter of discovering which graphs 
can't actually be contained within any $D_{2}$-regular planar graphs.
If $D_{2} \in \{ 3,5 \}$, it turns out that there are no such graphs:

\begin{Lemma} \label{bounded907}
Let $D_{2} \in \{ 3,5 \}$ be a fixed constant
and let $H$ be a planar graph with $\Delta(H) \leq D_{2}$.
Then there exists a $D_{2}$-regular planar graph $H^{*}$ that contains a copy of $H$.
\end{Lemma}
\textbf{Proof}
Let $L_{3}$ and $L_{5}$ denote the graphs shown in Figures~\ref{L3} and~\ref{L5}, respectively.
Then, for $D_{2} \in \{ 3,5 \}$,
\begin{figure} [ht] 
\setlength{\unitlength}{1cm}
\begin{picture}(10,1.5)(-5,-0.5)
\put(0,0){\line(1,0){2}}
\put(1,0){\line(0,1){1}}
\put(0,0){\line(1,1){1}}
\put(2,0){\line(-1,1){1}}
\put(0,0){\line(1,-1){1}}
\put(1,-1){\line(1,1){1}}

\put(0,0){\circle*{0.1}}
\put(1,0){\circle*{0.1}}
\put(2,0){\circle*{0.1}}
\put(1,1){\circle*{0.1}}
\put(1,-1){\circle*{0.1}}

\end{picture}
\caption{The graph $L_{3}$.} \label{L3}
\end{figure}
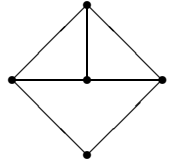 
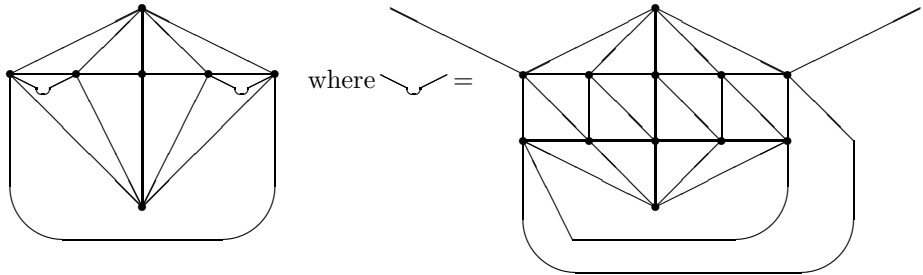
\begin{figure} [ht] 
\setlength{\unitlength}{0.88cm}
\begin{picture}(16,3.5)(0,-2.5)
\put(0,0){\line(1,0){4}}
\put(2,-2){\line(0,1){3}}
\put(0,0){\line(2,1){2}}
\put(1,0){\line(1,1){1}}
\put(3,0){\line(-1,1){1}}
\put(4,0){\line(-2,1){2}}
\put(0,0){\line(1,-1){2}}
\put(1,0){\line(1,-2){1}}
\put(3,0){\line(-1,-2){1}}
\put(4,0){\line(-1,-1){2}}

\put(2,0){\oval(4,5)[b]}

\put(0,0){\line(2,-1){0.4}}
\put(1,0){\line(-2,-1){0.4}}
\put(0.5,-0.2){\oval(0.2,0.2)[b]}
\put(3,0){\line(2,-1){0.4}}
\put(4,0){\line(-2,-1){0.4}}
\put(3.5,-0.2){\oval(0.2,0.2)[b]}

\put(0,0){\circle*{0.1}}
\put(1,0){\circle*{0.1}}
\put(2,0){\circle*{0.1}}
\put(3,0){\circle*{0.1}}
\put(4,0){\circle*{0.1}}
\put(2,-2){\circle*{0.1}}
\put(2,1){\circle*{0.1}}

\put(4.5,-0.2){where}

\put(5.6,0){\line(2,-1){0.4}}
\put(6.6,0){\line(-2,-1){0.4}}
\put(6.1,-0.2){\oval(0.2,0.2)[b]}

\put(6.7,-0.2){$=$}

\put(7.75,0){\line(1,0){4}}
\put(7.75,-1){\line(1,0){4}}

\put(7.75,0){\line(-2,1){2}}
\put(11.75,0){\line(2,1){2}}

\put(7.75,0){\line(0,-1){1}}
\put(8.75,0){\line(0,-1){1}}
\put(9.75,1){\line(0,-1){3}}
\put(10.75,0){\line(0,-1){1}}
\put(11.75,0){\line(0,-1){1}}

\put(9.75,1){\line(-2,-1){2}}
\put(9.75,1){\line(-1,-1){1}}
\put(9.75,1){\line(1,-1){1}}
\put(9.75,1){\line(2,-1){2}}

\put(9.75,-2){\line(-2,1){2}}
\put(9.75,-2){\line(-1,1){1}}
\put(9.75,-2){\line(1,1){1}}
\put(9.75,-2){\line(2,1){2}}

\put(9.75,-1){\oval(4,3)[br]}
\put(10.25,-1){\oval(5,4)[b]}
\put(8.5,-2.5){\line(-1,2){0.75}}
\put(8.5,-2.5){\line(1,0){1.25}}

\put(7.75,0){\line(1,-1){1}}
\put(8.75,0){\line(1,-1){1}}
\put(9.75,0){\line(1,-1){1}}
\put(10.75,0){\line(1,-1){1}}
\put(11.75,0){\line(1,-1){1}}

\put(7.75,0){\circle*{0.1}}
\put(8.75,0){\circle*{0.1}}
\put(9.75,0){\circle*{0.1}}
\put(10.75,0){\circle*{0.1}}
\put(11.75,0){\circle*{0.1}}

\put(7.75,-1){\circle*{0.1}}
\put(8.75,-1){\circle*{0.1}}
\put(9.75,-1){\circle*{0.1}}
\put(10.75,-1){\circle*{0.1}}
\put(11.75,-1){\circle*{0.1}}

\put(9.75,1){\circle*{0.1}}
\put(9.75,-2){\circle*{0.1}}

\end{picture}
\caption{The graph $L_{5}$.} \label{L5}
\end{figure} 
$L_{D_{2}}$ is a planar graph where all vertices have degree~$D_{2}$ except for exactly one vertex with degree $D_{2}-1$.
Thus, $H^{*}$ can be constructed by taking a copy of $H$
and attaching $D_{2} - \deg (v)$ copies of $L_{D_{2}}$ 
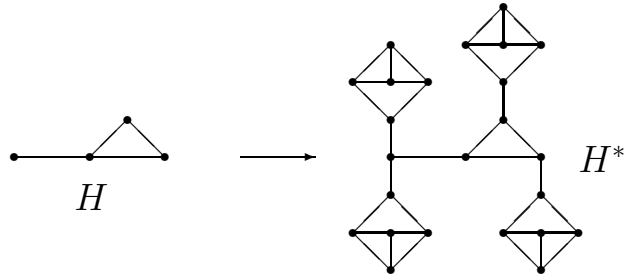
\begin{figure} [ht] 
\setlength{\unitlength}{1cm}
\begin{picture}(20,3)(-2.2,4)

\put(0,5){\line(1,0){2}}
\put(1.5,5.5){\line(1,-1){0.5}}
\put(1.5,5.5){\line(-1,-1){0.5}}

\put(0,5){\circle*{0.1}}
\put(1,5){\circle*{0.1}}
\put(2,5){\circle*{0.1}}
\put(1.5,5.5){\circle*{0.1}}

\put(0.8,4.3){\Large{$H$}}

\put(3,5){\vector(1,0){1}}

\put(5,5){\line(1,0){2}}
\put(6.5,5.5){\line(1,-1){0.5}}
\put(6.5,5.5){\line(-1,-1){0.5}}

\put(5,5){\circle*{0.1}}
\put(6,5){\circle*{0.1}}
\put(7,5){\circle*{0.1}}
\put(6.5,5.5){\circle*{0.1}}

\put(5,4.5){\line(0,1){1}}
\put(7,4.5){\line(0,1){0.5}}
\put(6.5,5.5){\line(0,1){0.5}}

\put(4.5,6){\line(1,0){1}}
\put(5,6){\line(0,1){0.5}}
\put(4.5,6){\line(1,1){0.5}}
\put(5.5,6){\line(-1,1){0.5}}
\put(4.5,6){\line(1,-1){0.5}}
\put(5,5.5){\line(1,1){0.5}}

\put(4.5,6){\circle*{0.1}}
\put(5,6){\circle*{0.1}}
\put(5.5,6){\circle*{0.1}}
\put(5,6.5){\circle*{0.1}}
\put(5,5.5){\circle*{0.1}}

\put(6,6.5){\line(1,0){1}}
\put(6.5,6.5){\line(0,1){0.5}}
\put(6,6.5){\line(1,1){0.5}}
\put(7,6.5){\line(-1,1){0.5}}
\put(6,6.5){\line(1,-1){0.5}}
\put(6.5,6){\line(1,1){0.5}}

\put(6,6.5){\circle*{0.1}}
\put(6.5,6.5){\circle*{0.1}}
\put(7,6.5){\circle*{0.1}}
\put(6.5,7){\circle*{0.1}}
\put(6.5,6){\circle*{0.1}}

\put(4.5,4){\line(1,0){1}}
\put(5,4){\line(0,-1){0.5}}
\put(4.5,4){\line(1,1){0.5}}
\put(5.5,4){\line(-1,1){0.5}}
\put(4.5,4){\line(1,-1){0.5}}
\put(5,3.5){\line(1,1){0.5}}

\put(4.5,4){\circle*{0.1}}
\put(5,4){\circle*{0.1}}
\put(5.5,4){\circle*{0.1}}
\put(5,4.5){\circle*{0.1}}
\put(5,3.5){\circle*{0.1}}

\put(6.5,4){\line(1,0){1}}
\put(7,4){\line(0,-1){0.5}}
\put(6.5,4){\line(1,1){0.5}}
\put(7.5,4){\line(-1,1){0.5}}
\put(6.5,4){\line(1,-1){0.5}}
\put(7,3.5){\line(1,1){0.5}}

\put(6.5,4){\circle*{0.1}}
\put(7,4){\circle*{0.1}}
\put(7.5,4){\circle*{0.1}}
\put(7,4.5){\circle*{0.1}}
\put(7,3.5){\circle*{0.1}}

\put(7.5,4.8){\Large{$H^{*}$}}

\end{picture}

\caption{Constructing a $3$-regular planar graph $H^{*} \supset H$.} \label{attachfig}
\end{figure}
to each vertex $v \in V(H)$ (see Figure~\ref{attachfig}).
$\phantom{qwerty}
\setlength{\unitlength}{.25cm}
\begin{picture}(1,1)
\put(0,0){\line(1,0){1}}
\put(0,0){\line(0,1){1}}
\put(1,1){\line(-1,0){1}}
\put(1,1){\line(0,-1){1}}
\end{picture}$ \\

\begin{figure} [ht] 
\setlength{\unitlength}{1cm}
\begin{picture}(10,2)(-5,0)

\put(0,0){\line(1,0){2}}
\put(1,0.666667){\line(0,1){1.333333}}
\put(0,0){\line(3,2){1}}
\put(0,0){\line(3,4){1}}
\put(0,0){\line(1,2){1}}
\put(2,0){\line(-3,2){1}}
\put(2,0){\line(-3,4){1}}
\put(2,0){\line(-1,2){1}}

\put(0,0){\circle*{0.1}}
\put(2,0){\circle*{0.1}}
\put(1,0.666667){\circle*{0.1}}
\put(1,1.333333){\circle*{0.1}}
\put(1,2){\circle*{0.1}}

\put(-0.1,-0.3){$x$}
\put(2,-0.3){$y$}
\put(0.9,0.366667){$w$}
\put(1.08,1.3){$v$}
\put(1,2.1){$u$}

\end{picture}
\caption{The unique planar embedding of $K_{5} \setminus \{ u,w \}$.} \label{PH}
\end{figure}
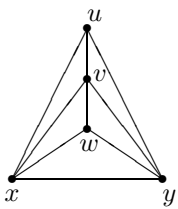

If $D_{2}=4$,
the following example shows that matters are more interesting:

\begin{Example} \label{bounded908}
No $4\textrm{-regular}$ planar graph contains a copy of the graph $K_{5}$ minus an edge.
\end{Example}
\textbf{Proof}
The graph $K_{5} \setminus \{ u,w \}$ is drawn with its \textit{unique} planar embedding (see~\cite{whi})
in Figure~\ref{PH}.

Consider any planar graph $G \supset K_{5} \setminus \{ u,w \}$ with $\Delta (G) = 4$.
Since we already have 
$\textrm{deg}_{H}(v) = \textrm{deg}_{H}(x) = \textrm{deg}_{H}(y) = 4$,
any new edge with at least one endpoint inside the triangle given by $vxy$
must have both endpoints inside.
Hence, the sum of degrees inside this triangle must remain odd,
and so this region must still contain a vertex of odd degree.
Thus, $G$ is not $4$-regular.
$\phantom{qwerty}
\setlength{\unitlength}{.25cm}
\begin{picture}(1,1)
\put(0,0){\line(1,0){1}}
\put(0,0){\line(0,1){1}}
\put(1,1){\line(-1,0){1}}
\put(1,1){\line(0,-1){1}}
\end{picture}$ \\
\\
\\

It is, in fact, possible to determine algorithmically whether or not a given graph $H$
can't be contained within any $4$-regular planar graphs,
since the following result (which is part of a joint paper with Louigi Addario-Berry \cite{add})
shows that it suffices just to check all $4$-regular planar \textit{multigraphs} 
with the same set of vertices as $H$:

\begin{Theorem} \label{bounded991}
Given a simple planar graph $H$ with $\Delta (H) \leq 4$,
there exists a $4$-regular simple planar graph $G \supset H$ 
if and only if there exists a $4$-regular planar multigraph $G^{\prime} \supset H$ 
with $V(G^{\prime}) = V(H)$.
\end{Theorem}
\textbf{Proof}
Suppose there exists a $4$-regular simple planar graph $G$ such that $H \subset G$
and let $G^{\prime}$ be a \textit{minimal} $4$-regular planar multigraph with $H \subset G^{\prime}$,
in the sense that $|V(G^{\prime}) \setminus V(H)|$ is as small as possible.

Suppose $|V(G^{\prime}) \setminus V(H)| \neq 0$
(hoping to obtain a contradiction)
and let $v \in V(G^{\prime}) \setminus V(H)$.
We shall show that we can obtain a $4$-regular planar multigraph $G^{*}$
such that~$H \subset~G^{*}$ and $V(G^{*}) = V(G^{\prime}) \setminus v$,
thus obtaining our desired contradiction: \\
Case (a): If $v$ has two loops to itself,
then we may simply take $G^{*}$ to be $G^{\prime} \setminus v$. \\
Case (b): If $v$ has exactly one loop to itself
and its other two neighbours are $v_{1}$ and $v_{2}$
(where we allow the possibility that $v_{1}=v_{2}$),
then we may take $G^{*}$ to be $G^{\prime} \setminus v + \{ (v_{1},v_{2}) \}$. \\
Case (c): If $v$ has no loops to itself,
then fix a plane drawing of $G^{\prime}$
and let~$e_{1},e_{2},e_{3}$ and $e_{4}$ be the edges incident to $v$ 
\textit{in clockwise order} in this drawing.
Let $v_{1},v_{2},v_{3}$ and $v_{4}$, respectively, denote the other endpoints of~$e_{1},e_{2},e_{3}$ and~$e_{4}$
(allowing the possibility that $v_{i}=v_{j}$ for some $i$ and $j$).
Then we may take $G^{*}$
to be $G^{\prime} \setminus v + \{ (v_{1},v_{2}), (v_{3},v_{4}) \}$,
since this can also be drawn in the plane.

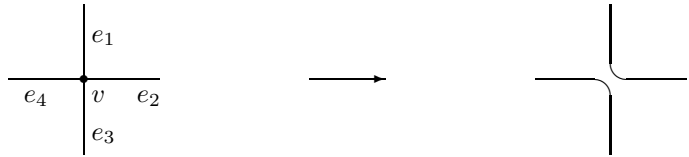
\begin{figure} [ht] 
\setlength{\unitlength}{1cm}
\begin{picture}(10,1.25)(-1.5,0.5)

\put(0,1){\line(1,0){2}}
\put(1,0){\line(0,1){2}}
\put(1,1){\circle*{0.1}}
\put(1.1,1.5){$e_{1}$}
\put(1.7,0.7){$e_{2}$}
\put(1.1,0.2){$e_{3}$}
\put(0.2,0.7){$e_{4}$}
\put(1.1,0.7){$v$}

\put(4,1){\vector(1,0){1}}

\put(7,1){\line(1,0){0.75}}
\put(8.25,1){\line(1,0){0.75}}
\put(8,0){\line(0,1){0.75}}
\put(8,1.25){\line(0,1){0.75}}
\put(8.25,1.25){\oval(0.5,0.5)[bl]}
\put(7.75,0.75){\oval(0.5,0.5)[tr]}

\end{picture}
\caption{$\textrm{Constructing a smaller $4$-regular planar multigraph in case (c).}$} \label{plane}
\end{figure}

For the converse direction,
suppose there exists a $4$-regular planar multigraph~$G^{\prime}$ such that $H \subset G^{\prime}$.
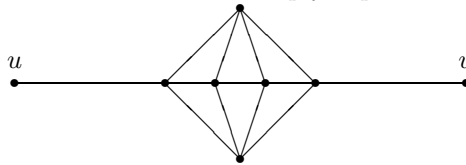
\begin{figure} [ht] 
\setlength{\unitlength}{1cm}
\begin{picture}(10,1.25)(-5,-0.5)

\put(-2,0){\line(1,0){6}}

\put(1,1){\line(-1,-1){1}}
\put(1,1){\line(-1,-3){0.333333}}
\put(1,1){\line(1,-3){0.333333}}
\put(1,1){\line(1,-1){1}}

\put(1,-1){\line(-1,1){1}}
\put(1,-1){\line(-1,3){0.333333}}
\put(1,-1){\line(1,3){0.333333}}
\put(1,-1){\line(1,1){1}}

\put(0,0){\circle*{0.1}}
\put(0.666667,0){\circle*{0.1}}
\put(1.333333,0){\circle*{0.1}}
\put(2,0){\circle*{0.1}}
\put(1,1){\circle*{0.1}}
\put(1,-1){\circle*{0.1}}

\put(-2,0){\circle*{0.1}}
\put(4,0){\circle*{0.1}}

\put(-2.1,0.2){$u$}
\put(3.9,0.2){$v$}

\end{picture}
\caption{$\textrm{Constructing a $4$-regular simple graph from a $4$-regular multigraph.}$}
\label{4regplanar}
\end{figure} 
Then simply replace every edge $e=uv$ of $E(G^{\prime}) \setminus~\!E(H)$~
by a copy of the graph shown in Figure~\ref{4regplanar}.
The resulting graph $G$ will be a $4$-regular simple planar graph with $H \subset G$.
$\phantom{qwerty}
\setlength{\unitlength}{.25cm}
\begin{picture}(1,1)
\put(0,0){\line(1,0){1}}
\put(0,0){\line(0,1){1}}
\put(1,1){\line(-1,0){1}}
\put(1,1){\line(0,-1){1}}
\end{picture}$

For those who are interested,
Section~\ref{4reg} shall be devoted to providing a \textit{polynomial-time} algorithm for determining
whether or not a given graph $H$ can ever be a subgraph of a $4$-regular planar graph. \\
\\
\\

We shall now return to the main thrust of this section
by proving an analogous result to Lemma~\ref{bounded901} 
for the case when $d_{1}(n) < D_{2}(n)$~$\forall n$.
Note that we still have to specify that no component of $H$ is $D_{2}(n)$-regular,
but that now the only other condition involving $H$ is that
$\Delta(H) \leq \liminf_{n \to \infty} D_{2}(n)$:

\begin{Lemma} \label{bounded904}
Let $H$ be a planar graph with components $H_{1}, H_{2}, \ldots, H_{k}$ for some~$\!k$.
Suppose $d_{1}(n)$ and $D_{2}(n)$ are integer-valued functions that for all large $n$ satisfy
(a) $d_{1}(n) < \min \{ 6, D_{2}(n) \}$, and
(b) $D_{2}(n) \geq \max \{ \Delta(H), \max_{i} (\delta(H_{i})+1), 3 \}$.
Then $\exists \beta > 0$ and $\exists N$ such that
\begin{eqnarray*}
& & \mathbf{P} \Big[
\textrm{$P_{n,d_{1},5,0,D_{2}}$
will \emph{not} have a set of $\beta n$ vertex-disjoint} \\
& & \phantom{wwwwwwwwwwq}\textrm{induced order-preserving copies of $H$}\Big] 
< e^{- \beta n}~
\forall n \geq N.
\end{eqnarray*}
\end{Lemma}
\textbf{Proof}
We shall show that there exists a \textit{connected} planar graph $H^{*}$ on 
the vertices $\{ 1,2,\ldots,|H^{*}| \}$
that
(i) contains an induced order-preserving copy of~$H$,
(ii) satisfies 
$\limsup_{n \to \infty} d_{1}(n) \leq \delta(H^{*}) \leq \Delta(H^{*}) \leq \liminf_{n \to \infty} D_{2}(n)$,
and (iii)~satisfies $\deg_{H^{*}}(1) \!<\! \liminf_{n \to \infty} D_{2}(n)$.
We can then use Theorem~\ref{bounded11}~on~$H^{*}$.

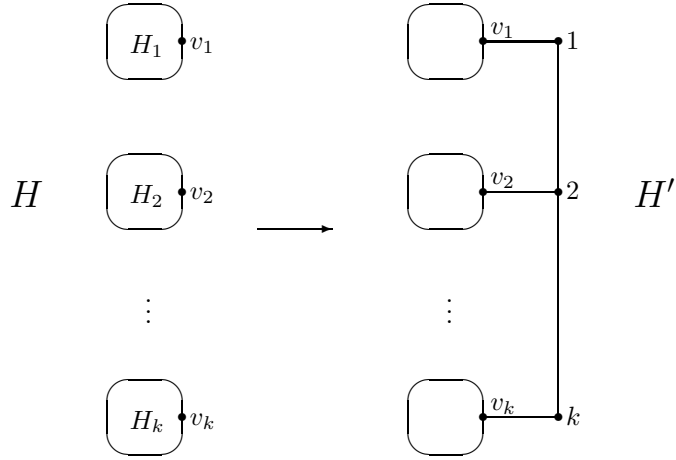
\begin{figure} [ht]
\setlength{\unitlength}{1cm}
\begin{picture}(20,5.75)(-2.9,5.25)

\put(-1.3,8.3){\Large{$H$}}

\put(0.3,5){\line(1,0){0.4}}
\put(0.3,6){\line(1,0){0.4}}
\put(0.3,8){\line(1,0){0.4}}
\put(0.3,9){\line(1,0){0.4}}
\put(0.3,10){\line(1,0){0.4}}
\put(0.3,11){\line(1,0){0.4}}

\put(0.3,5.5){\oval(0.6,1)[l]}
\put(0.3,8.5){\oval(0.6,1)[l]}
\put(0.3,10.5){\oval(0.6,1)[l]}

\put(0.7,5.5){\oval(0.6,1)[r]}
\put(0.7,8.5){\oval(0.6,1)[r]}
\put(0.7,10.5){\oval(0.6,1)[r]}

\put(1,10.5){\circle*{0.1}}
\put(1,8.5){\circle*{0.1}}
\put(1,5.5){\circle*{0.1}}

\put(0.5,6.75){$\vdots$}

\put(0.3,10.35){$H_{1}$}
\put(0.3,8.35){$H_{2}$}
\put(0.3,5.35){$H_{k}$}

\put(1.1,10.4){$v_{1}$}
\put(1.1,8.4){$v_{2}$}
\put(1.1,5.4){$v_{k}$}

\put(2,8){\vector(1,0){1}}

\put(7,8.3){\Large{$H^{\prime}$}}

\put(4.3,6){\line(1,0){0.4}}
\put(4.3,9){\line(1,0){0.4}}
\put(4.3,11){\line(1,0){0.4}}
\put(4.3,5){\line(1,0){0.4}}
\put(4.3,8){\line(1,0){0.4}}
\put(4.3,10){\line(1,0){0.4}}

\put(4.3,5.5){\oval(0.6,1)[l]}
\put(4.3,8.5){\oval(0.6,1)[l]}
\put(4.3,10.5){\oval(0.6,1)[l]}

\put(4.7,5.5){\oval(0.6,1)[r]}
\put(4.7,8.5){\oval(0.6,1)[r]}
\put(4.7,10.5){\oval(0.6,1)[r]}

\put(5,10.5){\circle*{0.1}}
\put(5,8.5){\circle*{0.1}}
\put(5,5.5){\circle*{0.1}}

\put(6,10.5){\circle*{0.1}}
\put(6,8.5){\circle*{0.1}}
\put(6,5.5){\circle*{0.1}}

\put(6,10.5){\line(0,-1){5}}

\put(5,10.5){\line(1,0){1}}
\put(5,8.5){\line(1,0){1}}
\put(5,5.5){\line(1,0){1}}

\put(4.5,6.75){$\vdots$}

\put(5.1,10.6){$v_{1}$}
\put(5.1,8.6){$v_{2}$}
\put(5.1,5.6){$v_{k}$}

\put(6.1,10.4){$1$}
\put(6.1,8.4){$2$}
\put(6.1,5.4){$k$}

\end{picture} 

\caption{Constructing $H^{\prime}$ from $H$.} \label{prime}
\end{figure}

Without loss of generality,
$V(H) = \{ k+1,k+2, \ldots, |H|+k \}$,
where~$\!k =~\!\kappa(H)$.
Let $S = \{ v_{1},v_{2}, \ldots, v_{k} \} \subset V(H)$
be such that we have $v_{i} \in V(H_{i})$ $\forall i$
and deg$(v_{i})<\liminf_{n \to \infty} D_{2}(n)$~$\forall i$,
and let us define $H^{\prime}$ to be the graph with
$V(H^{\prime}) = \{ 1,2, \ldots, |H|+k \}$
and $E(H^{\prime}) = E(H) \cup \bigcup_{i \leq k} (v_{i},i) \cup \bigcup_{i \leq k-1} (i, i+1)$
(see Figure~\ref{prime}).
Then $H^{\prime}$ is a \textit{connected} planar graph that
(i) contains an induced order-preserving copy of $H$,
(ii) satisfies $\Delta(H^{\prime}) \leq \liminf_{n \to \infty} D_{2}(n)$
(since $\deg_{H^{\prime}}(v_{i}) = \deg_{H}(v_{i})+1 \leq \liminf_{n \to \infty} D_{2}(n)$~$\forall i$,
and all the new vertices have degree at most $3 \leq \liminf_{n \to \infty} D_{2}(n)$),
and (iii) satisfies $\deg_{H^{\prime}}(1) <~\!3 \leq \liminf_{n \to \infty} D_{2}(n)$.
Thus, it remains only to extend~$H^{\prime}$ into a graph $H^{*}$ 
that also satisfies $\delta(H^{*}) \geq \limsup_{n \to \infty} d_{1}(n)$.

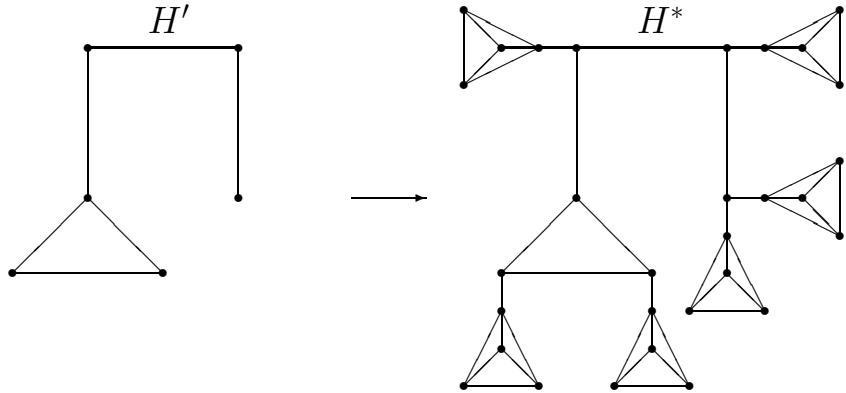
\begin{figure} [ht]
\setlength{\unitlength}{1cm}
\begin{picture}(20,4.75)(-0.5,3.75)

\put(1,8){\line(1,0){2}}
\put(1,6){\line(0,1){2}}
\put(3,6){\line(0,1){2}}
\put(0,5){\line(1,0){2}}
\put(1,6){\line(1,-1){1}}
\put(1,6){\line(-1,-1){1}}

\put(1,8){\circle*{0.1}}
\put(3,8){\circle*{0.1}}
\put(1,6){\circle*{0.1}}
\put(3,6){\circle*{0.1}}
\put(0,5){\circle*{0.1}}
\put(2,5){\circle*{0.1}}

\put(1.8,8.2){\Large{$H^{\prime}$}}

\put(4.5,6){\vector(1,0){1}}

\put(8.3,8.2){\Large{$H^{*}$}}

\put(7,8){\line(1,0){3}}
\put(7.5,6){\line(0,1){2}}
\put(9.5,5.5){\line(0,1){2.5}}
\put(6.5,5){\line(1,0){2}}
\put(7.5,6){\line(1,-1){1}}
\put(7.5,6){\line(-1,-1){1}}
\put(9.5,6){\line(1,0){0.5}}
\put(6.5,4.5){\line(0,1){0.5}}
\put(8.5,4.5){\line(0,1){0.5}}

\put(7.5,8){\circle*{0.1}}
\put(9.5,8){\circle*{0.1}}
\put(7.5,6){\circle*{0.1}}
\put(9.5,6){\circle*{0.1}}
\put(6.5,5){\circle*{0.1}}
\put(8.5,5){\circle*{0.1}}

\put(6,3.5){\line(1,0){1}}
\put(6.5,4){\line(0,1){0.5}}
\put(6,3.5){\line(1,1){0.5}}
\put(6,3.5){\line(1,2){0.5}}
\put(7,3.5){\line(-1,1){0.5}}
\put(7,3.5){\line(-1,2){0.5}}

\put(6.5,4.5){\circle*{0.1}}
\put(6.5,4){\circle*{0.1}}
\put(6,3.5){\circle*{0.1}}
\put(7,3.5){\circle*{0.1}}

\put(8,3.5){\line(1,0){1}}
\put(8.5,4){\line(0,1){0.5}}
\put(8,3.5){\line(1,1){0.5}}
\put(8,3.5){\line(1,2){0.5}}
\put(9,3.5){\line(-1,1){0.5}}
\put(9,3.5){\line(-1,2){0.5}}

\put(8.5,4.5){\circle*{0.1}}
\put(8.5,4){\circle*{0.1}}
\put(8,3.5){\circle*{0.1}}
\put(9,3.5){\circle*{0.1}}

\put(9,4.5){\line(1,0){1}}
\put(9.5,5){\line(0,1){0.5}}
\put(9,4.5){\line(1,1){0.5}}
\put(9,4.5){\line(1,2){0.5}}
\put(10,4.5){\line(-1,1){0.5}}
\put(10,4.5){\line(-1,2){0.5}}

\put(9.5,5.5){\circle*{0.1}}
\put(9.5,5){\circle*{0.1}}
\put(9,4.5){\circle*{0.1}}
\put(10,4.5){\circle*{0.1}}

\put(10,6){\line(1,0){0.5}}
\put(11,5.5){\line(0,1){1}}
\put(10,6){\line(2,1){1}}
\put(10,6){\line(2,-1){1}}
\put(10.5,6){\line(1,1){0.5}}
\put(10.5,6){\line(1,-1){0.5}}

\put(10,6){\circle*{0.1}}
\put(10.5,6){\circle*{0.1}}
\put(11,6.5){\circle*{0.1}}
\put(11,5.5){\circle*{0.1}}

\put(10,8){\line(1,0){0.5}}
\put(11,7.5){\line(0,1){1}}
\put(10,8){\line(2,1){1}}
\put(10,8){\line(2,-1){1}}
\put(10.5,8){\line(1,1){0.5}}
\put(10.5,8){\line(1,-1){0.5}}

\put(10,8){\circle*{0.1}}
\put(10.5,8){\circle*{0.1}}
\put(11,8.5){\circle*{0.1}}
\put(11,7.5){\circle*{0.1}}

\put(6,7.5){\line(0,1){1}}
\put(6.5,8){\line(1,0){0.5}}
\put(6.5,8){\line(-1,1){0.5}}
\put(6.5,8){\line(-1,-1){0.5}}
\put(7,8){\line(-2,1){1}}
\put(7,8){\line(-2,-1){1}}

\put(6,7.5){\circle*{0.1}}
\put(6,8.5){\circle*{0.1}}
\put(6.5,8){\circle*{0.1}}
\put(7,8){\circle*{0.1}}

\end{picture}

\caption{Constructing $H^{*}$ from $H^{\prime}$ in the case $\liminf_{n \to \infty} d_{1}(n) = 3$.} \label{star}
\end{figure}

Let $L$ be a planar connected $d$-regular graph,
where $d= \limsup_{n \to \infty} d_{1}(n)$
(it is clear from Section~\ref{growth} that such a graph must exist).
Then $H^{*}$ can be constructed from $H^{\prime}$ simply by attaching $d - \deg(v)$ copies of $L$
to each vertex $v \in~V(H^{\prime})$ with $\deg_{H^{\prime}}(v)<d$ (see Figure~\ref{star}).
Since $d < \liminf_{n \to \infty} D_{2}(n)$,
we still have
$\Delta(H^{*}) \!\leq~\!\!\liminf_{n \to \infty} D_{2}(n)$
(and $\deg_{H^{*}}(1) \!=\! d \!<\! \liminf_{n \to \infty} D_{2}(n)$).~
$\!\setlength{\unitlength}{.25cm}
\begin{picture}(1,1)
\put(0,0){\line(1,0){1}}
\put(0,0){\line(0,1){1}}
\put(1,1){\line(-1,0){1}}
\put(1,1){\line(0,-1){1}}
\end{picture}$ \\
\\
\\

We may now combine Lemmas~\ref{bounded901},~\ref{bounded907} and~\ref{bounded904}
to obtain our full result:

\begin{Theorem} \label{bounded1001}
Let $H$ be a planar graph with components $H_{1}, \ldots, H_{k}$, for some~$k$.
Suppose $d_{1}(n)$ and $D_{2}(n)$ are integer-valued functions that for all large $n$ satisfy
(a) $d_{1}(n) \leq \min \{ 5, D_{2}(n) \}$,
(b) $D_{2}(n) \geq \max \{ \Delta(H), \max_{i} (\delta(H_{i})+1), 3 \}$,
(c)~$(d_{1}(n),D_{2}(n)) \!\notin\! \{ (3,3),(5,5) \}$ for odd $n$,
and also (d) $(d_{1}(n),D_{2}(n)) \!\neq~\!\!(4,4)$ 
if $H$ happens to be a graph that can never be contained within a $4$-regular planar graph.
Then $\exists \beta > 0$ and $\exists N$ such that
\begin{eqnarray*}
& & \mathbf{P} \Big[
\textrm{$P_{n,d_{1},5,0,D_{2}}$
will \emph{not} have a set of $\beta n$ vertex-disjoint} \\
& & \phantom{wwwwwwwwwwq}\textrm{induced order-preserving copies of $H$}\Big] 
< e^{- \beta n}~
\forall n \geq N. 
\end{eqnarray*}
\end{Theorem}

\begin{displaymath}
\end{displaymath}

It now only remains to complete matters by dealing with the case when $H$ does have $D_{2}(n)$-regular components:

\begin{Theorem} \label{bounded1002}
Let $H$ be a planar graph with components $H_{1}, H_{2}, \ldots, H_{k},$ for some $k$,
and let $D_{2} \geq \max \{ \Delta(H),3 \}$ be a fixed constant equal to $\delta(H_{i})$ for some~$i$.
Suppose $d_{1}(n)$ is an integer-valued function that for all large $n$ satisfies
(a) $d_{1}(n) \leq \min \{ 5,D_{2} \}$,
(b) $(d_{1}(n),D_{2}) \notin \{ (3,3),(5,5) \}$ for odd $n$,
and (c)~$\!(d_{1}(n), D_{2}) \neq~\!(4,4)$ 
if $H$ happens to be a graph that can never be contained within a $4$-regular planar graph.
Then
\begin{displaymath}
\limsup_{n \to \infty} \mathbf{P} [P_{n,d_{1},5,0,D_{2}}
\textrm{ will have a copy of $H$}] < 1,
\end{displaymath}
but for any given constant $t$,
\begin{eqnarray*}
& & \liminf_{n \to \infty} \mathbf{P} \Big[
\textrm{$P_{n,d_{1},5,0,D_{2}}$
will have a set of $t$ vertex-disjoint} \\
& & \phantom{wwwwwwwwwww}\textrm{induced order-preserving copies of $H$}\Big]
>0.
\end{eqnarray*}
\end{Theorem} 
\textbf{Proof}
The first part follows from Theorem~\ref{bounded311},
since for $P_{n,d_{1},5,0,D_{2}}$ to have a copy of $H$
it would have to contain \textit{components} isomorphic to the $D_{2}$-regular components of $H$. 
To prove the second part,
note first that it suffices (by symmetry) to ignore the order-preserving condition.
Without loss of generality,
we may assume that $H_{1},H_{2}, \ldots, H_{l}$
are the $D_{2}$-regular components, for some $l$,
and that they are all distinct.
By Theorem~\ref{bounded404},
we know that
\begin{displaymath}
\liminf_{n \to \infty}
\textrm{$\mathbf{P}\left[ \bigcap_{i \leq l} \left(P_{n,d_{1},5,0,D_{2}} \right.\right.$ 
will have $\geq t$ components isomorphic to $H_{i}$} ) \Big] > 0,
\end{displaymath}
and, by Theorem~\ref{bounded1001},
we know that $\exists \beta > 0$ and $\exists N$ such that
\begin{eqnarray*}
\mathbf{P} [
\textrm{$P_{n,d_{1},5,0,D_{2}}$
\textit{won't} have a set of $\beta n$ vertex-disjoint induced copies of $H^{\prime}$}] \\
< e^{- \beta n} \textrm{ } \forall n \geq N,
\end{eqnarray*}
where $H^{\prime}$ is the graph with components $H_{l+1},H_{l+2}, \ldots, H_{k}$.
But clearly in any graph with at least $t$ components to $H_{i}$~$\forall i \leq l$
and with a set of $\beta n$ vertex-disjoint induced copies of $H^{\prime}$,
we can just fix exactly $t$ components isomorphic to $H_{i}$~$\forall i \leq l$
and then find a set of $t$ vertex-disjoint copies of $H^{\prime}$ in the rest of the graph
(assuming $n$ is sufficiently large).
Thus,
the result follows.
$\phantom{qwerty}
\setlength{\unitlength}{.25cm}
\begin{picture}(1,1)
\put(0,0){\line(1,0){1}}
\put(0,0){\line(0,1){1}}
\put(1,1){\line(-1,0){1}}
\put(1,1){\line(0,-1){1}}
\end{picture}$ \\
\\
\\

\newpage
\section{Subgraphs of \boldmath{$4$}-Regular Graphs} \label{4reg}

\textbf{(Joint work with Louigi Addario-Berry)} \\
\\
We have now completed our picture of $P_{n,d_{1},d_{2},D_{1},D_{2}}$ for the specific case when
$(d_{2}(n),D_{1}(n)) = (5,0)$~$\forall n$.
In Sections~\ref{general} and~\ref{unbounded},
we shall investigate what happens for general functions $d_{2}(n)$ and $D_{1}(n)$,
but first we shall take a break from the main thrust of Part II to 
return to the question of whether or not a given graph $H$ can ever be a subgraph of a $4$-regular planar graph.

We have already seen (in Theorem~\ref{bounded991}) that we can determine this matter algorithmically.
We shall now present a more efficient algorithm,
which has running time $O \left( |H|^{2.5} \right)$
and can be used to find an explicit $4$-regular planar graph $G \supset H$ if such a graph exists
(note that this improved algorithm will not be used anywhere in the rest of this thesis ---
we are interested in it purely for its own sake).

Recall (from Theorem~\ref{bounded991})
that a given simple planar graph~$H$ can be a subgraph of a $4$-regular simple planar graph if and only if
it can be a subgraph of a $4$-regular planar \textit{multigraph}.
Clearly, this second interpretation is just a special case of the more general problem of
determining whether or not a given planar \textit{multigraph} $H$ is a subgraph of some $4$-regular planar multigraph.
Hence, we will actually aim to produce an efficient algorithm for the latter problem.

Before we give the details of the algorithm itself,
we shall first note (in Lemma~\ref{lemma2})
that our problem is straightforward for graphs with a special structure.
We shall then give the algorithm itself,
which will essentially consist of breaking $H$ up
into more and more highly connected pieces until we can apply Lemma~\ref{lemma2}.

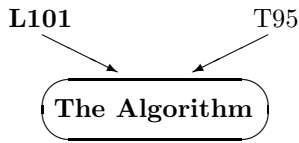
\begin{figure} [ht]
\setlength{\unitlength}{1cm}
\begin{picture}(20,0.9)(-2.5,0.5)

\put(1.55,1.5){\textbf{L\ref{lemma2}}}
\put(4.8,1.5){T\ref{bounded991}}
\put(2.15,0.3){\textbf{The Algorithm}}

\put(3.5,0.4){\oval(3,0.8)}
\put(2,1.4){\vector(2,-1){1}}
\put(5,1.4){\vector(-2,-1){1}}

\end{picture}

\caption{The structure of Section~\ref{4reg}.}
\end{figure}

Before we look at our key lemma,
we shall first pause to meet some important definitions.
The first concerns the idea of `discrepancy functions':

\begin{Definition}
Given a planar multigraph $H$,
we say
$f_{H}:V(H) \to \mathbf{N}$
is a \emph{\textbf{discrepancy function}} on $H$
if (a) $f_{H}(v) \leq 4-\deg_{H}(v)$~$\forall v \in V(H)$
(we call this \emph{\textbf{the discrepancy inequality}})
and (b) $\sum_{v \in V(H)} f_{H}(v)$ is even
(we call this \emph{\textbf{discrepancy parity}}).
If it is also the case that $f_{H}(v) + \deg_{H}(v)$ is even for all $v \in V(H)$,
we call $f_{H}$ an \emph{\textbf{even discrepancy function}} on $H$.

A plane multigraph $G$ \emph{\textbf{satisfies}} $(H,f_{H})$
if $V(G)=V(H)$, $E(G) \supset E(H)$ and $\deg_{G}(v) = \deg_{H}(v) + f_{H}(v)$~$\forall v$.
If such a plane multigraph $G$ exists,
we say that \emph{\boldmath{$f_{H}$}} \emph{\textbf{can be satisfied on}} \emph{\boldmath{$H$}},
or that \emph{\textbf{\boldmath{$(H,f_{H})$} can be satisfied}}.
\end{Definition}

Thus, using Theorem~\ref{bounded991},
it will suffice for us to determine whether or not
the discrepancy function $f_{H}$ defined by setting $f_{H}(v) = 4-\deg_{H}(v)$~$\forall v \in V(H)$
can be satisfied on $H$.
When we break $H$ up into pieces in our algorithm, however,
we shall often find it useful to also define discrepancy functions that are not equal to $4 - \deg$. \\

Our second definition concerns the concept of `augmentations',
which will later play a critical role in our algorithm:

\begin{Definition}
Given a multigraph $B$,
we define the operation of \emph{\textbf{placing a diamond}}
on an edge $uv \in E(B)$
to mean that we subdivide the edge with three vertices 
and then also add two other new vertices so that they are both adjacent to precisely these three vertices.
We define the operation of \emph{\textbf{placing a vertex}}
on an edge $xy \in E(B)$
to mean that we subdivide the edge with a single vertex.

Given multigraphs $B$ and $R$
and a discrepancy function $f_{R}$,
we say that $(R,f_{R})$ is an \emph{\textbf{augmentation}} of $B$ if $R$ can be formed from $B$ 
by placing vertices and diamonds on some of the edges of $B$
(in such a way that there is at most one vertex or diamond on each original edge)
and if $f_{R} = 4 - \deg_{R}$ for all vertices in the new diamonds
and $f_{R} \in \{1,2\}$ for the other new vertices.
\end{Definition}

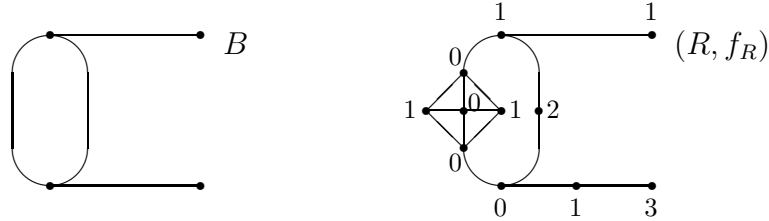
\begin{figure} [ht] 
\setlength{\unitlength}{1cm}
\begin{picture}(10,2.1)(-1.85,0.15)

\put(0,0){\line(1,0){2}}
\put(0,2){\line(1,0){2}}
\put(0,1){\oval(1,2)}

\put(6,0){\line(1,0){2}}
\put(6,2){\line(1,0){2}}
\put(6,1){\oval(1,2)}
\put(5,1){\line(1,1){0.5}}
\put(5,1){\line(1,-1){0.5}}
\put(5,1){\line(1,0){1}}
\put(6,1){\line(-1,1){0.5}}
\put(6,1){\line(-1,-1){0.5}}

\put(0,0){\circle*{0.1}}
\put(0,2){\circle*{0.1}}
\put(2,0){\circle*{0.1}}
\put(2,2){\circle*{0.1}}

\put(6,0){\circle*{0.1}}
\put(6,2){\circle*{0.1}}
\put(7,0){\circle*{0.1}}
\put(8,0){\circle*{0.1}}
\put(8,2){\circle*{0.1}}
\put(6.5,1){\circle*{0.1}}
\put(5,1){\circle*{0.1}}
\put(5.5,1){\circle*{0.1}}
\put(6,1){\circle*{0.1}}
\put(5.5,0.5){\circle*{0.1}}
\put(5.5,1.5){\circle*{0.1}}

\put(2.3,1.75){\large{$B$}}
\put(8.3,1.75){\large{$(R,f_{R})$}}

\put(5.9,-0.4){$0$}
\put(5.9,2.2){$1$}
\put(6.9,-0.4){$1$}
\put(7.9,-0.4){$3$}
\put(7.9,2.2){$1$}
\put(6.6,0.9){$2$}
\put(4.7,0.9){$1$}
\put(5.55,1){$0$}
\put(6.1,0.9){$1$}
\put(5.3,0.2){$0$}
\put(5.3,1.6){$0$}

\end{picture}
\caption{$\textrm{A planar multigraph and an augmentation of it.}$}
\label{D1}
\end{figure}

An example of an augmentation is given in Figure~\ref{D1}.
When we break $H$ up into pieces in our algorithm,
the augmentation of a piece will capture the key information about how it interacted with the rest of $H$.

\phantom{p}

We now come to our key lemma,
which essentially tells us that it will suffice if we can find an algorithm to reduce our problem to
trying to satisfy augmentations of $3$-vertex-connected graphs.

\begin{Lemma} \label{lemma2}
Let $B$ be a planar multigraph of maximum degree at most $4$ that contains no $2$-vertex-cuts
(so $|B| \leq 3$ or $B$ is $3$-vertex-connected),
and let $(R,f_{R})$ be an augmentation of $B$.
Suppose we know which parts of $R$ correspond to which edges of $B$.
Then $\exists \lambda$ such that we can determine in at most $\lambda|B|^{2.5}$ operations
whether or not $(R,f_{R})$ can be satisfied.
\end{Lemma}
\textbf{Proof}
Without loss of generality,
we may assume that $|B|>3$
(since if $B$ is bounded then there are only a finite number of possibilities for $(R,f_{R})$,
and the satisfiability of these can be determined in finite time).
Thus, $3$-vertex-connectivity implies that $B$ has no loops.
By a result of Whitney~\cite{whi} on $3$-vertex-connected simple graphs,
it then follows that $B$ has a unique planar embedding.
Thus, $R$ will also have a unique embedding,
apart from possibly at places where $B$ has multi-edges.

Note that all vertices in $B$ must have at least $3$ distinct neighbours,
since $B$ does not contain any $2\textrm{-vertex-cuts}$.
Hence (since $\deg_{B}(x) = \deg_{R}(x) \leq 4 -f_{R}(x)$~$\forall x \in V(B)$),
if vertices $u$ and $v$ have a multi-edge between them in~$B$,
then it must be only a double-edge
and it must be that $f_{R}(u) = f_{R}(v) = 0$.
We shall now use this information to find a pair $(R^{\prime},f_{R^{\prime}})$
such that $R^{\prime}$ has a unique planar embedding and
$(R,f_{R})$ can be satisfied if and only if
$(R^{\prime},f_{R^{\prime}})$ can be.

\begin{figure} [ht] 
\setlength{\unitlength}{1cm}
\begin{picture}(10,1.66)(-3.5,0.27)

\put(0,0){\line(0,1){2}}
\put(0.1,0){$A$}

\put(1.5,0){\line(0,1){2}}
\put(1.5,1){\circle*{0.1}}
\put(1.6,0.9){$1$}
\put(1.6,0){$B$}

\put(3,0){\line(0,1){2}}
\put(3,1){\circle*{0.1}}
\put(3.1,0.9){$2$}
\put(3.1,0){$C$}

\put(4.5,0){\line(0,1){2}}
\put(4,1){\line(1,0){1}}
\put(4,1){\line(1,1){0.5}}
\put(4,1){\line(1,-1){0.5}}
\put(5,1){\line(-1,1){0.5}}
\put(5,1){\line(-1,-1){0.5}}

\put(4,1){\circle*{0.1}}
\put(4.5,1){\circle*{0.1}}
\put(5,1){\circle*{0.1}}
\put(4.5,0.5){\circle*{0.1}}
\put(4.5,1.5){\circle*{0.1}}

\put(3.7,0.9){$1$}
\put(4.55,1){$0$}
\put(5.1,0.9){$1$}
\put(4.3,0.2){$0$}
\put(4.3,1.6){$0$}

\put(5.1,0){$D$}

\end{picture}
\caption{$\textrm{Augmented versions of an edge.}$}
\label{D2}
\end{figure}
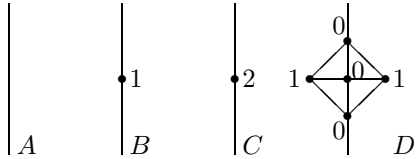

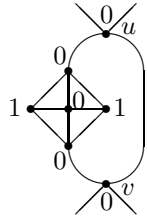
\begin{figure} [ht] 
\setlength{\unitlength}{1cm}
\begin{picture}(10,2.5)(-0.29,-0.13)

\put(6,1){\oval(1,2)}
\put(5,1){\line(1,1){0.5}}
\put(5,1){\line(1,-1){0.5}}
\put(5,1){\line(1,0){1}}
\put(6,1){\line(-1,1){0.5}}
\put(6,1){\line(-1,-1){0.5}}
\put(6,2){\line(1,1){0.4}}
\put(6,2){\line(-1,1){0.4}}
\put(6,0){\line(1,-1){0.4}}
\put(6,0){\line(-1,-1){0.4}}

\put(5,1){\circle*{0.1}}
\put(5.5,1){\circle*{0.1}}
\put(6,1){\circle*{0.1}}
\put(5.5,0.5){\circle*{0.1}}
\put(5.5,1.5){\circle*{0.1}}
\put(6,0){\circle*{0.1}}
\put(6,2){\circle*{0.1}}

\put(4.7,0.9){$1$}
\put(5.55,1){$0$}
\put(6.1,0.9){$1$}
\put(5.3,0.2){$0$}
\put(5.3,1.6){$0$}
\put(5.92,2.15){$0$}
\put(5.92,-0.35){$0$}

\put(6.2,2){$u$}
\put(6.2,-0.15){$v$}

\end{picture}
\caption{$\textrm{A Type A-Type D double edge.}$}
\label{D3}
\end{figure}

Let Type A, Type B, Type C and Type D denote the four possible `augmented versions' of an edge,
as shown in Figure~\ref{D2},
and recall that $R$ will have a unique embedding apart from at any places where $B$ has a double-edge.
If there exist vertices $u$ and $v$ with a Type A-Type D double-edge between them,
then it can be seen that it is impossible to satisfy $(R,f_{R})$,
since $f(u)=f(v)=0$ (see Figure~\ref{D3}).
If we have no Type A-Type D double-edges,
then let $R^{\prime}$ be formed from $R$ as follows: \\
(i) If the augmented versions of a double-edge are Type A and Type B,
then delete the Type A part; \\
(ii) If A and C, delete C; \\
(iii) If B and C, delete C; \\
(iv) If B and D, delete D. \\
(v) If C and D, delete C. \\
Let $f_{R^{\prime}}(v) = f_{R}(v)$~$\forall v \in V(R^{\prime})$.

Using the fact that the two ends of any double-edge must have $f_{R} = 0$,
it is easy to see that $(R,f_{R})$ can be satisfied if and only if $(R^{\prime},f_{R}^{\prime})$ can be satisfied.
It is also clear that $R^{\prime}$ will have a unique embedding.
Thus, to determine whether or not $(R^{\prime},f_{R}^{\prime})$ can be satisfied,
it suffices to see if we can satisfy $(R^{\prime},f_{R^{\prime}})$ in this embedding.

\begin{figure} [ht] 
\setlength{\unitlength}{1cm}
\begin{picture}(10,1.66)(-1.5,0.27)

\put(0,1){\line(1,0){2}}
\put(1,0){\line(0,1){2}}

\put(4,1){\vector(1,0){1}}

\put(7,1){\line(1,0){0.75}}
\put(8.25,1){\line(1,0){0.75}}
\put(8,0){\line(0,1){0.75}}
\put(8,1.25){\line(0,1){0.75}}
\put(8.25,1.25){\oval(0.5,0.5)[bl]}
\put(7.75,0.75){\oval(0.5,0.5)[tr]}

\end{picture}
\caption{$\textrm{Separating crossing edges of our matching.}$}
\label{D15}
\end{figure}
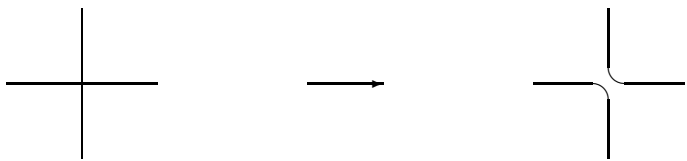

We shall now show that we can reduce this latter problem 
to finding a perfect matching in a suitably defined `auxiliary' graph.
We define this auxiliary graph
(which will not necessarily be planar)
to consist of the vertices of $R^{\prime}$
\textit{with a vertex $x$ appearing $f_{R}(x)$ times}
and with edges between two vertices if and only if 
our embedding has a face containing both of them.
If $(R^{\prime},f_{R}^{\prime})$ can be satisfied in our embedding,
say by a graph $M$,
then the edges in $M$ that are not edges of $R^{\prime}$ form a perfect matching in the auxiliary graph.
Conversely, if we can find a perfect matching,
then inserting the edges of this matching into our embedding will give us a (not necessarily plane) multigraph
satisfying $(R^{\prime},f_{R^{\prime}})$,
which can then be made into a plane multigraph satisfying $(R^{\prime},f_{R}^{\prime})$
simply by separating any crossing edges of our matching,
as in Figure~\ref{D15}.

$R^{\prime}$ can be generated from $R$ in $O \left( |B|^{2} \right)$ time,
since there are $O \left( |B|^{2} \right)$ pairs of vertices in $B$ to check for double-edges
and we know which parts of $R$ correspond to which edges of $B$.
Note that $|R^{\prime}|=O(|B|)$,
so a planar embedding of $R^{\prime}$ can then be found in $O(|B|)$ time (see~\cite{boo})
and we can obtain the auxiliary graph in $O \left( |B|^{2} \right)$
(since there are $O \left( |B|^{2} \right)$ possible edges).
This auxiliary graph will also have $O(|B|)$ vertices,
so we can then determine whether or not it has a perfect matching in $O \left( |B|^{2.5} \right)$ time (see~\cite{eve}).
$\phantom{qwerty}
\setlength{\unitlength}{.25cm}
\begin{picture}(1,1)
\put(0,0){\line(1,0){1}}
\put(0,0){\line(0,1){1}}
\put(1,1){\line(-1,0){1}}
\put(1,1){\line(0,-1){1}}
\end{picture}$ \\

Before we come to our algorithm,
let us first adopt one final convenient definition:

\begin{Definition}
We say that a planar multigraph is \emph{\textbf{4-embeddable}}
if it is a subgraph of some $4$-regular planar multigraph.
\end{Definition}

Thus, the aim of this section is to produce an efficient algorithm to determine whether or not
a given multigraph $H$ is $4$-embeddable. \\
\\

We now present our algorithm. 
We shall first provide a short sketch,
before then giving the details in full.
Afterwards,
we will investigate $\textrm{the running time}$. \\
\\
\textbf{Sketch of Algorithm} \\
The algorithm shall consist of four stages,
each of which will involve breaking $H$ up into more highly connected pieces,
until we can eventually apply Lemma~\ref{lemma2} to all of these.

We will start,
in Stages 1 and 2,
by straightforwardly showing that $H$ is $4$-embeddable
if and only if all its $2$-edge-connected components are.

In Stage 3, 
we will then break our $2$-edge-connected components into $2$-\textit{vertex}-connected blocks,
and show that the discrepancy function $f=4-\deg$ can be satisfied on our $2$-edge-connected components
if and only if certain specified discrepancy functions can be satisfied on all the $2$-vertex-connected blocks.

Stage 4 is where we will use the notion of augmentations.
We shall split our $2$-vertex-connected blocks into $3$-vertex-connected multigraphs
and define augmented versions of each of these.
There will be different cases depending on exactly how the $2$-vertex-cuts break up the graph,
and we will show that the discrepancy functions defined on our $2$-vertex-connected blocks can be satisfied
if and only if all these augmentations can be satisfied.
This can then be determined using Lemma~\ref{lemma2}.

\subsection*{FULL ALGORITHM}

\subsubsection*{STAGE 1}

Clearly, there exists a $4$-regular planar multigraph $G \supset H$ if and only if
there exist $4$-regular planar multigraphs $G_{i} \supset H_{i}$ for all components $H_{i}$ of $H$
(the `if' direction follows by taking $G$ to be the graph whose components are the $G_{i}$'s 
and the `only if' direction follows by taking $G_{i}=G$~$\forall i$).

Thus, the first stage of our algorithm will be to split $H$ into its components.

\subsubsection*{STAGE 2}

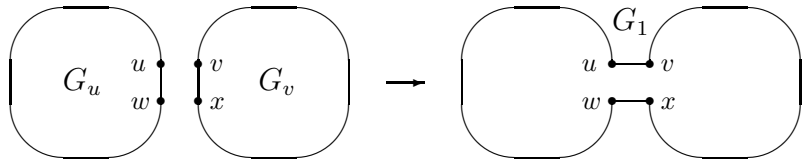
\begin{figure} [ht] 
\setlength{\unitlength}{1cm}
\begin{picture}(10,2)(-0.75,0)

\put(1,1){\oval(2,2)}
\put(3.5,1){\oval(2,2)}
\put(7,1.25){\oval(2,1.5)[t]}
\put(7,0.75){\oval(2,1.5)[b]}
\put(9.5,1.25){\oval(2,1.5)[t]}
\put(9.5,0.75){\oval(2,1.5)[b]}

\put(5,1){\vector(1,0){0.5}}

\put(8,0.75){\line(1,0){0.5}}
\put(8,1.25){\line(1,0){0.5}}
\put(6,0.75){\line(0,1){0.5}}
\put(10.5,0.75){\line(0,1){0.5}}

\put(2,0.75){\circle*{0.1}}
\put(2,1.25){\circle*{0.1}}
\put(2.5,0.75){\circle*{0.1}}
\put(2.5,1.25){\circle*{0.1}}
\put(8,0.75){\circle*{0.1}}
\put(8,1.25){\circle*{0.1}}
\put(8.5,0.75){\circle*{0.1}}
\put(8.5,1.25){\circle*{0.1}}

\put(1.6,1.15){$u$}
\put(1.6,0.65){$w$}
\put(2.65,1.15){$v$}
\put(2.65,0.65){$x$}

\put(7.6,1.15){$u$}
\put(7.6,0.65){$w$}
\put(8.65,1.15){$v$}
\put(8.65,0.65){$x$}

\put(0.7,0.9){\large{$G_{u}$}}
\put(3.3,0.9){\large{$G_{v}$}}
\put(8,1.7){\large{$G_{1}$}}

\end{picture}
\caption{$\textrm{Constructing a $4$-regular planar multigraph $G_{1}$ 
from $4$-regular planar}$}
\textrm{multigraphs $G_{u}$ and $G_{v}$.}
\label{D4}
\end{figure}

Let $H_{1}$ be a component of $H$ and suppose that $H_{1}$ has a cut-edge $e=uv$.
Let $H_{u}$ and $H_{v}$ denote the components of $H_{1} \setminus e$ containing $u$ and $v$, respectively.
Clearly, there exists a $4$-regular planar multigraph $G_{1} \supset H_{1}$ only if
there exist $4$-regular planar multigraphs $G_{u} \supset H_{u}$ and $G_{v} \supset H_{v}$
(this follows by taking $G_{u}=G_{v}=G_{1}$).
We shall now see that the converse is also true:

Suppose there exist $4$-regular planar multigraphs $G_{u} \supset H_{u}$ and $G_{v} \supset H_{v}$.
Note that $\deg_{H_{u}}(u) = \deg_{H_{1}}(u)-1 \leq 3$,
since $v \notin V(H_{u})$,
so $\exists w \in V(G_{u})$ such that
$uw \in E(G_{u}) \setminus E(H_{u})$.
Similarly,
$\exists x \in V(G_{v})$ such that
$vx \in E(G_{v}) \setminus E(H_{v})$.
Since $G_{u}$ and $G_{v}$ are both planar,
they can be drawn with the edges $uw$ and $vx$, respectively, in the outside face.
Thus, the graph $G_{1}$ formed by deleting these two edges and inserting edges $uv$ and $wx$ will also be planar,
as well as being a $4$-regular multigraph containing $H_{1}$ (see Figure~\ref{D4}).

We have shown that $H_{1}$ is $4$-embeddable if and only if $H_{u}$ and $H_{v}$ both are.
Thus, by repeated use of this result,
we find that $H_{1}$ is $4$-embeddable if and only if
all its $2$-edge-connected components are
(counting an isolated vertex as $2$-edge-connected).

Therefore, 
the second stage of our algorithm will be to split the components of $H$ into their $2$-edge-connected components.

\subsubsection*{STAGE 3}

Let $A$ be one of our $2$-edge-connected components.
We wish to determine whether or not there exists a $4$-regular planar multigraph $G_{A} \supset A$.
By Theorem~\ref{bounded991},
it suffices to discover whether or not there exists a $4$-regular planar multigraph $G^{\prime} \supset~\!A$
with $V(G^{\prime}) = V(A)$,
i.e.~to determine whether or not we can satisfy the even discrepancy function on $A$ given by
$f_{A}(v) =~\!4-~\!\deg_{A}(v)$~$\forall v \in~\!V(A)$.

Suppose that $A$ has a cut-vertex $v$.
Since $A$ contains no cut-edges,
it must be that $A \setminus v$ consists of exactly two components, $A_{1}$ and $A_{2}$,
with exactly two edges from $v$ to each of these components.
Thus, $\deg_{A}(v)=4$ and $f_{A}(v)=0$.

\begin{figure} [ht] 
\setlength{\unitlength}{0.75cm}
\begin{picture}(10,1.6)(-0.3,0.3)

\put(1,1){\oval(2,2)}
\put(4,1){\oval(2,2)}
\put(8,1){\oval(2,2)}
\put(12.5,1){\oval(2,2)}

\put(2,0.75){\line(2,1){1}}
\put(2,1.25){\line(2,-1){1}}
\put(9,0.75){\line(2,1){0.5}}
\put(9,1.25){\line(2,-1){0.5}}
\put(11,1){\line(2,1){0.5}}
\put(11,1){\line(2,-1){0.5}}

\put(5,1.5){$(A,f_{A})$}
\put(9,1.5){$\left( A_{1}^{*},f_{A_{1}^{*}} \right)$}
\put(13.5,1.5){$\left( A_{2}^{*},f_{A_{2}^{*}} \right)$}

\put(2.5,1){\circle*{0.1333}}
\put(9.5,1){\circle*{0.1333}}
\put(11,1){\circle*{0.1333}}

\put(2,0.75){\circle*{0.1333}}
\put(2,1.25){\circle*{0.1333}}
\put(3,0.75){\circle*{0.1333}}
\put(3,1.25){\circle*{0.1333}}

\put(9,0.75){\circle*{0.1333}}
\put(9,1.25){\circle*{0.1333}}
\put(11.5,0.75){\circle*{0.1333}}
\put(11.5,1.25){\circle*{0.1333}}

\put(2.4,0.5){$0$}
\put(2.4,1.2){$v$}
\put(9.4,0.5){$0$}
\put(9.65,0.9){$v$}
\put(10.9,0.5){$0$}
\put(10.65,0.9){$v$}

\put(0.75,0.9){$A_{1}$}
\put(3.75,0.9){$A_{2}$}

\end{picture}
\caption{$\textrm{The planar multigraphs $A,A_{1}^{*}$ and $A_{2}^{*}$.}$}
\label{D5}
\end{figure}
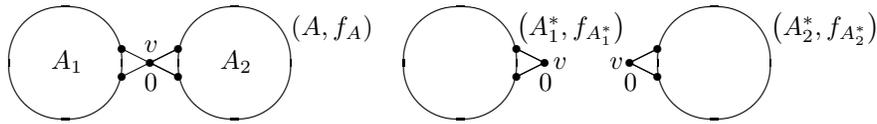

Let $A_{1}^{*}$ denote the planar multigraph induced by $V(A_{1}) \cup v$ 
and let $f_{A_{1}^{*}}$ denote the even discrepancy function on $A_{1}^{*}$ given by
$f_{A_{1}^{*}}(x) = f_{A}(x)$~$\forall x \in V(A_{1}^{*})$
(note $f_{A_{1}^{*}}(x) + \deg_{A_{1}^{*}}(x) = 4$~$\forall x \neq v$
and $f_{A_{1}^{*}}(v) + \deg_{A_{1}^{*}}(v) = 2$,
so $f_{A_{1}^{*}}$ is indeed an even discrepancy function).
Let $A_{2}^{*}$ and $f_{A_{2}^{*}} \textrm{be defined similarly (see Figure~\ref{D5}).}$

\begin{figure} [ht] 
\setlength{\unitlength}{1cm}
\begin{picture}(10,3.1)(-0.5,-0.7)

\put(0.75,0.75){\oval(1.5,1.5)}
\put(3.25,0.75){\oval(1.5,1.5)}

\put(1.5,0.5){\line(2,1){1}}
\put(1.5,1){\line(2,-1){1}}

\put(4.2,2.25){\large{$G^{\prime}$}}

\put(2,0.75){\circle*{0.1}}

\put(1.9,0.45){$v$}

\put(1,1.5){\line(1,0){2}}
\put(1,0){\line(1,0){2}}
\put(2,1.5){\oval(3,1)[t]}
\put(2,0.75){\oval(4,3.5)}
\put(2,0){\oval(3,1)[b]}

\put(0.6,0.6){\large{$A_{1}$}}
\put(3.1,0.6){\large{$A_{2}$}}

\put(0.5,1.5){\circle*{0.1}}
\put(1,1.5){\circle*{0.1}}
\put(0,0.5){\circle*{0.1}}
\put(0,1){\circle*{0.1}}
\put(0.5,0){\circle*{0.1}}
\put(1,0){\circle*{0.1}}
\put(1.5,0.5){\circle*{0.1}}
\put(1.5,1){\circle*{0.1}}

\put(3,1.5){\circle*{0.1}}
\put(3.5,1.5){\circle*{0.1}}
\put(2.5,0.5){\circle*{0.1}}
\put(2.5,1){\circle*{0.1}}
\put(3,0){\circle*{0.1}}
\put(3.5,0){\circle*{0.1}}
\put(4,0.5){\circle*{0.1}}
\put(4,1){\circle*{0.1}}

\put(5,0.75){\vector(1,0){1}}

\put(7.75,0.75){\oval(1.5,1.5)}
\put(10.25,0.75){\oval(1.5,1.5)}

\put(8.5,0.5){\line(2,1){1}}
\put(8.5,1){\line(2,-1){1}}

\put(11.2,1.5){\large{$G^{*}$}}

\put(9,0.75){\circle*{0.1}}

\put(7.5,1.5){\circle*{0.1}}
\put(8,1.5){\circle*{0.1}}
\put(7,0.5){\circle*{0.1}}
\put(7,1){\circle*{0.1}}
\put(7.5,0){\circle*{0.1}}
\put(8,0){\circle*{0.1}}
\put(8.5,0.5){\circle*{0.1}}
\put(8.5,1){\circle*{0.1}}

\put(10,1.5){\circle*{0.1}}
\put(10.5,1.5){\circle*{0.1}}
\put(9.5,0.5){\circle*{0.1}}
\put(9.5,1){\circle*{0.1}}
\put(10,0){\circle*{0.1}}
\put(10.5,0){\circle*{0.1}}
\put(11,0.5){\circle*{0.1}}
\put(11,1){\circle*{0.1}}

\put(8.9,0.45){$v$}

\put(7.75,1.5){\oval(0.5,0.5)[t]}
\put(7,0.75){\oval(0.5,0.5)[l]}
\put(7.75,0){\oval(0.5,0.5)[b]}

\put(10.25,1.5){\oval(0.5,0.5)[t]}
\put(11,0.75){\oval(0.5,0.5)[r]}
\put(10.25,0){\oval(0.5,0.5)[b]}

\end{picture}
\caption{$\textrm{Constructing the graph $G^{*}$ from $G^{\prime}$.}$}
\label{D6}
\end{figure}
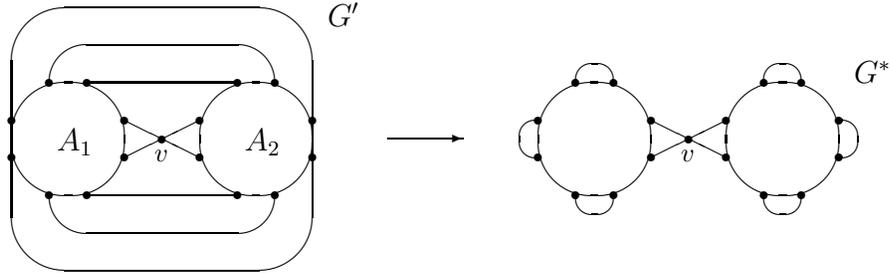

Clearly, we can satisfy $(A,f_{A})$ if we can satisfy both 
$\left( A_{1}^{*},f_{A_{1}^{*}} \right)$ and $\left( A_{2}^{*},f_{A_{2}^{*}} \right)$
(since if there exist plane multigraphs $G_{1}^{*}$ and $G_{2}^{*}$
satisfying $\left( A_{1}^{*},f_{A_{1}^{*}} \right)$ and $\left( A_{2}^{*},f_{A_{2}^{*}} \right)$, respectively,
then we may assume that $v$ is in the outside face of both of these,
and so we can then `glue' these two drawings together at $v$ to obtain a plane multigraph that satisfies $(A,f_{A})$).
We shall now see that the converse is also true:

Suppose $(A,f_{A})$ can be satisfied,
i.e.~there exists a plane multigraph $G^{\prime} \supset A$
with $V(G^{\prime})=V(A)$ and $\deg_{G^{\prime}}(x) = 4$~$\forall x$.
Let us consider the induced plane drawing of $A$.
Since $A_{2}$ is connected,
it must lie in a single face of~$A_{1}^{*}$.
Thus, we may assume that our plane drawing of $A$ is as shown in Figure~\ref{D5},
where without loss of generality we have drawn $A_{2}$ in the outside face of~$A_{1}^{*}$.
Note that the set of edges in $E(G^{\prime}) \setminus E(A)$
between $A_{1}$ and $A_{2}$ must all lie in a single face of our plane drawing
and that there must be an even number of such edges,
since $f_{A}$ is an even discrepancy function and $f_{A}(v)=0$.
Thus, we may `pair up' these edges,
as in Figure~\ref{D6},
to obtain a plane multigraph~$G^{*}$ satisfying $(f_{A},A)$ that has \textit{no} edges from $A_{1}$ to $A_{2}$.
It is then clear that $G_{1}^{*} = G^{*} \setminus A_{2}$ and $G_{2}^{*} = G^{*} \setminus A_{1}$
will satisfy $\left( A_{1}^{*},f_{A_{1}^{*}} \right)$ and $\left( A_{2}^{*},f_{A_{2}^{*}} \right)$,~respectively.

Thus, we have shown that the even discrepancy function $f_{A}$ can be satisfied on $A$ if and only if
the even discrepancy functions $f_{A_{1}^{*}}$ and $f_{A_{2}^{*}}$ 
can be satisfied on $A_{1}^{*}$ and $A_{2}^{*}$, respectively.
By repeatedly using this result, 
we may obtain a set of discrepancy functions 
defined on $2$-\textit{vertex}-connected planar multigraphs
such that $(A,f_{A})$ can be satisfied if and only if all these can be satisfied.

Therefore, the third stage of our algorithm will be to split our $2$-edge-connected components 
into $2$-vertex-connected blocks
(the decomposition is, in fact, unique),
and give each the appropriate discrepancy function.

\subsubsection*{STAGE 4}

Let $C$ be one of our $2$-vertex-connected blocks.
We wish to determine whether or not $(C,f_{C})$ can be satisfied. 
Analogously to Stages 1-3,
we shall split $C$ up into pieces at $2$-vertex-cuts.
However, unlike with these earlier stages,
this time if there exists a graph $M$ satisfying $(C,f_{C})$ 
there may be several different possibilities for how the edges of $M$ could interact with these pieces.
To keep track of this,
we shall define augmentations of the pieces
in such a way that $(C,f_{C})$ can be satisfied if and only if these augmentations can all be satisfied.

We will proceed iteratively.
At the start of each iteration,
we shall have a `blue' graph
(which will initially be $C$)
and an augmentation of it
(initially $(C,f_{C})$)
for which we want to determine satisfiability.
We will split our blue graph in two at a $2$-vertex-cut by breaking off a $3$-vertex-connected piece,
and we shall define augmentations of these two pieces
(in such a way that the augmentation of the blue graph can be satisfied if and only if
the augmentations of the pieces can).
Lemma~\ref{lemma2} can then be used to determine satisfiability of the augmentation of the $3$-vertex-connected piece,
while the other piece and its augmentation can be used as the inputs for the next iteration.
The iterative loop terminates when the blue graph is itself $3$-vertex-connected.

We shall now give the full details:

\subsubsection*{Initialising}

Let us define our initial `blue graph', $B$, to be $C$,
let us also define our initial `red graph', $R$, to be $C$,
and let $R$ have discrepancy function $f_{R}=f_{C}$.
Note that $(R,f_{R})$ is an augmentation of $B$.
At the start of each iteration,
we will always have a blue planar multigraph with no cut-vertex,
and an augmentation of this consisting of a red graph and a discrepancy function.

\subsubsection*{The Iterative Loop}

Check if $B$ has any $2$-vertex cuts.
If not,
then we are done, 
since we can simply use Lemma~\ref{lemma2}.
Otherwise, let us find a minimal $2$-vertex-cut $\{ u,v \}$,
where we use `minimal' to mean that the component of smallest order in $B \setminus \{ u,v \}$ is minimal
over all possible $2$-vertex-cuts.

We shall now proceed to define several graphs based on the pieces of $B \setminus~\{ u,v \}$.
Let $B_{1}$ denote a component of smallest order in $B \setminus \{ u,v \}$,
let $B_{1}^{*}$ denote the graph induced by $V(B_{1}) \cup \{ u,v \}$
and let $B_{1}^{\dag}$ denote the graph obtained from $B_{1}^{*}$ by deleting any edges from $u$ to~$v$.
Let $B_{2} = B \setminus B_{1}^{*}$,
let $B_{2}^{*} = B \setminus B_{1}$
and let $B_{2}^{\dag}$ denote the graph obtained from $B_{2}^{*}$ by deleting any edges from $u$ to $v$.
Let $R_{1}^{*},R_{2}^{*},R_{1}^{\dag}$ and $R_{2}^{\dag}$,
respectively,
denote the red versions of $B_{1}^{*},B_{2}^{*},B_{1}^{\dag}$ and~$B_{2}^{\dag}$
that follow `naturally' from $R$,
and let $R_{1} = R \setminus R_{2}^{*}$
and $R_{2} = R \setminus R_{1}^{*}$.

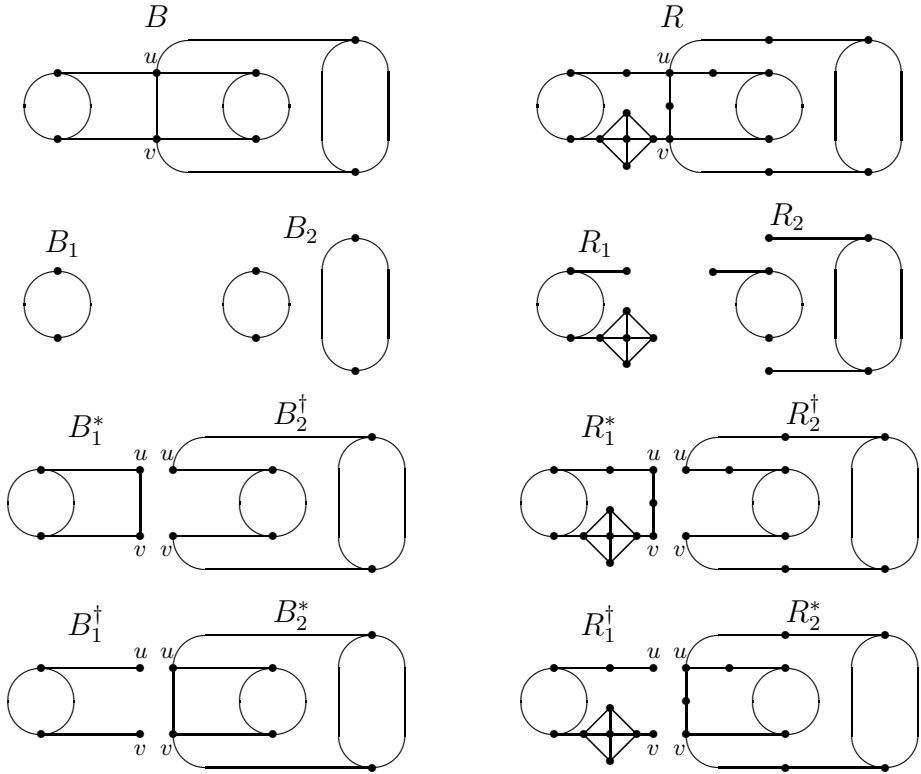
\begin{figure} [ht] 
\setlength{\unitlength}{0.88cm}
\begin{picture}(10,11)(0.75,-9.1)

\put(1.5,0.5){\oval(1,1)}
\put(4.5,0.5){\oval(1,1)}
\put(6,0.5){\oval(1,2)}
\put(1.5,0){\line(1,0){3}}
\put(1.5,1){\line(1,0){3}}
\put(3,0){\line(0,1){1}}
\put(6,1){\oval(6,1)[tl]}
\put(6,0){\oval(6,1)[bl]}
\put(3,0){\circle*{0.1}}
\put(3,1){\circle*{0.1}}
\put(2.8,1.15){\small{$u$}}
\put(2.8,-0.3){\small{$v$}}
\put(2.8,1.7){\large{$B$}}
\put(1.5,0){\circle*{0.1}}
\put(1.5,1){\circle*{0.1}}
\put(4.5,0){\circle*{0.1}}
\put(4.5,1){\circle*{0.1}}
\put(6,-0.5){\circle*{0.1}}
\put(6,1.5){\circle*{0.1}}

\put(1.5,-2.5){\oval(1,1)}
\put(4.5,-2.5){\oval(1,1)}
\put(6,-2.5){\oval(1,2)}
\put(1.3,-1.7){\large{$B_{1}$}}
\put(4.9,-1.5){\large{$B_{2}$}}
\put(1.5,-3){\circle*{0.1}}
\put(1.5,-2){\circle*{0.1}}
\put(4.5,-3){\circle*{0.1}}
\put(4.5,-2){\circle*{0.1}}
\put(6,-3.5){\circle*{0.1}}
\put(6,-1.5){\circle*{0.1}}

\put(9.25,0.5){\oval(1,1)}
\put(12.25,0.5){\oval(1,1)}
\put(13.75,0.5){\oval(1,2)}
\put(9.25,0){\line(1,0){3}}
\put(9.25,1){\line(1,0){3}}
\put(10.75,0){\line(0,1){1}}
\put(13.75,1){\oval(6,1)[tl]}
\put(13.75,0){\oval(6,1)[bl]}
\put(10.75,0){\circle*{0.1}}
\put(10.75,1){\circle*{0.1}}
\put(10.55,1.15){\small{$u$}}
\put(10.55,-0.3){\small{$v$}}
\put(10.6,1.7){\large{$R$}}
\put(10.1,1){\circle*{0.1}}
\put(10.1,-0.4){\line(0,1){0.8}}
\put(10.1,-0.4){\line(1,1){0.4}}
\put(10.1,-0.4){\line(-1,1){0.4}}
\put(10.1,0.4){\line(1,-1){0.4}}
\put(10.1,0.4){\line(-1,-1){0.4}}
\put(9.7,0){\circle*{0.1}}
\put(10.1,0){\circle*{0.1}}
\put(10.5,0){\circle*{0.1}}
\put(10.1,0.4){\circle*{0.1}}
\put(10.1,-0.4){\circle*{0.1}}
\put(10.75,0.5){\circle*{0.1}}
\put(12.25,1.5){\circle*{0.1}}
\put(12.25,-0.5){\circle*{0.1}}
\put(11.4,1){\circle*{0.1}}
\put(9.25,0){\circle*{0.1}}
\put(9.25,1){\circle*{0.1}}
\put(12.25,0){\circle*{0.1}}
\put(12.25,1){\circle*{0.1}}
\put(13.75,-0.5){\circle*{0.1}}
\put(13.75,1.5){\circle*{0.1}}

\put(9.25,-2.5){\oval(1,1)}
\put(12.25,-2.5){\oval(1,1)}
\put(13.75,-2.5){\oval(1,2)}
\put(9.35,-1.7){\large{$R_{1}$}}
\put(12.25,-1.3){\large{$R_{2}$}}
\put(10.1,-2){\circle*{0.1}}
\put(9.25,-2){\line(1,0){0.85}}
\put(9.25,-3){\line(1,0){1.25}}
\put(10.1,-3.4){\line(0,1){0.8}}
\put(10.1,-3.4){\line(1,1){0.4}}
\put(10.1,-3.4){\line(-1,1){0.4}}
\put(10.1,-2.6){\line(1,-1){0.4}}
\put(10.1,-2.6){\line(-1,-1){0.4}}
\put(9.7,-3){\circle*{0.1}}
\put(10.1,-3){\circle*{0.1}}
\put(10.5,-3){\circle*{0.1}}
\put(10.1,-2.6){\circle*{0.1}}
\put(10.1,-3.4){\circle*{0.1}}
\put(12.25,-1.5){\circle*{0.1}}
\put(12.25,-3.5){\circle*{0.1}}
\put(11.4,-2){\circle*{0.1}}
\put(11.4,-2){\line(1,0){0.85}}
\put(12.25,-1.5){\line(1,0){1.5}}
\put(12.25,-3.5){\line(1,0){1.5}}
\put(9.25,-3){\circle*{0.1}}
\put(9.25,-2){\circle*{0.1}}
\put(12.25,-3){\circle*{0.1}}
\put(12.25,-2){\circle*{0.1}}
\put(13.75,-3.5){\circle*{0.1}}
\put(13.75,-1.5){\circle*{0.1}}

\put(1.25,-5.5){\oval(1,1)}
\put(4.75,-5.5){\oval(1,1)}
\put(6.25,-5.5){\oval(1,2)}
\put(1.25,-6){\line(1,0){1.5}}
\put(1.25,-5){\line(1,0){1.5}}
\put(3.25,-6){\line(1,0){1.5}}
\put(3.25,-5){\line(1,0){1.5}}
\put(3.25,-6){\circle*{0.1}}
\put(3.25,-5){\circle*{0.1}}
\put(2.75,-6){\line(0,1){1}}
\put(6.25,-5){\oval(6,1)[tl]}
\put(6.25,-6){\oval(6,1)[bl]}
\put(2.75,-6){\circle*{0.1}}
\put(2.75,-5){\circle*{0.1}}
\put(3.05,-4.85){\small{$u$}}
\put(3.05,-6.3){\small{$v$}}
\put(2.65,-4.85){\small{$u$}}
\put(2.65,-6.3){\small{$v$}}
\put(4.75,-4.3){\large{$B_{2}^{\dag}$}}
\put(1.65,-4.5){\large{$B_{1}^{*}$}}
\put(1.25,-6){\circle*{0.1}}
\put(1.25,-5){\circle*{0.1}}
\put(4.75,-6){\circle*{0.1}}
\put(4.75,-5){\circle*{0.1}}
\put(6.25,-6.5){\circle*{0.1}}
\put(6.25,-4.5){\circle*{0.1}}

\put(9,-5.5){\oval(1,1)}
\put(12.5,-5.5){\oval(1,1)}
\put(14,-5.5){\oval(1,2)}
\put(9,-6){\line(1,0){1.5}}
\put(9,-5){\line(1,0){1.5}}
\put(11,-6){\line(1,0){1.5}}
\put(11,-5){\line(1,0){1.5}}
\put(10.5,-6){\line(0,1){1}}
\put(14,-5){\oval(6,1)[tl]}
\put(14,-6){\oval(6,1)[bl]}
\put(10.5,-6){\circle*{0.1}}
\put(10.5,-5){\circle*{0.1}}
\put(11,-6){\circle*{0.1}}
\put(11,-5){\circle*{0.1}}
\put(10.4,-4.85){\small{$u$}}
\put(10.4,-6.3){\small{$v$}}
\put(10.8,-4.85){\small{$u$}}
\put(10.8,-6.3){\small{$v$}}
\put(12.5,-4.3){\large{$R_{2}^{\dag}$}}
\put(9.85,-5){\circle*{0.1}}
\put(9.85,-6.4){\line(0,1){0.8}}
\put(9.85,-6.4){\line(1,1){0.4}}
\put(9.85,-6.4){\line(-1,1){0.4}}
\put(9.85,-5.6){\line(1,-1){0.4}}
\put(9.85,-5.6){\line(-1,-1){0.4}}
\put(9.45,-6){\circle*{0.1}}
\put(9.85,-6){\circle*{0.1}}
\put(10.25,-6){\circle*{0.1}}
\put(9.85,-5.6){\circle*{0.1}}
\put(9.85,-6.4){\circle*{0.1}}
\put(10.5,-5.5){\circle*{0.1}}
\put(12.5,-4.5){\circle*{0.1}}
\put(12.5,-6.5){\circle*{0.1}}
\put(11.65,-5){\circle*{0.1}}
\put(9.4,-4.5){\large{$R_{1}^{*}$}}
\put(9,-6){\circle*{0.1}}
\put(9,-5){\circle*{0.1}}
\put(12.5,-6){\circle*{0.1}}
\put(12.5,-5){\circle*{0.1}}
\put(14,-6.5){\circle*{0.1}}
\put(14,-4.5){\circle*{0.1}}

\put(9,-8.5){\oval(1,1)}
\put(12.5,-8.5){\oval(1,1)}
\put(14,-8.5){\oval(1,2)}
\put(9,-9){\line(1,0){1.5}}
\put(9,-8){\line(1,0){1.5}}
\put(11,-9){\line(1,0){1.5}}
\put(11,-8){\line(1,0){1.5}}
\put(11,-9){\line(0,1){1}}
\put(14,-8){\oval(6,1)[tl]}
\put(14,-9){\oval(6,1)[bl]}
\put(10.5,-9){\circle*{0.1}}
\put(10.5,-8){\circle*{0.1}}
\put(11,-9){\circle*{0.1}}
\put(11,-8){\circle*{0.1}}
\put(10.4,-7.85){\small{$u$}}
\put(10.4,-9.3){\small{$v$}}
\put(10.8,-7.85){\small{$u$}}
\put(10.8,-9.3){\small{$v$}}
\put(12.5,-7.3){\large{$R_{2}^{*}$}}
\put(9.85,-8){\circle*{0.1}}
\put(9.85,-9.4){\line(0,1){0.8}}
\put(9.85,-9.4){\line(1,1){0.4}}
\put(9.85,-9.4){\line(-1,1){0.4}}
\put(9.85,-8.6){\line(1,-1){0.4}}
\put(9.85,-8.6){\line(-1,-1){0.4}}
\put(9.45,-9){\circle*{0.1}}
\put(9.85,-9){\circle*{0.1}}
\put(10.25,-9){\circle*{0.1}}
\put(9.85,-8.6){\circle*{0.1}}
\put(9.85,-9.4){\circle*{0.1}}
\put(11,-8.5){\circle*{0.1}}
\put(12.5,-7.5){\circle*{0.1}}
\put(12.5,-9.5){\circle*{0.1}}
\put(11.65,-8){\circle*{0.1}}
\put(9.4,-7.5){\large{$R_{1}^{\dag}$}}
\put(9,-9){\circle*{0.1}}
\put(9,-8){\circle*{0.1}}
\put(12.5,-9){\circle*{0.1}}
\put(12.5,-8){\circle*{0.1}}
\put(14,-9.5){\circle*{0.1}}
\put(14,-7.5){\circle*{0.1}}

\put(1.25,-8.5){\oval(1,1)}
\put(4.75,-8.5){\oval(1,1)}
\put(6.25,-8.5){\oval(1,2)}
\put(1.25,-9){\line(1,0){1.5}}
\put(1.25,-8){\line(1,0){1.5}}
\put(3.25,-9){\line(1,0){1.5}}
\put(3.25,-8){\line(1,0){1.5}}
\put(3.25,-9){\circle*{0.1}}
\put(3.25,-8){\circle*{0.1}}
\put(3.25,-9){\line(0,1){1}}
\put(6.25,-8){\oval(6,1)[tl]}
\put(6.25,-9){\oval(6,1)[bl]}
\put(2.75,-9){\circle*{0.1}}
\put(2.75,-8){\circle*{0.1}}
\put(3.05,-7.85){\small{$u$}}
\put(3.05,-9.3){\small{$v$}}
\put(2.65,-7.85){\small{$u$}}
\put(2.65,-9.3){\small{$v$}}
\put(4.75,-7.3){\large{$B_{2}^{*}$}}
\put(1.65,-7.5){\large{$B_{1}^{\dag}$}}
\put(1.25,-9){\circle*{0.1}}
\put(1.25,-8){\circle*{0.1}}
\put(4.75,-9){\circle*{0.1}}
\put(4.75,-8){\circle*{0.1}}
\put(6.25,-9.5){\circle*{0.1}}
\put(6.25,-7.5){\circle*{0.1}}

\end{picture}
\caption{$\textrm{The planar multigraphs defined in the iterative loop of Stage 4.}$}
\label{D7}
\end{figure}

Let $u1$ denote the statement 
\begin{quote}
`$f_{R}(u)=0$ or there is only one edge in $B$ from $u$ to $B_{1}$' 
\end{quote}
(note that the latter implies $|B_{1}|=1$, by minimality,
but that it is not equivalent to this,
as we may have multi-edges).
It is important to note that the number of edges in $B$ from $u$ or $v$ to $B_{1}$
is exactly the same as the number of edges in $R$ from $u$ or $v$, respectively, to $R_{1}$
(and similarly for $B_{2}$ and $R_{2}$).
Thus, it would be equivalent to define $u1$ to denote
`$f_{R}(u)=0$ or there is only one edge in $R$ from $u$ to $R_{1}$'.
Let $v1$ denote the analogous statement to $u1$ for $v$,
and let $u2$ and $v2$ denote the analogous statements for $B_{2}$.
Let $\overline{u1},\overline{v1},\overline{u2}$ and $\overline{v2}$
denote the complements of $u1,v1,u2$ and $v2$. 

Recall that we wish to split our graph in two at each iteration.
Note that if we have $\overline{u2}$, for example,
then $f_{R}(u) \geq 1$
and there are at least two edges in $R$ from $u$ to $R_{2}$,
so there may be several possibilities for where a graph satisfying $(R,f_{R})$
could have a new $u-R_{2}$ edge.
This could complicate matters,
causing an exponential blow-up in the running time,
unless we choose to split the graph in such a way that only the edges from $u$ to $R_{1}$
are important to the analysis.
Thus,
our choice of how best to split the graph
depends on which of the statements $u1,v1,u2$ and $v2$ are true,
and hence our next step is to divide our iterative loop
into different cases based on this information.

\subsubsection*{Case (a): \boldmath{$u1 \land v1$}}

We shall now establish a couple of important facts,
before then splitting into two further subcases arising from parity issues.

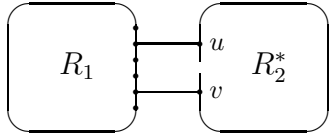
\begin{figure} [ht] 
\setlength{\unitlength}{0.85cm}
\begin{picture}(10,1.2)(-4.6,0.5)

\put(1,1.625){\oval(2,0.75)[t]}
\put(4,1.625){\oval(2,0.75)[t]}
\put(1,0.375){\oval(2,0.75)[b]}
\put(4,0.375){\oval(2,0.75)[b]}
\put(0,0.35){\line(0,1){1.3}}
\put(2,0.35){\line(0,1){1.3}}
\put(3,0.35){\line(0,1){0.55}}
\put(3,1.1){\line(0,1){0.55}}
\put(5,0.35){\line(0,1){1.3}}

\put(2,0.625){\line(1,0){1}}
\put(2,1.375){\line(1,0){1}}

\put(3,0.625){\circle*{0.1}}
\put(3,1.375){\circle*{0.1}}

\put(2,0.875){\circle*{0.1}}
\put(2,0.625){\circle*{0.1}}
\put(2,0.375){\circle*{0.1}}
\put(2,1.625){\circle*{0.1}}
\put(2,1.375){\circle*{0.1}}
\put(2,1.125){\circle*{0.1}}

\put(3.15,1.275){$u$}
\put(3.15,0.525){$v$}
\put(0.8,0.9){\large{$R_{1}$}}
\put(3.8,0.9){\large{$R_{2}^{*}$}}

\put(3,0.625){\line(-4,-1){0.4}}
\put(2.4,0.475){\line(-4,-1){0.4}}
\put(3,0.625){\line(-4,1){0.4}}
\put(2.4,0.775){\line(-4,1){0.4}}
\put(3,1.375){\line(-4,-1){0.4}}
\put(2.4,1.225){\line(-4,-1){0.4}}
\put(3,1.375){\line(-4,1){0.4}}
\put(2.4,1.525){\line(-4,1){0.4}}

\end{picture}
\caption{$\textrm{The planar multigraph $R$.}$}
\label{D8}
\end{figure}

By definition,
$B_{1}$ is connected.
Thus, $R_{1}$ must also be connected,
and so has to lie in a single face of $R_{2}^{*}$.
Hence, in any planar embedding $R$ must look as in Figure~\ref{D8},
where broken lines represent edges that may or may not exist
and where, without loss of generality,
we have drawn $R_{1}$ in the outside face of $R_{2}^{*}$.
Therefore, if a plane multigraph $M$ satisfies $(R,f_{R})$
then all edges in $E(M) \setminus~E(R)$ between $V(R_{1})$ and $V(R_{2}^{*})$
must lie within only two faces of the induced embedding of $R$
(since $u$ can have more than one edge to $R_{1}$ only if $f(u)=0$,
and similarly for $v$).

Secondly,
since $f_{R}$ satisfies discrepancy parity,
note that
$\sum_{x \in V(R_{1})} f_{R}(x)$ and $\sum_{x \in V \left( R_{2}^{*} \right) } f_{R}(x)$ 
must either both be odd or both be even.

\subsubsection*{Case (a)(i): 
\boldmath{$\sum_{x \in V(R_{1})} f_{R}(x)$} and \boldmath{$\sum_{x \in V \left( R_{2}^{*} \right) } f_{R}(x)$} both odd}

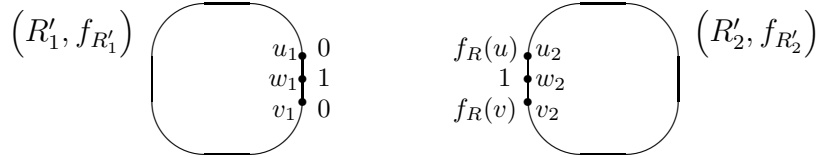
\begin{figure} [ht] 
\setlength{\unitlength}{1cm}
\begin{picture}(10,1.3)(-3.5,0.5)

\put(0,1){\oval(2,2)}
\put(5,1){\oval(2,2)}
\put(1,0.7){\circle*{0.1}}
\put(1,1){\circle*{0.1}}
\put(1,1.3){\circle*{0.1}}
\put(4,0.7){\circle*{0.1}}
\put(4,1){\circle*{0.1}}
\put(4,1.3){\circle*{0.1}}
\put(-2.9,1.5){\large{$\left( R_{1}^{\prime},f_{R_{1}^{\prime}} \right)$}}
\put(6.2,1.5){\large{$\left( R_{2}^{\prime},f_{R_{2}^{\prime}} \right)$}}
\put(0.6,1.3){$u_{1}$}
\put(0.55,0.9){$w_{1}$}
\put(0.6,0.5){$v_{1}$}
\put(4.1,1.3){$u_{2}$}
\put(4.1,0.9){$w_{2}$}
\put(4.1,0.5){$v_{2}$}
\put(1.2,1.3){$0$}
\put(1.2,0.9){$1$}
\put(1.2,0.5){$0$}
\put(3,1.3){$f_{R}(u)$}
\put(3.6,0.9){$1$}
\put(3,0.5){$f_{R}(v)$}

\end{picture}
\caption{$\textrm{The planar multigraphs $R_{1}^{\prime}$ and $R_{2}^{\prime}$, and their discrepancy functions.}$}
\label{D9}
\end{figure}

Let $B_{1}^{\prime} = B_{1}^{\dag} + uv$ and let $B_{2}^{\prime} = B_{2}^{*} + uv$
(so $uv$ will now be a multi-edge in $B_{2}^{\prime}$ if $uv \in E(B)$).
We shall now define an augmentation $(R_{1}^{\prime}, f_{R_{1}^{\prime}})$ of $B_{1}^{\prime}$
and an augmentation $(R_{2}^{\prime}, f_{R_{2}^{\prime}})$ of $B_{2}^{\prime}$
such that $(R,f_{R})$ can be satisfied if and only if
$\left( R_{1}^{\prime},f_{R_{1}^{\prime}} \right)$ and $\left( R_{2}^{\prime},f_{R_{2}^{\prime}} \right)$ 
can both be satisfied
(these new augmentations are illustrated in Figure~\ref{D9}).

Let $R_{1}^{\prime}$ be the graph formed from $R_{1}^{\dag}$
by relabelling $u$ and $v$ as $u_{1}$ and~$v_{1}$, respectively,
and introducing a new vertex $w_{1}$ with edges to both $u_{1}$ and $v_{1}$.
Similarly, let $R_{2}^{\prime}$ be the graph formed from $R_{2}^{*}$
by relabelling $u$ and $v$ as $u_{2}$ and $v_{2}$, respectively,
and introducing a new vertex $w_{2}$ with edges to both $u_{2}$ and $v_{2}$.
Let $f_{R_{1}^{\prime}}$ be the discrepancy function on $R_{1}^{\prime}$ defined by setting
$f_{R_{1}^{\prime}}(u_{1})=f_{R_{1}^{\prime}}(v_{1})=0$,
$f_{R_{1}^{\prime}}(w_{1})=1$, and $f_{R_{1}^{\prime}}(x)=f_{R}(x)$~$\forall x \in V(R_{1})$.
Let $f_{R_{2}^{\prime}}$ be the discrepancy function on $R_{2}^{\prime}$ defined by setting
$f_{R_{2}^{\prime}}(u_{2})=f_{R}(u)$,
$f_{R_{2}^{\prime}}(v_{2})=f_{R}(v)$,
$f_{R_{2}^{\prime}}(w_{2})=~1$, 
and $f_{R_{2}^{\prime}}(x)=f_{R}(x)$~$\forall x \in V(R_{2})$.
(Note that $f_{R_{1}^{\prime}}$ and $f_{R_{2}^{\prime}}$ are both valid discrepancy functions,
since the discrepancy inequality is clearly satisfied by both
and discrepancy parity follows from the fact that
$\sum_{x \in V \left( R_{1}^{\prime} \right) } \!f_{R_{1}^{\prime}}(x) \!=~\!\!\sum_{x \in V(R_{1})} \!f_{R}(x)+~\!1$,
$\sum_{x \in V \left( R_{2}^{\prime} \right) } \!f_{R_{2}^{\prime}}(x) 
\!=~\!\!\sum_{x \in V \left( R_{2}^{*} \right) } \!f_{R}(x)+~\!1$
and 
$\sum_{x \in V(R_{1})} f_{R}(x)$ and $\sum_{x \in V \left( R_{2}^{*} \right) } f_{R}(x)$ are both odd).

\begin{Claim}
$(R,f_{R})$ can be satisfied if and only if
$\left( R_{1}^{\prime},f_{R_{1}^{\prime}} \right)$ and $\left( R_{2}^{\prime},f_{R_{2}^{\prime}} \right)$ 
can both be satisfied.
\end{Claim}
\textbf{Proof}
Suppose first that there exists a plane multigraph $M$ satisfying $(R,f_{R})$.
Since $\sum_{x \in V(R_{1})} f_{R}(x)$ and $\sum_{x \in V \left( R_{2}^{*} \right) } f_{R}(x)$ are both odd,
there must be an odd number of edges in $E(M) \setminus E(R)$ between $V(R_{1})$ and $V(R_{2}^{*})$.
As already noted,
these edges must all lie within two faces of the embedding of $R$ induced from $M$.
Thus, one of these faces must have an odd number of new edges 
and the other must have an even number.
By pairing edges up,
as in the second half of Stage $3$,
we can hence obtain a planar multigraph satisfying $(R,f_{R})$ 
that has exactly one new edge between $V(R_{1})$ and $V(R_{2}^{*})$.
It is then easy $\textrm{to see that we can satisfy both 
$\left( R_{1}^{\prime},f_{R_{1}^{\prime}} \right)$ and $\left( R_{2}^{\prime},f_{R_{2}^{\prime}} \right)$.}$

Suppose now that 
$\left( R_{1}^{\prime},f_{R_{1}^{\prime}} \right)$ and $\left( R_{2}^{\prime},f_{R_{2}^{\prime}} \right)$ 
can both be satisfied,
by plane multigraphs $M_{R_{1}^{\prime}}$ and $M_{R_{2}^{\prime}}$ respectively,
and let the edges adjacent to $w_{1}$ in $E \left( M_{R_{1}^{\prime}} \right) \setminus E(R_{1}^{\prime})$ and 
$w_{2}$ in $E \left( M_{R_{2}^{\prime}} \right) \setminus E(R)$
be denoted by $e_{1}=z_{1}w_{1}$ and $e_{2}=z_{2}w_{2}$ respectively.
We may assume that $e_{1}$ is in the outside face of $M_{R_{1}^{\prime}}$.
Note that the edges $u_{1}w_{1}$ and $v_{1}w_{1}$ must then be in the outside face of $M_{R_{1}^{\prime}} \setminus e_{1}$,
since these are the only edges incident to $w_{1}$ in $M_{R_{1}^{\prime}} \setminus e_{1}$.
Hence,
by turning our drawing upside-down if necessary,
we may assume that $u_{1},w_{1}$ and $v_{1}$ 
are in clockwise order around this outer face of $M_{R_{1}^{\prime}} \setminus e_{1}$,
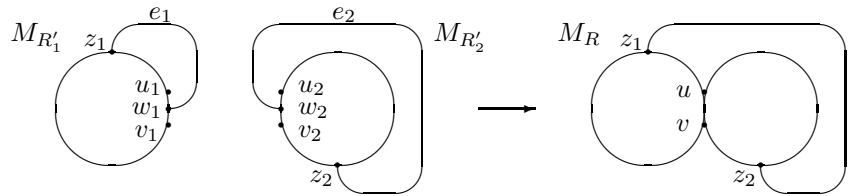
\begin{figure} [ht] 
\setlength{\unitlength}{0.75cm}
\begin{picture}(10,2.2)(-1.5,0.2)

\put(1,1){\oval(2,2)}
\put(5,1){\oval(2,2)}
\put(2,0.7){\circle*{0.1}}
\put(2,1){\circle*{0.1}}
\put(2,1.3){\circle*{0.1}}
\put(4,0.7){\circle*{0.1}}
\put(4,1){\circle*{0.1}}
\put(4,1.3){\circle*{0.1}}
\put(-0.8,2.2){$M_{R_{1}^{\prime}}$}
\put(6.7,2.2){$M_{R_{2}^{\prime}}$}
\put(1.4,1.3){$u_{1}$}
\put(1.35,0.9){$w_{1}$}
\put(1.4,0.5){$v_{1}$}
\put(4.3,1.3){$u_{2}$}
\put(4.3,0.9){$w_{2}$}
\put(4.3,0.5){$v_{2}$}
\put(1,2){\circle*{0.1}}
\put(0.5,2.1){$z_{1}$}
\put(2,1.75){\oval(1,1.5)[r]}
\put(1.75,2){\oval(1.5,1)[tl]}
\put(4,1.75){\oval(1,1.5)[l]}
\put(5.75,0){\oval(1.5,1)[b]}
\put(6,1.75){\oval(1,1.5)[tr]}
\put(4,2.5){\line(1,0){2}}
\put(6.5,0){\line(0,1){2}}
\put(1.75,2.5){\line(1,0){0.25}}
\put(5,0){\circle*{0.1}}
\put(4.5,-0.3){$z_{2}$}
\put(1.65,2.6){$e_{1}$}
\put(4.9,2.6){$e_{2}$}

\put(7.5,1){\vector(1,0){1}}

\put(10.5,1){\oval(2,2)}
\put(11.5,0.7){\circle*{0.1}}
\put(11.5,1.3){\circle*{0.1}}
\put(10.5,2){\circle*{0.1}}
\put(11,1.2){$u$}
\put(11,0.6){$v$}
\put(10,2.1){$z_{1}$}
\put(11.25,2){\oval(1.5,1)[tl]}
\put(12.5,1){\oval(2,2)}
\put(13.25,0){\oval(1.5,1)[b]}
\put(13.5,1.75){\oval(1,1.5)[tr]}
\put(11.25,2.5){\line(1,0){2.25}}
\put(14,0){\line(0,1){2}}
\put(12.5,0){\circle*{0.1}}
\put(12,-0.3){$z_{2}$}
\put(8.9,2.2){$M_{R}$}

\end{picture}
\caption{$\textrm{Constructing a planar multigraph $M_{R}$ satisfying $(R,f_{R})$.}$}
\label{D10}
\end{figure} 
and so $M_{R_{1}^{\prime}}$ is as shown in Figure~\ref{D10}
(where, without loss of generality, 
we have drawn $e_{1}$ so that $v_{1}$ is also in the outside face of $M_{R_{1}^{\prime}}$).
Similarly, we may assume that $M_{R_{2}^{\prime}}$ is also as shown in Figure~\ref{D10}.
It is then clear that we can delete $w_{1}$ and $w_{2}$,
`glue' $u_{1}$ to $u_{2}$ and $v_{1}$ to $v_{2}$
(i.e.~identify $u_{1}$ and $u_{2}$ and, separately, $v_{1}$ and $v_{2}$),
and insert the edge $z_{1}z_{2}$  
to obtain a plane multigraph $M_{R}$ that will satisfy $(R,f_{R})$
(note that it doesn't matter whether $z_{2} \in \{ u_{2},v_{2} \}$).
$\phantom{qwerty}
\setlength{\unitlength}{.25cm}
\begin{picture}(1,1)
\put(0,0){\line(1,0){1}}
\put(0,0){\line(0,1){1}}
\put(1,1){\line(-1,0){1}}
\put(1,1){\line(0,-1){1}}
\end{picture}$ \\

Recall that $B_{1}^{\prime} = B_{1}^{\dag} + uv$ 
and note that $B_{1}^{\prime}$ must not contain any $2$-vertex-cuts, by the minimality of $B_{1}$.
Thus, by Lemma~\ref{lemma2},
in $O \left( |B_{1}^{\prime}|^{2.5} \right)$ time
we can determine whether or not $\left( R_{1}^{\prime},f_{R_{1}^{\prime}} \right)$ can be satisfied.
If it cannot, we terminate the algorithm.
If it can, we return to the start of the iterative loop with $B_{2}^{\prime}$ as our new blue graph,
$R_{2}^{\prime}$ as our new red graph and $f_{R_{2}^{\prime}}$ as our new discrepancy function
(note that, as required,
$B_{2}^{\prime}$ does not contain a cut-vertex
since otherwise this would also be a cut-vertex in $B$
--- this property will be required for case~(b)).

\subsubsection*{Case (a)(ii): 
\boldmath{$\sum_{x \in V(R_{1})} f_{R}(x)$} and \boldmath{$\sum_{x \in V \left( R_{2}^{*} \right) } f_{R}(x)$} both even}

\begin{figure} [ht] 
\setlength{\unitlength}{1cm}
\begin{picture}(10,7)(-3,1)

\put(1,7.625){\oval(2,0.75)[t]}
\put(1,6.375){\oval(2,0.75)[b]}
\put(0,6.35){\line(0,1){1.3}}
\put(2,6.35){\line(0,1){1.3}}
\put(-1.9,7.5){\large{$\left( R_{1}^{\prime},f_{R_{1}^{\prime}} \right)$}}
\put(2,6.375){\circle*{0.1}}
\put(2,7.625){\circle*{0.1}}
\put(1.6,7.6){$u_{1}$}
\put(1.6,6.3){$v_{1}$}
\put(2.2,7.55){\footnotesize{$0$}}
\put(2.2,6.3){\footnotesize{$0$}}

\put(5,7.625){\oval(2,0.75)[t]}
\put(5,6.375){\oval(2,0.75)[b]}
\put(4,6.35){\line(0,1){1.3}}
\put(6,6.35){\line(0,1){1.3}}
\put(6.2,7.5){\large{$\left( R_{2}^{\prime},f_{R_{2}^{\prime}} \right)$}}
\put(4,6.375){\circle*{0.1}}
\put(4,7.625){\circle*{0.1}}
\put(4.2,7.6){$u_{2}$}
\put(4.2,6.3){$v_{2}$}
\put(3.1,7.6){\footnotesize{$f_{R}(u)$}}
\put(3.1,6.3){\footnotesize{$f_{R}(v)$}}

\put(1,5.125){\oval(2,0.75)[t]}
\put(1,3.875){\oval(2,0.75)[b]}
\put(0,3.85){\line(0,1){1.3}}
\put(2,3.85){\line(0,1){1.3}}
\put(-1.9,5){\large{$\left( R_{1}^{\prime\prime},f_{R_{1}^{\prime\prime}} \right)$}}
\put(2,3.875){\circle*{0.1}}
\put(2,5.125){\circle*{0.1}}
\put(1.6,5.1){$u_{1}$}
\put(1.6,3.8){$v_{1}$}
\put(1.6,4.5){\line(1,0){0.8}}
\put(1.6,4.5){\line(1,1){0.4}}
\put(1.6,4.5){\line(1,-1){0.4}}
\put(2.4,4.5){\line(-1,1){0.4}}
\put(2.4,4.5){\line(-1,-1){0.4}}
\put(2,4.5){\circle*{0.1}}
\put(1.6,4.5){\circle*{0.1}}
\put(2.4,4.5){\circle*{0.1}}
\put(2,4.1){\circle*{0.1}}
\put(2,4.9){\circle*{0.1}}
\put(1.3,4.4){\footnotesize{$1$}}
\put(2.5,4.4){\footnotesize{$1$}}
\put(2.1,4.5){\footnotesize{$0$}}
\put(2.1,3.7){\footnotesize{$0$}}
\put(2.1,5.1){\footnotesize{$0$}}
\put(2.1,4){\footnotesize{$0$}}
\put(2.1,4.85){\footnotesize{$0$}}

\put(5,5.125){\oval(2,0.75)[t]}
\put(5,3.875){\oval(2,0.75)[b]}
\put(4,3.85){\line(0,1){1.3}}
\put(6,3.85){\line(0,1){1.3}}
\put(6.2,5){\large{$\left( R_{2}^{\prime\prime},f_{R_{2}^{\prime\prime}} \right)$}}
\put(4,3.875){\circle*{0.1}}
\put(4,5.125){\circle*{0.1}}
\put(4.2,5.1){$u_{2}$}
\put(4.2,3.8){$v_{2}$}
\put(3.6,4.5){\line(1,0){0.8}}
\put(3.6,4.5){\line(1,1){0.4}}
\put(3.6,4.5){\line(1,-1){0.4}}
\put(4.4,4.5){\line(-1,1){0.4}}
\put(4.4,4.5){\line(-1,-1){0.4}}
\put(4,4.5){\circle*{0.1}}
\put(3.6,4.5){\circle*{0.1}}
\put(4.4,4.5){\circle*{0.1}}
\put(4,4.1){\circle*{0.1}}
\put(4,4.9){\circle*{0.1}}
\put(3.3,4.4){\footnotesize{$1$}}
\put(4.5,4.4){\footnotesize{$1$}}
\put(3.75,4.5){\footnotesize{$0$}}
\put(3.15,3.75){\footnotesize{$f_{R}(v)$}}
\put(3.15,5.15){\footnotesize{$f_{R}(u)$}}
\put(3.75,4){\footnotesize{$0$}}
\put(3.75,4.85){\footnotesize{$0$}}

\put(5,2.625){\oval(2,0.75)[t]}
\put(5,1.375){\oval(2,0.75)[b]}
\put(4,1.35){\line(0,1){1.3}}
\put(6,1.35){\line(0,1){1.3}}
\put(6.2,2.5){\large{$\left( R_{2}^{\prime\prime\prime},f_{R_{2}^{\prime\prime\prime}} \right)$}}
\put(4,1.375){\circle*{0.1}}
\put(4,2.625){\circle*{0.1}}
\put(4.2,2.6){$u_{2}$}
\put(4.2,1.3){$v_{2}$}
\put(4,2){\circle*{0.1}}
\put(4.2,1.9){$w$}
\put(3.1,2.6){\footnotesize{$f_{R}(u)$}}
\put(3.1,1.3){\footnotesize{$f_{R}(v)$}}
\put(3.7,1.9){\footnotesize{$2$}}

\end{picture}
\caption{$\textrm{The planar multigraphs 
$R_{1}^{\prime}$, $R_{1}^{\prime\prime}$,
$R_{2}^{\prime}$, $R_{2}^{\prime\prime}$ and $R_{2}^{\prime\prime\prime}$,
and their}$} 
\textrm{\phantom{qqq} discrepancy functions.}
\label{D11}
\end{figure}
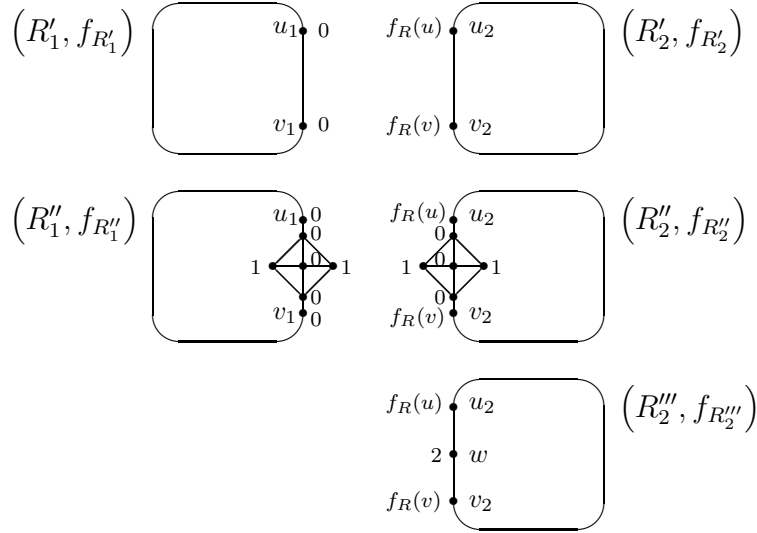

Again, we let $B_{1}^{\prime} = B_{1}^{\dag} + uv$ and $B_{2}^{\prime} = B_{2}^{*} + uv$.
This time,
we shall define augmentations 
$\left( R_{1}^{\prime},f_{R_{1}^{\prime}} \right)$ 
and $\left( R_{1}^{\prime\prime},f_{R_{1}^{\prime\prime}} \right)$ of $B_{1}^{\prime}$
and augmentations 
$\left( R_{2}^{\prime},f_{R_{2}^{\prime}} \right)$, $\left( R_{2}^{\prime\prime},f_{R_{2}^{\prime\prime}} \right)$
and $\left( R_{2}^{\prime\prime\prime},f_{R_{2}^{\prime\prime\prime}} \right)$ of $B_{2}^{\prime}$
(see Figure~\ref{D11})
such that
$(R,f_{R})$ can be satisfied if and only if: 
\begin{quote}
(1) $\left( R_{1}^{\prime},f_{R_{1}^{\prime}} \right)$ and $\left( R_{2}^{\prime},f_{R_{2}^{\prime}} \right)$ 
can both be satisfied, 
but $\left( R_{1}^{\prime\prime},f_{R_{1}^{\prime\prime}} \right)$~can't; \\ 
(2) $\left( R_{1}^{\prime\prime},f_{R_{1}^{\prime\prime}} \right)$ 
and $\left( R_{2}^{\prime\prime},f_{R_{2}^{\prime\prime}} \right)$ 
can both be satisfied, 
but $\left( R_{1}^{\prime},f_{R_{1}^{\prime}} \right)$~can't;~or \\
(3) $\left( R_{1}^{\prime},f_{R_{1}^{\prime}} \right),\left( R_{1}^{\prime\prime},f_{R_{1}^{\prime\prime}} \right)$ 
and $\left( R_{2}^{\prime\prime\prime},f_{R_{2}^{\prime\prime\prime}} \right)$ can all be satisfied. 
\end{quote}

Let $R_{1}^{\prime}$ be the graph formed from $R_{1}^{\dag}$ by relabelling $u$ and $v$ as $u_{1}$ and $v_{1}$,
respectively,
and inserting an edge between $u_{1}$ and $v_{1}$.
Let $f_{R_{1}^{\prime}}$ be the discrepancy function on $R_{1}^{\prime}$ defined by setting
$f_{R_{1}^{\prime}}(u_{1})=f_{R_{1}^{\prime}}(v_{1})=0$
and $f_{R_{1}^{\prime}}(x) = f_{R}(x)$ otherwise.
Let $R_{1}^{\prime\prime}$ be the graph formed from $R_{1}^{\prime}$ by placing a diamond on the $u_{1}v_{1}$ edge,
and let $f_{R_{1}^{\prime\prime}}$ be defined by setting
$f_{R_{1}^{\prime\prime}}(x)=f_{R_{1}^{\prime}}(x)$~$\forall x \in V(R_{1}^{\prime})$
and $f_{R_{1}^{\prime\prime}}(x) = 4 - \deg_{R_{1}^{\prime\prime}}(x)$~$\forall x \notin V(R_{1}^{\prime})$.

Let $R_{2}^{\prime}$ be the graph formed from $R_{2}^{*}$ by relabelling $u$ and $v$ as $u_{2}$ and $v_{2}$,
respectively,
and inserting a new edge between $u_{1}$ and $v_{1}$
(so $u_{1}v_{1}$ will now be a multi-edge if $uv \in E(R)$).
Let $f_{R_{2}^{\prime}} = f_{R_{2}^{*}}$.
Let $R_{2}^{\prime\prime}$ be the graph formed from $R_{2}^{\prime}$ by
placing a diamond on the new $u_{2}v_{2}$ edge,
and let $f_{R_{2}^{\prime\prime}}$ be defined by setting
$f_{R_{2}^{\prime\prime}}(x)=f_{R_{2}^{\prime}}(x)$~$\forall x \in V(R_{2}^{\prime})$
and $f_{R_{2}^{\prime\prime}} = 4 - \deg_{R_{2}^{\prime\prime}}(x)$~$\forall x \notin V(R_{2}^{\prime})$.
Let $R_{2}^{\prime\prime\prime}$ be the graph formed from $R_{2}^{\prime}$
by instead subdividing the new $u_{2}v_{2}$ edge with a vertex $w$,
and let $f_{2}^{\prime\prime\prime}$ be defined by
$f_{2}^{\prime\prime\prime}(w)=2$
and $f_{2}^{\prime\prime\prime}(x) = f_{2}^{\prime}(x)$~$\forall x \in V(R_{2}^{\prime})$.

\begin{Claim}
 $(R,f_{R})$ can be satisfied if and only if
one of (1),(2) or (3) holds.
\end{Claim}
\textbf{Proof}
The `if' direction follows from a similar `gluing' argument to case~(a)(i),
since we can again assume that the appropriate parts of our graphs are drawn in the outside face,
so we shall now proceed with proving the `only if' direction:

Suppose that a plane multigraph $M$ satisfies $(R,f_{R})$.
Since $\sum_{x \in V(R_{1})} f_{R}(x)$ and $\sum_{x \in V \left( R_{2}^{*} \right)} f_{R}(x)$ are both even,
there must be an even number of edges in $E(M) \setminus E(R)$ between $V(R_{1})$ and $V(R_{2}^{*})$.
As in case (a)(i),
these edges must all lie in two faces,
so we must either have an even number in both of these faces or an odd number in both.
By the same argument as with (a)(i),
we may in fact without loss of generality assume that 
there are either no new edges in both faces or exactly one in both.
In the former,
it is clear that we can satisfy both 
$\left( R_{1}^{\prime},f_{R_{1}^{\prime}} \right)$ and $\left( R_{2}^{\prime},f_{R_{2}^{\prime}} \right)$,
and in the latter it is clear that we can satisfy both 
$\left( R_{1}^{\prime\prime},f_{R_{1}^{\prime\prime}} \right)$ 
and $\left( R_{2}^{\prime\prime},f_{R_{2}^{\prime\prime}} \right)$.
Note that we can satisfy $\left( R_{2}^{\prime\prime\prime},f_{R_{2}^{\prime\prime\prime}} \right)$ if 
we can satisfy 
$\left( R_{2}^{\prime},f_{R_{2}^{\prime}} \right)$ or $\left( R_{2}^{\prime\prime},f_{R_{2}^{\prime\prime}} \right)$.
Thus, we can either satisfy 
$\left( R_{1}^{\prime},f_{R_{1}^{\prime}} \right)$,$\left( R_{2}^{\prime},f_{R_{2}^{\prime}} \right)$
and $\left( R_{2}^{\prime\prime\prime},f_{R_{2}^{\prime\prime\prime}} \right)$,
or $\left( R_{1}^{\prime\prime},f_{R_{1}^{\prime\prime}} \right)$,
$\left( R_{2}^{\prime\prime},f_{R_{2}^{\prime\prime}} \right)$ 
and $\left( R_{2}^{\prime\prime\prime},f_{R_{2}^{\prime\prime\prime}} \right)$.
In the first case, either (1) or (3) must hold,
and in the second case either (2) or (3) must hold.
$\phantom{qwerty}
\setlength{\unitlength}{.25cm}
\begin{picture}(1,1)
\put(0,0){\line(1,0){1}}
\put(0,0){\line(0,1){1}}
\put(1,1){\line(-1,0){1}}
\put(1,1){\line(0,-1){1}}
\end{picture}$ \\

We have now shown that $(R,f_{R})$ can be satisfied if and only if (1),(2) or (3) hold.
As in case (a)(i),
we can use Lemma~\ref{lemma2} to determine in $O \left( |B_{1}|^{2.5} \right)$ time whether 
$\left( R_{1}^{\prime},f_{R_{1}^{\prime}} \right)$ and $\left( R_{1}^{\prime\prime},f_{R_{1}^{\prime\prime}} \right)$ 
can be satisfied.
If neither can be satisfied,
we terminate the algorithm.
If at least one can be satisfied,
then we return to the start of the iterative loop with $B_{2}^{\prime}$ as our new blue graph 
and either 
$(R_{2}^{\prime\prime\prime},f_{R_{2}^{\prime\prime\prime}})$,
$(R_{2}^{\prime},f_{R_{2}^{\prime}})$ or $(R_{2}^{\prime\prime},f_{R_{2}^{\prime\prime}})$
as our augmentation,
according to whether both
$\left( R_{1}^{\prime},f_{R_{1}^{\prime}} \right)$ and $\left( R_{1}^{\prime\prime},f_{R_{1}^{\prime\prime}} \right)$,
just $\left( R_{1}^{\prime},f_{R_{1}^{\prime}} \right)$,
or just $\left( R_{1}^{\prime\prime},f_{R_{1}^{\prime\prime}} \right)$
can be satisfied, respectively.

\subsubsection*{Case (b): \boldmath{$(\overline{u1} \lor \overline{v1}) \land u2 \land v2$}}

We shall again start with some groundwork on the structure of $R$,
analogously to case (a),
before splitting into subcases.

Since $\overline{u1} \lor \overline{v1}$ holds,
we can't have $f(u) = f(v) =0$.
Thus, since $u2 \land v2$ also holds,
it must be that either $u$ or $v$ has only one edge to $B_{2}$.
Hence, since $B$ contains no cut-vertices,
it must be that $B_{2}$ is connected.
Therefore, $R_{2}$ must also be connected 
and so must lie in a single face of $R_{1}^{*}$.
Hence, we may proceed in a similar way to case (a),
but this time we will split into subcases depending on the parity of $R_{1}^{*}$ and $R_{2}$,
rather than $R_{1}$ and $R_{2}^{*}$.

\subsubsection*{Case (b)(i): 
\boldmath{$\sum_{x \in V \left( R_{1}^{*} \right) } f_{R}(x)$} and \boldmath{$\sum_{x \in V(R_{2})} f_{R}(x)$} both odd}

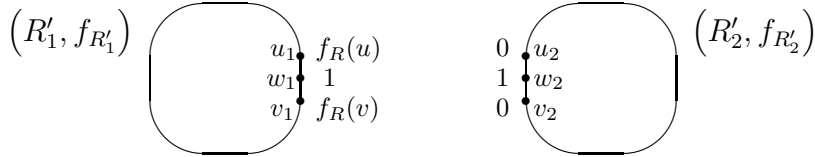
\begin{figure} [ht] 
\setlength{\unitlength}{1cm}
\begin{picture}(10,2)(-3.5,0)

\put(0,1){\oval(2,2)}
\put(5,1){\oval(2,2)}
\put(1,0.7){\circle*{0.1}}
\put(1,1){\circle*{0.1}}
\put(1,1.3){\circle*{0.1}}
\put(4,0.7){\circle*{0.1}}
\put(4,1){\circle*{0.1}}
\put(4,1.3){\circle*{0.1}}
\put(-2.9,1.5){\large{$\left( R_{1}^{\prime},f_{R_{1}^{\prime}} \right)$}}
\put(6.2,1.5){\large{$\left( R_{2}^{\prime},f_{R_{2}^{\prime}} \right)$}}
\put(0.6,1.3){$u_{1}$}
\put(0.55,0.9){$w_{1}$}
\put(0.6,0.5){$v_{1}$}
\put(4.1,1.3){$u_{2}$}
\put(4.1,0.9){$w_{2}$}
\put(4.1,0.5){$v_{2}$}
\put(1.2,1.3){$f_{R}(u)$}
\put(1.3,0.9){$1$}
\put(1.2,0.5){$f_{R}(v)$}
\put(3.6,1.3){$0$}
\put(3.6,0.9){$1$}
\put(3.6,0.5){$0$}

\end{picture}
\caption{$\textrm{The planar multigraphs $R_{1}^{\prime}$ and $R_{2}^{\prime}$, and their discrepancy functions.}$}
\label{D12}
\end{figure}

This time,
we let $B_{1}^{\prime} = B_{1}^{*} + uv$
(so $uv$ will be a multi-edge in $B_{1}^{\prime}$ if $uv \in E(B)$)
and let $B_{2}^{\prime} = B_{2}^{\dag} + uv$.
We will define augmentations
$(R_{1}^{\prime}, f_{R_{1}^{\prime}})$ of $B_{1}^{\prime}$
and $(R_{2}^{\prime}, f_{R_{2}^{\prime}})$ of $B_{2}^{\prime}$
(see Figure~\ref{D12})
such that $(R,f_{R})$ can be satisfied if and only if
$\left( R_{1}^{\prime},f_{R_{1}^{\prime}} \right)$ and $\left( R_{2}^{\prime},f_{R_{2}^{\prime}} \right)$ 
can both be satisfied.

Let $R_{1}^{\prime}$ be the graph formed from $R_{1}^{*}$
by relabelling $u$ and $v$ as $u_{1}$ and $v_{1}$, respectively,
and introducing a new vertex $w_{1}$ with edges to both $u_{1}$ and $v_{1}$.
Similarly, let $R_{2}^{\prime}$ be the graph formed from $R_{2}^{\dag}$
by relabelling $u$ and $v$ as $u_{2}$ and~$v_{2}$, respectively,
and introducing a new vertex $w_{2}$ with edges to both $u_{2}$ and $v_{2}$.
Let $f_{R_{1}^{\prime}}$ be the discrepancy function on $R_{1}^{\prime}$ defined by setting
$f_{R_{1}^{\prime}}(u_{1})=f_{R}(u)$,
$f_{R_{1}^{\prime}}(v_{1})=f_{R}(v)$,
$f_{R_{1}^{\prime}}(w_{1})=1$, 
and $f_{R_{1}^{\prime}}(x)=f_{R}(x)$~$\forall x \in V(R_{1})$.
Let $f_{R_{2}^{\prime}}$ be the discrepancy function on $R_{2}^{\prime}$ defined by setting
$f_{R_{2}^{\prime}}(u_{2})=f_{R_{2}^{\prime}}(v_{2})=0$,
$f_{R_{2}^{\prime}}(w_{2})=1$, and $f_{R_{2}^{\prime}}(x)=f_{R}(x)$~$\forall x \in V(R_{2})$.

The proof that $(R,f_{R})$ may be satisfied if and only if 
both $\left( R_{1}^{\prime},f_{R_{1}^{\prime}} \right)$ and $\left( R_{2}^{\prime},f_{R_{2}^{\prime}} \right)$ 
may be satisfied is as with case (a)(i).
Again,
we can determine in $O \left( |B_{1}^{\prime}|^{2.5} \right)$ time 
whether or not $\left( R_{1}^{\prime},f_{R_{1}^{\prime}} \right)$ can be satisfied,
and if so we return to the start of the iterative loop with $B_{2}^{\prime}$ as our new blue graph,
$R_{2}^{\prime}$ as our new red graph 
and $f_{R_{2}^{\prime}}$ as our new discrepancy function. 
Otherwise, we terminate the algorithm.

\subsubsection*{Case (b)(ii): 
\boldmath{$\sum_{x \in V \left( R_{1}^{*} \right)} f_{R}(x)$} and \boldmath{$\sum_{x \in V(R_{2})} f_{R}(x)$} both even}

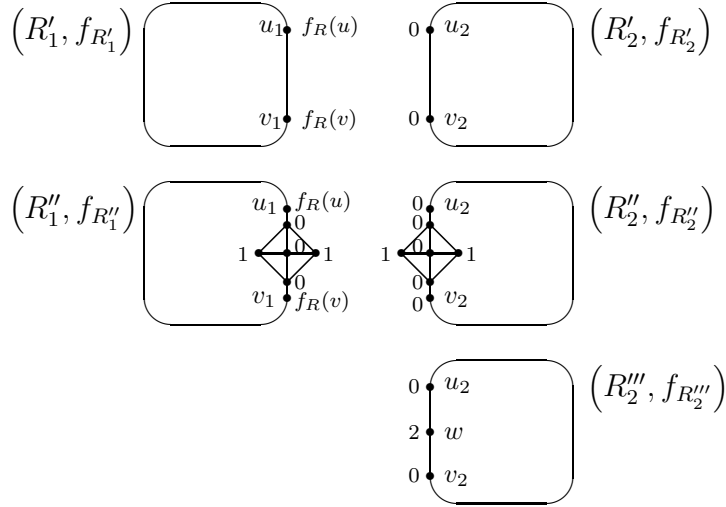
\begin{figure} [ht] 
\setlength{\unitlength}{0.95cm}
\begin{picture}(10,7)(-3,1)

\put(1,7.625){\oval(2,0.75)[t]}
\put(1,6.375){\oval(2,0.75)[b]}
\put(0,6.35){\line(0,1){1.3}}
\put(2,6.35){\line(0,1){1.3}}
\put(-1.9,7.5){\large{$\left( R_{1}^{\prime},f_{R_{1}^{\prime}} \right)$}}
\put(2,6.375){\circle*{0.1}}
\put(2,7.625){\circle*{0.1}}
\put(1.6,7.6){$u_{1}$}
\put(1.6,6.3){$v_{1}$}
\put(2.2,7.6){\footnotesize{$f_{R}(u)$}}
\put(2.2,6.3){\footnotesize{$f_{R}(v)$}}

\put(5,7.625){\oval(2,0.75)[t]}
\put(5,6.375){\oval(2,0.75)[b]}
\put(4,6.35){\line(0,1){1.3}}
\put(6,6.35){\line(0,1){1.3}}
\put(6.2,7.5){\large{$\left( R_{2}^{\prime},f_{R_{2}^{\prime}} \right)$}}
\put(4,6.375){\circle*{0.1}}
\put(4,7.625){\circle*{0.1}}
\put(4.2,7.6){$u_{2}$}
\put(4.2,6.3){$v_{2}$}
\put(3.7,7.55){\footnotesize{$0$}}
\put(3.7,6.3){\footnotesize{$0$}}

\put(1,5.125){\oval(2,0.75)[t]}
\put(1,3.875){\oval(2,0.75)[b]}
\put(0,3.85){\line(0,1){1.3}}
\put(2,3.85){\line(0,1){1.3}}
\put(-1.9,5){\large{$\left( R_{1}^{\prime\prime},f_{R_{1}^{\prime\prime}} \right)$}}
\put(2,3.875){\circle*{0.1}}
\put(2,5.125){\circle*{0.1}}
\put(1.5,5.1){$u_{1}$}
\put(1.5,3.8){$v_{1}$}
\put(1.6,4.5){\line(1,0){0.8}}
\put(1.6,4.5){\line(1,1){0.4}}
\put(1.6,4.5){\line(1,-1){0.4}}
\put(2.4,4.5){\line(-1,1){0.4}}
\put(2.4,4.5){\line(-1,-1){0.4}}
\put(2,4.5){\circle*{0.1}}
\put(1.6,4.5){\circle*{0.1}}
\put(2.4,4.5){\circle*{0.1}}
\put(2,4.1){\circle*{0.1}}
\put(2,4.9){\circle*{0.1}}
\put(1.3,4.4){\footnotesize{$1$}}
\put(2.5,4.4){\footnotesize{$1$}}
\put(2.1,4.5){\footnotesize{$0$}}
\put(2.1,3.75){\footnotesize{$f_{R}(v)$}}
\put(2.1,5.15){\footnotesize{$f_{R}(u)$}}
\put(2.1,4){\footnotesize{$0$}}
\put(2.1,4.85){\footnotesize{$0$}}

\put(5,5.125){\oval(2,0.75)[t]}
\put(5,3.875){\oval(2,0.75)[b]}
\put(4,3.85){\line(0,1){1.3}}
\put(6,3.85){\line(0,1){1.3}}
\put(6.2,5){\large{$\left( R_{2}^{\prime\prime},f_{R_{2}^{\prime\prime}} \right)$}}
\put(4,3.875){\circle*{0.1}}
\put(4,5.125){\circle*{0.1}}
\put(4.2,5.1){$u_{2}$}
\put(4.2,3.8){$v_{2}$}
\put(3.6,4.5){\line(1,0){0.8}}
\put(3.6,4.5){\line(1,1){0.4}}
\put(3.6,4.5){\line(1,-1){0.4}}
\put(4.4,4.5){\line(-1,1){0.4}}
\put(4.4,4.5){\line(-1,-1){0.4}}
\put(4,4.5){\circle*{0.1}}
\put(3.6,4.5){\circle*{0.1}}
\put(4.4,4.5){\circle*{0.1}}
\put(4,4.1){\circle*{0.1}}
\put(4,4.9){\circle*{0.1}}
\put(3.3,4.4){\footnotesize{$1$}}
\put(4.5,4.4){\footnotesize{$1$}}
\put(3.75,4.5){\footnotesize{$0$}}
\put(3.75,3.7){\footnotesize{$0$}}
\put(3.75,5.1){\footnotesize{$0$}}
\put(3.75,4){\footnotesize{$0$}}
\put(3.75,4.85){\footnotesize{$0$}}

\put(5,2.625){\oval(2,0.75)[t]}
\put(5,1.375){\oval(2,0.75)[b]}
\put(4,1.35){\line(0,1){1.3}}
\put(6,1.35){\line(0,1){1.3}}
\put(6.2,2.5){\large{$\left( R_{2}^{\prime\prime\prime},f_{R_{2}^{\prime\prime\prime}} \right)$}}
\put(4,1.375){\circle*{0.1}}
\put(4,2.625){\circle*{0.1}}
\put(4.2,2.6){$u_{2}$}
\put(4.2,1.3){$v_{2}$}
\put(4,2){\circle*{0.1}}
\put(4.2,1.9){$w$}
\put(3.7,2.55){\footnotesize{$0$}}
\put(3.7,1.3){\footnotesize{$0$}}
\put(3.7,1.9){\footnotesize{$2$}}

\end{picture}
\caption{$\textrm{The planar multigraphs 
$R_{1}^{\prime}$, $R_{1}^{\prime\prime}$, $R_{2}^{\prime}$, $R_{2}^{\prime\prime}$
and $R_{2}^{\prime\prime\prime}$, and their}$} 
\textrm{\phantom{qqq} discrepancy functions.}
\label{D13}
\end{figure}

Again, we let $B_{1}^{\prime} = B_{1}^{*} + uv$
and $B_{2}^{\prime} = B_{2}^{\dag} + uv$.
This time, as with case (a)(ii),
we shall define augmentations 
$\left( R_{1}^{\prime},f_{R_{1}^{\prime}} \right)$ 
and $\left( R_{1}^{\prime\prime},f_{R_{1}^{\prime\prime}} \right)$ of $B_{1}^{\prime}$
and augmentations 
$\left( R_{2}^{\prime},f_{R_{2}^{\prime}} \right), \left( R_{2}^{\prime\prime},f_{R_{2}^{\prime\prime}} \right)$
and $\left( R_{2}^{\prime\prime\prime},f_{R_{2}^{\prime\prime\prime}} \right)$ of $B_{2}^{\prime}$
(see Figure~\ref{D13})
such that
$(R,f_{R})$ can be satisfied if and only if: 
\begin{quote}
(1) $\left( R_{1}^{\prime},f_{R_{1}^{\prime}} \right)$ and $\left( R_{2}^{\prime},f_{R_{2}^{\prime}} \right)$ 
can both be satisfied, 
but $\left( R_{1}^{\prime\prime},f_{R_{1}^{\prime\prime}} \right)$~can't; \\ 
(2) $\left( R_{1}^{\prime\prime},f_{R_{1}^{\prime\prime}} \right)$ 
and $\left( R_{2}^{\prime\prime},f_{R_{2}^{\prime\prime}} \right)$ 
can both be satisfied, 
but $\left( R_{1}^{\prime},f_{R_{1}^{\prime}} \right)$~can't;~or \\
(3) $\left( R_{1}^{\prime},f_{R_{1}^{\prime}} \right),\left( R_{1}^{\prime\prime},f_{R_{1}^{\prime\prime}} \right)$ 
and $\left( R_{2}^{\prime\prime\prime},f_{R_{2}^{\prime\prime\prime}} \right)$ can all be satisfied.
\end{quote}

Let $R_{1}^{\prime}$ be the graph formed from $R_{1}^{*}$ by relabelling $u$ and $v$ as $u_{1}$ and~$v_{1}$,
respectively,
and inserting an edge between $u_{1}$ and $v_{1}$
(so $u_{1}v_{1}$ will now be a multi-edge if $uv \in E(R)$).
Let $f_{R_{1}^{\prime}}$ be defined by setting $f_{R_{1}^{\prime}}(x) = f_{R}(x)$~$\forall x \in~\!V(R_{1}^{\prime})$.
Let $R_{1}^{\prime\prime}$ be the graph formed from $R_{1}^{\prime}$ by placing a diamond on the $u_{1}v_{1}$ edge,
and let $f_{R_{1}^{\prime\prime}}$ be the function defined by setting
$f_{R_{1}^{\prime\prime}}(x)=f_{R_{1}^{\prime}}(x)$~$\forall x \in V(R_{1}^{\prime})$
and $f_{R_{1}^{\prime\prime}}(x) = 4-\deg_{R_{1}^{\prime\prime}}(x)$~$\forall x \notin V(R_{1}^{\prime})$.

Let $R_{2}^{\prime}$ be the graph formed from $R_{2}^{\dag}$ by relabelling $u$ and $v$ as $u_{2}$ and~$v_{2}$,
respectively,
and inserting a new edge between $u_{2}$ and $v_{2}$.
Let $f_{R_{2}^{\prime}}$ be the discrepancy function on $R_{2}^{\prime}$ defined by setting
$f_{R_{2}^{\prime}}(u_{2})=f_{R_{2}^{\prime}}(v_{2})=0$
and $f_{R_{2}^{\prime}}(x) =~\!f_{R}(x)$~$\forall x \in V(R_{2})$.
Let $R_{2}^{\prime\prime}$ be the graph formed from $R_{2}^{\prime}$ by
placing a diamond on the new $u_{2}v_{2}$ edge,
and let $f_{R_{2}^{\prime\prime}}$ be defined by setting
$f_{R_{2}^{\prime\prime}}(x)=~f_{R_{2}^{\prime}}(x)$~$\forall x \in V(R_{2}^{\prime})$
and $f_{R_{2}^{\prime\prime}}(x) = 4 - \deg_{R_{2}^{\prime\prime}}(x)$~$\forall x \notin V(R_{2}^{\prime})$.
Let $R_{2}^{\prime\prime\prime}$ be the graph formed from $R_{2}^{\prime}$
by instead subdividing the new $u_{2}v_{2}$ edge with a vertex $w$,
and let $f_{R_{2}}^{\prime\prime\prime}$ be defined by
$f_{R_{2}}^{\prime\prime\prime}(w)=2$
and $f_{R_{2}}^{\prime\prime\prime}(x) =~\!f_{R_{2}}^{\prime}(x)$~$\forall x \in~\!V(R_{2}^{\prime})$.

The proof that $(R,f_{R})$ can be satisfied if and only if (1),(2) or (3) hold is as with case (a)(ii).
Again, we can determine in $O \left( |B_{1}^{\prime}|^{2.5} \right)$ time whether or not 
$\left( R_{1}^{\prime},f_{R_{1}^{\prime}} \right)$ and $\left( R_{1}^{\prime\prime},f_{R_{1}^{\prime\prime}} \right)$ 
can be satisfied,
and if at least one can
then we return to the start of the iterative loop with $B_{2}^{\prime}$ as our new blue graph 
and either 
$(R_{2}^{\prime\prime\prime},f_{R_{2}^{\prime\prime\prime}})$,
$(R_{2}^{\prime},f_{R_{2}^{\prime}})$ or $(R_{2}^{\prime\prime},f_{R_{2}^{\prime\prime}})$ as our augmentation,
according to whether both 
$\left( R_{1}^{\prime},f_{R_{1}^{\prime}} \right)$ and $\left( R_{1}^{\prime\prime},f_{R_{1}^{\prime\prime}} \right)$,
just $\left( R_{1}^{\prime},f_{R_{1}^{\prime}} \right)$,
or just $\left( R_{1}^{\prime\prime},f_{R_{1}^{\prime\prime}} \right)$
can be satisfied, respectively. 
If neither $\left( R_{1}^{\prime},f_{R_{1}^{\prime}} \right)$
nor $\left( R_{1}^{\prime\prime},f_{R_{1}^{\prime\prime}} \right)$ can be satisfied,
we terminate the algorithm.

\subsubsection*{Case (c): \boldmath{$(\overline{u1} \lor \overline{v1}) \land (\overline{u2} \lor \overline{v2})$}}

We will now deal with the remaining case,
which will follow from a detailed investigated of the properties that are forced upon us 
if $(\overline{u1} \lor \overline{v1}) \land (\overline{u2} \lor \overline{v2})$~holds.

Recall that if we have $\overline{u1}$,
then by definition $f_{R}(u) \geq 1$ and $u$ has at least two edges to $R_{1}$,
so $u$ must have only one edge to $R_{2}$,
and hence we have $u2$.
Similarly,
$\overline{v1} \Rightarrow v2, \overline{u2} \Rightarrow u1$ and $\overline{v2} \Rightarrow v1$.
Thus, the only possibilities are $u1 \land \overline{u2} \land \overline{v1} \land v2$
and $\overline{u1} \land u2 \land v1 \land \overline{v2}$.
By swapping $u$ and $v$ if necessary,
we can without loss of generality assume that we have the former.

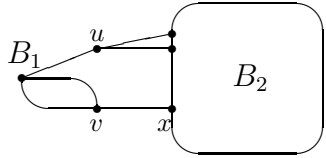
\begin{figure} [ht] 
\setlength{\unitlength}{1cm}
\begin{picture}(10,2)(-3,2.5)

\put(4,4.125){\oval(2,0.75)[t]}
\put(4,2.875){\oval(2,0.75)[b]}
\put(3,2.85){\line(0,1){1.3}}
\put(5,2.85){\line(0,1){1.3}}
\put(2,3.9){\line(1,0){1}}
\put(2,3.9){\line(5,1){1}}
\put(2,3.1){\line(1,0){1}}
\put(1,3.5){\circle*{0.1}}
\put(2,3.1){\circle*{0.1}}
\put(2,3.9){\circle*{0.1}}
\put(3,3.1){\circle*{0.1}}
\put(1,3.5){\line(5,2){1}}
\put(2,3.5){\oval(2,0.8)[bl]}
\put(1,3.1){\oval(2,0.8)[tr]}
\put(3.8,3.4){\large{$B_{2}$}}
\put(0.8,3.7){\large{$B_{1}$}}
\put(1.9,4){$u$}
\put(1.9,2.8){$v$}
\put(2.8,2.8){$x$}
\put(3,3.9){\circle*{0.1}}
\put(3,4.1){\circle*{0.1}}

\end{picture}
\caption{$\textrm{The structure of $B$ in case (c).}$}
\label{D14a}
\end{figure}

Note that the only way to obtain $u1 \land \overline{u2}$
is to have exactly one edge in $B$ from $u$ to $B_{1}$
(or, equivalently, exactly one edge in $R$ from $u$ to $R_{1}$),
exactly two edges in $B$ from $u$ to $B_{2}$,
no edges in $B$ from $u$ to $v$,
and $f_{R}(u)=1$.
Similarly, we must have exactly one edge in $B$ from $v$ to $B_{2}$,
exactly two edges in $B$ from $v$ to $B_{1}$,
and $f_{R}(v)=1$.
Note also that we must have $|B_{1}|=1$,
since otherwise the minimality of $B_{1}$ would imply that $u$ and $v$ 
would both have to have at least two edges in $B$ to $B_{1}$,
which would in turn imply that we would have to have $u2 \land v2$.
Thus, $B$ must be as shown in Figure~\ref{D14a}.

If $|B_{2}|=1$,
then $|R|$ is bounded by a constant
and so we can determine the satisfiabilty of $(R,f_{R})$ in $O(1)$ time
(simply by checking all graphs with $|R|$ vertices to see if any of these do satisfy $(R, f_{R})$).

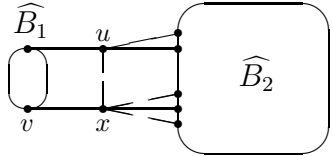
\begin{figure} [ht] 
\setlength{\unitlength}{1cm}
\begin{picture}(10,2)(-3,0)

\put(4,1.625){\oval(2,0.75)[t]}
\put(4,0.375){\oval(2,0.75)[b]}
\put(3,0.35){\line(0,1){1.3}}
\put(5,0.35){\line(0,1){1.3}}
\put(2,1.4){\line(1,0){1}}
\put(2,1.4){\line(5,1){1}}
\put(2,0.6){\line(1,0){1}}
\put(2,0.6){\line(5,1){0.4}}
\put(2.6,0.72){\line(5,1){0.4}}
\put(2,0.6){\line(5,-1){0.4}}
\put(2.6,0.48){\line(5,-1){0.4}}
\put(2,0.6){\line(0,1){0.36}}
\put(2,1.08){\line(0,1){0.36}}
\put(1,1.4){\line(1,0){1}}
\put(1,0.6){\line(1,0){1}}
\put(1,1){\oval(0.5,0.8)}
\put(1,0.6){\circle*{0.1}}
\put(1,1.4){\circle*{0.1}}
\put(2,0.6){\circle*{0.1}}
\put(2,1.4){\circle*{0.1}}
\put(3.8,0.9){\large{$\widehat{B_{2}}$}}
\put(0.8,1.6){\large{$\widehat{B_{1}}$}}
\put(1.9,1.5){$u$}
\put(1.9,0.3){$x$}
\put(0.9,0.3){$v$}
\put(3,0.4){\circle*{0.1}}
\put(3,0.6){\circle*{0.1}}
\put(3,0.8){\circle*{0.1}}
\put(3,1.4){\circle*{0.1}}
\put(3,1.6){\circle*{0.1}}

\end{picture}
\caption{$\textrm{The $2$-vertex-cut $\{ u,x \}$.}$}
\label{D14}
\end{figure}

If $|B_{2}|>1$, then let $x$ denote the neighbour of $v$ in $B_{2}$,
let $\widehat{B_{1}} = B_{1} \cup v$
and let $\widehat{B_{2}} = B_{2} \setminus x$.
Note that $ux$ forms a $2$-vertex-cut where $u$ and $x$ both have just one edge to $\widehat{B_{1}}$
(see Figure~\ref{D14}).
Hence, we can copy case (a) with $B_{1}$ and $B_{2}$ replaced by $\widehat{B_{1}}$ and $\widehat{B_{2}}$, respectively,
to again obtain graphs $B_{1}^{\prime}$ and~$B_{2}^{\prime}$ and appropriate augmentations.
It may be that the graph $B_{1}^{\prime}$ will have a $2$-vertex-cut,
so this time we won't be able to use Lemma~\ref{lemma2}
to determine the satisfiability of augmentations of it.
However, we know that we will have $|B_{1}^{\prime}|=4$,
so the number of vertices in any augmentation of $B_{1}^{\prime}$ will be bounded by a constant,
and hence we will be able to determine satisfiability of these augmentations in $O(1)$ time.

\subsection*{Running Time}

We shall now show that the algorithm takes $O \left( |H|^{2.5} \right)$ time.
It is fairly easy to see that the first three stages can be accomplished within this limit
(in fact, they take only $O \left( |H|^{2} \right)$),
so we will proceed straight to an examination of Stage~4.

We apply Stage 4 to each of the $2$-vertex-connected blocks derived from Stage~3.
It is easy to see that the total number of vertices in all these blocks is at most $2|H|$,
since each vertex of $H$ will only appear in at most two of these,
so it will actually suffice just to deal with the case when $H$ is itself a $2$-vertex-connected block,
i.e.~when we start Stage 4 with only one $2$-vertex-connected block,
and it has $|H|$ vertices.

During Stage 4,
we take a graph $B$
and use it to construct graphs $B_{1}^{\prime}$ and~$B_{2}^{\prime}$,
where $|B_{1}^{\prime}| + |B_{2}^{\prime}| = |B|+2$
and $|B_{2}^{\prime}| < |B|$,
before replacing $B$ with $B_{2}^{\prime}$ and iterating.
Let $B_{1,1}^{\prime},B_{1,2}^{\prime}, \ldots, B_{1,l}^{\prime}$, for some $l$, 
denote the various graphs that take the role of $B_{1}^{\prime}$ during our algorithm.
Since $|B_{2}^{\prime}| < |B|$,
we can only have at most $|H|$ iterations,
and so we must have $\sum_{i} |B_{1,i}^{\prime}| \leq 3|H|$
(by telescoping, since we always have $|B_{1}^{\prime}| + |B_{2}^{\prime}| = |B|+2$).
We need to apply the algorithm given by Lemma~\ref{lemma2}
to at most three augmentations of each $B_{1,i}^{\prime}$,
so the total time taken by all such applications will be at most
$3\lambda \sum_{i} \left( |B_{1,i}^{\prime}|^{2.5} \right) 
\leq 3\lambda \left( \sum_{i}|B_{1,i}^{\prime}| \right)^{2.5}=~O \left( |H|^{2.5} \right)$.

At the start of each iteration,
we wish to determine whether $B$ has any $2$-vertex-cuts and,
if so, find a minimal one.
Using an algorithm from~\cite{hop} for decomposing a
graph into its so-called `triconnected components',
this takes $O(|B|)=O(|H|)$ time.
It is fairly clear that all other operations involved in an iteration of Stage 4,
aside from applications of Lemma~\ref{lemma2},
can also be accomplished within $O(|H|)$ time,
so (since we recall that there are at most $|H|$ iterations)
this all takes $O \left( |H|^{2} \right)$ time in total
(in fact, by careful bookkeeping, this could be reduced to $O(|H|)$).
Hence,
it follows that the whole algorithm takes $O \left( |H|^{2.5} \right)$ time.

\subsection*{Comments}

By keeping track of all the operations,
the algorithm can be used to find an explicit $4$-regular planar multigraph $G \supset H$ if such a graph exists,
also in $O \left( |H|^{2.5} \right)$ time.
If $H$ is simple,
then we can also obtain a $4$-regular simple planar graph $G^{\prime} \supset H$
without affecting the order of the overall running time,
using the proof of Theorem~\ref{bounded991}.

\newpage
\section{General Bounds} \label{general}

In Section~\ref{subs}, we completed our picture of $P_{n,d_{1},d_{2},D_{1},D_{2}}$
for the specific case when $(d_{2}(n),D_{1}(n))=(5,0)$ $\forall n$.
In this section,
we shall deduce that the same results actually also hold for general $d_{2}(n)$ and \textit{bounded} $D_{1}(n)$,
apart from for some trivial differences.

We will start by showing (in Theorem~\ref{bounded791}) that if 
$D_{1}(n) \leq K$~$\forall n$, for some~$K$, and $d_{2}(n)>0$~$\forall n$,
then $P_{n,d_{1},d_{2},D_{1},D_{2}}$ behaves in exactly the same way as $P_{n,d_{1},5,0,D_{2}}$
(in terms of whether or not the probabilities 
of being connected or of containing given components or subgraphs are bounded away from $0$ or $1$).
This is simply because we may use our appearance results to see that $P_{n,d_{1},5,0,D_{2}}$
will a.a.s.~satisfy these more restrictive bounds on the minimum and maximum degrees anyway.

We shall then look at what happens if $D_{1}(n) \leq K$~$\forall n$ and $d_{2}=0$~$\forall n$.
We will first show (in Theorem~\ref{bounded147}) that 
$P_{n,0,0,D_{1},D_{2}}$ behaves in exactly the same way as $P_{n,0,0,0,D_{2}}$,
and then see (in Theorems~\ref{bounded803}-\ref{bounded806})
that this latter graph actually has all the standard characteristics too,
except that obviously
$\mathbf{P}[P_{n,0,0,0,D_{2}} \textrm{ will be connected}]=0$
and $\mathbf{P}[P_{n,0,0,0,D_{2}}$ 
will have a component isomorphic to $H]=1$ if $|H|=1$.

This leaves the case when $D_{1}(n)$ is unbounded,
which we shall look at in Section~\ref{unbounded}.

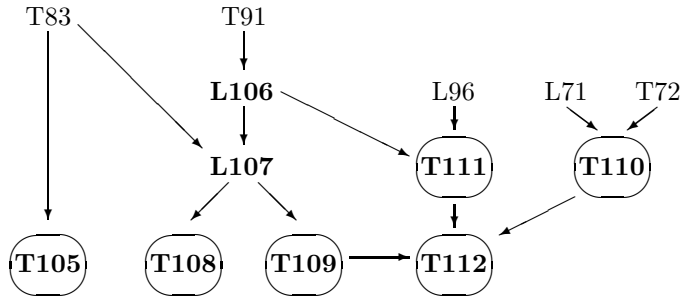
\begin{figure} [ht]
\setlength{\unitlength}{1cm}
\begin{picture}(20,4)(-1.65,0)

\put(0.1,3.5){T\ref{bounded11}}
\put(2.7,3.5){T\ref{bounded404}}

\put(2.55,2.5){\textbf{L\ref{bounded811}}}
\put(5.5,2.5){L\ref{bounded904}}

\put(2.55,1.5){\textbf{L\ref{bounded812}}}
\put(5.35,1.5){\textbf{T\ref{bounded802}}}
\put(7.45,1.5){\textbf{T\ref{bounded807}}}

\put(-0.05,0.2){\textbf{T\ref{bounded791}}}
\put(1.75,0.2){\textbf{T\ref{bounded147}}}
\put(3.35,0.2){\textbf{T\ref{bounded803}}}
\put(5.35,0.2){\textbf{T\ref{bounded806}}}

\put(0.4,0.3){\oval(1,0.8)}
\put(2.2,0.3){\oval(1,0.8)}
\put(3.8,0.3){\oval(1,0.8)}
\put(5.8,0.3){\oval(1,0.8)}

\put(5.8,1.6){\oval(1,0.8)}
\put(7.9,1.6){\oval(1,0.8)}

\put(0.4,3.4){\vector(0,-1){2.5}}
\put(0.8,3.5){\vector(1,-1){1.65}}
\put(3,3.4){\vector(0,-1){0.5}}
\put(3,2.4){\vector(0,-1){0.5}}
\put(2.8,1.4){\vector(-1,-1){0.5}}
\put(3.2,1.4){\vector(1,-1){0.5}}
\put(3.5,2.6){\vector(2,-1){1.7}}
\put(4.4,0.4){\vector(1,0){0.8}}
\put(5.8,2.4){\vector(0,-1){0.3}}
\put(5.8,1.1){\vector(0,-1){0.3}}
\put(7.4,1.2){\vector(-2,-1){1}}

\put(7,2.5){L\ref{bounded1}}
\put(8.2,2.5){T\ref{bounded311}}
\put(8.5,2.4){\vector(-4,-3){0.4}}
\put(7.3,2.4){\vector(4,-3){0.4}}

\end{picture}

\caption{The structure of Section~\ref{general}.}
\end{figure}

We start with our aforementioned result for when $d_{2}(n)$ is strictly positive and $D_{1}(n)$ is bounded
(these form condition (d) in the statement of the theorem,
while conditions (a) and (b) are just trivial necessities 
and (c) is our ever-present condition that $D_{2}(n) \geq 3$):

\begin{Theorem} \label{bounded791}
Let $K$ be fixed
and let $d_{1}(n)$, $d_{2}(n)$, $D_{1}(n)$ and $D_{2}(n)$ be integer-valued functions 
that for all large $n$ satisfy
(a) $d_{1}(n) \leq \min \{ 5,d_{2}(n),D_{2}(n) \}$ and $D_{1}(n) \leq D_{2}(n)$,
(b) $(d_{1}(n), D_{2}(n)) \notin \{ (3,3),(5,5) \}$ for odd $n$,
(c) $D_{2}(n) \geq 3$,
and (d) $d_{2}(n) > 0$ and $D_{1}(n) \leq K$.
Then
\begin{displaymath}
\mathbf{P}[P_{n,d_{1},5,0,D_{2}} \in \mathcal{P}(n,d_{1},d_{2},D_{1},D_{2})] \to 1 \textrm{ as } n \to \infty.
\end{displaymath} 
\end{Theorem} 
\textbf{Proof}
Without loss of generality, $K \geq 5$.
Let us split the proof into two different cases for
(i) the values of $n$ for which $D_{2}(n) \geq K$
and (ii) the values of~$n$ for which $D_{2}(n) < K$.

For case (i),
it clearly suffices to prove the result for when $d_{1},d_{2}$ and $D_{1}$
are arbitrary \textit{fixed constants} in
$\{ 0,1, \ldots, 5 \}, \{ 1,2, \ldots, 5 \}$ and $\{ 0,1, \ldots, K \}$,
respectively, satisfying $d_{1} \leq d_{2}$
(if we choose $d_{1} \in \{ 3,5 \}$,
then we may ignore any odd values of $n$ for which $D_{2}(n) = d_{1}$).
But note that Theorem~\ref{bounded11},
on appearances,
then implies
$\mathbf{P}
[\delta(P_{n,d_{1},5,0,D_{2}}) \leq \max \{ d_{1},1 \} \textrm{ and } \Delta(P_{n,d_{1},5,0,D_{2}}) \geq K] \to 1$
as $n \to \infty$,
and so the result follows.

For case (ii),
the proof is the same except that we also take $D_{2}$ to be an arbitrary fixed constant in $\{ 3,4, \ldots, K-1 \}$
satisfying $D_{2} \geq \max \{ d_{1}, D_{1} \}$.
This~time,
Theorem~\ref{bounded11} gives
$\mathbf{P}
[\delta(P_{n,d_{1},5,0,D_{2}}) \leq \max \{ d_{1},1 \} \textrm{ and } \Delta(P_{n,d_{1},5,0,D_{2}}) = D_{2}] \to~\!1$
as $n \to \infty$,
and again the result follows.
$\phantom{qwerty} 
\setlength{\unitlength}{.25cm}
\begin{picture}(1,1)
\put(0,0){\line(1,0){1}}
\put(0,0){\line(0,1){1}}
\put(1,1){\line(-1,0){1}}
\put(1,1){\line(0,-1){1}}
\end{picture}$ \\
\\

Hence, it follows that
all our results for $P_{n,d_{1},5,0,D_{2}}$ also hold for $P_{n,d_{1},d_{2},D_{1},D_{2}}$
if $d_{2}(n) > 0$~$\forall n$ and $D_{1}(n) \leq K$ $\forall n$. \\
\\

In the remainder of this section,
we shall deal with the case $P_{n,0,0,0,D_{2}}$,
i.e.~a graph with at least one isolated vertex and with maximum degree at most $D_{2}(n)$.
We shall see (in Theorem~\ref{bounded147}) that for all fixed $K \leq \liminf_{n \to \infty} D_{2}(n)$ we have
$\Delta(P_{n,0,0,0,D_{2}}) \geq K$ a.a.s.,
and so results for $P_{n,0,0,D_{1},D_{2}}$ when $D_{1}(n) \leq K$ $\forall n$ will just follow automatically
from results for $P_{n,0,0,0,D_{2}}$,
which we shall then investigate in Theorems~\ref{bounded803}--\ref{bounded806}. \\
\\

Analogously to Theorem~\ref{bounded791},
we shall prove Theorem~\ref{bounded147} via a result on appearances in $P_{n,0,0,0,D_{2}}$.
In order to obtain this appearance result,
we first note the following simple lemma on
$\mathbf{P}[P_{n,0,5,0,D_{2}} \in \mathcal{P}(n,0,0,0,D_{2})]$,
which will also be useful later on in this section:

\begin{Lemma} \label{bounded811}
Let $D_{2}(n)$ be an integer-valued function satisfying $D_{2}(n) \geq 3$~$\forall n$.
Then $\liminf_{n \to \infty} \mathbf{P}[P_{n,0,5,0,D_{2}} \in \mathcal{P}(n,0,0,0,D_{2})]>0.$
\end{Lemma}
\textbf{Proof}
By Theorem~\ref{bounded404} on components,
$\liminf_{n \to \infty} \mathbf{P}[\delta(P_{n,0,5,0,D_{2}})=~0] >~0.$
Thus, the result follows.~$\phantom{qwerty} 
\setlength{\unitlength}{.25cm}
\begin{picture}(1,1)
\put(0,0){\line(1,0){1}}
\put(0,0){\line(0,1){1}}
\put(1,1){\line(-1,0){1}}
\put(1,1){\line(0,-1){1}}
\end{picture}$ \\

Recall from Definition~\ref{maxdisjoint} that
$\widehat{f_{H}^{0}}(G)$ denotes the maximum size of a set of totally edge-disjoint appearances of $H$ in $G$.
It now follows from Lemma~\ref{bounded811} that we have the following appearance result for~$P_{n,0,0,0,D_{2}}$:

\begin{Lemma} \label{bounded812}
Let $H$ be a fixed connected planar graph on $\{ 1,2, \ldots, h \}$.
Then there exists $\beta(h)>0$ such that,
given any integer-valued function $D_{2}(n)$ satisfying
$\liminf_{n \to \infty} D_{2}(n) \geq \max \{ \Delta(H), \deg_{H}(1)+1, 3 \}$,
we have
\begin{displaymath}
\mathbf{P}[\widehat{f_{H}^{0}} (P_{n,0,0,0,D_{2}}) \leq \beta n] < e^{- \beta n} 
\textrm{ for all sufficiently large } n.
\end{displaymath}
\end{Lemma}
\textbf{Proof}
This follows from Theorem~\ref{bounded11} and Lemma~\ref{bounded811}.
$\phantom{qwerty}
\setlength{\unitlength}{.25cm}
\begin{picture}(1,1)
\put(0,0){\line(1,0){1}}
\put(0,0){\line(0,1){1}}
\put(1,1){\line(-1,0){1}}
\put(1,1){\line(0,-1){1}}
\end{picture}$ 

\phantom{p}

As mentioned, an important consequence of this last result 
is that for any fixed $K \leq \liminf_{n \to \infty} D_{2}(n)$,
we have $\mathbf{P}[\Delta(P_{n,0,0,0,D_{2}}) \geq K] \to 1$ as $n \to \infty$,
so

\begin{Theorem} \label{bounded147}
Let $K$ be a fixed constant and let $D_{1}(n)$ and $D_{2}(n)$ be integer-valued functions satisfying
$D_{1}(n) \leq K$ $\forall n$ and $D_{2}(n) \geq 3$ $\forall n$.
Then
\begin{displaymath}
\mathbf{P}[P_{n,0,0,0,D_{2}} \in \mathcal{P}(n,0,0,D_{1},D_{2})] \to 1 \textrm{ as } n \to \infty.
\end{displaymath} 
\end{Theorem}

\phantom{p}

Thus, results for $P_{n,0,0,D_{1},D_{2}}$ when $D_{1}(n) \leq K$ will actually just be the same as for $P_{n,0,0,0,D_{2}}$.
Therefore, to complete our picture of $P_{n,d_{1},d_{2},D_{1},D_{2}}$ for the case when $D_{1}(n) \leq K$ $\forall n$,
it will suffice just to deal with the case $d_{1}(n)=d_{2}(n)=D_{1}(n)=0$~$\forall n$.
Clearly, $\mathbf{P}[P_{n,0,0,0,D_{2}} \textrm{ will be connected}]=0$,
and so we are only left with looking at the limiting probabilities for
$P_{n,0,0,0,D_{2}}$ having a component isomorphic to $H$
and for $P_{n,0,0,0,D_{2}}$ having a copy of $H$,
for given $H$. \\
\\

A lower bound for $\mathbf{P}[P_{n,0,0,0,D_{2}}$ will have a component isomorphic to $H]$
may in fact be obtained exactly as in Section~\ref{cpts}:

\begin{Theorem} \label{bounded803}
Let $D_{2}(n)$ be an integer-valued function satisfying $D_{2}(n) \geq\!~\!\!3$~$\forall n$,
and let $t$ be a constant.
Then, given any connected planar graphs
$H_{1},H_{2}, \ldots, H_{k}$ 
with
$\Delta(H_{i}) \leq \liminf_{n \to \infty} D_{2}(n)$~$\forall i$,
we have
\begin{eqnarray*}
& & \liminf_{n \to \infty}
\mathbf{P}
\Big[\bigcap_{i \leq k} \left(\textrm{$P_{n,0,0,0,D_{2}}$ 
will have $\geq t$ components}\right. \\
& & \textrm{\phantom{wwwwwwwwwww}with order-preserving isomorphisms to $H_{i}$} ) ] 
> 0.
\end{eqnarray*}
\end{Theorem}
\textbf{Proof}
We may use the same proof as for Lemma~\ref{bounded7}
(by deleting the associated cut-edges of some appearances),
with Theorem~\ref{bounded11} replaced by Lemma~\ref{bounded812}.~$
\setlength{\unitlength}{.25cm}
\begin{picture}(1,1)
\put(0,0){\line(1,0){1}}
\put(0,0){\line(0,1){1}}
\put(1,1){\line(-1,0){1}}
\put(1,1){\line(0,-1){1}}
\end{picture}$ 

\phantom{p}

Clearly, 
$\mathbf{P}[P_{n,0,0,0,D_{2}} \textrm{ will have a component isomorphic to }H]=1$
if $H$ is an isolated vertex,
by definition.
However, if $H$ is not an isolated vertex then
we are able to bound the probability away from $1$ by the following result:

\begin{Theorem} \label{bounded807}
There exists a constant $\epsilon>0$ such that,
given any integer-valued function $D_{2}(n)$ with $D_{2}(n) \geq 3$ $\forall n$,
we have
\begin{eqnarray*}
\mathbf{P}[P_{n,0,0,0,D_{2}} 
\textrm{ will consist of \emph{exactly} one isolated vertex plus a connected graph}] \\
> \epsilon 
\textrm{ } \forall n.
\end{eqnarray*}
\end{Theorem}
\textbf{Proof}
Clearly,
the result holds for $n \leq 2$.
Now let us choose any $n \geq 3$ and any $D_{2}(n)$.
We shall find a constant $\epsilon$,
independent of these choices of $n$ and $D_{2}(n)$,
satisfying the conditions of the theorem.

Let $\mathcal{I}(n,D_{2},k)$ denote the set of all graphs in $\mathcal{P}(n,0,0,0,D_{2})$
with exactly~$k$ isolated vertices.
Note that 
$|\mathcal{I}(n,D_{2},1)| = n|\mathcal{P}(n-1,1,5,0,D_{2})|$
and that 
$|\mathcal{I}_{c}(n,D_{2},1)| = n|\mathcal{P}_{c}(n-1,1,5,0,D_{2})|$,
where $\mathcal{P}_{c}(n-1,1,5,0,D_{2})$ denotes the set of connected graphs in $\mathcal{P}(n-1,1,5,0,D_{2})$
and $\mathcal{I}_{c}(n,D_{2},1)$ denotes the set of graphs in $\mathcal{I}(n,D_{2},1)$
that are connected apart from the isolated vertex.
By Theorem~\ref{bounded311},
there exists a strictly positive constant $c$
(independent of $n$ and $D_{2}(n)$)
such that
$\frac{ |\mathcal{P}_{c}(n-1,1,5,0,D_{2})| }{ |\mathcal{P}(n-1,1,5,0,D_{2})| } \!>\! c$.
Thus, $|\mathcal{I}_{c}(n,D_{2},1)| \!>\! c |\mathcal{I}(n,D_{2},1)|$,
and so it suffices to find a strictly positive constant $\epsilon^{\prime}$,
independent of $n$ and $D_{2}(n)$,
such that
$\frac{ |\mathcal{I}(n,D_{2},1)| }{ |\mathcal{P}(n,0,0,0,D_{2})| } \!>\! \epsilon^{\prime}$,
i.e.~$\mathbf{P}[P_{n,0,0,0,D_{2}} \textrm{ will have exactly one isolated vertex}] \!>~\!\!\epsilon^{\prime}$. \\

The remainder of the proof 
will now be a `downwards cascade' argument similar to that of Theorem~\ref{bounded311}. 

Let $k \in \{ 2,3, \ldots, n-2 \}$, let $G \in \mathcal{I}(n,D_{2},k)$,
and let $G^{*}$ denote the graph of order $n-k$ obtained by deleting all $k$ isolated vertices from $G$.
Starting from~$G$, we shall create a new graph $G^{\prime} \in \mathcal{I}(n,D_{2},k-1)$
by considering different cases depending on $G^{*}$: \\
\\
Case (a) If $G^{*}$ has $> \frac{n-k}{43}$ vertices of degree $<D_{2}(n)$ \\
Starting with $G$,
simply insert an edge between an isolated vertex (we have $k$ choices)
and a non-isolated vertex of degree $<D_{2}(n)$
(we have $> \frac{n-k}{43}$ choices). 
\begin{figure} [ht]
\setlength{\unitlength}{1cm}
\begin{picture}(20,1.5)(-2,0)

\put(1,0.5){\line(0,1){0.5}}
\put(2.5,0.5){\line(0,1){0.5}}
\put(6.5,0.5){\line(0,1){0.5}}
\put(8,0.5){\line(0,1){0.5}}

\put(5.5,0.75){\line(1,0){1}}

\put(1.75,1){\oval(1.5,1)[t]}
\put(1.75,0.5){\oval(1.5,1)[b]}
\put(7.25,1){\oval(1.5,1)[t]}
\put(7.25,0.5){\oval(1.5,1)[b]}

\put(0,0.75){\circle*{0.1}}
\put(1,0.75){\circle*{0.1}}
\put(5.5,0.75){\circle*{0.1}}
\put(6.5,0.75){\circle*{0.1}}

\put(3.5,0.75){\vector(1,0){1}}

\end{picture}

\caption{Reducing the number of isolated vertices in case (a).}
\end{figure}
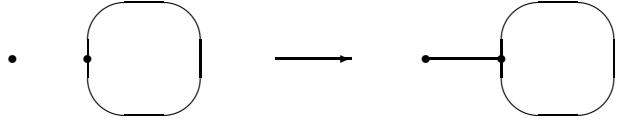
\\
Case (b) If $G^{*}$ has $\leq \frac{n-k}{43}$ vertices of degree $<D_{2}(n)$ 
(in which case $D_{2}(n) \leq 6$,
since $e(G^{*}) < 3(n-k)$
and so we can only have at most $\frac{6(n-k)}{7}$ vertices of degree~$\geq~\!7$) \\
By Lemma~\ref{bounded1},
$G$ must contain at least $\frac{n-k}{43}$ cycles of size $\leq 6$.
Delete an edge $uv$ in one of these cycles
(we have at least $\frac{3}{D_{2}(n)+D_{2}(n)^{2}+D_{2}(n)^{3}+D_{2}(n)^{4}}\frac{n-k}{43}
\geq \frac{3}{6+6^{2}+6^{3}+6^{4}}\frac{n-k}{43}$ choices for this edge,
since each cycle must contain at least $3$ edges and each edge is in at most 
$(D_{2}(n)-1)^{m-2}<D_{2}(n)^{m-2}$ cycles of size~$m$),
and insert an edge between $u$ and an isolated vertex (we have at least $k$ choices for this). \\

\begin{figure} [ht]
\setlength{\unitlength}{1cm}
\begin{picture}(20,1)(-2,0)

\put(1,0.5){\line(0,1){0.5}}
\put(2.5,0.5){\line(0,1){0.5}}
\put(8,0.5){\line(0,1){0.5}}

\put(5.5,0.75){\line(4,1){1}}

\put(1.75,1){\oval(1.5,1)[t]}
\put(1.75,0.5){\oval(1.5,1)[b]}
\put(7.25,1){\oval(1.5,1)[t]}
\put(7.25,0.5){\oval(1.5,1)[b]}

\put(0,0.75){\circle*{0.1}}
\put(1,0.5){\circle*{0.1}}
\put(1,1){\circle*{0.1}}
\put(5.5,0.75){\circle*{0.1}}
\put(6.5,0.5){\circle*{0.1}}
\put(6.5,1){\circle*{0.1}}

\put(3.5,0.75){\vector(1,0){1}}

\put(1,0.75){\oval(1,0.5)[r]}
\put(6.5,0.75){\oval(1,0.5)[r]}

\put(0.7,0.95){$u$}
\put(0.7,0.4){$v$}

\put(6.2,1.05){$u$}
\put(6.2,0.4){$v$}

\end{picture}

\caption{Reducing the number of isolated vertices in case (b).}
\end{figure}
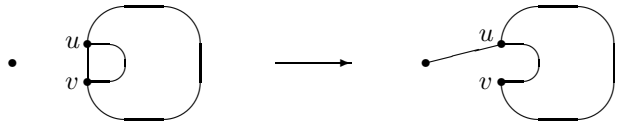

\phantom{p}

In both cases,
we have at least $\frac{3}{6+6^{2}+6^{3}+6^{4}}\frac{k(n-k)}{43}$
choices when constructing $G^{\prime}$ from $G$.
Thus, we have at least $\frac{3}{6+6^{2}+6^{3}+6^{4}}\frac{k(n-k)}{43}|\mathcal{I}(n,D_{2}(n),k)|$
ways to construct a graph in $\mathcal{I}(n,D_{2}(n),k-1)$. \\

Given one of our constructed graphs,
there are at most $2$ possibilities for how the graph was obtained (case (a) or (b)).

If case (a) was used,
then we can re-obtain the original graph simply by deleting the inserted edge,
for which there are at most $3(n-(k-1))-3 = 3(n-k)$
possibilities 
(including the case when $n-(k-1)<3$).

If case (b) was used,
then we can re-obtain the original graph by deleting the inserted edge ($\leq 3(n-k)$ possibilities)
and re-inserting the deleted edge
(at most $D_{2}(n)^{2} + D_{2}(n)^{3} + D_{2}(n)^{4} + D_{2}(n)^{5}$ possibilities,
once the inserted edge is located,
since the deleted edge $uv$ was originally part of a cycle of size $\leq 6$,
and so $v$ is still at distance at most $5$ from $u$).

Thus, 
recalling that $D_{2} \leq 6$ if case (b) was used,
we find that there are at most
$3(n-k) + 3(n-k)(6^{2}+6^{3}+6^{4}+6^{5}) = 
3(1+6^{2}+6^{3}+6^{4}+6^{5})(n-k)$ 
possibilities for the original graph in total.

Let $\alpha = \frac{1}{43(6+6^{2}+6^{3}+6^{4})
(1+6^{2}+6^{3}+6^{4}+6^{5})}$.
Then, for all $k \in \{ 2,3, \ldots, n-2 \}$, we have
$|\mathcal{I}(n,D_{2},k-1)| \geq \alpha k|\mathcal{I}(n,D_{2},k)|
\geq \alpha(k-1) |\mathcal{I}(n,D_{2},k)|$. 

Let
$p_{k} = \frac{|\mathcal{I}(n,D_{2},k+1)|}{|\mathcal{P}(n,0,0,0,D_{2})|}$
and note that
$p_{0}$
is the probability that 
$P_{n,0,0,0,D_{2}}$ will have exactly one isolated vertex.
Since $|\mathcal{I}(n,D_{2},k+1)| \leq \frac{|\mathcal{I}(n,D_{2},k)|}{\alpha k}$
for all $k \in \{ 1,2, \ldots, n-3 \}$,
we have $p_{k} \leq \frac{p_{0}}{\alpha^{k}k!}$
for all $k \in \{ 0,1, \ldots, n-3 \}$.
Note that we also have $\sum_{k \geq 0}^{n-3} p_{k} = 1 - \frac{1}{|\mathcal{P}(n,0,0,0,D_{2})|} \geq \frac{7}{8}$
for $n \geq 3$
(since every graph in $\mathcal{P}(n,0,0,0,D_{2})$ has at least one isolated vertex
and can only have more than $n-2$ if it is the empty graph $E_{n}$),
so $\sum_{k \geq 0} \frac{p_{0}}{\alpha^{k}k!} \geq \frac{7}{8}$
and hence 
$p_{0} \geq \frac{7}{8} \left( \sum_{k \geq 0} \frac{\left( \frac{1}{\alpha} \right)^{k}}{k!} \right)^{-1} 
= \frac{7 e^{-\frac{1}{\alpha}}}{8}$.
\phantom{qwerty}
\setlength{\unitlength}{0.25cm}
\begin{picture}(1,1)
\put(0,0){\line(1,0){1}}
\put(0,0){\line(0,1){1}}
\put(1,1){\line(-1,0){1}}
\put(1,1){\line(0,-1){1}}
\end{picture} \\
\\

It now only remains to look at $\mathbf{P}[P_{n,0,0,0,D_{2}} \textrm{ will have a copy of } H]$.
As in Section~\ref{subs},
the following two results show that the behaviour of this probability
(in terms of whether or not it is bounded away from $1$)
depends only on whether there are arbitrarily many $n$ for which $H$ has a $D_{2}(n)$-regular component:

\begin{Theorem} \label{bounded802}
Let $D_{2}(n)$ be an integer-valued function satisfying $D_{2}(n) \!\geq~\!\!3$~$\forall n$,
and let $H$ be a planar graph with components $H_{1},H_{2}, \ldots, H_{k}$, for some $k$.
Suppose for all $i$ we have 
$\Delta(H_{i}) \leq \liminf_{n \to \infty} D_{2}(n)$ and $\delta(H_{i}) < \liminf_{n \to \infty} D_{2}(n)$.
Then $\exists \beta > 0$ and $\exists N$ such that
\begin{eqnarray*}
& & \mathbf{P} \Big[
\textrm{$P_{n,0,0,0,D_{2}}$
will \emph{not} have a set of $\beta n$ vertex-disjoint} \\ 
& & \phantom{wwwwwwwwww}\textrm{induced order-preserving copies of $H$}\Big] 
< e^{- \beta n}~ 
\forall n \geq N.
\end{eqnarray*}
\end{Theorem}
\textbf{Proof}
This follows immediately when Lemma~\ref{bounded904}
(the analogous result for $P_{n,0,5,0,D_{2}}$) is combined with Lemma~\ref{bounded811}.
$\phantom{qwerty}
\setlength{\unitlength}{.25cm}
\begin{picture}(1,1)
\put(0,0){\line(1,0){1}}
\put(0,0){\line(0,1){1}}
\put(1,1){\line(-1,0){1}}
\put(1,1){\line(0,-1){1}}
\end{picture}$ \\

\begin{Theorem} \label{bounded806}
Let $D_{2} \geq 3$ be a fixed integer,
and let $H$ be a planar graph with components $H_{1},H_{2}, \ldots, H_{k}$, for some $k$.
Suppose for all $i$ we have $\Delta(H_{i}) \leq D_{2}$
and that for some $i$ $H_{i}$ is $D_{2}$-regular.
Then
\begin{displaymath}
\limsup_{n \to \infty} \mathbf{P} [P_{n,0,0,0,D_{2}}
\textrm{ will have a copy of $H$}] < 1,
\end{displaymath}
but for any given constant $t$,
\begin{eqnarray*}
& & \liminf_{n \to \infty} \mathbf{P} \Big[
\textrm{$P_{n,0,0,0,D_{2}}$
will have a set of $t$ vertex-disjoint} \\
& & \phantom{wwwwwwwwwwi}\textrm{induced order-preserving copies of $H$}\Big]
>0.
\end{eqnarray*}
\end{Theorem}
\textbf{Proof}
The upper bound follows from Theorem~\ref{bounded807}
and the lower bound from Theorems~\ref{bounded803} and~\ref{bounded802},
exactly as in Theorem~\ref{bounded1002}.
$\phantom{qwerty}
\setlength{\unitlength}{.25cm}
\begin{picture}(1,1)
\put(0,0){\line(1,0){1}}
\put(0,0){\line(0,1){1}}
\put(1,1){\line(-1,0){1}}
\put(1,1){\line(0,-1){1}}
\end{picture}$

\newpage
\section{\boldmath{$D_{1}(n) \to \infty$}} \label{unbounded}

We now have a complete picture of $P_{n,d_{1},d_{2},D_{1},D_{2}}$ for the case when $D_{1}(n)$ is bounded
above by an arbitrary constant (see page~\pageref{sum}),
i.e.~if we are given any functions at all for $d_{1}(n)$, $d_{2}(n)$ and $D_{2}(n)$
and we are given a function $D_{1}(n)$ with $\limsup_{n \to \infty} D_{1}(n) < \infty$,
then we can tell how likely it is
(in terms of whether the probabilities are bounded away from $0$ and/or $1$)
that $P_{n,d_{1},d_{2},D_{1},D_{2}}$
will be connected or contain any particular component/subgraph.
This leaves the matter of what happens when $\limsup_{n \to \infty} D_{1}(n) = \infty$,
which we shall now discuss very briefly in this section. \\

Recall that our picture of $P_{n,d_{1},d_{2},D_{1},D_{2}}$
for the case when $D_{1}(n)$ is bounded followed immediately from our results for $P_{n,d_{1},d_{2},0,D_{2}}$,
since we were able to show that the maximum degree in this latter graph will a.a.s.~be larger than any given constant.
Hence,
if we could obtain a higher bound for $\Delta(P_{n,d_{1},d_{2},0,D_{2}})$,
then this would automatically enable us to extend our current results.
In fact,
it has very recently been shown in \cite{mcd2} that (a.a.s.)~the standard random planar graph
$P_{n,0,5,0,n-1}$ has maximum degree of the order of $\log n$,
and it seems likely that such a result should also hold for $P_{n,d_{1},d_{2},0,D_{2}}$ in general.
Hence, it is probable that the description of $P_{n,d_{1},d_{2},D_{1},D_{2}}$
given on page~\pageref{sum} will still hold even if $D_{1}(n)$ is allowed to grow slowly with $n$.

If we allow $D_{1}(n)$ to become very large, 
then of course eventually our picture of $P_{n,d_{1},d_{2},D_{1},D_{2}}$ will have to change,
since (for example) $P_{n,d_{1},d_{2},n-1,n-1}$ will clearly be connected!
We now conclude this thesis with one final nice result
to show that this change will a.a.s.~have happened by the time $D_{1}(n)=n-o(n)$:

\begin{Theorem} \label{bounded315}
Let $d_{1}(n)$, $d_{2}(n)$, $D_{1}(n)$ and $D_{2}(n)$ 
be integer-valued functions that for all large $n$ satisfy
(a) $d_{1}(n) \leq \min \{ 5, d_{2}(n), D_{2}(n) \}$ and $D_{1}(n) \leq D_{2}(n)$,
and (b) $d_{2}(n)>0$ and $D_{1}(n)=n-o(n)$.
Then
\begin{displaymath}
\mathbf{P}[P_{n,d_{1},d_{2},D_{1},D_{2}} \textrm{ will be connected}] \to 1 \textrm{ as } n \to \infty.
\end{displaymath}
\end{Theorem}
\textbf{Proof}
Let $d_{1}(n)$, $d_{2}(n)$, $D_{1}(n)$ and $D_{2}(n)$ 
be as in the statement of the theorem,
and let $\mathcal{G}_{n}$ denote the set of graphs in $\mathcal{P}(n,d_{1},d_{2},D_{1},D_{2})$ 
that are \textit{not} connected.
We shall use $\mathcal{G}_{n}$ to construct so many graphs in $\mathcal{P}(n,d_{1},d_{2},D_{1},D_{2})$
that we must have $\frac{|\mathcal{G}_{n}|}{|\mathcal{P}(n,d_{1},d_{2},D_{1},D_{2})|} \to 0$ as $n \to \infty$.

Let $G$ be an arbitrary graph in $\mathcal{G}_{n}$, for some $n$.
Since $D_{1}(n) = n-o(n)$,
we may assume that $n$ is large enough that 
only one component of $G$ contains a vertex with degree $\Delta(G)$.
Let us call this component $C$.
Let $x$ be a vertex with deg$(x) = \delta(C)$,
and let us choose a vertex $u \in V(C) \setminus x$ with degree at most~$\Delta(G)-2$.
Note that $e(G) \leq 3n-6$ implies $\sum_{v \in V(G)} \deg(v) \leq 6n-12$,
so we can assume that $n$ is large enough that there
only exist at most $6$ vertices with degree greater than $\Delta(G)-2$.
Thus, we have at least $|C \setminus x|-6 =|C|-1-6 \geq D_{1}(n)-6$
ways to choose a vertex $u$.

Before we continue with our argument,
let us choose a vertex $w \in G \setminus C$ in one of two ways,
depending on whether $|G \setminus C|=1$ or $|G \setminus C| \geq 2$.
If $|G \setminus C|=1$,
then we let $w$ be the unique vertex in~$G \setminus C$.
If $|G \setminus C| \geq 2$,
then let $y$ be a vertex in $G \setminus C$ with deg$(y)=\delta(G \setminus C)$
and let $w$ be any vertex in $(G \setminus C) \setminus y$.

Let $G^{*}$ denote the graph formed from $G$ by inserting the edge $uw$.
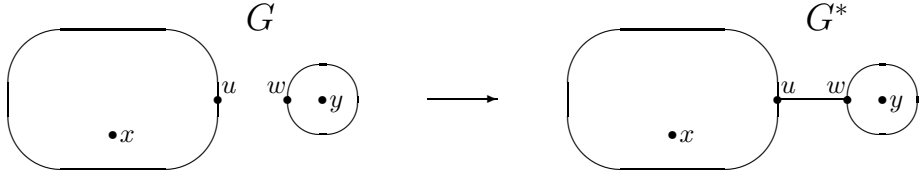
\begin{figure} [ht]
\setlength{\unitlength}{0.93cm}
\begin{picture}(20,2.6)(0,0)

\put(1.5,1){\oval(3,2)}
\put(4.5,1){\oval(1,1)}

\put(3,1){\circle*{0.1}}
\put(1.5,0.5){\circle*{0.1}}
\put(4,1){\circle*{0.1}}
\put(4.5,1){\circle*{0.1}}

\put(3.05,1.1){$u$}
\put(1.6,0.4){$x$}
\put(3.7,1.1){$w$}
\put(4.6,0.9){$y$}

\put(3.4,2){\Large${G}$}

\put(6,1){\vector(1,0){1}}

\put(9.5,1){\oval(3,2)}
\put(12.5,1){\oval(1,1)}

\put(11,1){\circle*{0.1}}
\put(9.5,0.5){\circle*{0.1}}
\put(12,1){\circle*{0.1}}
\put(12.5,1){\circle*{0.1}}

\put(11.05,1.1){$u$}
\put(9.6,0.4){$x$}
\put(11.7,1.1){$w$}
\put(12.6,0.9){$y$}

\put(11.4,2){\Large${G^{*}}$}

\put(11,1){\line(1,0){1}}

\end{picture}

\caption{Forming the graph $G^{*}$.}
\end{figure}
Note that $G^{*}$ is still planar,
since $u$ and $w$ were in separate components,
and that we still have $\delta(G^{*}) \geq d_{1}(n)$ and $\Delta(G^{*}) \geq D_{1}(n)$,
since we have not deleted any edges.
Note also that we still have $\Delta(G^{*}) \leq D_{2}(n)$,
since deg$_{G^{*}}(w)=o(n)$ and deg$_{G^{*}}(u) =$ deg$_{G}(u)+1 \leq \Delta (G) -1$,
and that we still have
$\delta(G^{*}) \leq d_{2}(n)$,
since if $|G \setminus C|$ was $1$ then deg$_{G^{*}}(w)=1 \leq d_{2}(n)$
and if $|G \setminus C|$ was at least $2$ then
deg$_{G^{*}}(y)=$ deg$_{G}(y) = \delta_{G} (G \setminus C)$
and deg$_{G^{*}}(x)=$ deg$_{G}(x) = \delta_{G}(C)$
and so $\delta(G^{*}) = \delta(G) \leq d_{2}(n)$.

Hence, since we had $|\mathcal{G}_{n}|$ choices for $G$ and at least $D_{1}(n)-6$ choices for $u$,
we have at least $(D_{1}(n)-6)|\mathcal{G}_{n}|$ ways to construct a graph in $\mathcal{P}(n,d_{1},d_{2},D_{1},D_{2})$.

Let us now consider how many times a graph $G^{*} \in \mathcal{P}(n,d_{1},d_{2},D_{1},d_{2})$
will be constructed.
Clearly, we can re-obtain the original graph $G$ by deleting the edge~$uw$.
Since deg$_{G^{*}}(u)=$ deg$_{G}(u)+1 \leq \Delta(G)-1$,
$w$ is one of at most $n-D_{1}(n)$ vertices that is not adjacent to a vertex of degree $\Delta(G)$.
Once we have located $w$,
we then have only one possibility for the edge $uw$,
since all paths between $w$ and any vertices of degree $\Delta(G)$ must use this edge.

Thus, we have constructed at least $\frac{(D_{1}(n)-6)|\mathcal{G}_{n}|}{(n-D_{1}(n))}$ \textit{distinct}
graphs in the set $\mathcal{P}(n,d_{1},d_{2},D_{1},D_{2})$ and, therefore,
\begin{eqnarray*}
\frac{|\mathcal{G}_{n}|}{|\mathcal{P}(n,d_{1},d_{2},D_{1},D_{2})|} & \leq & \frac{(n-D_{1}(n))}{D_{1}(n)-6} \\
& = & \frac{o(n)}{n-o(n)} \\
& \to & 0 \textrm{ as } n \to \infty.
\phantom{qwerty}
\setlength{\unitlength}{.25cm}
\begin{picture}(1,1)
\put(0,0){\line(1,0){1}}
\put(0,0){\line(0,1){1}}
\put(1,1){\line(-1,0){1}}
\put(1,1){\line(0,-1){1}}
\end{picture}
\end{eqnarray*}

\end{document}